\documentclass[hidelinks,onefignum,onetabnum]{siamart220329}

\usepackage{lipsum}
\usepackage{amsfonts}
\usepackage{graphicx}
\usepackage{epstopdf}
\usepackage{algorithmic}
\ifpdf
  \DeclareGraphicsExtensions{.eps,.pdf,.png,.jpg}
\else
  \DeclareGraphicsExtensions{.eps}
\fi

\newsiamremark{remark}{Remark}
\newsiamremark{hypothesis}{Hypothesis}
\crefname{hypothesis}{Hypothesis}{Hypotheses}
\newsiamthm{claim}{Claim}

\headers{Domain Decomposition Learning Methods}{Q. Sun, X. Xu, and H. Yi}

\title{Domain Decomposition Learning Methods for Solving Elliptic Problems}

\author{Qi Sun\thanks{School of Mathematical Sciences, Tongji University, Shanghai 200092, China
(\email{qsun\_irl@tongji.edu.cn}, \email{2111166@tongji.edu.cn}).}
\and Xuejun Xu\footnotemark[1]$\ ^,$\thanks{Institute of Computational Mathematics, AMSS, Chinese Academy of Sciences, Beijing 100190, China (\email{xxj@lsec.cc.ac.cn}).}
\and Haotian Yi\footnotemark[1]}

\usepackage{amsopn}

\usepackage{adjustbox} 
\usepackage{tikz} 
\usetikzlibrary{arrows,shapes,positioning,automata}
\usetikzlibrary{calc,decorations.markings}
\usetikzlibrary{decorations.pathreplacing}
\usetikzlibrary{arrows.meta, positioning, quotes}
\usetikzlibrary{matrix}
\usepackage{makecell} 
\usepackage{booktabs} 
\usepackage{multirow}
\usepackage{subcaption}
\usepackage{diagbox}
\usepackage[most]{tcolorbox}
\usepackage[symbol]{footmisc}
\usepackage{footnote}
\usepackage{refcount}
\usepackage{rotating}
\usepackage{bm}
\usetikzlibrary{snakes}

\newlength{\commentindent}
\makeatletter
\renewcommand{\algorithmiccomment}[1]{\unskip\hfill\makebox[\commentindent][l]{//~#1}\par}
\LetLtxMacro{\oldalgorithmic}{\algorithmic}
\renewcommand{\algorithmic}[1][0]{%
  \oldalgorithmic[#1]%
  \renewcommand{\ALC@com}[1]{%
    \ifnum\pdfstrcmp{##1}{default}=0\else\algorithmiccomment{##1}\fi}%
}
\makeatother

\begin{document}

\maketitle

\begin{abstract}
With recent advancements in computer hardware and software platforms, there has been a surge of interest in solving partial differential equations with deep learning-based methods, and the integration with domain decomposition strategies has attracted considerable attention owing to its enhanced representation and parallelization capacities of the network solution. While there are already several works that substitute the subproblem solver with neural networks for overlapping Schwarz methods, the non-overlapping counterpart has not been extensively explored because of the inaccurate flux estimation at interface that would propagate errors to neighbouring subdomains and eventually hinder the convergence of outer iterations. In this study, a novel learning approach for solving elliptic boundary value problems, i.e., the compensated deep Ritz method using neural network extension operators, is proposed to enable reliable flux transmission across subdomain interfaces, thereby allowing us to construct effective learning algorithms for realizing non-overlapping domain decomposition methods (DDMs) in the presence of erroneous interface conditions. Numerical experiments on a variety of elliptic problems, including regular and irregular interfaces, low and high dimensions, two and four subdomains, and smooth and high-contrast coefficients are carried out to validate the effectiveness of our proposed algorithms.
\end{abstract}

\begin{keywords}
elliptic boundary value problem, domain decomposition methods, artificial neural networks, erroneous Dirichlet-to-Neumann map, compensated deep Ritz method.
\end{keywords}

\begin{MSCcodes}
65M55, 92B20, 49S05.
\end{MSCcodes}

\section{Introduction}
Many problems of interest in science and engineering are modeled by partial differential equations, which help us to understand and control complex systems across a broad range of real-world applications \cite{evans2010partial,lions1971optimal}. Unfortunately, finding the analytic solution for many problems is often difficult or even impossible, therefore various numerical techniques such as finite difference, finite volume, and finite element methods \cite{leveque2007finite,leveque2002finite,brenner2008mathematical} are developed to obtain the approximate solution. Based on a discretization of the solution space by dividing the computational domain into a polygon mesh, these mesh-based numerical methods are highly accurate and efficient for low-dimensional problems on regular domains. However, there are still many challenging issues to be addressed, e.g., mesh generation remains complex when the boundary is geometrically complicated or dynamically changing, computation of high-dimensional problems is often infeasible due to the curse of dimensionality, and others \cite{karniadakis2021physics}. As classical methods are continuously being improved, it also raises the need for new methods and tools to tackle the difficulties mentioned above.

With recent advancements in computer hardware and software platforms, deep learning-based approaches \cite{lecun2015deep} have emerged as an attractive alternative for solving different types of partial differential equations in both forward and inverse problems \cite{sirignano2018dgm,han2017deep,brunton2019machine}. Thanks to its universal approximation capabilities \cite{scarselli1998universal}, the use of neural networks as an ansatz to the solution function or operator mapping has achieved remarkable success in diverse disciplines. One noteworthy work is the physics-informed neural networks (PINNs) \cite{raissi2019physics,lagaris2000neural,lagaris1998artificial} that incorporate the residual of underlying equations into the training loss function, where integer-order differential operators can be directly calculated through automatic differentiation \cite{paszke2017automatic}. Another important work is the deep Ritz method \cite{yu2018deep}, which resorts to the Ritz formulation and performs better than PINNs for problems with low-regularity solutions \cite{chen2020comparison}. It is also possible to design learning tasks according to the Galerkin formulation \cite{zang2020weak}, but the training process often struggles to converge due to the imbalance between generator and discriminator models. In addition, to improve the boundary condition satisfaction, several techniques including, but not limited to the deep Nitsche method \cite{liao2019deep}, augmented Lagrangian relaxation \cite{huang2021augmented}, and auxiliary network with distance functions \cite{mcfall2009artificial,berg2018unified} have been developed. Compared to traditional numerical methods \cite{leveque2007finite,brenner2008mathematical}, deep learning-based approaches offer advantages of flexible and meshless implementation, strong ability to handle non-linearity and to break the curse of dimensionality \cite{karniadakis2021physics}. However, they may exhibit poor performance when handling problems with multi-scale phenomena \cite{wight2020solving,jagtap2020extended}, and the large training cost is also a major drawback that limits their use in large-scale scientific computing. To address these challenges as well as enhance the representation and parallelization capacity of network solutions, integrating deep learning with domain decomposition strategies \cite{heinlein2021combining,heinlein2021combining2,heinlein2023predicting,heinlein2020machine} has attracted increasing attention in recent years.

One way is to incorporate the distributed training techniques \cite{ben2019demystifying}, e.g., the data and module parallelization, into the original PINNs approach \cite{jagtap2020extended,jagtap2020conservative,hu2021extended,hu2022augmented}, where the learning task is split into multiple training sections through a non-overlapping partition of the domain and various continuity conditions are enforced on subdomain interfaces. Although this combination is quite general and parallelizable, it differs from the conventional way of splitting a partial differential equation \cite{toselli2004domain,quarteroni1999domain}. Besides, the averaging of solution, flux, and residual on the interface \cite{hu2021extended,hu2022augmented} may be problematic for solutions with jump conditions\footnote{For instance, the solution of elliptic interface problem with high-contrast coefficients \cite{li2006immersed} lies in the Sobolev space $H^{1+\epsilon}(\Omega)$ with $\epsilon>0$ possibly close to zero \cite{mercier2003minimal}, thereby enforcing the residual continuity condition on the interface \cite{jagtap2020extended,hu2022augmented} cannot be directly applied due to the lack of regularity.}. On the other hand, conventional DDMs \cite{toselli2004domain} can be formulated at the continuous or variational level, which also allows deep learning-based methods to be employed for solving the decomposed subproblem. As a result, the machine learning analogue of overlapping Schwarz methods have emerged recently and successfully handled many elliptic problems \cite{liao2019deep,li2019d3m,mercier2021coarse,sheng2022pfnn}, however, the non-overlapping counterpart has not been systematically studied yet. A major challenge is that the local network solution is prone to returning erroneous flux prediction along the subdomain interface, which would propagate errors to neighbouring subproblems and eventually hamper the convergence of outer iterations. In other words, the low accuracy of flux estimation is a key threat to the integration of deep learning and non-overlapping DDMs, especially for those based on a direct flux exchange across subdomain interfaces, but has not been fully addressed or resolved in the existing literature.

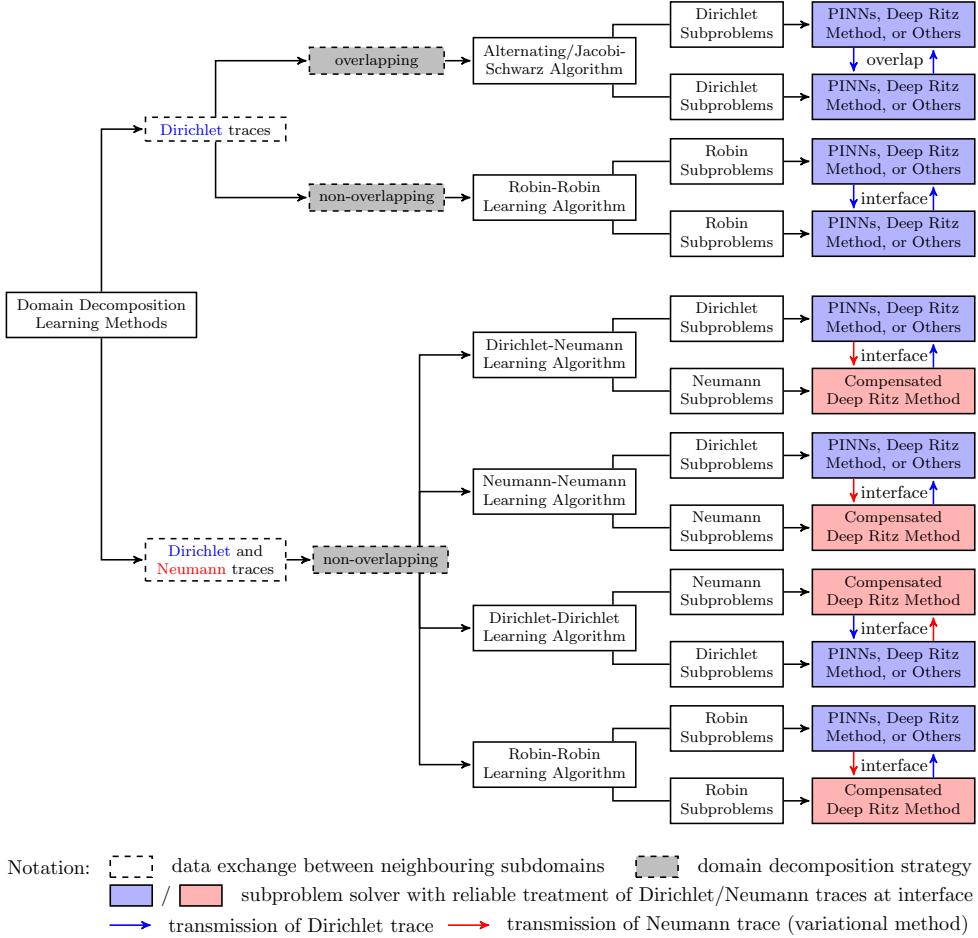
\begin{figure}[t!]
\begin{adjustbox}{max totalsize={0.99\textwidth}{0.99\textheight},center}
\begin{tikzpicture}[shorten >=1pt,auto,node distance=2.5cm,thick,main node/.style={rectangle,draw,font=\footnotesize},decoration={brace,mirror,amplitude=7},dimension node/.style={draw,dashed,fill=lightgray,font=\footnotesize},decoration={brace,mirror,amplitude=7}]

\node[main node] (1) {\parbox{3.2cm}{\centering Domain Decomposition \\ Learning Methods}};

\node[draw=none,fill=none] (2) [right = 0.2cm of 1] {};
\node[dimension node, fill=white] (3) [above  = 3.cm of 2] {\parbox{2.3cm}{\centering \textcolor{blue}{Dirichlet} traces}};
\node[dimension node, fill=white] (4) [below  = 3.9cm of 2] {\parbox{2.3cm}{\centering \textcolor{blue}{Dirichlet} and \textcolor{red}{Neumann} traces}};

\node[draw=none,fill=none] (5) [right = 1.5cm of 3] {};
\node[dimension node] (6) [above = 0.85cm of 5] {\parbox{2.2cm}{\centering overlapping }};
\node[dimension node] (7) [below = 0.85cm of 5] {\parbox{2.2cm}{\centering non-overlapping}};

\node[main node] (8) [right = 0.5cm of 6] {\parbox{2.7cm}{\centering Alternating/Jacobi- \\ Schwarz Algorithm}};
\node[main node] (9) [right = 0.5cm of 7] {\parbox{2.7cm}{\centering Robin-Robin \\ Learning Algorithm}};

\node[draw=none,fill=none] (10) [right = 1.5cm of 8] {};
\node[main node] (11) [above = 0.1cm of 10] {\parbox{1.8cm}{\centering Dirichlet \\ Subproblems}};
\node[main node] (12) [below = 0.1cm of 10] {\parbox{1.8cm}{\centering Dirichlet \\ Subproblems}};
\node[main node,fill = white!70!blue] (13) [right = 0.5cm of 11] {\parbox{2.7cm}{\centering PINNs, Deep Ritz \\ Method, or Others}};
\node[main node,fill = white!70!blue] (14) [right = 0.5cm of 12] {\parbox{2.7cm}{\centering PINNs, Deep Ritz \\ Method, or Others}};

\node[draw=none,fill=none] (15) [right = 1.5cm of 9] {};
\node[main node] (16) [above = 0.1cm of 15] {\parbox{1.8cm}{\centering Robin \\ Subproblems}};
\node[main node] (17) [below = 0.1cm of 15] {\parbox{1.8cm}{\centering Robin \\ Subproblems}};

\node[main node,fill = white!70!blue] (18) [right = 0.5cm of 16] {\parbox{2.7cm}{\centering PINNs, Deep Ritz \\ Method, or Others}};
\node[main node,fill = white!70!blue] (19) [right = 0.5cm of 17] {\parbox{2.7cm}{\centering PINNs, Deep Ritz \\ Method, or Others}};

\path[every node/.style={font=\sffamily\small,sloped},->,>=stealth',color=blue]
	 (13) edge [transform canvas={xshift=-2.05em}] node[draw=none, rotate=90, right] {\ \textcolor{black}{\textnormal{overlap}}} (14)
	 (14) edge [transform canvas={xshift=2.05em}] (13)
	 (18) edge [transform canvas={xshift=-2.05em}] node[draw=none, rotate=90, right] {\textcolor{black}{\textnormal{interface}}} (19)
	 (19) edge [transform canvas={xshift=2.05em}] (18);	 	
	  
\path[every node/.style={font=\sffamily\small,sloped},->,>=stealth']	 
	 (11) edge (13)
	 (12) edge (14)
	 (16.east) edge (18.west)
	 (17.east) edge (19.west)
	 (6) edge (8)
	 (7) edge (9);
\draw[->,>=stealth'] (3.north) |- (6.west);
\draw[->,>=stealth'] (3.south) |- (7.west);
\draw[->,>=stealth'] (1.north) |- (3.west);
\draw ([xshift=3.em]8.north) |- ([xshift=.1em]11.west);
\draw ([xshift=3.em]8.south) |- ([xshift=.1em]12.west);
\draw ([xshift=3.em]9.north) |- ([xshift=.1em]16.west);
\draw ([xshift=3.em]9.south) |- ([xshift=.1em]17.west);

\node[dimension node] (20) [right = 0.47cm of 4] {\parbox{2.2cm}{\centering non-overlapping}};

\node[draw=none,fill=none] (21) [right = 1.78 cm of 20] {};
\node[main node] (22) [above = 3.15cm of 21] {\parbox{2.7cm}{\centering Dirichlet-Neumann \\ Learning Algorithm}};
\node[main node] (23) [above = 0.68cm of 21] {\parbox{2.7cm}{\centering Neumann-Neumann \\ Learning Algorithm}};
\node[main node] (24) [below = 0.68cm of 21] {\parbox{2.7cm}{\centering Dirichlet-Dirichlet \\ Learning Algorithm}};
\node[main node] (25) [below = 3.15cm of 21] {\parbox{2.7cm}{\centering Robin-Robin \\ Learning Algorithm}};

\node[draw=none,fill=none] (26) [right = 1.5cm of 22] {};
\node[main node] (27) [above = 0.1cm of 26] {\parbox{1.8cm}{\centering Dirichlet \\ Subproblems}};
\node[main node] (28) [below = 0.1cm of 26] {\parbox{1.8cm}{\centering Neumann \\ Subproblems}};	
\node[main node,fill = white!70!blue] (29) [right = 0.5cm of 27] {\parbox{2.7cm}{\centering PINNs, Deep Ritz \\ Method, or Others}};
\node[main node,fill = white!70!red] (30) [right = 0.5cm of 28] {\parbox{2.7cm}{\centering Compensated \\ Deep Ritz Method}};

\path[every node/.style={font=\sffamily\small,sloped},->,>=stealth']
	 (29) edge [transform canvas={xshift=-2.05em},color=red] node[draw=none, rotate=90, right] {\textcolor{black}{\textnormal{interface}}} (30)
	 (30) edge [transform canvas={xshift=2.05em},color=blue] (29);	  	 	 	
\path[every node/.style={font=\sffamily\small,sloped},->,>=stealth']	 
	 (27.east) edge (29.west)	
	 (28.east) edge (30.west);
\draw ([xshift=3.em]22.north) |- ([xshift=.1em]27.west);
\draw ([xshift=3.em]22.south) |- ([xshift=.1em]28.west);	 

\node[draw=none,fill=none] (31) [right = 1.5cm of 23] {};
\node[main node] (32) [above = 0.1cm of 31] {\parbox{1.8cm}{\centering Dirichlet \\ Subproblems}};
\node[main node] (33) [below = 0.1cm of 31] {\parbox{1.8cm}{\centering Neumann \\ Subproblems}};	
\node[main node,fill = white!70!blue] (34) [right = 0.5cm of 32] {\parbox{2.7cm}{\centering PINNs, Deep Ritz \\ Method, or Others}};
\node[main node,fill = white!70!red] (35) [right = 0.5cm of 33] {\parbox{2.7cm}{\centering Compensated \\ Deep Ritz Method}};

\path[every node/.style={font=\sffamily\small,sloped},->,>=stealth']
	 (34) edge [transform canvas={xshift=-2.05em},color=red] node[draw=none, rotate=90, right] {\textcolor{black}{\textnormal{interface}}} (35)
	 (35) edge [transform canvas={xshift=2.05em},color=blue] (34);	  	 	 	
\path[every node/.style={font=\sffamily\small,sloped},->,>=stealth']	 
	 (32.east) edge (34.west)	
	 (33.east) edge (35.west);
\draw ([xshift=3.em]23.north) |- ([xshift=.1em]32.west);
\draw ([xshift=3.em]23.south) |- ([xshift=.1em]33.west);	 

\node[draw=none,fill=none] (36) [right = 1.5cm of 24] {};
\node[main node] (37) [above = 0.1cm of 36] {\parbox{1.8cm}{\centering Neumann \\ Subproblems}};
\node[main node] (38) [below = 0.1cm of 36] {\parbox{1.8cm}{\centering Dirichlet \\ Subproblems}};	
\node[main node,fill = white!70!red] (39) [right = 0.5cm of 37] {\parbox{2.7cm}{\centering Compensated \\ Deep Ritz Method}};
\node[main node,fill = white!70!blue] (40) [right = 0.5cm of 38] {\parbox{2.7cm}{\centering PINNs, Deep Ritz \\ Method, or Others}};

\path[every node/.style={font=\sffamily\small,sloped},->,>=stealth']
	 (40) edge [transform canvas={xshift=2.05em},color=red] (39)
	 (39) edge [transform canvas={xshift=-2.05em},color=blue] node[draw=none, rotate=90, right] {\textcolor{black}{\textnormal{interface}}} (40);	  	 	 	
\path[every node/.style={font=\sffamily\small,sloped},->,>=stealth']	 
	 (37.east) edge (39.west)	
	 (38.east) edge (40.west);
\draw ([xshift=3.em]24.north) |- ([xshift=.1em]37.west);
\draw ([xshift=3.em]24.south) |- ([xshift=.1em]38.west);	 

\node[draw=none,fill=none] (41) [right = 1.5cm of 25] {};
\node[main node] (42) [above = 0.1cm of 41] {\parbox{1.8cm}{\centering Robin \\ Subproblems}};
\node[main node] (43) [below = 0.1cm of 41] {\parbox{1.8cm}{\centering Robin \\ Subproblems}};	
\node[main node,fill = white!70!blue] (44) [right = 0.5cm of 42] {\parbox{2.7cm}{\centering PINNs, Deep Ritz \\ Method, or Others}};
\node[main node,fill = white!70!red] (45) [right = 0.5cm of 43] {\parbox{2.7cm}{\centering Compensated \\ Deep Ritz Method}};

\path[every node/.style={font=\sffamily\small,sloped},->,>=stealth']
	 (44) edge [transform canvas={xshift=-2.05em},color=red] node[draw=none, rotate=90, right] {\textcolor{black}{\textnormal{interface}}} (45)
	 (45) edge [transform canvas={xshift=2.05em},color=blue] (44);	  	 	 	
\path[every node/.style={font=\sffamily\small,sloped},->,>=stealth']	 
	 (42.east) edge (44.west)	
	 (43.east) edge (45.west)
	 (4.east) edge (20.west);
\draw ([xshift=3.em]25.north) |- ([xshift=.1em]42.west);
\draw ([xshift=3.em]25.south) |- ([xshift=.1em]43.west);	 

\draw[->,>=stealth'] ([xshift=2.em]20.north) |- (22.west);
\draw[->,>=stealth'] ([xshift=2.em]20.north) |- (23.west);
\draw[->,>=stealth'] ([xshift=2.em]20.south) |- (24.west);
\draw[->,>=stealth'] ([xshift=2.em]20.south) |- (25.west);
\draw[->,>=stealth'] (1.south) |- (4.west);
\end{tikzpicture}
\end{adjustbox}
\begin{adjustbox}{max totalsize={1\textwidth}{0.5\textheight},left}
\begin{tikzpicture}[shorten >=1pt,auto,node distance=2.5cm,thick,main node/.style={rectangle,draw,font=\footnotesize},decoration={brace,mirror,amplitude=7},dimension node/.style={draw,dashed,fill=lightgray,font=\footnotesize},decoration={brace,mirror,amplitude=7}]

\node[dimension node, fill=white] (1)  {\parbox{0.5cm}{ \textcolor{white}{em} }}; 
\node[draw=none] (2) [right = 0.2cm of 1] {\centering data exchange between neighbouring subdomains };
\node[dimension node] (3) [right = 0.4cm of 2] {\parbox{.5cm}{ \textcolor{lightgray}{em} }};
\node[draw=none] (4) [right = 0.2cm of 3] {\centering domain decomposition strategy };

\node[main node, fill = white!70!blue] (5) [below = 0.11cm of 1] {\parbox{.5cm}{ \textcolor{white!70!blue}{em} }};
\node[draw=none] (6) [right = 0.01cm of 5] {\centering $/$ };
\node[main node, fill = white!70!red] (7) [right = 0.01cm of 6] {\parbox{.5cm}{ \textcolor{white!70!red}{em} }};
\node[draw=none] (8) [right = 0.2cm of 7] {\centering subproblem solver with reliable treatment of Dirichlet/Neumann traces at interface};

\node[draw=none] (9) [below = 0.11cm of 5] {\parbox{0.5cm}{ \textcolor{white}{em} }};
\draw[->,>=stealth',color=blue] (9.west) -- (9.east);	
\node[draw=none] (11) [right = 0.14cm of 9] {\centering transmission of Dirichlet trace};

\node[draw=none] (12) [right = 0.11cm of 11] {\parbox{0.5cm}{ \textcolor{white}{em} }};
\draw[->,>=stealth',color=red] (12.west) -- (12.east);	
\node[draw=none] (13) [right = 0.14cm of 12] {\centering transmission of Neumann trace (variational method) };
 
\node[draw=none] (14) [left = 0.2cm of 1] {\centering Notation: };
\node[draw=none] (15) [above = 0.2cm of 14] {};
\end{tikzpicture}
\end{adjustbox}
\vspace{-0.7cm}
\caption{Our proposed framework of domain decomposition learning methods for solving second-order elliptic boundary value problems.}
\label{fig-big-picture}
\vspace{-0.7cm}
\end{figure}

This study mainly focuses on the benchmark Poisson's equation that serves as a necessary prerequisite to validate the effectiveness of deep learning-based domain decomposition approaches \cite{heinlein2021combining,li2019d3m,li2020deep,mercier2021coarse}, namely,
\begin{equation}
\begin{array}{cl}
- \Delta u(x) = f(x)\ \ & \text{in}\ \Omega,\\
u(x)=0\ \ & \text{on}\ \partial \Omega,
\end{array}
\label{Poisson-StrongForm}
\end{equation}
where $\Omega\subset\mathbb{R}^d$ is a bounded Lipschitz domain, $d\in\mathbb{N}_+$ the dimension, and $f(x)\in L^2(D)$ a given function. The classification of DDMs for solving problem \eqref{Poisson-StrongForm} are typically categorized as either an overlapping or a non-overlapping approach \cite{toselli2004domain,quarteroni1999domain,mathew2008domain}, which can be further refined according to the information exchange between neighbouring subdomains (see \autoref{fig-big-picture}). Here, the refined category is adopted throughout this work to distinguish between various deep learning-based domain decomposition algorithms. Note that the trained solution of Dirichlet subproblem using PINNs \cite{raissi2019physics}, deep Ritz \cite{yu2018deep}, or other similar methods \cite{karniadakis2021physics} is often found to exhibit erroneous Neumann traces on the interface (see Remark \ref{Remark-DtN-Map}). Accordingly, the flux transmission between neighbouring subdomains is thus of low accuracy, which would hinder the convergence of out iterations. To deal with this issue, we propose a novel learning approach, i.e., the compensated deep Ritz method using neural network extension operators, that allows reliable flux transmission between neighbouring subdomains but without explicitly involving the computation of Dirichlet-to-Neumann map on subdomain interfaces. This enables us to construct effective learning approaches for realizing classical Dirichlet-Neumann, Neumann-Neumann, Dirichlet-Dirichlet, and Robin-Robin algorithms in the non-overlapping regime (see \autoref{fig-big-picture}). It is noteworthy that although the Robin-Robin algorithm only requires the exchange of Dirichlet traces \cite{chen2014optimal}, two additional parameters within the interface conditions need to be determined, which may lead to incorrect network solutions. Fortunately, our compensated deep Ritz method can also help alleviate this issue (see section \ref{Section-RRLM}). 

The remainder of this paper is organized as follows. In section 2, we provide a brief review of classical DDMs for solving elliptic boundary value problems, as well as several deep learning approaches that can be employed as our subproblem solver. Next, following the most straightforward idea (see \autoref{fig-big-picture}), we introduce the machine learning analogue of Robin-Robin algorithm in section 3. To realize other non-overlapping DDMs using neural networks, a detailed illustration of our compensated deep Ritz method is presented in section 4. Experimental results on a series of benchmark problems are reported in section 5 to validate the effectiveness of our methods, as well as an interface problem with high-contrast coefficients. Finally, in section 6, we conclude the paper and outline some directions for future work.

\section{Preliminaries}
This section is devoted to provide a concise overview of classical DDMs \cite{toselli2004domain,quarteroni1999domain} for solving the Poisson's equation, together with the widely used deep learning approaches \cite{raissi2019physics,yu2018deep} that can be adopted as our subproblem solver.

\subsection{Domain Decomposition Methods}
The idea of domain decomposition for solving the Poisson's equation has a long history dating back to the 18th century \cite{Schwarz1870alternative}, and there is an extensive literature on DDMs owing to the emergence and improvement of parallel computers (see \cite{toselli2004domain,quarteroni1999domain,mathew2008domain} and references cited therein). For illustrative purposes, let us assume that the computational domain $\Omega\subset\mathbb{R}^d$ is partitioned into two subdomains $\{\Omega_i\}_{i=1}^2$ (see \autoref{fig-domain-decomposition} for example), while the case of multiple subdomains can be treated in a similar fashion. Depending on the partition strategy being employed, DDMs are usually categorized into two groups: overlapping and non-overlapping approaches, in which the decomposed subproblem is typically solved through mesh-based finite difference or finite element methods \cite{toselli2004domain,quarteroni1999domain}.

\begin{figure}[t!]
\begin{adjustbox}{max totalsize={0.95\textwidth}{0.95\textheight},center}
\begin{tikzpicture}

  \draw[color=gray, fill=lightgray, thick] (1.7,0) rectangle (2.3,4);   
  \draw[thick] (0,0) rectangle (4,4);
  \node[draw=none] at (1,2) {$\Omega_1$};
  \node[draw=none] at (3,2) {$\Omega_2$};
  \node [draw=none] at (2,2) {\rotatebox{90}{$\Omega_1\cap\Omega_2$}};
  \node[draw=none] at (1.3,3.5) {$\Gamma_2$};  
  \node[draw=none] at (2.9,3.5) {$\Gamma_1$};    
  
  \draw[color=gray, thick] (7,0) -- (7,4);   
  \draw[thick] (5,0) rectangle (9,4);
  \node[draw=none] at (6,2) {$\Omega_1$};
  \node[draw=none] at (8,2) {$\Omega_2$};
  \node[draw=none] at (7.3,3.5) {$\Gamma$};
   
  \draw[decorate,decoration={zigzag,segment length=8.mm, amplitude=2.mm, pre=lineto,pre length=0pt,post=lineto,post length=0pt},color=gray, thick] (12,0) -- (12,4);      
  \draw[thick] (10,0) rectangle (14,4);   
  \node[draw=none] at (11,2) {$\Omega_1$};
  \node[draw=none] at (13,2) {$\Omega_2$};
  \node[draw=none] at (12.5,3.5) {$\Gamma$};
  
   \draw[color=gray, thick] (17,0) -- (17,4);     
   \draw[color=gray, thick] (15,2) -- (19,2);        
   \draw[thick] (15,0) rectangle (19,4);
   \draw[fill=lightgray] (15,0) rectangle (17,2);  
   \draw[fill=lightgray] (17,2) rectangle (19,4);   
   \node[draw=none] at (16,1) {$\Omega_B$};
   \node[draw=none] at (16,3) {$\Omega_R$};
   \node[draw=none] at (18,1) {$\Omega_R$};
   \node[draw=none] at (18,3) {$\Omega_B$};
   \draw[black,fill=black] (17,2) circle (.3ex);  
 
\end{tikzpicture}
\end{adjustbox}
\vspace{-0.4cm}
\caption{Decomposition of a bounded domain $\Omega\subset\mathbb{R}^2$ into two subdomains. \textbf{Left:} Overlapping partition. \textbf{Middle:} Non-overlapping partition that is separated by curved regular or irregular interfaces. \textbf{Right:} Red-Black partition with the intersection between two subdomains in the same class marked in bold. }
\label{fig-domain-decomposition}
\vspace{-0.7cm}
\end{figure}
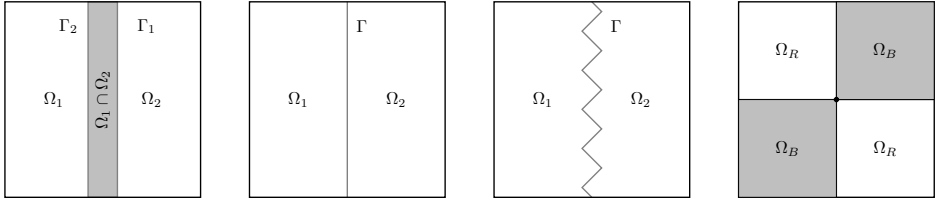

\begin{figure}[htp]
\begin{algorithm}[H]
\caption{DDMs with Exchange of Dirichlet Traces}
\begin{algorithmic}
\STATE{Start with the initial interface condition, $h_1^{[0]}$ and $h_2^{[0]}$, for each subsolution;}
\FOR{$k \gets 0$ to $K$ (maximum number of outer iterations)}
\WHILE{stopping criteria are not satisfied}
\STATE{\% \textit{Subproblem-Solving} }
\STATE{
\vspace{-0.72cm}
\begingroup
\renewcommand*{\arraystretch}{1.3}
\begin{equation*}
\left\{
\begin{array}{cl}
- \Delta u_1^{[k]} \!=\! f\ \ & \text{in}\ \Omega_1\\
u_1^{[k]} \!=\! 0\ \ & \text{on}\ \partial \Omega_1 \!\setminus\! \Gamma_1\\
\mathcal{B}_1 u_1^{[k]} \!=\! h_1^{[k]}\ \ & \text{on}\ \Gamma_1
\end{array} \right.
\text{with\ \ }
\mathcal{B}_1 u_1^{[k]} \!=\! \left\{
\begin{array}{cl}
u_1^{[k]} & \text{(SAM)} \\
\nabla u_1^{[k]}\cdot\bm{n}_1 \!+\! \kappa_1  u_1^{[k]} & \text{(RRA)} \\
\end{array}\right.
\end{equation*}
\endgroup
\vspace{-0.38cm}
}
\STATE{\% \textit{Exchange of Dirichlet Trace}}  
\STATE{
\vspace{-0.72cm}
\begingroup
\renewcommand*{\arraystretch}{1.3}
\begin{equation*}
h_2^{[k]} \!=\! \left\{
\begin{array}{cll}
u_1^{[k]}\ \ & \text{on}\ \Gamma_2 & \text{(SAM)} \\
\displaystyle - h_1^{[k]} + (\kappa_1 + \kappa_2) u_1^{[k]}\ \ & \text{on}\ \Gamma_2 & \text{(RRA)} \\
\end{array}\right.
\end{equation*}
\endgroup
\vspace{-0.38cm}
}
\STATE{\% \textit{Subproblem-Solving}}
\STATE{
\vspace{-0.72cm}
\begingroup
\renewcommand*{\arraystretch}{1.3}
\begin{equation*}
\left\{
\begin{array}{cl}
- \Delta u_2^{[k]} \!=\! f\ \ & \text{in}\ \Omega_2\\
u_2^{[k]} \!=\! 0\ \ & \text{on}\ \partial \Omega_2 \!\setminus\! \Gamma_2\\
\mathcal{B}_2 u_2^{[k]} \!=\! h_2^{[k]}\ \ & \text{on}\ \Gamma_2
\end{array} \right.
\text{with\ \ }
\mathcal{B}_2 u_2^{[k]} \!=\! \left\{
\begin{array}{cl}
u_2^{[k]} & \text{(SAM)} \\
\nabla u_2^{[k]}\cdot\bm{n}_2 \!+\! \kappa_2  u_2^{[k]} & \text{(RRA)} \\
\end{array}\right.
\end{equation*}
\endgroup
\vspace{-0.38cm}
}
\STATE{\% \textit{Exchange of Dirichlet Trace}}
\STATE{
\vspace{-0.72cm}
\begingroup
\renewcommand*{\arraystretch}{1.3}
\begin{equation*}
h_1^{[k+1]} \!=\! \left\{
\begin{array}{cll}
u_2^{[k]}\ \ & \text{on}\ \Gamma_1 & \text{(SAM)} \\
\displaystyle \rho h_1^{[k]} + (1-\rho)   (- h_2^{[k]} + (\kappa_1 + \kappa_2) u_2^{[k]})\ \ & \text{on}\ \Gamma_1 & \text{(RRA)}
\end{array}\right.
\end{equation*}
\endgroup
\vspace{-0.38cm}
}
\ENDWHILE
\ENDFOR
\STATE{(\textbf{Remark:} RRA is defined in the non-overlapping regime, i.e., $\Gamma_1=\Gamma_2$.) }
\end{algorithmic}
\label{DDM-Solution-Exchange}
\end{algorithm}
\vspace{-0.7cm}
\end{figure}

\begin{figure}[htp]
\begin{algorithm}[H]
\caption{DDMs with Exchange of Dirichlet and Neumann Traces}
\begin{algorithmic}
\STATE{Start with the initial interface condition, $h_1^{[0]}$ and $h_2^{[0]}$, for each subsolution;}
\FOR{$k \gets 0$ to $K$ (maximum number of outer iterations)}
\WHILE{stopping criteria are not satisfied}
\STATE{\% \textit{Subproblem-Solving} }
\STATE{
\vspace{-0.72cm}
\begingroup
\renewcommand*{\arraystretch}{1.3}
\begin{equation*}
\left\{\!\!
\begin{array}{cl}
- \Delta u_i^{[k]} \!=\! f & \text{in}\ \Omega_i\\
u_i^{[k]} \!=\! 0 & \text{on}\ \partial \Omega_i \!\setminus\! \Gamma\\
\mathcal{B}_i u_i^{[k]} \!=\! h_i^{[k]} & \text{on}\ \Gamma
\end{array}\right.
\text{with\ \ }
\mathcal{B}_i u_i^{[k]} \!=\! \left\{\!\!
\begin{array}{cll}
u_i^{[k]} & \!\!\text{for}\ i=1, & \!\!\text{(DNA)} \\
\displaystyle u_i^{[k]} & \!\!\text{for}\ i=1,2, & \!\!\text{(NNA)} \\
\displaystyle \nabla u_i^{[k]} \cdot \bm{n}_i & \!\!\text{for}\ i=1,2, & \!\!\text{(DDA)} 
\end{array}\right.
\end{equation*}
\endgroup
\vspace{-0.38cm}
}
\STATE{\% \textit{Exchange of Dirichlet or Neumann Trace}}
\STATE{
\vspace{-0.72cm}
\begingroup
\renewcommand*{\arraystretch}{1.3}
\begin{equation*}
h_i^{[k]} \!=\! \left\{\!\!
\begin{array}{clll}
-\nabla u_{1}^{[k]} \cdot \bm{n}_{1} & \text{on}\ \Gamma & \text{for}\ i=2, & \text{(DNA)} \\ 
\displaystyle \nabla u_{1}^{[k]} \cdot \bm{n}_{1} + \nabla u_{2}^{[k]} \cdot \bm{n}_{2} & \text{on}\ \Gamma & \text{for}\ i=1,2, & \text{(NNA)} \\ 
\displaystyle u_{1}^{[k]} - u_{2}^{[k]} & \text{on}\ \Gamma & \text{for}\ i=1,2, & \text{(DDA)}
\end{array}\right.
\end{equation*}
\endgroup
\vspace{-0.38cm}
}
\STATE{\% \textit{Subproblem-Solving}}
\STATE{
\vspace{-0.72cm}
\begingroup
\renewcommand*{\arraystretch}{1.3}
\begin{equation*}
\left\{\!\!
\begin{array}{cl}
- \Delta u_i^{[k]} \!=\! f & \text{in}\ \Omega_i\\
u_i^{[k]} \!=\! 0 & \text{on}\ \partial \Omega_i \!\setminus\! \Gamma\\
\mathcal{B}_i u_i^{[k]} \!=\! h_i^{[k]} & \text{on}\ \Gamma
\end{array} \right.
\text{with\ \ }
\mathcal{B}_i u_i^{[k]} \!=\! \left\{\!\!
\begin{array}{cll}
\nabla u_i^{[k]} \cdot \bm{n}_i & \!\!\text{for}\ i=2, & \!\!\text{(DNA)} \\
\displaystyle \nabla u_i^{[k]} \cdot \bm{n}_i & \!\!\text{for}\ i=1,2, & \!\!\text{(NNA)} \\
\displaystyle u_i^{[k]} & \!\!\text{for}\ i=1,2, & \!\!\text{(DDA)} \\
\end{array}\right.
\end{equation*}
\endgroup
\vspace{-0.38cm}
}
\STATE{\% \textit{Exchange of Dirichlet or Neumann Trace}}
\STATE{
\vspace{-0.72cm}
\begingroup
\renewcommand*{\arraystretch}{1.3}
\begin{equation*}
h_i^{[k+1]} = \left\{\!\!
\begin{array}{clll}
\rho u_{2}^{[k]} + (1-\rho) u_{1}^{[k]} & \text{on}\ \Gamma & \text{for}\ i=1, \ & \text{(DNA)} \\
h_i^{[k]} - \rho ( u_{1}^{[k]} + u_{2}^{[k]} ) & \text{on}\ \Gamma & \text{for}\ i=1,2, \ & \text{(NNA)} \\
h_i^{[k]} - \rho ( \nabla u_{1}^{[k]} \cdot \bm{n}_1 + \nabla u_{2}^{[k]} \cdot \bm{n}_2 )  & \text{on}\ \Gamma & \text{for}\ i=1,2,\ & \text{(DDA)} \\ 
\end{array}\right.
\end{equation*}
\endgroup
\vspace{-0.38cm}
}
\ENDWHILE
\ENDFOR
\end{algorithmic}
\label{DDM-Flux-Exchange}
\end{algorithm}
\vspace{-0.7cm}
\end{figure}

When using a neural network as the solution ansatz to boundary value problem \eqref{Poisson-StrongForm}, it has been observed that the trained model often agrees with the Dirichlet boundary condition but exhibits erroneous Neumann traces \cite{dockhorn2019discussion,bajaj2021robust}, which sets it apart from traditional mesh-based numerical methods \cite{toselli2004domain,quarteroni1999domain}. In this regard, the classification of DDMs adopted in this paper is based on the information exchange between neighbouring subdomains rather than the partition strategy as depicted in \autoref{fig-big-picture}. To be specific, we summarize some representative decomposition-based approaches in the literature \cite{toselli2004domain}, referring to the Schwarz alternating method and the Robin-Robin algorithm as SAM and RRA, respectively, in Algorithm \ref{DDM-Solution-Exchange}. On the other hand, the Dirichlet-Neumann, Neumann-Neumann, and Dirichlet-Dirichlet algorithms are abbreviated as DNA, NNA, and DDA, respectively, in Algorithm \ref{DDM-Flux-Exchange}\footnote{Abbreviations SAM, RRA, DNA, NNA, and DDA are only used in Algorithm \ref{DDM-Solution-Exchange} and \ref{DDM-Flux-Exchange}.}. In addition, the relaxation parameter $\rho$ should lie between $(0,\rho_{max})$ in order to achieve convergence \cite{na2022domain,funaro1988iterative}. Notably, overlapping methods with small overlap are cheap and easy to implement, it usually comes at the price of slower convergence than non-overlapping ones. Besides, non-overlapping DDMs are more nature and efficient in handling elliptic problems with large jumps in the coefficient \cite{xu1998some}.

\subsection{Deep Learning Solvers}\label{Section-DL-Solvers}
As can be concluded from the previous discussion, the decomposed subproblem on each subdomain takes on the form
\begin{equation}
\begin{array}{cl}
- \Delta u_i(x) = f(x)\ \ & \text{in}\ \Omega_i,\\
u_i(x)=0\ \ & \text{on}\ \partial \Omega_i\setminus\Gamma,\\
\mathcal{B}_i u_i(x)=h_i(x)\ \ & \text{on}\ \Gamma,
\end{array}
\label{Subproblem-StrongForm}
\end{equation}
where $\mathcal{B}_i$ is a boundary operator on the interface that may represent the Dirichlet, Neumann, or Robin boundary condition, namely,
\begin{equation*}
\begin{array}{rl}
\textnormal{Dirichlet boundary condition:}\ &\ \mathcal{B}_i u_i(x) = u_i(x), \\
\textnormal{Neumann boundary condition:}\ &\ \mathcal{B}_i u_i(x) = \nabla u_i(x) \cdot \bm{n}_i, \\
\textnormal{Robin boundary condition:}\ &\ \mathcal{B}_i u_i(x) = \nabla u_i(x)\cdot \bm{n}_i + \kappa_i u_i(x),
\end{array}
\end{equation*}
while the function $h_i(x)$ is iteratively determined along the outer iteration \cite{toselli2004domain}. 

When deep learning-based approaches are utilized to solve \eqref{Subproblem-StrongForm}, the hypothesis space of local solution is first built using a neural network. If not otherwise stated, we shall use the fully-connected neural network of depth $L\in\mathbb{N}_+$ \cite{hornik1989multilayer}, in which the $\ell$-th hidden layer receives an input $x^{\ell-1}\in\mathbb{R}^{n_{\ell-1}}$ from its previous layer and transforms it to $T^{\ell}(x^{\ell-1}) = \bm{W}^{\ell} x^{\ell-1} + \bm{b}^{\ell}$. Here, $\bm{W}^{\ell} \in \mathbb{R}^{ n_{\ell} \times n_{\ell - 1} }$ and $\bm{b}^{\ell} \in \mathbb{R}^{n_{\ell}}$ are the weights and biases to be learned, and $\theta=\{ \bm{W}^{\ell}, \bm{b}^{\ell} \}_{\ell = 1}^L$ denotes the collection of all trainable parameters. By choosing an activation function $\sigma(\cdot)$ for each hidden layer, the solution ansatz can then be expressed as $\hat{u}_i(x;\theta) = \big(T^L \circ \sigma \circ T^{L-1} \cdots \circ \sigma \circ T^1\big)(x)$, where $\circ$ represents the composition operator. One can also employ other network architectures, e.g., residual neural network and its variants \cite{he2016deep,han2017deep}, for the parametrization of unknown solutions.

To update trainable parameters using the celebrated backpropagation algorithm \cite{hecht1992theory}, various training loss functions (before applying numerical integration) have been proposed, e.g., PINNs \cite{raissi2019physics} that are based on the strong formulation of \eqref{Subproblem-StrongForm}, i.e.,
\begin{equation*}
	        \hat{u}_i(x;\theta) = \operatorname*{arg\,min}_{\theta} \int_{\Omega_i} |\Delta \hat{u}_i + f|^2 dx + \beta \Big( \int_{\partial\Omega_i\backslash \Gamma}|\hat{u}_i|^2ds+\int_{\Gamma}|\mathcal{B}_i \hat{u}_i-h_i|^2ds\Big),
\end{equation*}
where $\beta>0$ is a user-defined penalty coefficient. Alternatively, the deep Ritz method \cite{yu2018deep} resorts to the Ritz formulation of \eqref{Subproblem-StrongForm}, namely,
\begin{equation*}
	\hat{u}_i(x;\theta) = \operatorname*{arg\,min}_{\theta} \int_{\Omega_i} \Big( \frac12 |\nabla \hat{u}_i|^2 - f\hat{u}_i \Big) dx + \beta \int_{\partial\Omega_i\backslash \Gamma}|\hat{u}_i|^2ds + L_\Gamma(\hat{u}_i)
\end{equation*}
where the last term depends on the interface condition being imposed
\begingroup
\renewcommand*{\arraystretch}{2.1}
\begin{equation*}
L_\Gamma(\hat{u}_i) = \left\{ 
\begin{array}{ll}
\displaystyle \beta \int_{\Gamma}|\hat{u}_i(x;\theta) - h_i(x)|^2ds\ \ &\ \text{(Dirichlet condition),}  \\
\displaystyle  - \int_{\Gamma} h_i(x) \hat{u}_i(x;\theta)\,ds\ \ &\ \text{(Neumann condition),} \\
\displaystyle \int_\Gamma  \left( \frac{\kappa_i}{2} |\hat{u}_i(x;\theta)|^2 - h_i(x) \hat{u}_i(x;\theta) \right) ds\ \ &\ \text{(Robin condition).}
\end{array}\right.
\end{equation*}
\endgroup
In addition to these two widely-used techniques, the weak adversarial network \cite{zang2020weak} is based on the Galerkin formulation of \eqref{Subproblem-StrongForm}, while another series of learning tasks is designed to use separate networks to fit the interior and boundary equations respectively \cite{mcfall2009artificial,berg2018unified}. We refer the readers to \cite{karniadakis2021physics,sheng2022pfnn,heinlein2021combining} for a more detailed review of deep learning-based numerical methods. Notably, with the interface conditions being included as penalty terms in the training loss function and the number of interface points being small compare to that of interior ones, the trained model of \eqref{Subproblem-StrongForm} is often prone to returning erroneous Dirichlet-to-Neumann map on the interface \cite{dockhorn2019discussion,bajaj2021robust} (see Remark \ref{Remark-DtN-Map}). This emerges as a key threat to the integration of deep learning and flux exchange-based DDMs but has not been fully addressed in the literature. 

\begin{remark}\label{Remark-DtN-Map}
To validate our statements, we first study the Dirichlet-to-Neumann map, i.e., $\mathcal{B}_1u_1=u_1$ in \eqref{Subproblem-StrongForm}, that sends boundary value data to normal derivative data through the trained network solution. Here, the PINNs approach \cite{raissi2019physics} is adopted for network training, where $\Omega_1 = [0,0.5]\times [0,1]$, $\Gamma=\{0.5\}\times [0,1]$, $f(x,y)$ and $h_1(x,y)$ are derived from the exact solution $u_1(x,y) = \sin(2 \pi x)  (\cos(2 \pi y) - 1)$. As can be seen from \autoref{table-DtN-map-DL-section}, the network solution using a fully-connected neural network ($\textnormal{depth}=8$, $\textnormal{width}=50$, also known as the multilayer perceptron) agrees with our true solution, and the performance can be further improved through the use of a more sophisticated network architecture \cite{wang2021understanding}. However, the corresponding Neumann traces are of unsatisfactory low accuracy, which is often unacceptable for flux transmission between neighbouring subdomains. Fortunately, the prediction of $\nabla u_1$ using fully-connected neural networks performs well inside the subdomain $\Omega_1$.

\begin{table}[ht]
\vspace{-0.3cm}
\caption{Trained network solutions $\hat{u}_1$ of Dirichlet subproblem \eqref{Subproblem-StrongForm} using different architectures, together with their error profiles $| u_1-\hat{u}_1|$ and $|\partial_x (u_1-\hat{u}_1)|$.}
\centering
\begin{tabular}{ c || c | c || c | c}
\toprule
& $\hat{u}_1$ & $\vert u_1-\hat{u}_1\vert$ & $\partial_x \hat{u}_1$ & $|\partial_x (u_1-\hat{u}_1)|$\\ 
\midrule
\!\!\makecell{multilayer \\ perceptron \\ ($\beta \equiv 400$) \\ \cite{raissi2019physics} }\!\! &
\begin{minipage}{.17\textwidth}
\centering
\includegraphics[width=\linewidth]{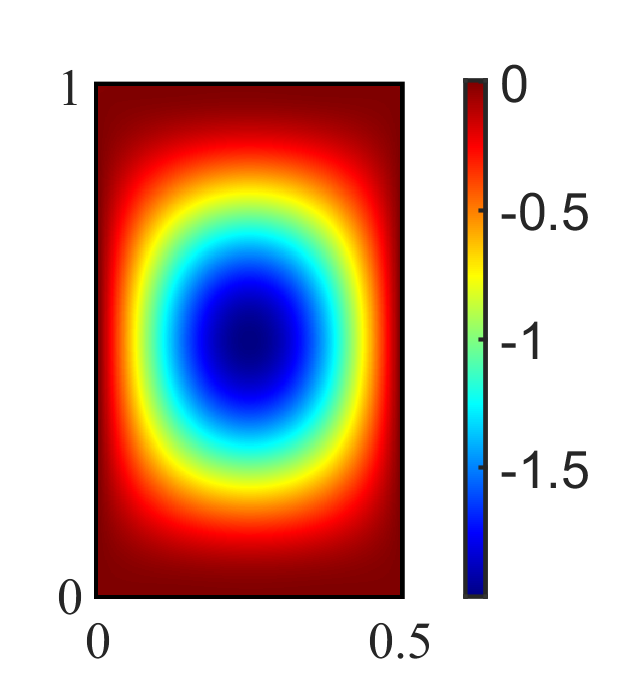}
\end{minipage}
& 
\begin{minipage}{.18\textwidth}
\centering
\includegraphics[width=0.98\linewidth]{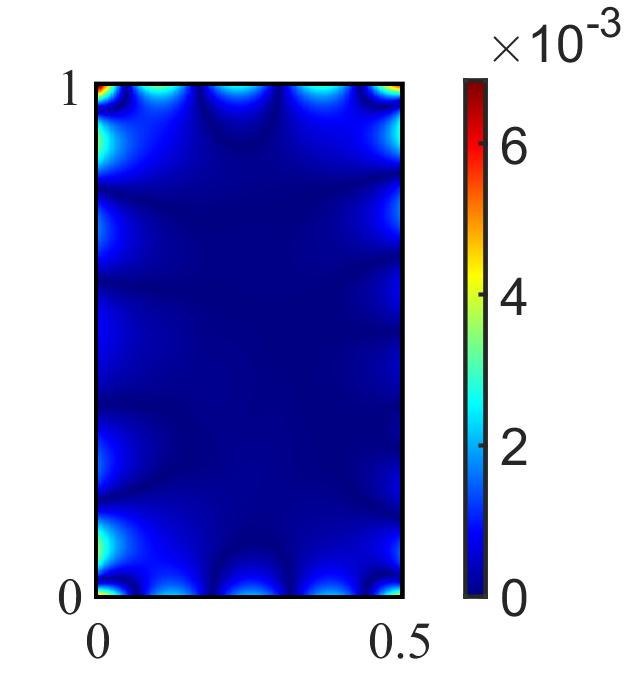}\vspace*{0.16cm}
\end{minipage}
& 
\begin{minipage}{.17\textwidth}
\centering
\includegraphics[width=0.96\linewidth]{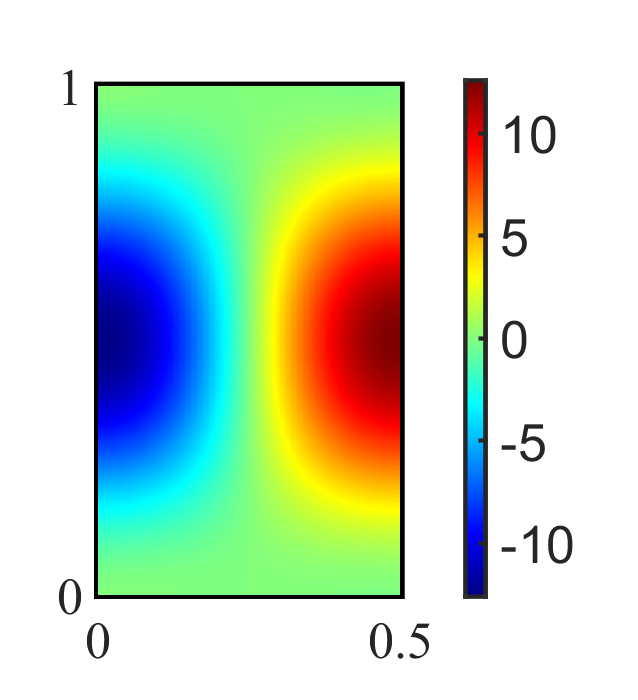}
\end{minipage}
& 
\begin{minipage}{.17\textwidth}
\centering
\includegraphics[width=\linewidth]{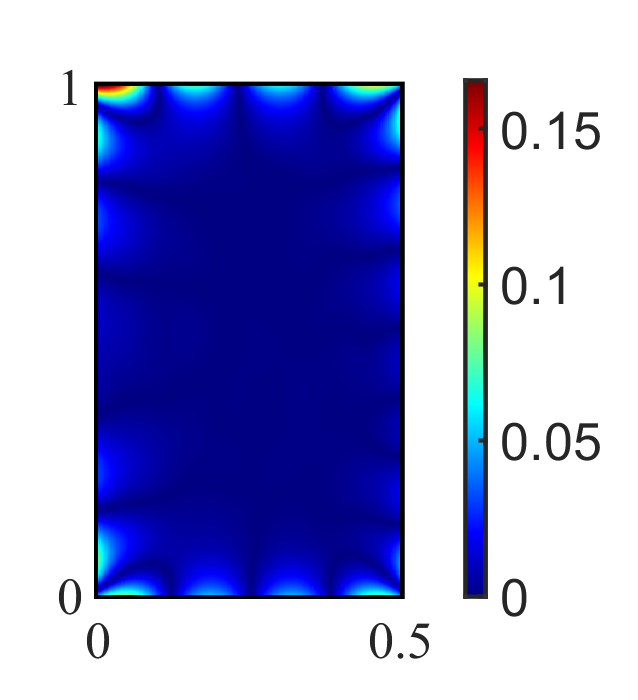}
\end{minipage}
\\ 
\midrule
\!\!\makecell{transformer \\ network \\ (adaptive $\beta$) \\ \cite{wang2021understanding}}\!\! &
\begin{minipage}{.17\textwidth}
\centering
\includegraphics[width=\linewidth]{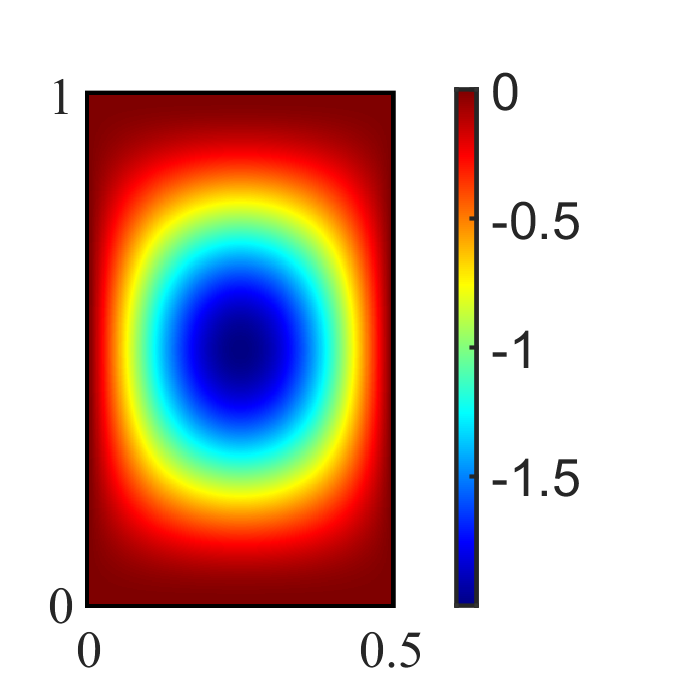}
\end{minipage}
& 
\begin{minipage}{.18\textwidth}
\centering
\includegraphics[width=0.98\linewidth]{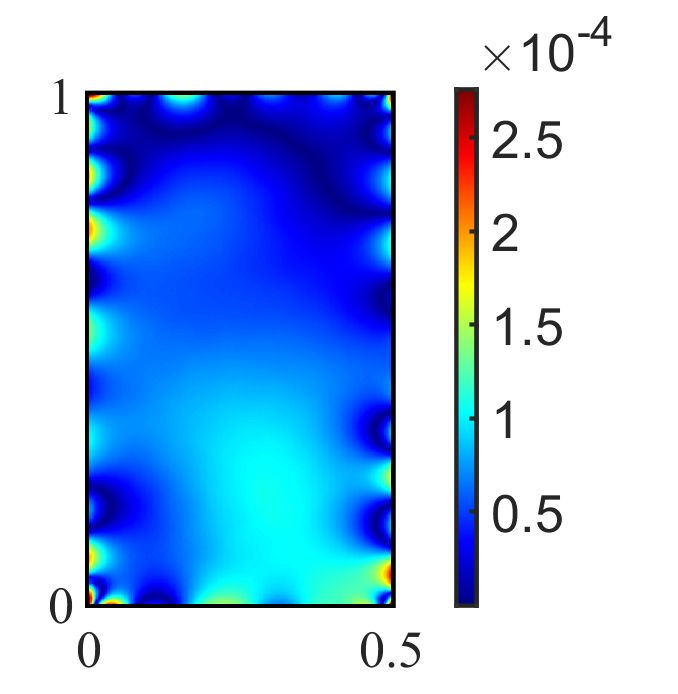}\vspace*{0.19cm}
\end{minipage}
& 
\begin{minipage}{.17\textwidth}
\centering
\includegraphics[width=0.9\linewidth]{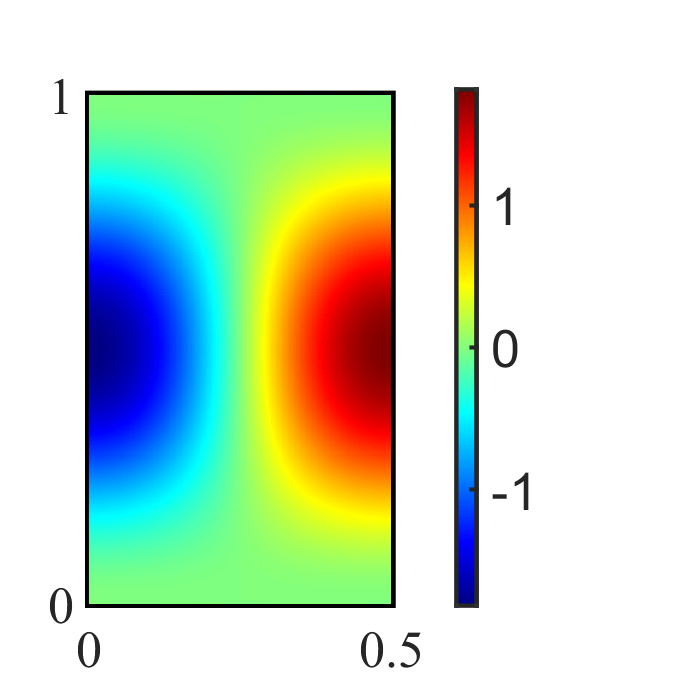}
\end{minipage}
& 
\begin{minipage}{.17\textwidth}
\centering
\includegraphics[width=0.9\linewidth]{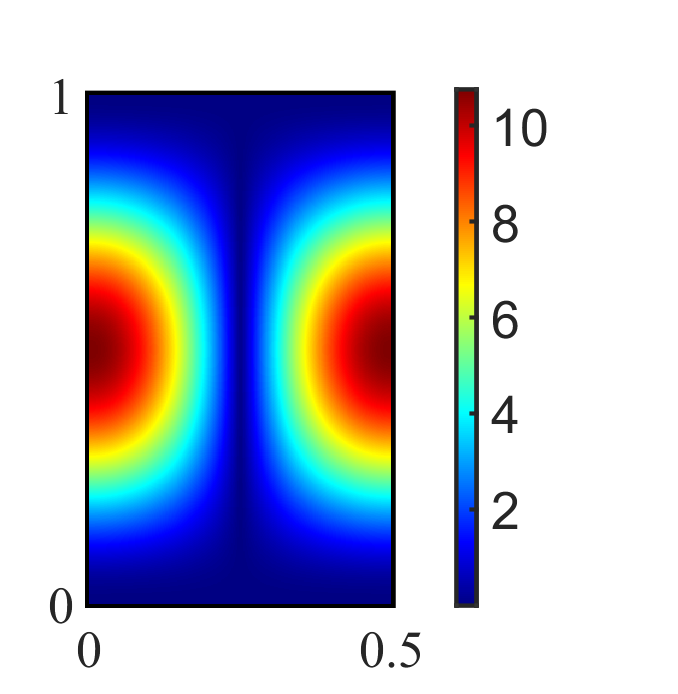}
\end{minipage}
\\ 
\bottomrule
\end{tabular}
\label{table-DtN-map-DL-section}
\vspace{-0.8cm}
\end{table}
\end{remark}

\begin{remark}\label{Remark-Robin}
Next, the Robin subproblem \eqref{Subproblem-StrongForm}, i.e., $\mathcal{B}_1u_1= \nabla u_1\cdot\bm{n}_1 + \kappa_1 u_1$, with the same exact solution $u_1(x,y) = \sin(2 \pi x)  (\cos(2 \pi y) - 1)$ is studied numerically. By employing the standard PINNs approach using fully-connected neural networks ($\textnormal{depth}=8$, $\textnormal{width}=50$) and $\beta=400$, the numerical results with different values of coefficient $\kappa_1$ are reported in \autoref{table-RobinProb-DL-section}. Clearly, the trained model fails to recover the true solution in the case of $\kappa_1=10^4$, which is due to the weight imbalance between $\hat{u}_1$ and $\nabla \hat{u}_1\cdot \bm{n}_1$ in the boundary penalty loss term 
$$\displaystyle L_\Gamma( \hat{u}_i, h_1 ) = \beta \int_{\Gamma} \Big| \kappa_1 \hat{u}_1 + \frac{\partial \hat{u}_1}{\partial \bm{n}_1} - h_1 \Big|^2ds.$$ 
In addition, this issue is inherent within the Robin boundary condition and cannot be fixed by fine-tuning the penalty coefficient $\beta>0$.

\begin{table}[htp]
\caption{Trained network solutions $\hat{u}_1$ of Robin subproblem \eqref{Subproblem-StrongForm} using different architectures, together with their error profiles $| u_1-\hat{u}_1|$ and $|\partial_x (u_1-\hat{u}_1)|$.}
\centering
\begin{tabular}{ c || c | c || c | c}
\toprule
& $\hat{u}_1$ & $\vert u_1-\hat{u}_1\vert$ & $\partial_x \hat{u}_1$ & $|\partial_x (u_1-\hat{u}_1)|$ \\
\midrule
\!\!\makecell{multilayer \\ perceptron \\ ($\kappa_1=1$) }\!\! &
\begin{minipage}{.17\textwidth}
\centering
\includegraphics[width=\linewidth]{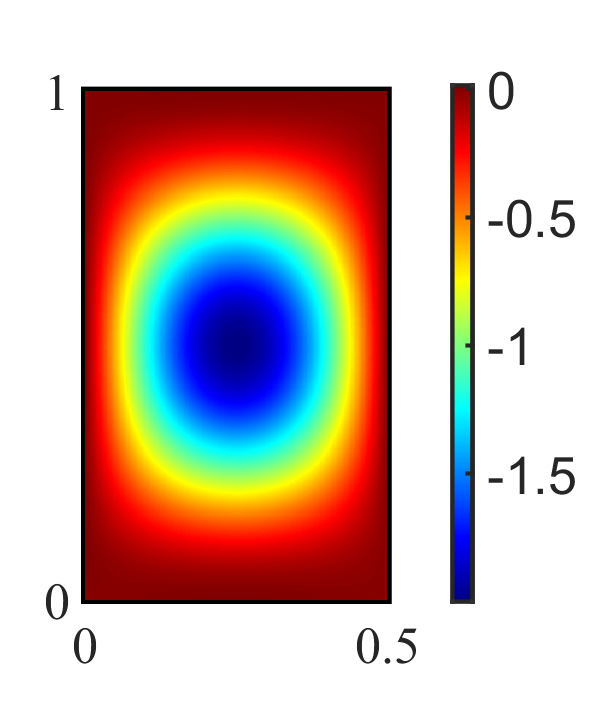}
\end{minipage}
& 
\begin{minipage}{.18\textwidth}
\centering
\includegraphics[width=\linewidth]{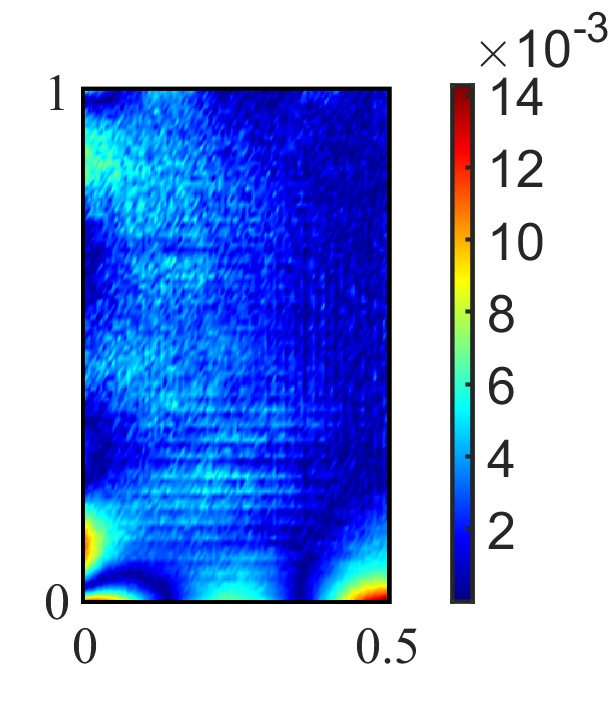}\vspace*{0.165cm}
\end{minipage}
& 
\begin{minipage}{.17\textwidth}
\centering
\includegraphics[width=.98\linewidth]{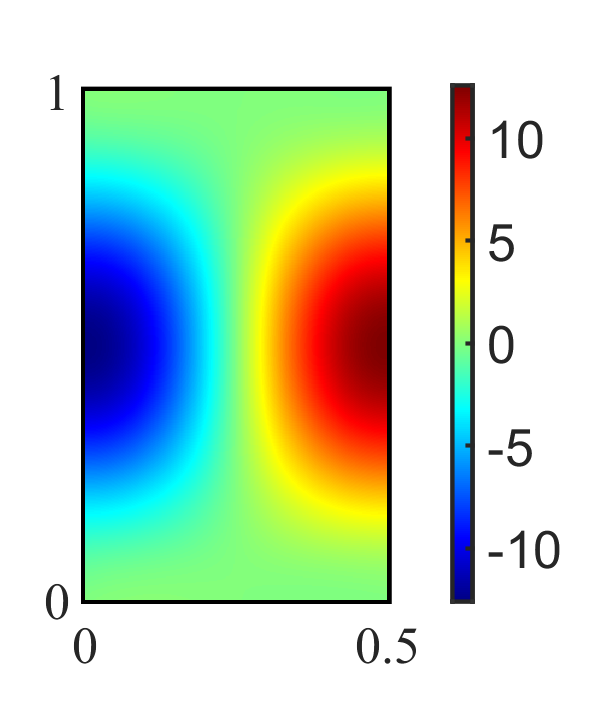}
\end{minipage}
& 
\begin{minipage}{.17\textwidth}
\centering
\includegraphics[width=\linewidth]{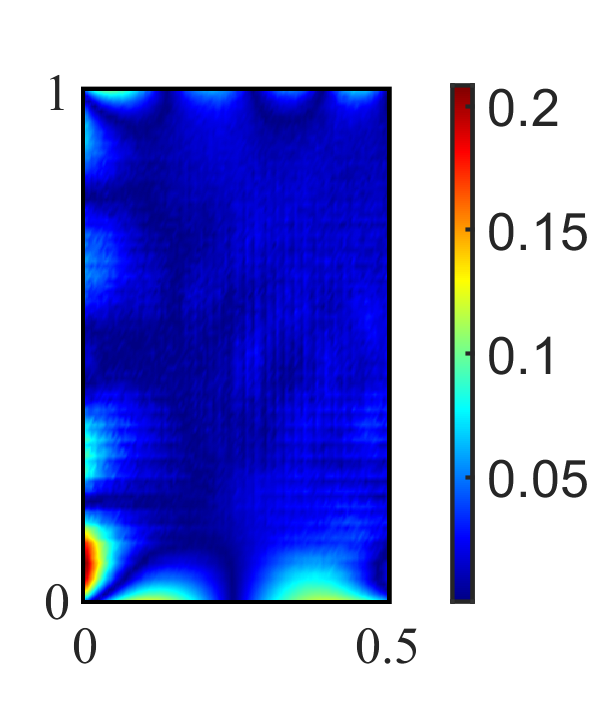}
\end{minipage}
\\ 
\midrule
\!\!\makecell{multilayer \\ perceptron \\ ($\kappa_1=10^4$) }\!\! &
\begin{minipage}{.17\textwidth}
\centering
\includegraphics[width=.9\linewidth]{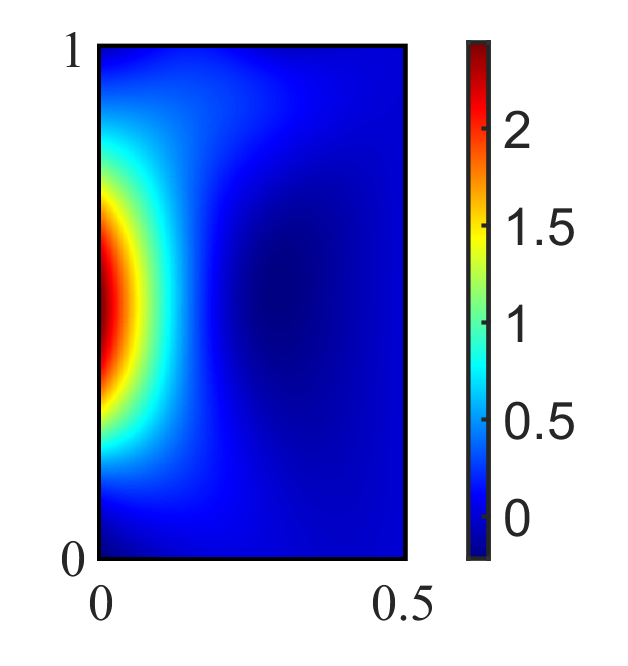}
\end{minipage}
& 
\begin{minipage}{.17\textwidth}
\centering
\includegraphics[width=.9\linewidth]{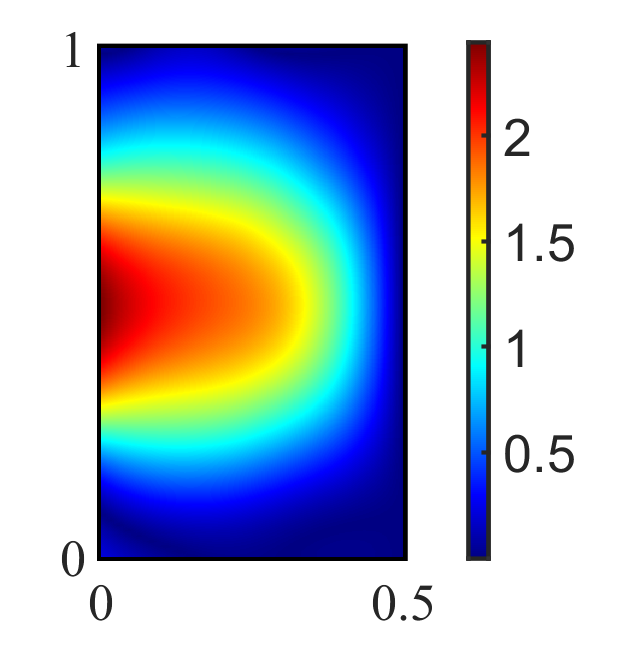}
\end{minipage}
& 
\begin{minipage}{.17\textwidth}
\centering
\includegraphics[width=.9\linewidth]{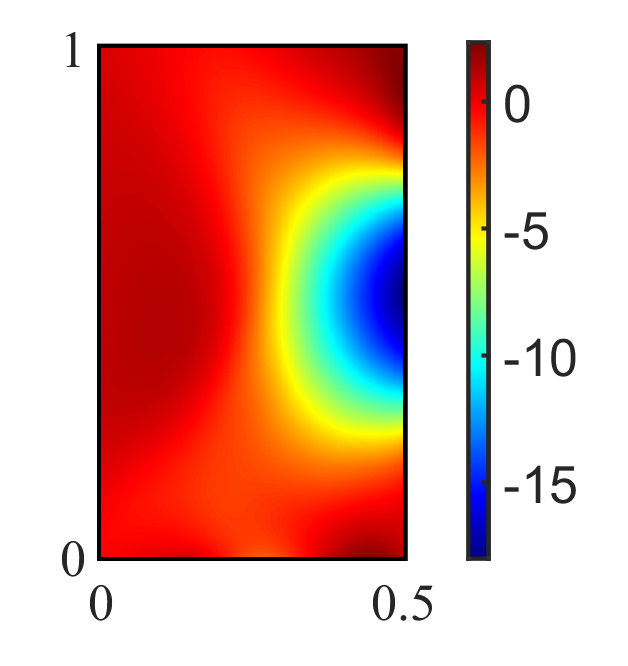}
\end{minipage}
& 
\begin{minipage}{.17\textwidth}
\centering
\includegraphics[width=.85\linewidth]{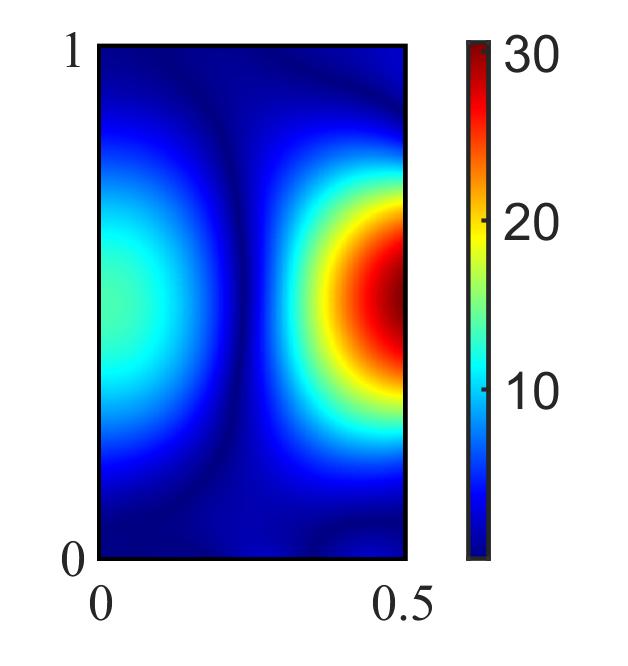}
\end{minipage}
\\
\bottomrule
\end{tabular}
\label{table-RobinProb-DL-section}
\vspace{-0.3cm}
\end{table}

\end{remark}


\section{Robin-Robin Algorithm using Physics-Informed Neural Networks}\label{Section-RR-PINNs}
In addition to the deep learning analogue of overlapping Schwarz methods \cite{li2019d3m,mercier2021coarse,li2020deep,sheng2022pfnn,taghibakhshi2022learning}, the non-overlapping Robin-Robin algorithm \cite{toselli2004domain,quarteroni1999domain} is also based on the exchange of Dirichlet traces between neighbouring subproblems (see \autoref{fig-big-picture} or Algorithm \autoref{DDM-Solution-Exchange}). As the decomposition leads to simpler functions to be learned on each subdomain, the PINNs approach \cite{raissi2019physics}, rather than the deep Ritz method, is employed here as the subproblem solver since it is known to empirically work better for problems with smooth solutions \cite{chen2020comparison}. However, a major drawback is the determination of two additional parameters, i.e., $\kappa_1$ and $\kappa_2$ within the Robin boundary conditions, which may require more outer iterations to converge or cause difficulties for the optimization process (see Remark \ref{Remark-Robin}).

For ease of illustration, we consider the case of two non-overlapping subdomains in what follows (see \autoref{fig-domain-decomposition} for example), where the interface conditions are invariable of the Robin type \cite{toselli2004domain,quarteroni1999domain,chen2014optimal}. The detailed iterative process (in terms of differential operators) is presented in Algorithm \autoref{DDM-Solution-Exchange}, from which it can be observed that the update of interface conditions only involves the Dirichlet traces.

\begin{figure}[htp]
\begin{algorithm}[H]
\caption{Robin-Robin Algorithm using PINNs (2 Subdomains)}
\begin{algorithmic}
\STATE{\% \textit{Initialization} }
\STATE{-- divide domain $\Omega\subset\mathbb{R}^d$ into two non-overlapping subdomains $\Omega_1$ and $\Omega_2$;}
\STATE{-- specify network structures $\hat{u}_1(x;\theta_1)$ and $\hat{u}_2(x;\theta_2)$ for each subproblem;}
\STATE{-- generate Monte Carlo training samples $X_\Gamma$, $X_{\Omega_i}$, and $X_{D_i}$ for $i=1$, 2; }
\STATE{\% \textit{Outer Iteration Loop} }
\STATE{Start with the initial guess $h^{[0]}$ along the interface $\Gamma$;}
\FOR{$k \gets 0$ to $K$ (maximum number of outer iterations)}
\WHILE{stopping criteria are not satisfied}
\STATE{\% \textit{Subproblem-Solving using PINNs} }
\STATE{
\vspace{-0.65cm}
\begingroup
\renewcommand*{\arraystretch}{1.1}
\begin{equation*}
\hat{u}_1^{[k]} = \operatorname*{arg\,min}_{\theta_1} L_{\Omega_1}( \hat{u}_1 ) + \beta \big( L_{D_1}( \hat{u}_1 ) + L_\Gamma( \hat{u}_1, h_1^{[k]} ) \big)
\end{equation*}
\endgroup
\vspace{-0.45cm}
}
\STATE{\% \textit{Exchange of Dirichlet Trace} }
\STATE{
\vspace{-0.65cm}
\begingroup
\renewcommand*{\arraystretch}{1.1}
\begin{equation*}
h_2^{[k]}(x_n^\Gamma) = - h_1^{[k]}(x_n^\Gamma) + (\kappa_1 + \kappa_2) \hat{u}_1^{[k]}(x_n^\Gamma), \ \ n=1,\cdots,N_\Gamma,
\end{equation*}
\endgroup
\vspace{-0.6cm}
}
\STATE{\% \textit{Subproblem-Solving using PINNs} }
\STATE{
\vspace{-0.65cm}
\begingroup
\renewcommand*{\arraystretch}{1.1}
\begin{equation*}
\hat{u}_2^{[k]} = \operatorname*{arg\,min}_{\theta_2} L_{\Omega_2}( \hat{u}_2 ) + \beta \big( L_{D_2}( \hat{u}_2 ) + L_\Gamma( \hat{u}_2, h_2^{[k]} ) \big)
\end{equation*}
\endgroup
\vspace{-0.45cm}
}
\STATE{\% \textit{Exchange of Dirichlet Trace} }
\STATE{
\vspace{-0.65cm}
\begingroup
\renewcommand*{\arraystretch}{1.1}
\begin{equation*}
h_1^{[k+1]}(x_n^\Gamma) \!=\! \rho  h_1^{[k]}(x_n^\Gamma) + (1-\rho)\big( (\kappa_1 + \kappa_2) u_2^{[k]}(x_n^\Gamma) - h_2^{[k]}(x_n^\Gamma)  \big), \ \ n=1,\cdots,N_\Gamma,
\end{equation*}
\endgroup
\vspace{-0.6cm}
}
\ENDWHILE
\ENDFOR
\end{algorithmic}
\label{Algorithm-RR-Learning-2Subdomains}
\end{algorithm}
\vspace{-0.7cm}
\end{figure}

To realize the Robin-Robin algorithm using PINNs, the decomposed subproblem is first rewritten as an optimization problem through the residual of equations, i.e.,
\begin{equation}
	u_i^{[k]} = \operatorname*{arg\,min}_{u_i} \int_{\Omega_i} | \Delta u_i + f |^2\,dx + \beta\left( \int_{\Gamma} \Big| \kappa_i u_i + \frac{\partial u_i}{\partial \bm{n}_i} - h_i^{[k]} \Big|^2\,ds + \int_{\partial\Omega\cap\partial\Omega_i} |u_i|^2\,ds \right),
	\label{RRLM-PINN-Subprob-Functional}
\end{equation}
for $i=1$, 2,  where the boundary and interface conditions are included as penalty terms during training. Then, by introducing the neural network parametrization\footnote{For notational simplicity, $\hat{u}_i(x,\theta_i)$ and $\hat{u}_i(x,\theta_i^{[k]})$ are sometimes abbreviated as $\hat{u}_i$ and $\hat{u}_i^{[k]}$.}
\begin{equation*}
	u_1^{[k]}(x) \approx \hat{u}_1^{[k]}(x) := \hat{u}_1(x;\theta_1^{[k]}) \ \ \ \text{and}\ \ \ u_2^{[k]}(x) \approx \hat{u}_2^{[k]}(x) := \hat{u}_2(x;\theta_2^{[k]}),
\end{equation*}
and generating the training sample points inside each subdomain and at its boundary
\begin{equation*}
	X_{\Omega_i} = \big\{ x_n^{\Omega_i} \big\}_{n=1}^{N_{\Omega_i}},\ \ \ X_{D_i} = \big\{ x_n^{D_i} \big\}_{n=1}^{N_{D_i}},\ \ \ \text{and}\ \ \ X_{\Gamma} = \big\{ x_n^{\Gamma} \big\}_{n=1}^{N_{\Gamma}},
\end{equation*}
for $i=1$ and 2, the stochastic tools \cite{bottou2018optimization} can be applied for fulfilling the corresponding optimization problems. Specifically, the learning tasks at the $k$-th outer iteration are
\begin{equation}
	\hat{u}_i^{[k]} = \operatorname*{arg\,min}_{\theta_i} L_{\Omega_i}( \hat{u}_i ) + \beta \big( L_{D_i}( \hat{u}_i ) + L_{\Gamma}( \hat{u}_i, h_i^{[k]} ) \big),\ \ \ i=1, 2,
	\label{RRLM-Subprob-Discrete}
\end{equation}
where the loss functions (not relabelled) are defined as
\begingroup
\vspace{-0.25cm}
\begin{equation*}
\begin{array}{c}
\displaystyle L_{\Omega_i}( \hat{u}_i) = \frac{1}{N_{\Omega_i}} \sum_{n=1}^{N_{\Omega_i}} \big| \Delta \hat{u}_i(x_n^{\Omega_i};\theta_i) + f(x_n^{\Omega_i}) \big|^2, \ \ \ L_{D_i}( \hat{u}_i ) = \frac{1}{N_{D_i}} \sum_{n=1}^{N_{D_i}} | \hat{u}_i(x_n^{D_i};\theta_i) |^2, \\
\displaystyle L_{\Gamma}( \hat{u}_i, h_i^{[k]}) = \frac{1}{N_\Gamma} \sum_{n=1}^{N_\Gamma} \big| \kappa_i \hat{u}_i(x_n^{\Gamma_i};\theta_i) + \frac{\partial \hat{u}_i}{\partial \bm{n}_i}(x_n^{\Gamma_i};\theta_i) - h_i^{[k]}(x_n^{\Gamma_i}) \big|^2.
\end{array}
\end{equation*}
\endgroup
Here, and in what follows, the sampling points are drawn uniformly at random from their corresponding domains. One can also use adaptive or adversarial sampling strategies \cite{gao2022failure,he2022mesh} to reduce the training cost.

To sum up, by employing PINNs as subproblem solvers, the deep learning analogue of Robin-Robin algorithm is presented in Algorithm \ref{Algorithm-RR-Learning-2Subdomains}, where $\kappa_1$, $\kappa_2>0$ are two additional user-defined parameters. We can assume, without loss of generality, that $\kappa_1=1$ and leave the other parameter to be tuned. In fact, as the number of interface points is typically much smaller than that of the interior of the subdomains, too large (or small) value of $\kappa_2$ may cause weight imbalance in the interface penalty term (see Remark \ref{Remark-Robin}), while a moderate value of $\kappa_2$ can guarantee convergence but at the cost of extra outer iterations. Such an imbalance issue greatly differs from the conventional finite element setting \cite{chen2014optimal}, and is further demonstrated through numerical experiments in \autoref{Section-Experiments}. Fortunately, this problem can be tackled through the use of our compensated deep Ritz method, which is theoretically and numerically studied in the following sections.

\section{Compensated Deep Ritz Method}\label{Section-Compensated-DeepRitz}
This section begins by studying non-overlapping DDMs that rely on a direct flux exchange across subdomain interfaces, then the Robin-Robin algorithm is revisited from a variational viewpoint. 

Note that when the Dirichlet subproblem is solved using the PINNs or deep Ritz approach \cite{raissi2019physics,yu2018deep}, it is common for the trained model to converge to a local minimizer that nearly satisfies the given Dirichlet boundary condition but with inaccurate Neumann traces. As a result, the flux transmission between neighbouring subdomains may be of low accuracy that would hinder the convergence of outer iterations. To address this issue, we propose in this section the compensated deep Ritz method that enables reliable flux transmission in the presence of erroneous interface conditions. Moreover, our proposed learning algorithm can also help with the network training when realizing the Robin-Robin algorithm with large coefficients.

\subsection{Dirichlet-Neumann Learning Algorithm}
Here, we focus on the classical Dirichlet-Neumann algorithm \cite{toselli2004domain,quarteroni1999domain}, where the detailed iterative process (in terms of differential operators) is presented in Algorithm \ref{DDM-Flux-Exchange}. To avoid the explicit computation and transmission of Dirichlet-to-Neumann maps at interface, the variational formulation of multidomain problem is taken into consideration. More precisely, the Galerkin formulation of problem \eqref{Poisson-StrongForm} reads: find $u\in H_0^1(\Omega)$ such that
\begin{equation}
	a(u,v) = (f,v)\ \ \ \text{for any}\ v\in H_0^1(\Omega),
	\label{Poisson-WeakForm}
\end{equation}
where the bilinear forms are defined as
\begin{equation*}
	a(u,v) = \int_\Omega \nabla u\cdot \nabla v\,dx \ \ \ \text{and}\ \ \ (f,v) = \int_\Omega fv\,dx.
\end{equation*}

Here, we consider a two subdomain decomposition of \eqref{Poisson-WeakForm}, while similar results can be obtained for multidomain cases \cite{toselli2004domain}. Let $\Omega_1$ and $\Omega_2$ denote a non-overlapping decomposition of the computational domain $\Omega\subset\mathbb{R}^d$, with the interface $\Gamma=\partial \Omega_1\cap \partial \Omega_2$ separating our subdomains as shown in \autoref{fig-domain-decomposition}. Moreover, we set $V_i = \big\{ v_i\in H^1(\Omega_i)\, \big|\, v_i|_{\partial\Omega\cap\partial\Omega_i} = 0 \big\}$, $V_i^0=H_0^1(\Omega_i)$, and define the bilinear terms
\begin{equation*}
	a_i(u_i, v_i) = \int_{\Omega_i} \nabla u_i \cdot \nabla v_i\,dx \ \ \ \text{and}\ \ \ (f, v_i)_i = \int_{\Omega_i} f v_i\, dx
\end{equation*}
for $i=1$, 2. Then, the Green's formula implies that \eqref{Poisson-WeakForm} can be reformulated as: find $u_1\in V_1$ and $u_2\in V_2$ such that
\begingroup
\renewcommand*{\arraystretch}{1.1}
\begin{equation*}
\begin{array}{cl}
a_1(u_1, v_1) = (f, v_1)_1\ \ & \text{for any}\ v_1\in V_1^0,\\
u_1 = u_2\ \ & \text{on}\ \Gamma, \\
a_2(u_2, v_2) = (f, v_2)_2 + (f, R_1\gamma_0 v_2)_1 - a_1(u_1, R_1\gamma_0 v_2)\ \ & \text{for any}\ v_2\in V_2,
\end{array}
\end{equation*}
\endgroup
where $\gamma_0 v= v|_\Gamma$ indicates the restriction of $v\in H^1(\Omega_i)$ on the interface $\Gamma$, and $R_i: H_{00}^{\frac12}(\Gamma)\to V_i$ any differentiable extension operator \cite{quarteroni1999domain,toselli2004domain}. Based on the minimum total potential energy principle \cite{evans2010partial}, we obtain its equivalent Ritz formulation, i.e.,
\begingroup
\renewcommand*{\arraystretch}{1.3}
\begin{equation*}
\begin{array}{c}
\displaystyle \operatorname*{arg\,min}_{u_1\in V_1,\ u_1|_\Gamma = u_2}  \frac12  a_1(u_1, u_1) - (f, u_1)_1, \\
\displaystyle \operatorname*{arg\,min}_{u_2\in V_2} \frac12 a_2(u_2, u_2) - (f, u_2)_2 - (f, R_1\gamma_0 u_2)_1 + a_1(u_1, R_1\gamma_0 u_2).
\end{array}
\end{equation*}
\endgroup
Therefore, the Dirichlet-Neumann algorithm \cite{toselli2004domain,quarteroni1999domain} can be written in terms of their energy functionals: given the initial guess $h^{[0]}\in H_{00}^{\frac12}(\Gamma)$, then solve for $k\geq 0$, 
\begin{itemize}
\vspace{0.4em}
\setlength\itemsep{0.4em}
\item[1)] $\displaystyle u_1^{[k]} = \operatorname*{arg\,min}_{u_1\in V_1,\ u_1|_\Gamma = h^{[k]}}  \frac12  a_1(u_1, u_1) - (f, u_1)_1,$
\item[2)] $\displaystyle u_2^{[k]} = \operatorname*{arg\,min}_{u_2\in V_2} \frac12 a_2(u_2, u_2) - (f, u_2)_2 - (f, R_1\gamma_0 u_2)_1 + a_1(u_1^{[k]}, R_1\gamma_0 u_2),$
\item[3)] $\displaystyle h^{[k+1]} = \rho u_2^{[k]} + (1-\rho) h^{[k]}\ \ \ \text{on}\ \Gamma,$
\vspace{0.4em}
\end{itemize}
with $\rho\in(0,\rho_{\textnormal{max}})$ being the acceleration parameter \cite{funaro1988iterative}. Notably, the flux continuity across the subdomain interface is now guaranteed without explicitly calculating and exchanging the Neumann trace of our Dirichlet subproblem.

Such a variational method also makes it possible to integrate with deep learning approaches. Next, the unknown solutions are parametrized by neural networks
\begin{equation}
	u_1^{[k]}(x) \approx \hat{u}_1^{[k]}(x) := \hat{u}_1(x;\theta_1^{[k]}) \ \ \ \text{and}\ \ \ u_2^{[k]}(x) \approx \hat{u}_2^{[k]}(x) := \hat{u}_2(x;\theta_2^{[k]})
	\label{Compensated-DeepRitz-Network-Parametrization}
\end{equation}
where $\hat{u}_i(x;\theta_i^{[k]})$ denotes the solution ansatz with trainable parameters $\theta_i^{[k]}$ for $i=1$ and 2. We note that in contrast to the standard finite element method \cite{toselli2004domain} where the approximate solution of Neumann subproblem is locally defined and the extension operation is mesh-dependent, the neural network parametrization $\hat{u}_2^{[k]}$ in \eqref{Compensated-DeepRitz-Network-Parametrization} is meshless and thus can extend itself to neighbouring subdomains. Therefore, we obtain a natural extension operator, i.e., the neural network extension operator,
\begin{equation}
	R_1\gamma_0 \hat{u}_2(x;\theta_2) = \hat{u}_2(x;\theta_2) \in V_1
	\label{DNLA-Extension}
\end{equation}
which extends the restriction of $\hat{u}_2(x,\theta_2)$ on interface $\Gamma$ to the subdomain $\Omega_1$ with zero boundary value on $\partial\Omega_1\cap\partial\Omega$. Here, the requirement of homogeneous boundary condition on $\partial\Omega_1\cap\partial\Omega$ is dealt by introducing an additional penalty term into the loss function of our extended Neumann subproblem \eqref{Compensated-DeepRitz-Neumann-Subprob-Functional}. In addition, as the extension function is required to be weakly differentiable and the solution of Neumann subproblem is typically regular enough in its subdomain, the hyperbolic tangent or sigmoid activation function is preferred rather than the ReLU activation function. 

Accordingly, by introducing penalty terms for enforcing essential boundary conditions, the Dirichlet subproblem on $\Omega_1$ can be formulated as
\begin{equation}
	\theta_1^{[k]} = \operatorname*{arg\,min}_{\theta_1} \int_{\Omega_1} \Big( \frac12 | \nabla \hat{u}_1 |^2 - f \hat{u}_1\Big) dx + \beta \left( \int_{\partial\Omega_1\cap \partial\Omega} |\hat{u}_1|^2\,ds + \int_{\Gamma} |\hat{u}_1 - h^{[k]}|^2\,ds \right),
	\label{Compensated-DeepRitz-Dirichlet-Subprob-Functional}
\end{equation}
where $\beta>0$ is the penalty coefficient. In fact, as the decomposition usually leads to simpler functions to be learned on each subdomain, the second-order derivatives can thus be involved during the training. As such, the residual form 
\begin{equation}
	\theta_1^{[k]} = \operatorname*{arg\,min}_{\theta_1} \int_{\Omega_1}  | \Delta \hat{u}_1 + f |^2 dx + \beta \left( \int_{\partial\Omega_1\cap \partial\Omega} |\hat{u}_1|^2\,ds + \int_{\Gamma} |\hat{u}_1 - h^{[k]}|^2\,ds \right)
	\label{Compensated-PINN-Dirichlet-Subprob-Functional}
\end{equation}
is then preferred to the Ritz energy \eqref{Compensated-DeepRitz-Dirichlet-Subprob-Functional} since PINNs \eqref{Compensated-PINN-Dirichlet-Subprob-Functional} are empirically found to be capable of offering more accurate estimation of $\nabla u_1$ inside $\Omega_1$. On the other hand, the learning task associated with our Neumann subproblem gives
\begin{equation}
	\theta_2^{[k]} = \operatorname*{arg\,min}_{\theta_2} \int_{\Omega_2} \Big( \frac12 | \nabla \hat{u}_2 |^2  - f \hat{u}_2 \Big) dx + \int_{\Omega_1} \Big( \nabla \hat{u}_1^{[k]} \cdot \nabla \hat{u}_2  - f \hat{u}_2 \Big) dx + \beta \int_{\partial\Omega} |\hat{u}_2|^2\,ds,
	\label{Compensated-DeepRitz-Neumann-Subprob-Functional}
\end{equation}
which relies on the precision of $\nabla \hat{u}_1^{[k]}$ and therefore benefits from \eqref{Compensated-PINN-Dirichlet-Subprob-Functional}.

Now we are ready to discretize functional integrals (\ref{Compensated-DeepRitz-Dirichlet-Subprob-Functional}, \ref{Compensated-PINN-Dirichlet-Subprob-Functional}) and \eqref{Compensated-DeepRitz-Neumann-Subprob-Functional}, where the Monte Carlo method is adopted to overcome the curse of dimensionality \cite{metropolis1949monte}. To be specific, the training sample points are generated uniformly at random inside each subdomain and at its boundary, i.e.,
\begin{equation*}
	X_{\Omega_i} = \big\{ x_n^{\Omega_i} \big\}_{n=1}^{N_{\Omega_i}},\ \ \ X_{D_i} = \big\{ x_n^{D_i} \big\}_{n=1}^{N_{D_i}},\ \ \ \text{and}\ \ \ X_{\Gamma} = \big\{ x_n^{\Gamma} \big\}_{n=1}^{N_{\Gamma}},
\end{equation*}
where $D_i = \partial\Omega_i\cap\partial\Omega$, $N_{\Omega_i}$, $N_{D_i}$, and $N_\Gamma$ represent the sample size of training datasets $X_{\Omega_i}$, $X_{D_i}$, and $X_\Gamma$, respectively. Consequently, by defining the following loss functions
\begingroup
\renewcommand*{\arraystretch}{2.6}
\begin{equation*}
\begin{array}{c}
L_{\Omega_i}( \hat{u}_i ) = \left\{\!
\begin{array}{lc}
\displaystyle \frac{1}{N_{\Omega_i}} \sum_{n=1}^{N_{\Omega_i}} \big| \Delta \hat{u}_i(x_n^{\Omega_i};\theta_i) + f(x_n^{\Omega_i}) \big|^2 & \!\!\!\text{(residual form),}  \\
\displaystyle \frac{1}{N_{\Omega_i}} \sum_{n=1}^{N_{\Omega_i}} \left( \frac12 | \nabla \hat{u}_i(x_n^{\Omega_i};\theta_i) |^2 - f(x_n^{\Omega_i}) \hat{u}_i(x_n^{\Omega_i};\theta_i)\right) & \!\!\!\text{(Ritz form),}
\end{array}\right. \\
\displaystyle L_{D_i}( \hat{u}_j ) = \frac{1}{N_{D_i}} \sum_{n=1}^{N_{D_i}} | \hat{u}_j(x_n^{D_i};\theta_j) |^2, \ \ \ L_{\Gamma}( \hat{u}_1, h^{[k]} ) = \frac{1}{N_\Gamma} \sum_{n=1}^{N_\Gamma} | \hat{u}_1(x_n^\Gamma;\theta_1) - h^{[k]}( x_n^\Gamma ) |^2, \\
\displaystyle L_N( \hat{u}_2, \hat{u}_1^{[k]} ) = \frac{1}{N_{\Omega_1}} \sum_{n=1}^{N_{\Omega_1}} \left( \nabla \hat{u}_1(x_n^{\Omega_1};\theta_1^{[k]}) \cdot \nabla \hat{u}_2(x_n^{\Omega_1};\theta_2) - f(x_n^{\Omega_1}) \hat{u}_2( x_n^{\Omega_1}; \theta_2 ) \right),
\end{array}
\end{equation*}
\endgroup
the learning task associated with (\ref{Compensated-DeepRitz-Dirichlet-Subprob-Functional}, \ref{Compensated-PINN-Dirichlet-Subprob-Functional}) is defined as
\begin{equation}
	\theta_1^{[k]} = \operatorname*{arg\,min}_{\theta_1} L_{\Omega_1}( \hat{u}_1 ) + \beta \big( L_{D_1}( \hat{u}_1 ) + L_{\Gamma}( \hat{u}_1, h^{[k]} ) \big),
	\label{Compensated-DeepRitz-Dirichlet-Subprob-Discrete}
\end{equation}
while that of the functional integral \eqref{Compensated-DeepRitz-Neumann-Subprob-Functional} is given by 
\begin{equation}
	\theta_2^{[k]} = \operatorname*{arg\,min}_{\theta_2} L_{\Omega_2}( \hat{u}_2 ) + L_N( \hat{u}_2,\hat{u}_1^{[k]} ) + \beta \big( L_{D_1}( \hat{u}_2 ) + L_{D_2}( \hat{u}_2 ) \big).
	\label{Compensated-DeepRitz-Neumann-Subprob-Discrete}
\end{equation}
Although the solution of Dirichlet subproblem is often prone to returning erroneous Neumann traces along the interface \cite{dockhorn2019discussion,bajaj2021robust}, it is evident from \eqref{Compensated-DeepRitz-Neumann-Subprob-Discrete} that our extended Neumann subproblem can be numerically solved without involving the issue of erroneous Dirichlet-to-Neumann map. Moreover, with the second-order differential operator being explicitly involved during the network training of Dirichlet subproblem \eqref{Compensated-PINN-Dirichlet-Subprob-Functional}, the resulting solution's gradient $\nabla \hat{u}_1$ is rather accurate inside the subdomain $\Omega_1$ (see Remark \ref{Remark-DtN-Map}), which is highly desirable for solving our extended Neumann subproblem \eqref{Compensated-DeepRitz-Neumann-Subprob-Functional}. 

\begin{figure}[t!]
\begin{algorithm}[H]
\caption{Dirichlet-Neumann Learning Algorithm (Two Subdomains)}
\begin{algorithmic}
\STATE{\% \textit{Initialization} }
\STATE{-- divide domain $\Omega\subset\mathbb{R}^d$ into two non-overlapping subdomains $\Omega_1$ and $\Omega_2$;}
\STATE{-- specify network structures $\hat{u}_1(x;\theta_1)$ and $\hat{u}_2(x;\theta_2)$ for each subproblem;}
\STATE{-- generate Monte Carlo training samples $X_\Gamma$, $X_{\Omega_i}$, and $X_{D_i}$ for $i=1$, 2; }
\STATE{\% \textit{Outer Iteration Loop} }
\STATE{Start with the initial guess $h^{[0]}$ along the subdomain interface $\Gamma$;}
\FOR{$k \gets 0$ to $K$ (maximum number of outer iterations)}
\WHILE{stopping criteria are not satisfied}
\STATE{\% \textit{Dirichlet Subproblem-Solving using PINNs or Deep Ritz Method} }
\STATE{
\vspace{-.65cm}
\begingroup
\renewcommand*{\arraystretch}{1.1}
\begin{equation*}
u_1^{[k]} = \operatorname*{arg\,min}_{\theta_1} L_{\Omega_1}( \hat{u}_1 ) + \beta \big( L_{D_1}( \hat{u}_1 ) + L_{\Gamma}( \hat{u}_1,, h^{[k]} ) \big)
\end{equation*}
\endgroup
\vspace{-0.45cm}
}
\STATE{\% \textit{Neumann Subproblem-Solving using Compensated Deep Ritz Method} }
\STATE{
\vspace{-.6cm}
\begingroup
\renewcommand*{\arraystretch}{1.1}
\begin{equation*}
u_2^{[k]} = \operatorname*{arg\,min}_{\theta_2} L_{\Omega_2}( \hat{u}_2 ) + L_N( \hat{u}_2, \hat{u}_1^{[k]} ) + \beta \big( L_{D_1}( \hat{u}_2 ) + L_{D_2}( \hat{u}_2 ) \big)
\end{equation*}
\endgroup
\vspace{-0.45cm}
}
\STATE{\% \textit{Update of Interface Condition using Dirichlet Trace} }
\STATE{
\vspace{-.65cm}
\begingroup
\renewcommand*{\arraystretch}{1.1}
\begin{equation*}
h^{[k+1]}(x_n^\Gamma) = \rho \hat{u}_2^{[k]}(x_n^\Gamma) + (1-\rho) h^{[k]}(x_n^\Gamma), \ \ n=1,\cdots,N_\Gamma,
\end{equation*}
\endgroup
\vspace{-0.7cm}
}
\ENDWHILE
\ENDFOR
\end{algorithmic}
\label{Algorithm-DN-Learning-2Subdomains}
\end{algorithm}
\vspace{-0.2cm}
\end{figure}

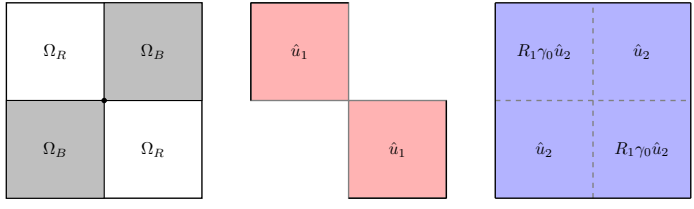
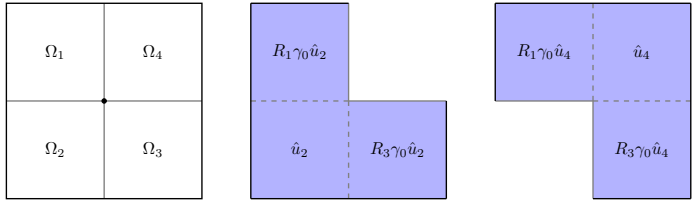
\begin{figure}[htp]
\begin{subfigure}[t]{\textwidth}
\centering
\begin{adjustbox}{max totalsize={0.7\textwidth}{0.95\textheight},center}
\begin{tikzpicture}
   \draw[color=gray, thick] (2,0) -- (2,4); 
   \draw[color=gray, thick] (0,2) -- (4,2); 
   \draw[thick] (0,0) rectangle (4,4);
   \draw[fill=lightgray] (0,0) rectangle (2,2);  
   \draw[fill=lightgray] (2,2) rectangle (4,4);   
   \node[draw=none] at (1,1) {$\Omega_B$};
   \node[draw=none] at (1,3) {$\Omega_R$};
   \node[draw=none] at (3,1) {$\Omega_R$};
   \node[draw=none] at (3,3) {$\Omega_B$};
   \draw[black,fill=black] (2,2) circle (.3ex);  
   
   \fill [red!30] (5,2) rectangle (7,4);   
   \fill [red!30] (7,0) rectangle (9,2);      
   
   \draw[color=gray, thick] (5,2) rectangle (7,4);
   \draw[color=gray, thick] (7,0) rectangle (9,2);
      
   \draw[thick] (5,2) -- (5,4);     
   \draw[thick] (5,4) -- (7,4);   
   \draw[thick] (7,0) -- (9,0);
   \draw[thick] (9,0) -- (9,2);    
      
   \node[draw=none] at (6,3) {$\hat{u}_1$};
   \node[draw=none] at (8,1) {$\hat{u}_1$};

	\fill [blue!30] (10,0) rectangle (12,2);
	\fill [blue!30] (10,2) rectangle (12,4);
	\fill [blue!30] (12,0) rectangle (14,2);
	\fill [blue!30] (12,2) rectangle (14,4);
	   
	\draw[color=gray, dashed, thick] (10,2) -- (14,2); 
	\draw[color=gray, dashed, thick] (12,0) -- (12,4); 
	
	\draw[thick] (10,0) -- (14,0);     
   \draw[thick] (10,0) -- (10,4);   
   \draw[thick] (10,4) -- (14,4);
   \draw[thick] (14,0) -- (14,4); 
	
   \node[draw=none] at (11,3) {$R_{1}\gamma_0\hat{u}_{2}$};
   \node[draw=none] at (13,3) {$\hat{u}_{2}$};
   \node[draw=none] at (13,1) {$R_{1}\gamma_0\hat{u}_{2}$};
    \node[draw=none] at (11,1) {$\hat{u}_{2}$};
   
\end{tikzpicture}
\end{adjustbox}
\vspace{-0.4cm}
\caption{\textbf{Left:} Red-Black partition into two sets. \textbf{Middle:} Dirichlet subproblem solved on $\Omega_R$. \textbf{Right:} Extended Neumann subproblem solved on $\Omega_B \cup \Omega_R$ due to the operation \eqref{DNLA-Extension}.}
\label{fig-extended-red-black}
\end{subfigure}

\begin{subfigure}[htp]{\textwidth}
\centering
\begin{adjustbox}{max totalsize={0.7\textwidth}{0.95\textheight},center}
\begin{tikzpicture}
   \draw[color=gray, thick] (2,0) -- (2,4); 
   \draw[color=gray, thick] (0,2) -- (4,2); 
   \draw[thick] (0,0) rectangle (4,4);
   \node[draw=none] at (1,1) {$\Omega_{2}$};
   \node[draw=none] at (1,3) {$\Omega_{1}$};
   \node[draw=none] at (3,1) {$\Omega_{3}$};
   \node[draw=none] at (3,3) {$\Omega_{4}$};
   \draw[black,fill=black] (2,2) circle (.3ex);  

    \fill [blue!30] (5,2) rectangle (7,4);  
    \fill [blue!30] (7,0) rectangle (9,2);   
    \fill [blue!30] (5,0) rectangle (7,2);   
   
   \draw[color=gray, dashed, thick] (7,0) -- (7,2); 
   \draw[color=gray, dashed, thick] (5,2) -- (7,2);
   \draw[color=gray, thick] (7,2) -- (7,4);
   \draw[color=gray, thick] (7,2) -- (9,2);
   
 	\draw[thick] (5,0) -- (9,0);     
   \draw[thick] (5,0) -- (5,4);   
   \draw[thick] (5,4) -- (7,4);
   \draw[thick] (9,0) -- (9,2);   
   
   \node[draw=none] at (6,1) {$\hat{u}_{2}$};
   \node[draw=none] at (6,3) {$R_{1}\gamma_0\hat{u}_{2}$};
   \node[draw=none] at (8,1) {$R_{3}\gamma_0\hat{u}_{2}$};
   
   \fill [blue!30] (10,2) rectangle (12,4);  
   \fill [blue!30] (12,0) rectangle (14,2);   
   \fill [blue!30] (12,2) rectangle (14,4);  
      
   \draw[color=gray, dashed, thick] (12,2) -- (12,4); 
   \draw[color=gray, dashed, thick] (12,2) -- (14,2); 
    \draw[color=gray, thick] (12,0) -- (12,2);
   \draw[color=gray, thick] (10,2) -- (12,2);

 	\draw[thick] (12,0) -- (14,0);     
   \draw[thick] (14,0) -- (14,4);   
   \draw[thick] (10,4) -- (14,4);
   \draw[thick] (10,4) -- (10,2);      
   
   \node[draw=none] at (11,3) {$R_{1}\gamma_0\hat{u}_{4}$};
   \node[draw=none] at (13,3) {$\hat{u}_{4}$};
   \node[draw=none] at (13,1) {$R_{3}\gamma_0\hat{u}_{4}$};
   
\end{tikzpicture}
\end{adjustbox}
\caption{\textbf{Left:} Partition of domain into four sets, with Dirichlet (or Neumann) subproblems defined on $\Omega_1$ and $\Omega_3$ (or $\Omega_2$ and $\Omega_4$). \textbf{Middle:} Extended Neumann subproblem solved on $\Omega_2\cup\Omega_1\cup\Omega_3$. \textbf{Right:} Extended Neumann subproblem solved on $\Omega_4\cup\Omega_1\cup\Omega_3$.
}
\label{fig-extended-4-area}
\end{subfigure}
\vspace{-0.5cm}
\caption{Illustration of neural network extension operators for 4 subdomains.}
\vspace{-0.8cm}
\end{figure}

In summary, our proposed Dirichelt-Neumann learning algorithm is presented in Algorithm \ref{Algorithm-DN-Learning-2Subdomains}, where the mini-batch data are not relabelled for notational simplicity and the stopping criteria can be constructed by measuring the difference between two consecutive iterations \cite{li2019d3m}. We also note that our Dirichlet-Neumann learning algorithm has sequential steps that inherited from the original scheme \cite{toselli2004domain,quarteroni1999domain}, and various techniques have been developed to solve subproblems in parallel (see \cite{mathew2008domain} and references cited therein).

\begin{remark}
Note that in the case of two subdomains, e.g., red-black partition in \autoref{fig-extended-red-black}, the solution $\hat{u}_2(x;\theta_2)$ of our extended Neumann subproblem \eqref{Compensated-DeepRitz-Neumann-Subprob-Functional} is defined over the entire domain, which seems to incur enormous cost at the first glance. In fact, the extension operation \eqref{DNLA-Extension} only involves subdomains that have a common interface with the underlying subproblem. Therefore, the computational domain of our extended subproblem \eqref{Compensated-DeepRitz-Neumann-Subprob-Functional} can be locally defined (see \autoref{fig-extended-4-area} for example).
\end{remark}

\begin{remark}\label{Remark-High-Contrast}
Our proposed method can also be used to solve the elliptic interface problem with high-contrast coefficients \cite{li2006immersed,he2022mesh,sun2023dirichlet}, which is formally written as
\begingroup
\renewcommand*{\arraystretch}{1.1}
\begin{equation*}
\begin{array}{cl}
-\nabla \cdot \left( c(x) \nabla u(x)  \right) = f(x)\ \ & \text{in}\ \Omega,\\
u(x) = 0\ \ & \text{on}\ \partial \Omega, 
\end{array}
\end{equation*}
\endgroup
where $\Gamma=\partial \Omega_1\cap\partial\Omega_2$ is an immersed interface (see \autoref{fig-domain-decomposition} for example), the coefficient function $c(x)$ is piecewise constant with respect to the decomposition of domain
\begingroup
\renewcommand*{\arraystretch}{1.1}
\begin{equation*}
c(x) = \left\{
\begin{array}{cl}
c_1>0\ \ & \text{in}\ \Omega_1,\\
c_2\gg c_1\ \ & \text{in}\ \Omega_2,
\end{array}\right.
\end{equation*}
\endgroup
and natural jump conditions \cite{li2006immersed} are given by
\begin{equation*}
	[u] = 0\ \ \ \text{and}\ \ \ \left[ c\frac{\partial u}{\partial \bm{n}} \right] = q\ \ \ \text{on}\ \Gamma.
\end{equation*}
Applying Green's formula in each subdomain and then adding them together, we obtain the Galerkin formulation: find $u_1\in V_1$ and $u_2\in V_2$ such that
\begingroup
\renewcommand*{\arraystretch}{1.1}
\begin{equation*}
\begin{array}{cl}
b_1(u_1, v_1) = (f, v_1)_1\ \ & \text{for any}\ v_1\in V_1^0,\\
u_1 = u_2\ \ & \text{on}\ \Gamma, \\
b_2(u_2, v_2) = (f, v_2)_2 + (f, R_1\gamma_0 v_2)_1 - b_1(u_1, R_1\gamma_0 v_2) - (q,v_2)_\Gamma,\ \ & \text{for any}\ v_2\in V_2,
\end{array}
\end{equation*}
\endgroup
where the bilinear forms are defined as
\begin{equation*}
	b_i(u_i, v_i) = \int_{\Omega_i} c_i \nabla u_i \cdot \nabla v_i\,dx,\ \ \ (f, v_i)_i = \int_{\Omega_i} f v_i\, dx,\ \ \ \text{and}\ \ \ (q,v_2)_\Gamma = \int_\Gamma qv\,ds,
\end{equation*}
for $i=1$, 2. By parametrizing solutions as neural networks, i.e., $u_i(x)\approx \hat{u}_i(x;\theta_i)$, and employing the neural network extension operator $R_1\gamma_0 \hat{u}_2 = \hat{u}_2$, the learning task associated with the Dirichlet subproblem\footnote{Here, the residual form is used instead since the solution on each subdomain can be assumed regular enough, which would result in a good approximation of $\nabla u_1$ inside the subdomain $\Omega_1$.} on $\Omega_1$ gives
\begin{equation}
	\theta_1^{[k]} = \operatorname*{arg\,min}_{\theta_1} \int_{\Omega_1} | \nabla\cdot(c_1 \nabla \hat{u}_1) + f  |^2dx + \beta \left( \int_{\partial\Omega_1\cap \partial\Omega} |\hat{u}_1|^2\,ds + \int_{\Gamma} |\hat{u}_1 - \hat{u}_2|^2\,ds \right),
	\label{High-Contrast-Dirichlet-Subprob}
\end{equation}
while that of the Neumann subproblem takes on the form
\begingroup
\vspace{-0.2cm}
\renewcommand*{\arraystretch}{2}
\begin{equation}
\begin{array}{cl}
\displaystyle \theta_2^{[k]} = \operatorname*{arg\,min}_{\theta_2} \!\!\!\!\! & \displaystyle \int_{\Omega_2} \Big( \frac{c_2}{2} | \nabla \hat{u}_2 |^2  - f \hat{u}_2 \Big) dx + \int_{\Omega_1} \Big( c_1 \nabla \hat{u}_1^{[k]} \cdot \nabla \hat{u}_2  - f \hat{u}_2 \Big) dx \\
& \displaystyle  + \int_\Gamma q\hat{u}_2\,ds + \beta \int_{\partial\Omega} |\hat{u}_2|^2\,ds.
\end{array}
\label{High-Contrast-Neumann-Subprob}
\end{equation}
\endgroup
Accordingly, an iterative learning approach for solving the elliptic interface problem with high-contrast coefficients can be immediately constructed from \eqref{High-Contrast-Dirichlet-Subprob} and \eqref{High-Contrast-Neumann-Subprob}, while a further theoretical investigation can be found in \cite{sun2023dirichlet}.
\end{remark}

\subsection{Neumann-Neumann  Learning Algorithm}
Similar in spirit, the compensated deep Ritz method can be applied to construct the Neumann-Neumann learning algorithm (see \autoref{fig-big-picture}). Using same notations as before, the Neumann-Neumann scheme (see Algorithm \ref{DDM-Flux-Exchange}) can be written in an equivalent Ritz formulation: given the initial guess $h^{[0]}\in H_{00}^{\frac12}(\Gamma)$, then solve for $k\geq 0$ and $i=1$, 2,
\begingroup
\vspace{-0.2cm}
\renewcommand*{\arraystretch}{2}
\begin{equation*}
\begin{array}{l}
\displaystyle 1)\ u_i^{[k]} = \operatorname*{arg\,min}_{u_i\in V_i,\ u_i|_\Gamma = h^{[k]}}  \frac12  a_i(u_i, u_i) - (f, u_i)_i, \\
\displaystyle 2)\ \psi_i^{[k]} = \operatorname*{arg\,min}_{\psi_i \in V_i} \ \frac12 a_i(\psi_i, \psi_i) + (f, \psi_i)_i + (f, R_{3-i}\gamma_0 \psi_i)_{3-i} - a_i(u_i^{[k]}, \psi_i) \\
\displaystyle \qquad\qquad\qquad\qquad - a_{3-i}(u_{3-i}^{[k]}, R_{3-i}\gamma_0 \psi_i), \\
\displaystyle 3)\ h^{[k+1]} = h^{[k]} - \rho (\psi_1^{[k]} + \psi_2^{[k]}) \ \ \ \text{on}\ \Gamma,
\end{array}
\end{equation*}
\endgroup
with $\rho\in(0,\rho_{\text{max}})$ denoting the acceleration parameter. Next, by parametrizing the unknown solutions as neural networks, that is, for $i=1$, 2,
\begin{equation*}
	u_i^{[k]}(x) \approx \hat{u}_i^{[k]}(x) := \hat{u}_i(x;\theta_i^{[k]}) \ \ \ \text{and}\ \ \ \psi_i^{[k]}(x) \approx \hat{\psi}_i^{[k]}(x) := \hat{\psi}_i(x;\eta_i^{[k]})
\end{equation*}
and by employing extension operators $R_1\gamma_0 \hat{\psi}_2(x;\eta_2) = \hat{\psi}_2(x;\eta_2)$ and $R_2\gamma_0 \hat{\psi}_1(x;\eta_1) = \hat{\psi}_1(x;\eta_1)$, the learning tasks associated with the Neumann-Neumann algorithm are given by, for $i=1$, 2, 
\begingroup
\vspace{-0.25cm}
\renewcommand*{\arraystretch}{2.5}
\begin{equation*}
\begin{array}{c}
\displaystyle \theta_i^{[k]} = \operatorname*{arg\,min}_{\theta_i} \int_{\Omega_i} | \Delta \hat{u}_i + f |^2 dx + \beta \left( \int_{\partial\Omega_i\cap \partial\Omega} |\hat{u}_i|^2\,ds + \int_{\Gamma} |\hat{u}_i - h^{[k]}|^2\,ds \right), \\
\displaystyle \eta_i^{[k]} = \operatorname*{arg\,min}_{\eta_i} \int_{\Omega_i} \Big( \frac12 | \nabla \hat{\psi}_i |^2 + f \hat{\psi}_i - \nabla \hat{u}_i^{[k]} \cdot \nabla \hat{\psi}_i \Big) dx + \beta \int_{\partial\Omega} |\hat{\psi}_i|^2ds \\
\displaystyle + \int_{\Omega_{3-i}} \Big( f\hat{\psi}_i - \nabla \hat{u}_{3-i}^{[k]} \cdot \nabla \hat{\psi}_i \Big) dx,\qquad\qquad
\end{array}
\end{equation*}
\endgroup
where $\beta>0$ is the penalty coefficient and training tasks associated with Dirichlet subproblems are defined in a residual form as before. Therefore, the iterative learning approach can be constructed after applying numerical integration.

\subsection{Dirichlet-Dirichlet Learning Algorithm}
Next, to build the Dirichlet-Dirichlet learning algorithm, we first rewrite the iterative process (see Algorithm \ref{DDM-Flux-Exchange}) as: given the initial guess $h^{[0]}\in H_{00}^{\frac12}(\Gamma)$ and $\psi_1^{[0]}=\psi_2^{[0]}=0$, then solve for $k\geq 0$,
\begingroup
\vspace{-0.2cm}
\renewcommand*{\arraystretch}{3.8}
\begin{equation*}
\begin{array}{l}
\displaystyle 1)\ \textnormal{solve for}\ u_i^{[k]}: 
\begingroup
\renewcommand*{\arraystretch}{1.3}
\left\{
\begin{array}{cl}
- \Delta u_i^{[k]} = f \ \ & \text{in}\ \Omega_i,\\
u_i^{[k]} = 0\ \ & \text{on}\ \partial\Omega\cap\partial\Omega_i, \\
\nabla u_i^{[k]} \cdot \bm{n}_i = h^{[k]} - \rho ( \nabla \psi_1^{[k]} \cdot \bm{n}_1 + \nabla \psi_2^{[k]} \cdot \bm{n}_2 )  \ \ & \text{on}\ \Gamma,
\end{array}\right.
\endgroup \\
\displaystyle 2)\ \textnormal{solve for}\ \psi_i^{[k+1]}: 
\begingroup
\renewcommand*{\arraystretch}{1.3}
\left\{
\begin{array}{cl}
- \Delta \psi_i^{[k+1]} = 0 \ \ & \text{in}\ \Omega_i,\\
\psi_i^{[k+1]} = 0\ \ & \text{on}\ \partial\Omega\cap\partial\Omega_i, \\
\psi_i^{[k+1]} = u_1^{[k]} - u_2^{[k]} \ \ & \text{on}\ \Gamma,
\end{array}\right.
\endgroup
\end{array}
\end{equation*}
\endgroup
where $i=1$, 2. Using same notations as before, the Green's theorem indicates that the variational formulation of our Neumann subproblem reads: 
\begin{equation*}
u_i^{[k]} \!=\! \operatorname*{arg\,min}_{u_i\in V_i}  \frac12  a_i(u_i, u_i) \!-\! (f, u_i)_i \!-\! (h^{[k]},u_i)_\Gamma \!+\! \rho \big( a_i(\psi_i^{[k]}, u_i) \!+\! a_{3-i}( \psi_{3-i}^{[k]}, R_{3-i} \gamma_0 u_i ) \big) 
\end{equation*}
where $(\cdot,\cdot)_\Gamma$ denotes the $L_2$ inner product on $\Gamma$ and the Dirichlet subproblem gives
\begingroup
\renewcommand*{\arraystretch}{1.3}
\begin{equation*}
\begin{array}{c}
\displaystyle \psi_i^{[k]} = \operatorname*{arg\,min}_{\psi_i\in V_i,\ \psi_i|_\Gamma = u_1^{[k]} - u_2^{[k]}}  \frac12 a_i(\psi_i, \psi_i)  
\end{array}
\end{equation*}
\endgroup
for $i=1$, 2. Consequently, by parametrizing unknown solutions as neural networks
\begin{equation*}
	u_i^{[k]}(x) \approx \hat{u}_i^{[k]}(x) := \hat{u}_i(x;\theta_i^{[k]}) \ \ \ \text{and}\ \ \ \psi_i^{[k]}(x) \approx \hat{\psi}_i^{[k]}(x) := \hat{\psi}_i(x;\eta_i^{[k]})
\end{equation*}
for $i=1$, 2, the neural network extension operators $R_1\gamma_0 \hat{u}_2(x,\theta_2) = \hat{u}_2(x,\theta_2)$ and $R_2\gamma_0 \hat{u}_1(x,\theta_1) = \hat{u}_1(x,\theta_1)$ are used to construct the learning tasks, i.e., for $i=1$, 2,
\begingroup
\vspace{-.25cm}
\renewcommand*{\arraystretch}{2.5}
\begin{equation*}
\begin{array}{c}
\displaystyle \theta_i^{[k]} = \operatorname*{arg\,min}_{\theta_i} \int_{\Omega_i} \Big( \frac12 |\nabla \hat{u}_i  |^2 - f \hat{u}_i + \rho \nabla \hat{\psi}_i^{[k]}\cdot \nabla \hat{u}_i\Big) dx + \rho \int_{\Omega_{3-i}} \nabla \hat{\psi}^{[k]}_{3-i} \cdot \nabla \hat{u}_i\,dx \\
\displaystyle - \int_\Gamma h^{[k]}\hat{u}_i\,ds + \beta \int_{\partial\Omega_i\cap \partial\Omega} |\hat{u}_i|^2\,ds,\qquad\qquad\qquad\ \\
\displaystyle \eta_i^{[k]} = \operatorname*{arg\,min}_{\eta_i} \int_{\Omega_i}  |\Delta \hat{\psi}_i  |^2 \,dx + \beta \left( \int_{\partial\Omega_i\cap \partial\Omega} |\hat{\psi}_i|^2\,ds + \int_{\Gamma} |\hat{\psi}_i - \hat{u}_1^{[k]} + \hat{u}_2^{[k]}|^2\,ds \right), \\
\end{array}
\end{equation*}
\endgroup
where the loss functions of Dirichlet subproblems are defined in a residual form.

\subsection{Robin-Robin Learning Algorithm}\label{Section-RRLM}
As mentioned before, the Robin-Robin algorithm only requires the exchange of Dirichlet traces between neighbouring subproblems, however, it may suffer from the issue of weight imbalance (see Remark \ref{Remark-Robin}). More specifically, let $\kappa_1=1$ in what follows, then a relatively large value of $\kappa_2\gg\kappa_1$ is typically required in order to achieve fast convergence along the outer iteration \cite{chen2014optimal}. To alleviate the negative influence of $\kappa_2\gg\kappa_1$, our compensated deep Ritz method is a promising alternative for realizing the Robin-Robin algorithm. 

Note that in terms of differential operator, the decomposed subproblem with parameter $\kappa_2\gg \kappa_1=1$ in the Robin-Robin algorithm \cite{quarteroni1999domain} can be rewritten as
\begin{equation*}
\begingroup
\renewcommand*{\arraystretch}{1.3}
\left\{
\begin{array}{ll}
- \Delta u_2^{[k]} = f \ \ & \text{in}\ \Omega_2,\\
u_2^{[k]} = 0\ \ & \text{on}\ \partial\Omega\cap\partial\Omega_2, \\
\kappa_2 u_2^{[k]} + \nabla u_2^{[k]} \cdot \bm{n}_2 = \kappa_2 u_1^{[k]} - \nabla u_1^{[k]} \cdot \bm{n}_1 \ \ & \text{on}\ \Gamma.
\end{array}\right.
\endgroup
\end{equation*}
Using same notations as before, it is equivalent to find $u_2^{[k]}\in V_2$ such that
\begin{equation}
	a_2(u_2^{[k]}, v_2) = (f, v_2)_2 + ( \kappa_2 (u_1^{[k]} - u_2^{[k]}) - \nabla u_1^{[k]} \cdot \bm{n}_1, v_2 )_\Gamma\ \ \ \textnormal{for any}\ v_2\in V_2.
	\label{RRLM-2nd-RobinProb-Weak}
\end{equation}
Next, by using the Green's formula, we arrive at another form of \eqref{RRLM-2nd-RobinProb-Weak}, that is,
\begin{equation*}
	a_2(u_2^{[k]}, v_2) = (f, v_2)_2 + \kappa_2 ( u_1^{[k]} - u_2^{[k]}, v_2 )_\Gamma - a_1( u_1^{[k]}, R_1\gamma_0 v_2 ) + (f, R_1\gamma_0 v_2)_1
\end{equation*}
for any $v_2\in V_2$. Therefore, the energy formulation of \eqref{RRLM-2nd-RobinProb-Weak} is given by
\begingroup
\renewcommand*{\arraystretch}{1.3}
\begin{equation*}
\begin{array}{c}
\displaystyle u_2^{[k]} = \operatorname*{arg\,min}_{u_2\in V_2}  \frac12  a_2(u_2, u_2) - (f, u_2)_2 - \kappa_2( u_1^{[k]} - \frac12 u_2, u_2)_\Gamma +  a_1( u_1^{[k]},  R_1\gamma_0 u_2) \\
\displaystyle - (f, R_1\gamma_0u_2)_1 \qquad\qquad\qquad\qquad\qquad\qquad\qquad\qquad
\end{array}
\end{equation*}
\endgroup
which completely differs from the original PINNs approach \eqref{RRLM-PINN-Subprob-Functional}. Next, by parametrizing unknown solutions as neural networks, i.e.,
\begin{equation*}
	u_1^{[k]}(x) \approx \hat{u}_1^{[k]}(x) := \hat{u}_1(x;\theta_1^{[k]})\ \ \ \text{and}\ \ \ u_2^{[k]}(x) \approx \hat{u}_2^{[k]}(x) := \hat{u}_2(x;\theta_2^{[k]}),
\end{equation*}
and by employing our neural network extension operator
\begin{equation*}
	R_1\gamma_0 \hat{u}_2(x;\theta_2) = \hat{u}_2(x;\theta_2)\in V_1,
\end{equation*}
the learning task associated with the second Robin problem takes on the form:
\begingroup
\renewcommand*{\arraystretch}{2.5}
\begin{equation*}
\begin{array}{c}
\displaystyle \theta_2^{[k]} = \operatorname*{arg\,min}_{\theta_2} \int_{\Omega_2} \Big( \frac12 |\nabla \hat{u}_2  |^2 - f \hat{u}_2 \Big) dx - \kappa_2 \int_\Gamma \Big( \hat{u}_1^{[k]}\hat{u}_2 - \frac12 |\hat{u}_2|^2 \Big) ds  \\
\displaystyle \qquad\qquad +\ \int_{\Omega_1} \Big( \nabla \hat{u}_1^{[k]} \cdot \nabla \hat{u}_2 - f\hat{u}_2 \Big) dx + \beta \int_{\partial\Omega} |\hat{u}_2|^2\,ds,
\end{array}
\end{equation*}
\endgroup
which removes the issue of weight imbalance within the Robin boundary condition.

\section{Numerical Experiments}\label{Section-Experiments}
To validate the effectiveness of our proposed domain decomposition learning algorithms, we conduct experiments using Dirichlet-Neumann and Robin-Robin learning algorithms on a wide range of elliptic boundary value problems in this section. Here, the Neumann-Neumann and Dirichlet-Dirichlet learning algorithms are omitted for space considerations. For brevity, we refer to our Dirichlet-Neumann learning algorithm as DNLA (PINNs/deep Ritz), with bracket indicating the type of deep learning method used for solving the Dirichlet subproblem. In contrast to our proposed algorithms, the existing learning approach \cite{li2020deep} for realizing the Dirichlet-Neumann algorithm is based on a direct substitution of local solvers with PINNs, which we refer to as DN-PINNs in what follows. On the other hand, the update of interface conditions in the Robin-Robin algorithm only relies on the exchange of Dirichlet traces, however, the subproblem-solving may suffer from the issue of weight imbalance as discussed in Remark \ref{Remark-Robin}. To further investigate its influence on the convergence of outer iterations, the Robin-Robin algorithm is realized using PINNs and compensated deep Ritz methods after the empirical study of DNLA, which is referred to as RR-PINNs and RRLA (PINN/deep Ritz) in a similar fashion.

For practical implementation \cite{li2020deep,karniadakis2021physics}, the network architecture deployed for each subproblem is a fully-connected neural network with $8$ hidden layers of $50$ neurons each \cite{goodfellow2016deep}. The hyperbolic tangent activation function is assigned to each neuron, which is differentiable and smooth enough to capture our local solutions. During training, we randomly sample $N_{\Omega_i}=20k$ points from each interior subdomain $\Omega_i$,  $N_{\Gamma}=5k$ points from the interface $\Gamma$, and $N_D=5k$ points from each boundary $\partial\Omega_i\setminus\Gamma$ of length equal to the interface. The trained models are then evaluated on the test dataset, i.e., $N_\Omega = 10k$ points that are uniformly distributed over the entire domain, and compared with the true solution to assess their performances. The penalty coefficient is set to $\beta=400$ and the number of mini-batches is chosen as $5$ for all simulations. When executing the learning task on each subdomain, the initial learning rate of the Adam optimizer is set to be $0.01$, which is divided by 10 at the 600-th and 800-th epoch. The training process terminates after $1k$ epochs for each decomposed subproblem, and we choose the model with minimum training loss for subsequent operations. The stopping criterion we set here is that either the relative-$L_2$ error between two consecutive iterations is less than 0.01 or the number of outer iteration reaches 30. All experiments are implemented using PyTorch 1.8.1 and trained on the NVIDIA GeForce RTX 3090.

\subsection{Dirichlet-Neumann Learning Algorithm}
As a representative benchmark, we consider deep learning-based approaches for realizing the non-overlapping Dirichlet-Neumann algorithm in this subsection. More precisely, a comparative study between DN-PINNs, DNLA (PINN), and DNLA (deep Ritz) is presented, with experiments conducted on a wide variety of elliptic boundary value problems to demonstrate the effectiveness and flexibility of our proposed methods. 

\subsubsection{Poisson's Equation with Simple Interface}
\begin{figure}[t!]
\centering
\includegraphics[width=0.17\textwidth]{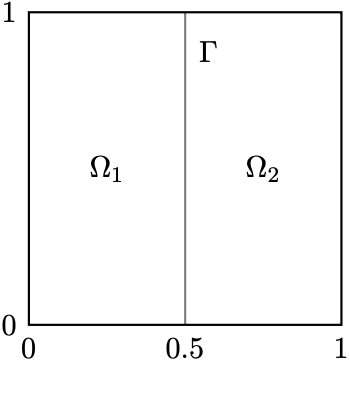}
\hspace*{0.1cm}
\includegraphics[width=0.192\textwidth]{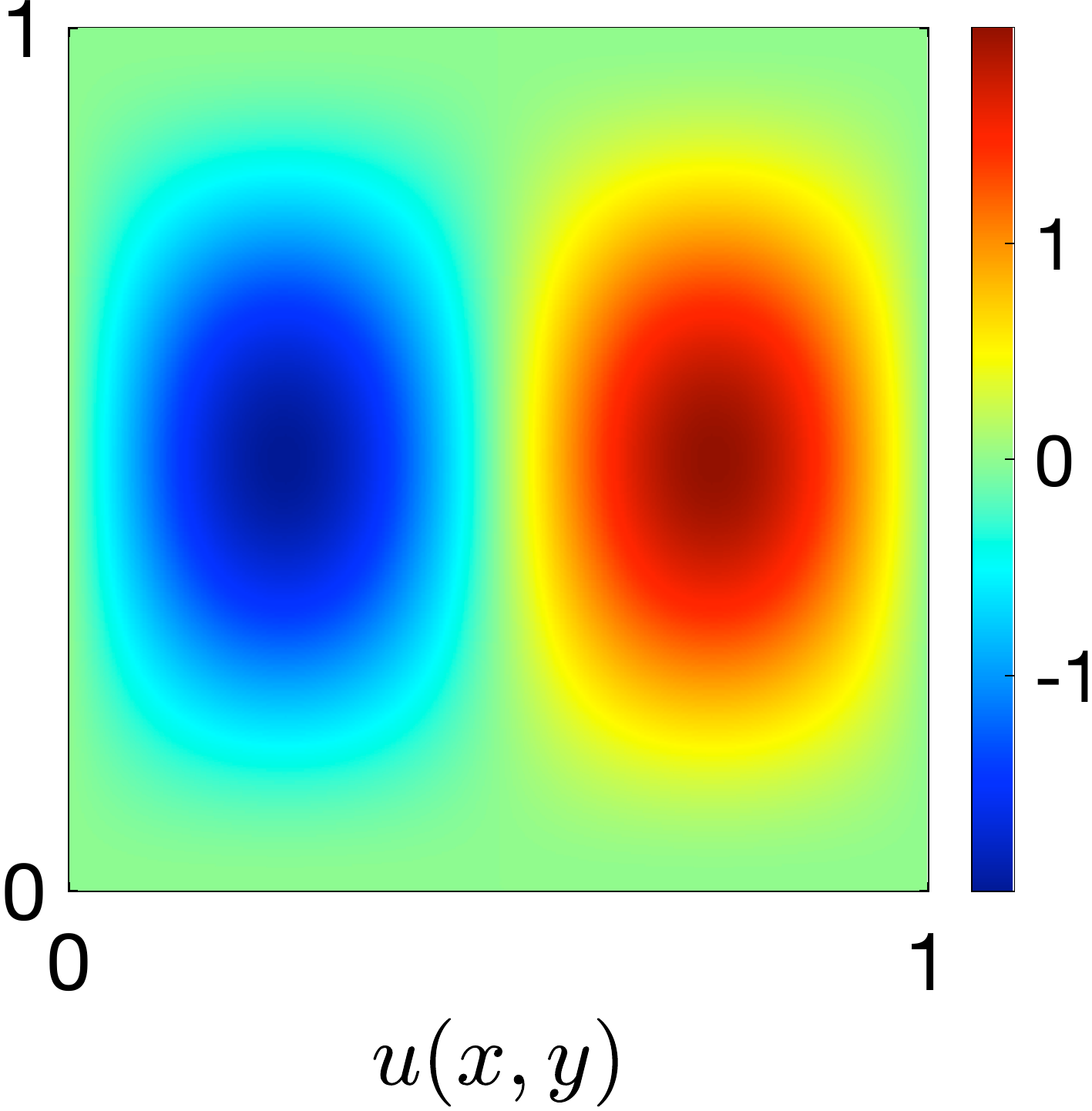}
\includegraphics[width=0.2\textwidth]{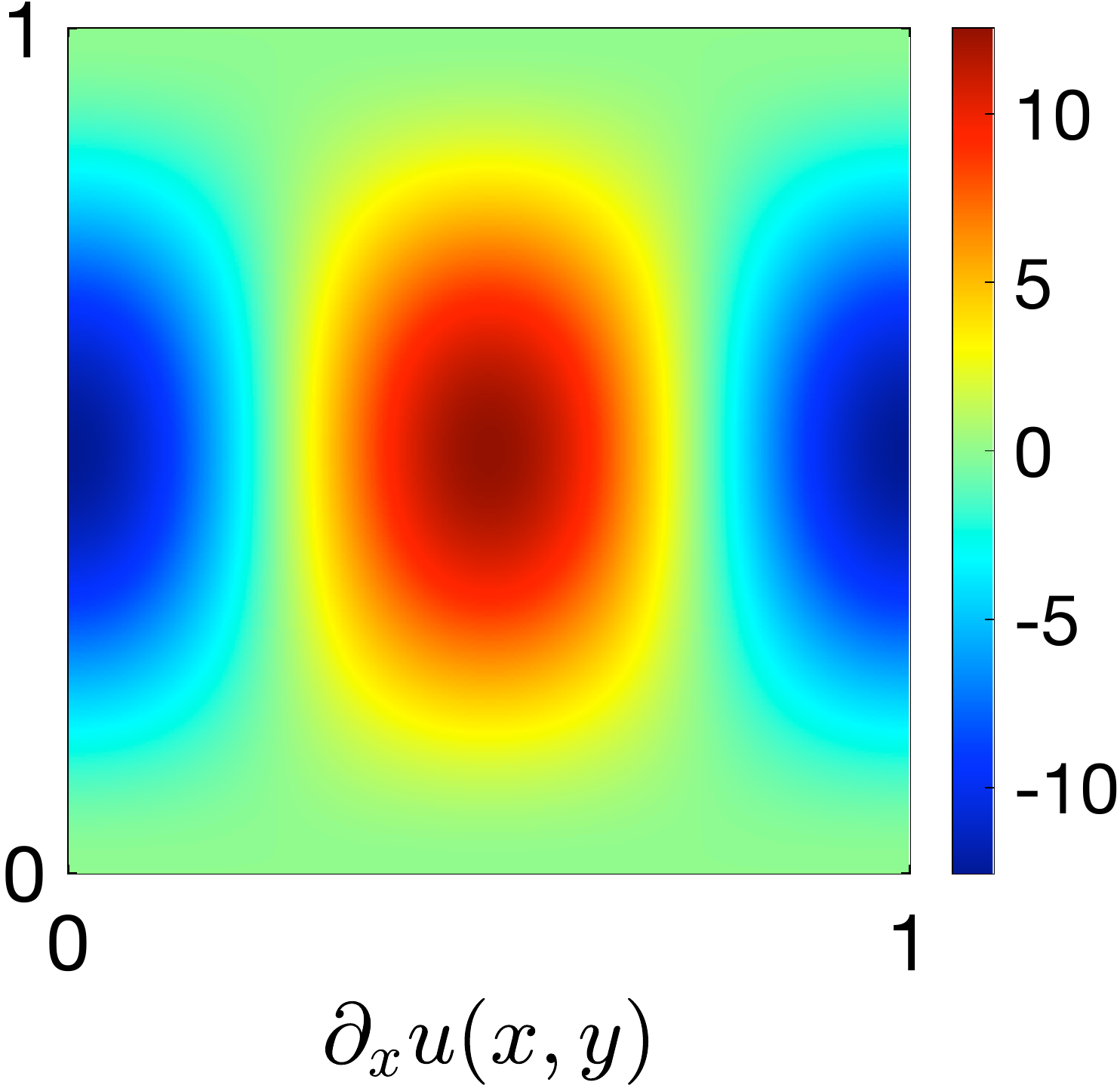}
\includegraphics[width=0.191\textwidth]{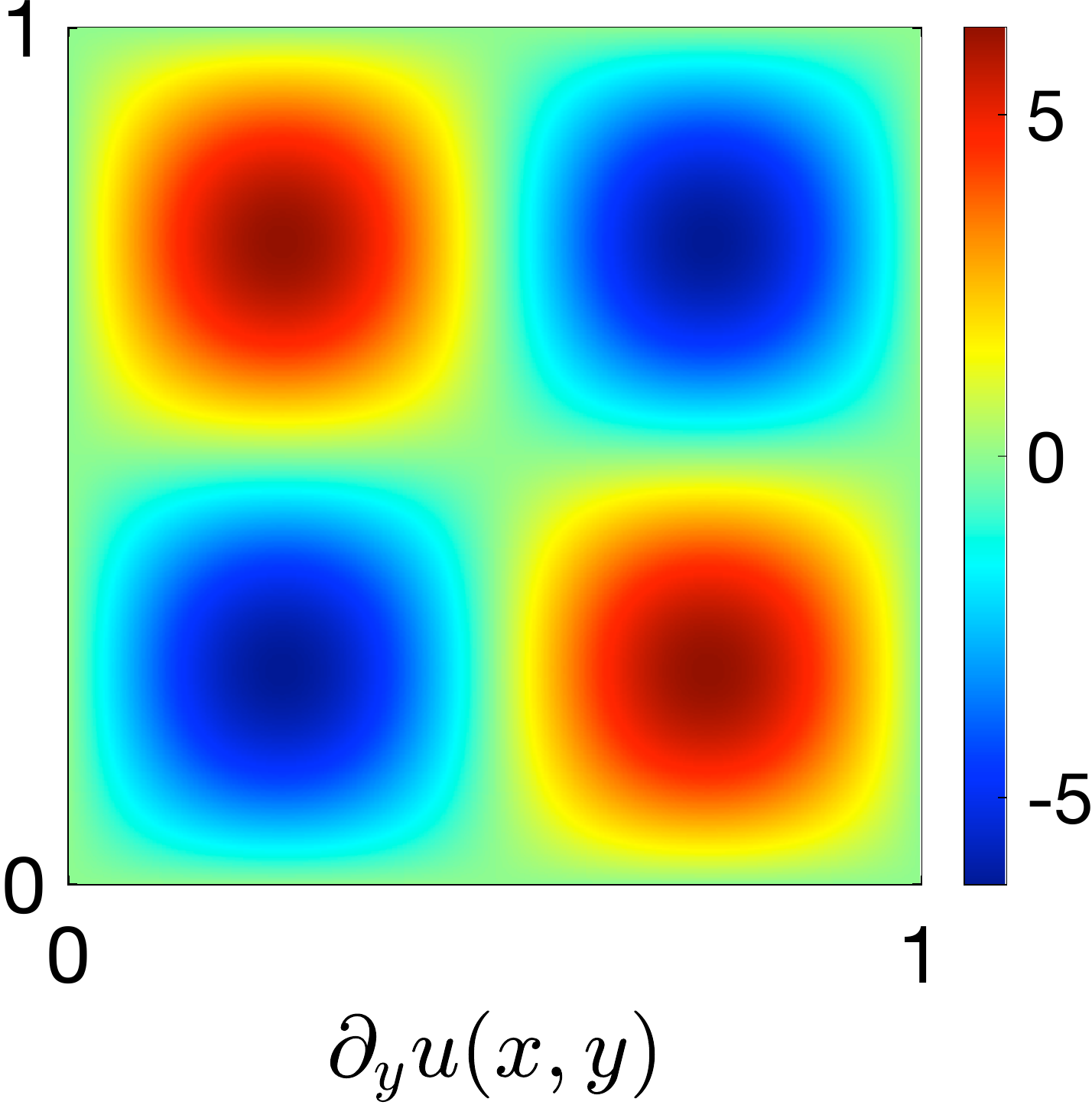}
\vspace{-0.1cm}
\caption{From left to right: decomposition into two subdomains, true solution $u(x,y)$, and its partial derivatives $\partial_x u(x,y)$, $\partial_y u(x,y)$ for the numerical example \eqref{Experiments-DNLA-ex1}.}
\label{Experiments-DNLA-ex1-exact-solution}
\vspace{-0.45cm}
\end{figure}

First, we consider the benchmark Poisson problem in two dimension, that is,
\begin{equation}
\begin{array}{cl}
-\Delta u(x,y)  = 4 \pi^2 \sin(2 \pi x)  (2 \cos(2 \pi y) - 1)  \ & \text{in}\ \Omega=(0,1)^2,\\
u(x,y) = 0\ \ & \text{on}\ \partial \Omega,
\end{array}
\label{Experiments-DNLA-ex1}
\end{equation}
where the true solution is given by $u(x,y) = \sin(2\pi x)(\cos(2\pi y)-1)$ and the interface $\Gamma=\partial\Omega_1\cap\partial\Omega_2$ is a straight line segment from $(0.5,0)$ to $(0.5,1)$ as shown in \autoref{Experiments-DNLA-ex1-exact-solution}. It is noteworthy that our exact solution reaches local extrema at $(0.5,0.5)$, thereby deviations in estimating the Neumann trace at and near the extreme point \cite{bajaj2021robust} can create a cascading effect in the convergence of outer iterations, which differs from other examples that have simple gradients on the interface \cite{li2020deep}.
 
\begin{figure}[t!]
\centering
\includegraphics[width=0.192\textwidth]{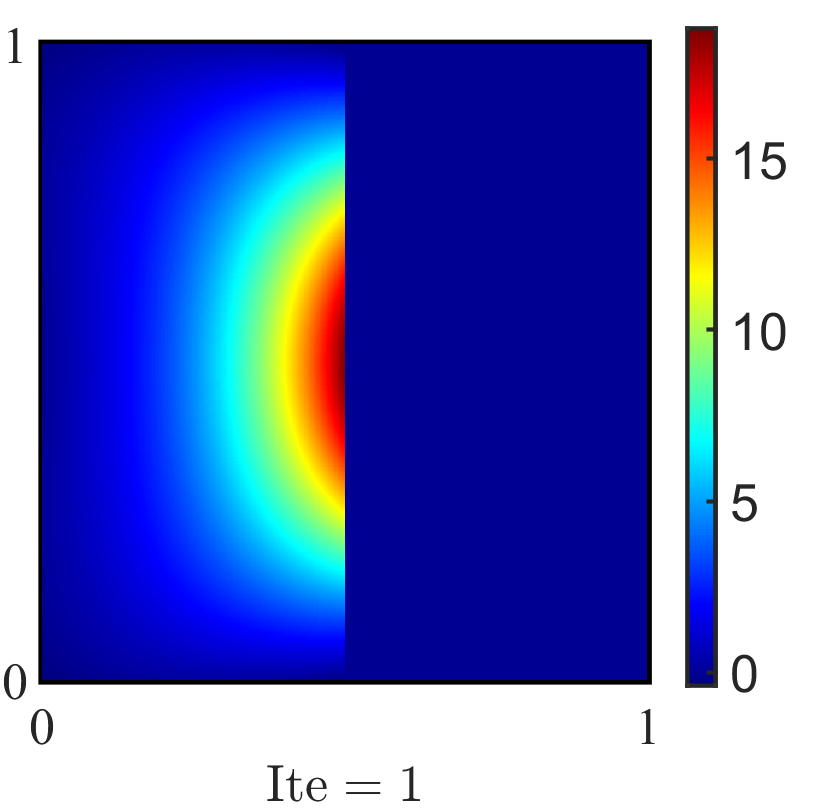}
\includegraphics[width=0.192\textwidth]{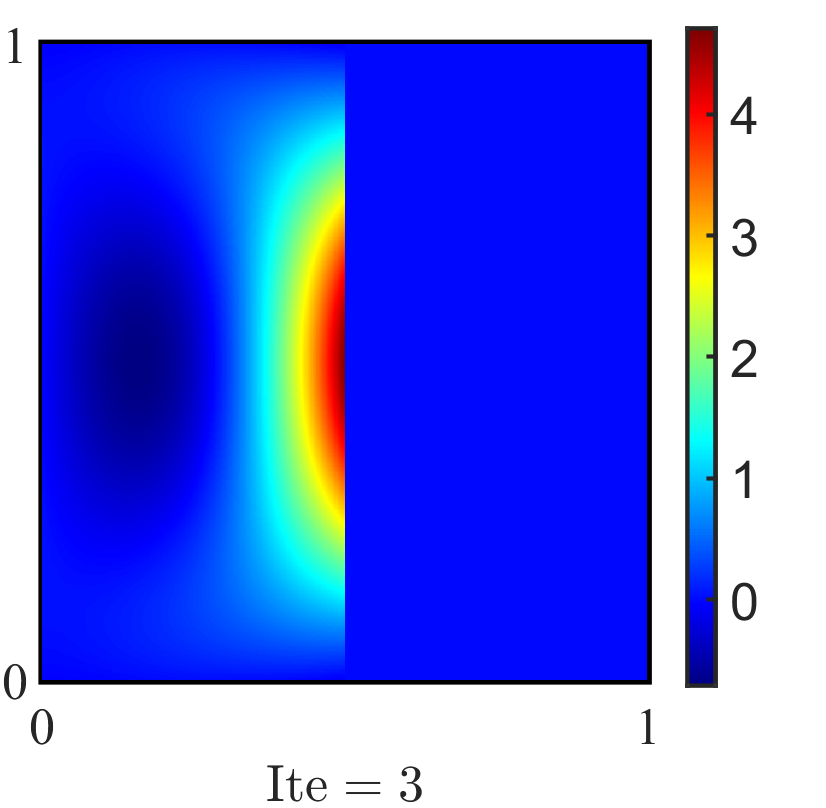}
\includegraphics[width=0.192\textwidth]{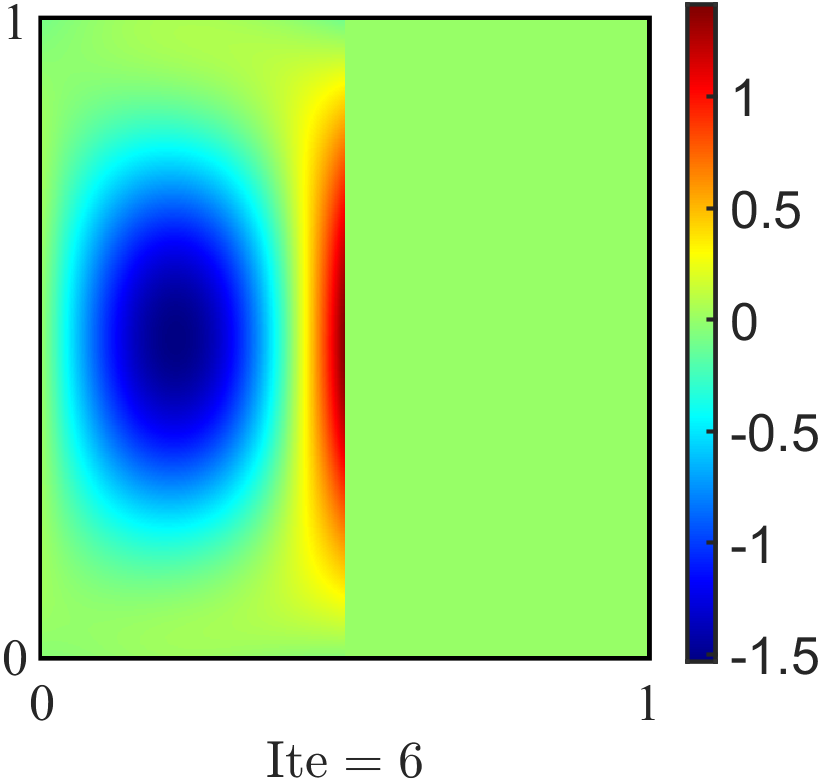}
\includegraphics[width=0.192\textwidth]{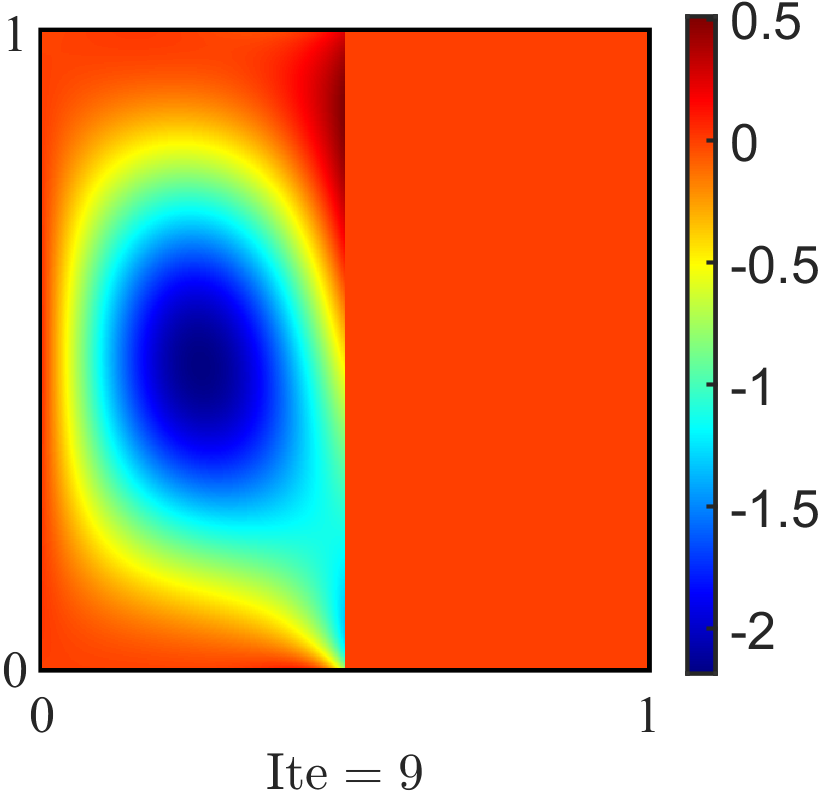}
\vspace{-0.1cm}
\caption{Iterative solutions $\hat{u}^{[k]}(x,y)$ of \eqref{Experiments-DNLA-ex1} using DN-PINNs on the test dataset.}
\label{Experiments-DNLA-ex1-DN-PINNs}
\vspace{-0.4cm}
\end{figure}

\begin{table}[t!]
\vspace{-0.2cm}
\small
\caption{Relative-$L_2$ errors of the network solution along the outer iteration for \eqref{Experiments-DNLA-ex1}, with mean value ($\pm$ standard deviation) being reported over 5 independent runs.}
\vspace{-0.1cm}
\centering
\renewcommand{\arraystretch}{1.1}
\begin{tabular}{ | c || c | c | c | c | c | c |  }
\hline
\multicolumn{2}{|c|}{ \diagbox[width=16em]{Relative Errors}{Outer Iterations} } & 1  & 3 & 5 & 7 & 9  \\
\hline	
\hline
\multirow{5}{*}{$ \displaystyle \!\! \frac{ \lVert \hat{u}^{[k]} \!-\! u \rVert_{L_2} } { \lVert u \rVert_{L_2} }\!\!$} & DN-PINNs & \makecell{24.16 \\ \!($\pm$\! 48.04)\!} & \makecell{6.70 \\ \!($\pm$\! 12.10)\!} &  \makecell{2.28 \\ \!($\pm$\! 2.70)\!} &  \makecell{0.91 \\ \!($\pm$\! 0.52)\!} &  \makecell{1.15 \\ \!($\pm$\! 0.49)\!}  \\ 
\cline{2-7}
& DNLA (PINNs) &  \makecell{9.16 \\ \!($\pm$\! 2.34)\!} &  \makecell{2.81 \\ \!($\pm$\! 0.56)\!} &  \makecell{1.17 \\ \!($\pm$\! 0.71)\!} &  \makecell{0.08 \\ \!($\pm$\! 0.07)\!} &  - \\ 
\cline{2-7}
& \!\!\! DNLA (deep Ritz)\! &  \makecell{9.50 \\ \!($\pm$\! 2.46)\!} &  \makecell{1.83 \\ \!($\pm$\! 0.86)\!} &  \makecell{0.41 \\ \!($\pm$\! 0.35)\!} &  \makecell{ 0.62 \\ \!($\pm$\! 0.55)\!} &  \makecell{0.40 \\ \!($\pm$\! 0.49)\!} \\ 
\hline		                                                     
\end{tabular}
\label{Experiments-DNLA-ex1-Err-Table}
\vspace{-0.2cm}
\end{table}

We first conduct experiments using the DN-PINNs approach, i.e., PINNs \cite{raissi2019physics,karniadakis2021physics} are used as the numerical solver for both Dirichlet and Neumann subproblems. The iterative solutions over the entire domain in a typical simulation are depicted in \autoref{Experiments-DNLA-ex1-DN-PINNs}, with the initial guess for the interface value data given by
\begin{equation*}
	h^{[0]}(x,y) = \left(2\pi \cos(2\pi x)+ \sin(2\pi x)\right) (\cos(2\pi y) - 1) - 50xy(x-1)(y-1)\ \ \ \text{on}\ \Gamma,
\end{equation*}
which remains unchanged for other methods tested below. As the trained networks tend to provide erroneous Neumann trace on the interface even when the training loss is very small (see Remark \ref{Remark-DtN-Map} or \cite{dockhorn2019discussion,bajaj2021robust}), DN-PINNs fails to converge to the correct solution of \eqref{Experiments-DNLA-ex1} as shown in \autoref{Experiments-DNLA-ex1-DN-PINNs}. Such an inaccurate flux prediction would hamper the convergence of outer iterations but is perhaps inevitable in practice for problems with complex interface conditions. In fact, a straightforward replacement of the numerical solver by other learning strategies, e.g., the deep Ritz method \cite{yu2018deep}, also suffer from the same issue.  

\begin{figure}[t!]
\centering
\begin{subfigure}[htp]{\textwidth}
\centering
\includegraphics[width=0.192\textwidth]{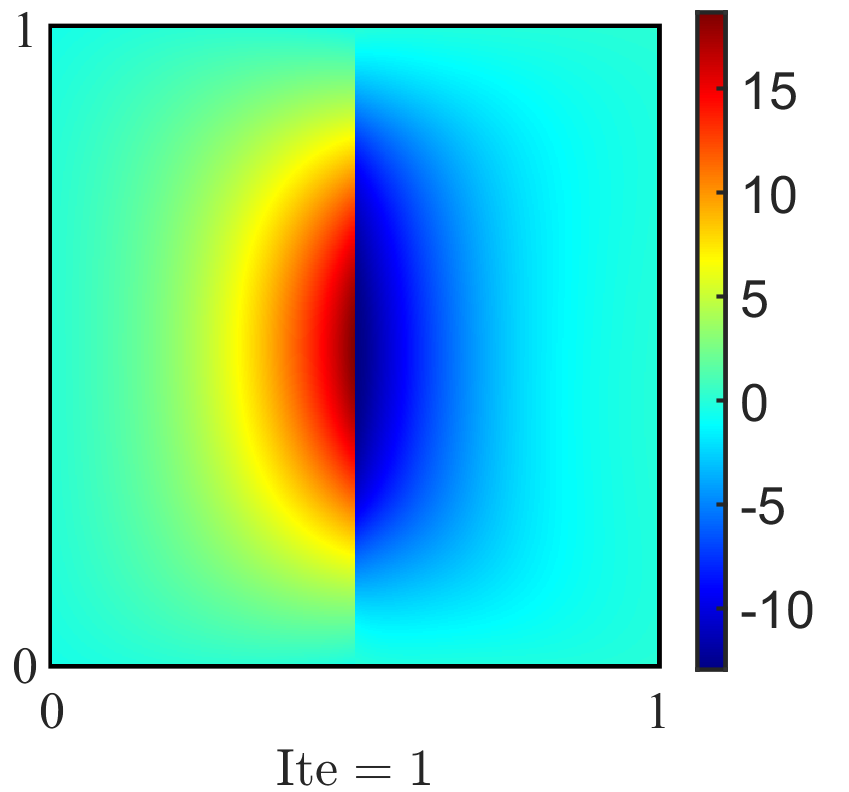}
\includegraphics[width=0.192\textwidth]{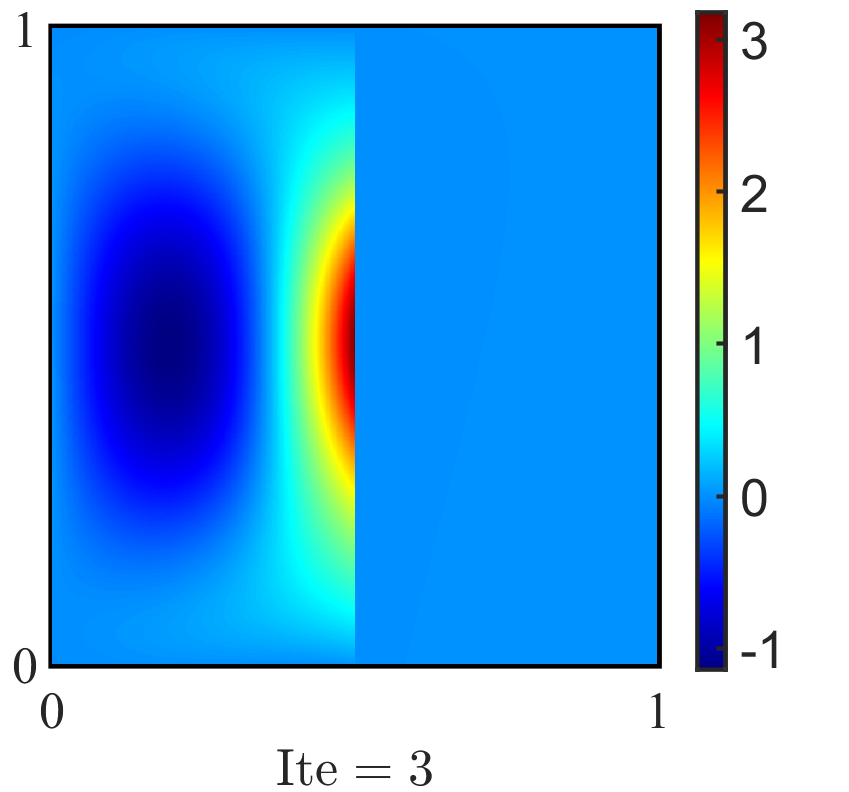}
\includegraphics[width=0.192\textwidth]{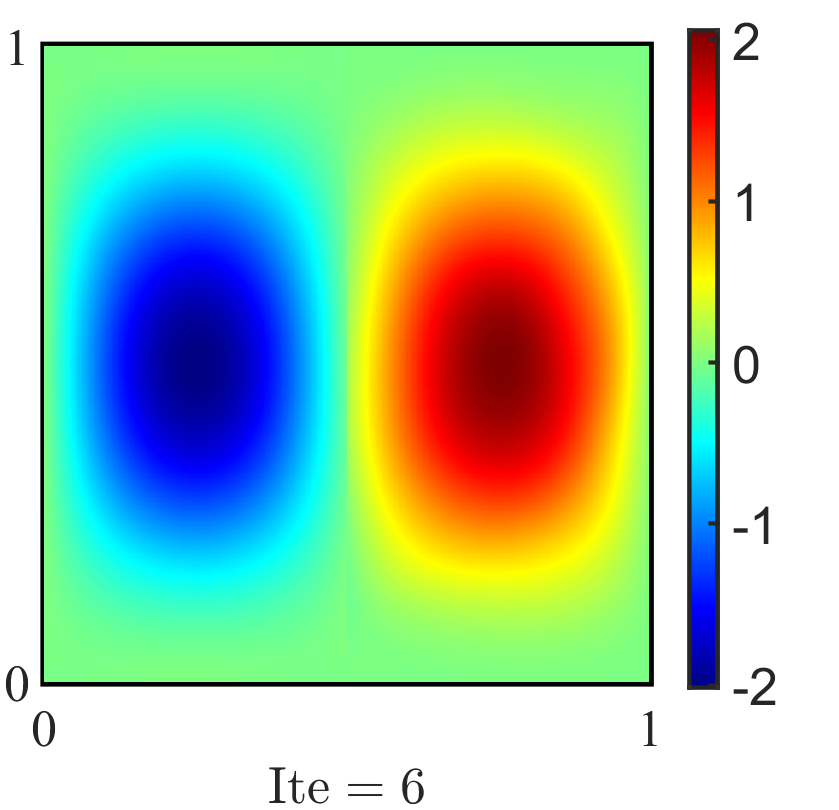}
\includegraphics[width=0.192\textwidth]{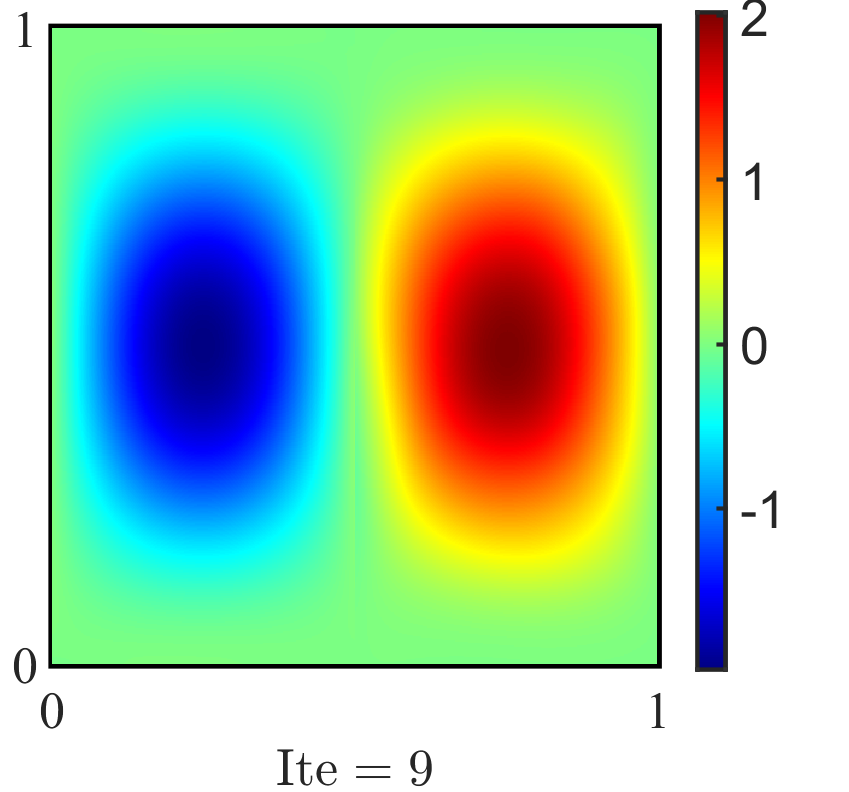}
\vspace{-0.1cm}
\caption{Iterative solutions $\hat{u}^{[k]}(x,y)$ along the outer iteration. }
\label{Experiments-DNLA-ex1-DNLA-PINN-solution}
\vspace{-0.2cm}
\end{subfigure}
\begin{subfigure}[htp]{\textwidth}
\centering
\includegraphics[width=0.192\textwidth]{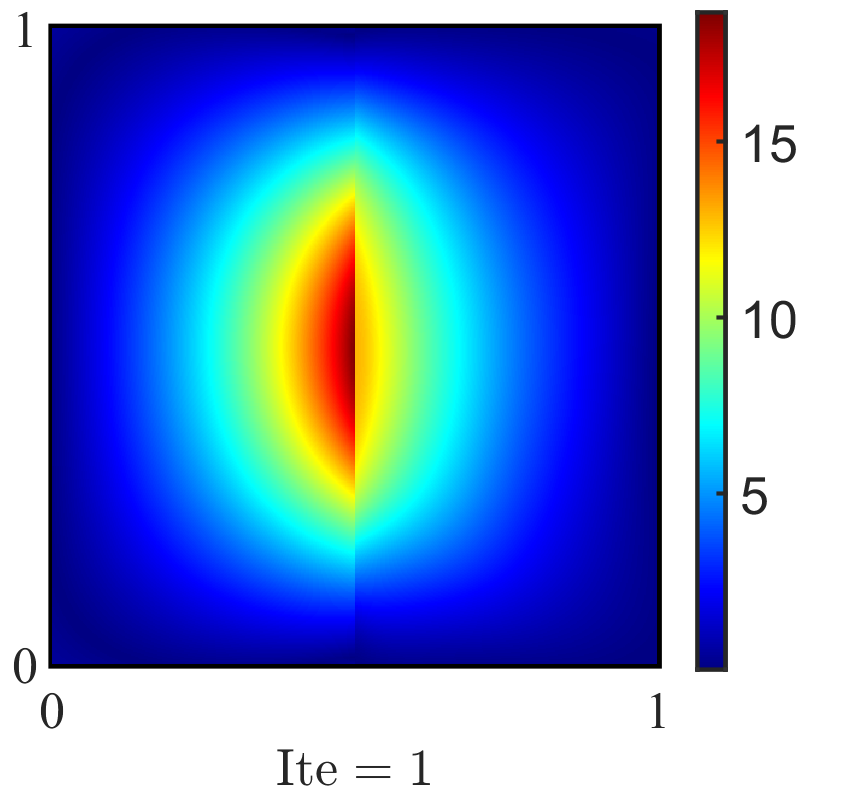}
\includegraphics[width=0.192\textwidth]{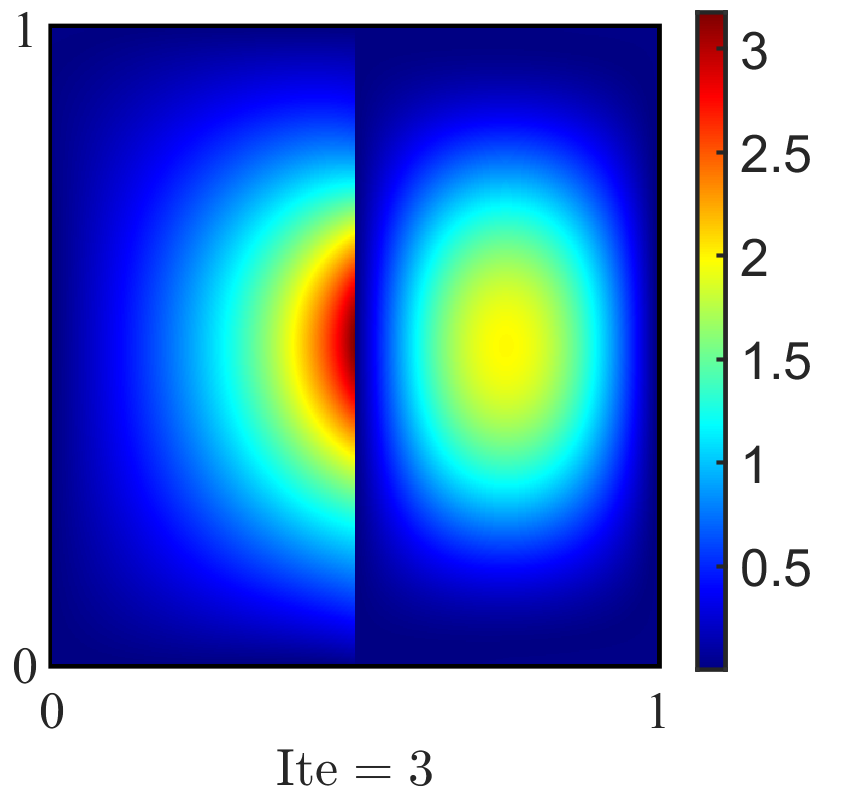}
\includegraphics[width=0.192\textwidth]{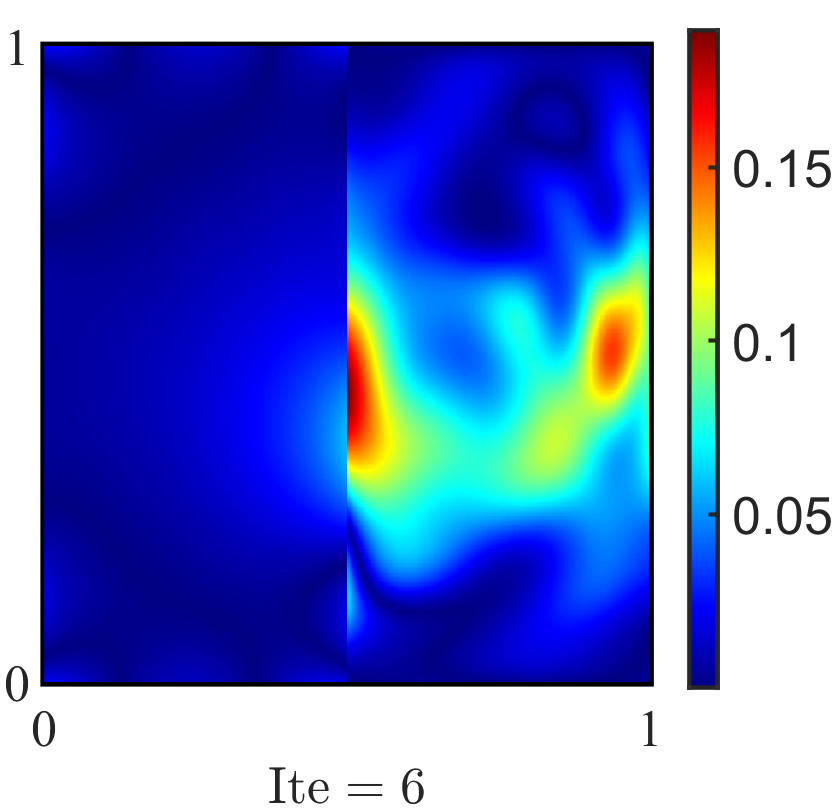}
\includegraphics[width=0.192\textwidth]{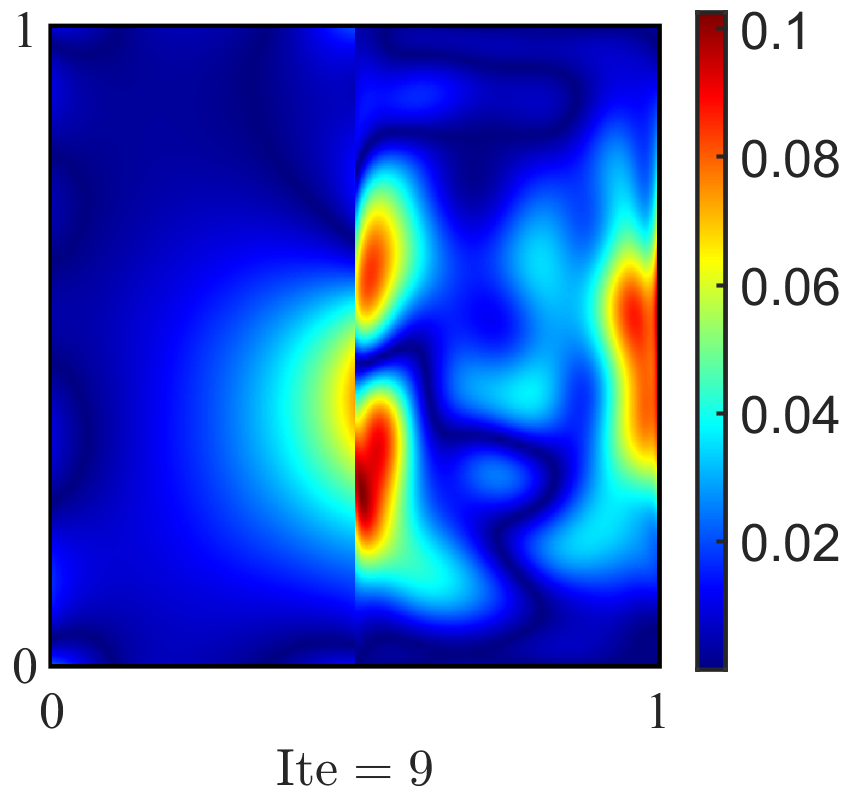}
\vspace{-0.1cm}
\caption{Error profiles $|\hat{u}^{[k]}(x,y) - u(x,y)|$ along the outer iteration. }
\label{Experiments-DNLA-ex1-DNLA-PINN-error}
\end{subfigure}
\vspace{-0.5cm}
\caption{Numerical results of \eqref{Experiments-DNLA-ex1} using DNLA (PINNs) on the test dataset.}
\label{Experiments-DNLA-ex1-DNLA-PINN}
\vspace{-0.4cm}
\end{figure}

\begin{figure}[t!]
\centering
\begin{subfigure}[htp]{\textwidth}
\centering
\includegraphics[width=0.192\textwidth]{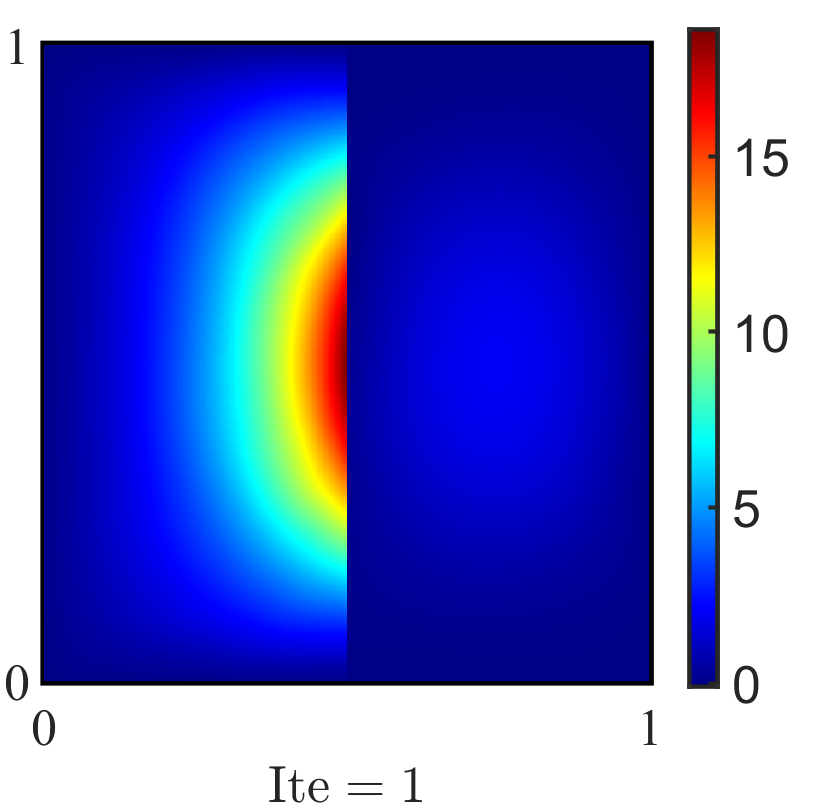}
\includegraphics[width=0.192\textwidth]{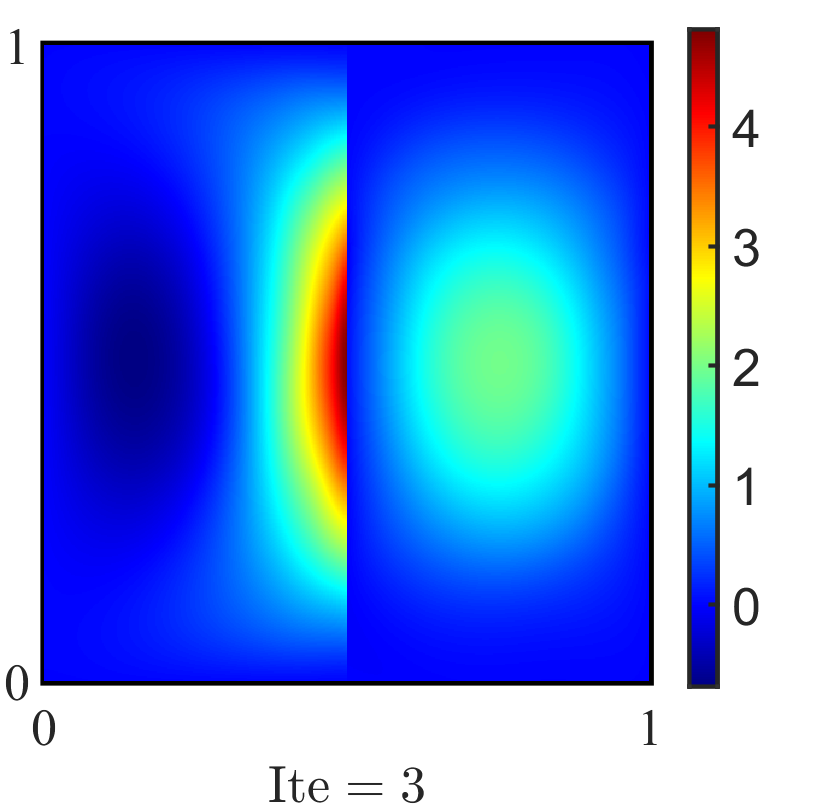}
\includegraphics[width=0.192\textwidth]{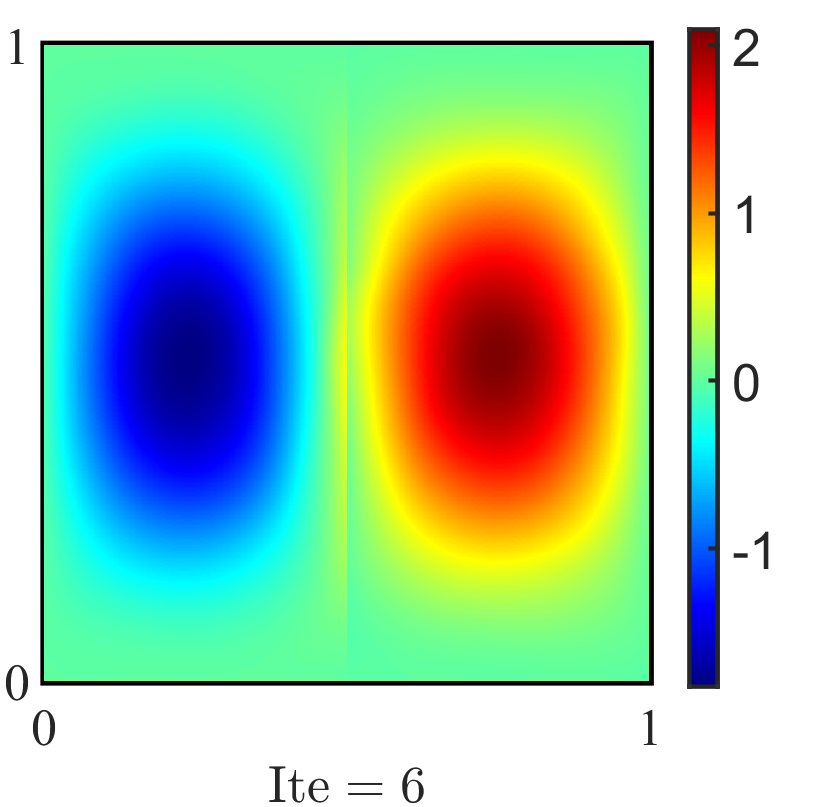}
\includegraphics[width=0.192\textwidth]{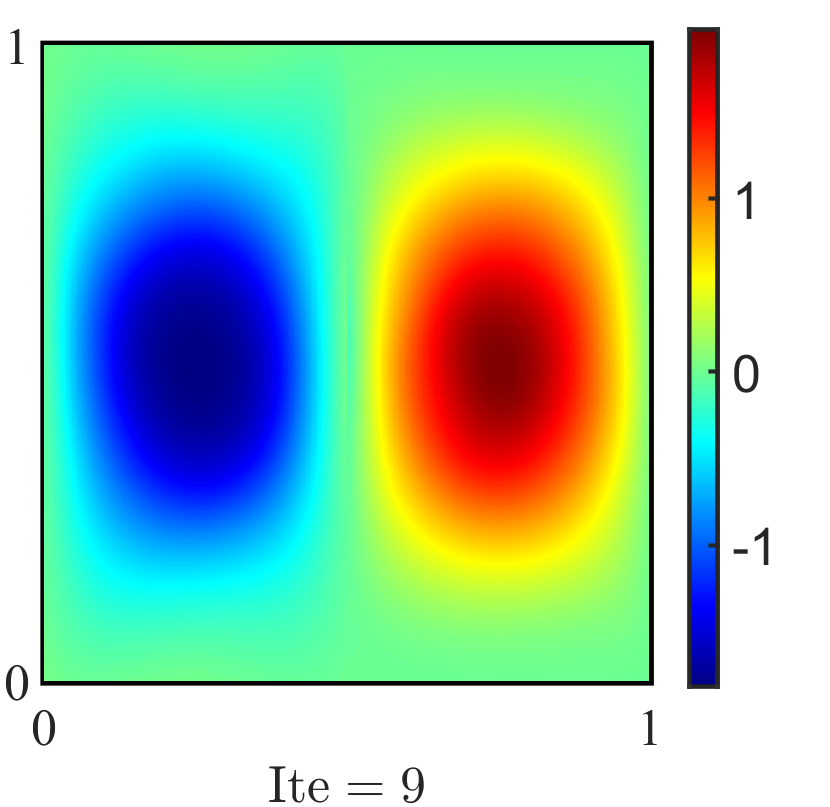}
\vspace{-0.1cm}
\caption{Iterative solutions $\hat{u}^{[k]}(x,y)$ along the outer iteration. }
\label{Experiments-DNLA-ex1-DNLA-DeepRitz-solution}
\vspace{-0.2cm}
\end{subfigure}
\begin{subfigure}[htp]{\textwidth}
\centering
\includegraphics[width=0.192\textwidth]{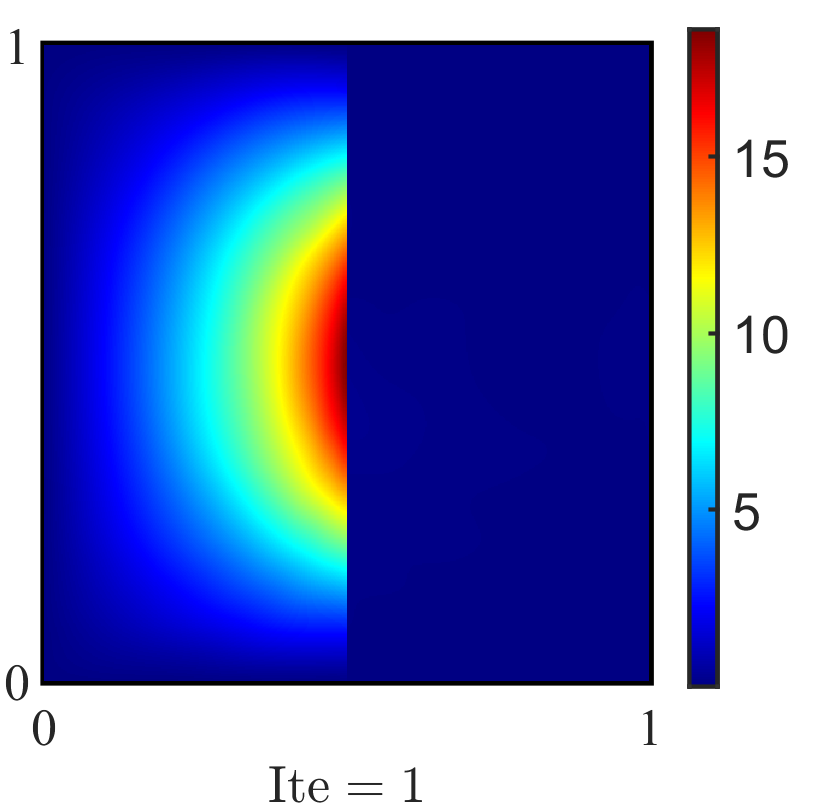}
\includegraphics[width=0.192\textwidth]{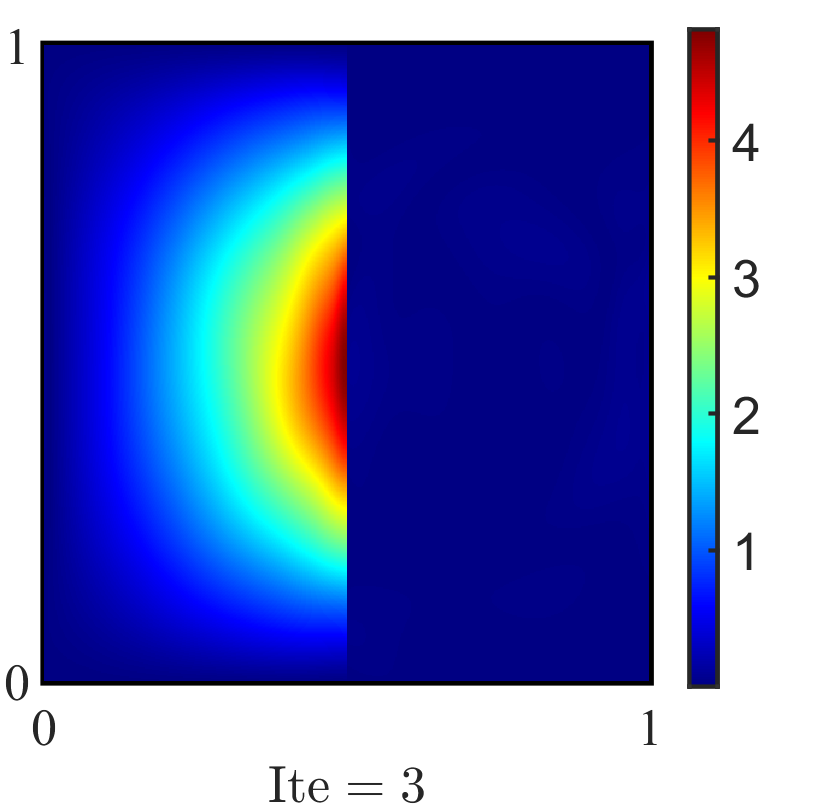}
\includegraphics[width=0.192\textwidth]{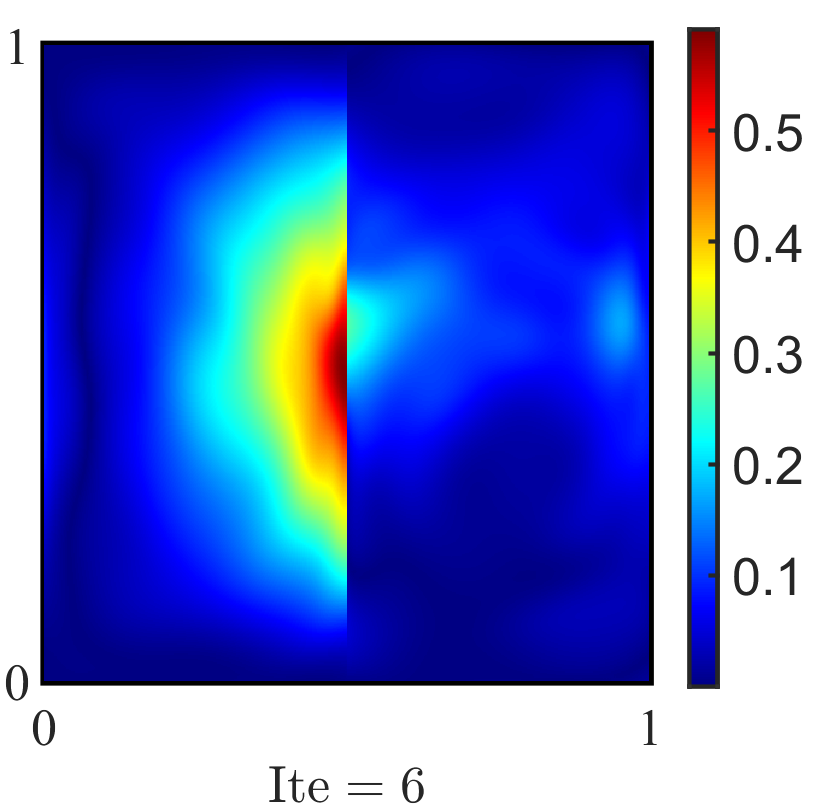}
\includegraphics[width=0.192\textwidth]{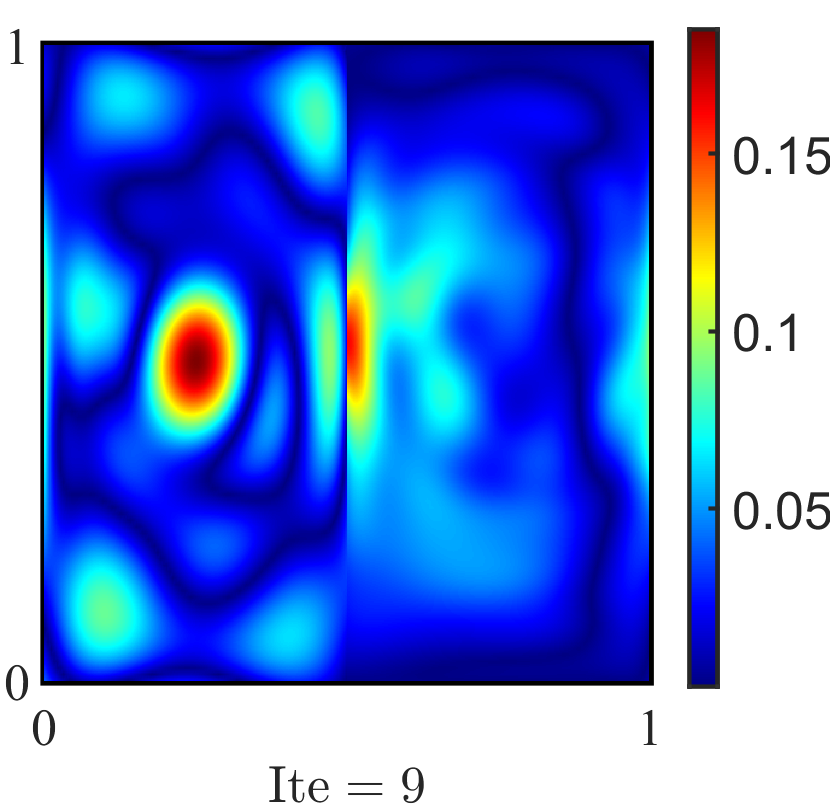}
\vspace{-0.1cm}
\caption{Error profiles $|\hat{u}^{[k]}(x,y) - u(x,y)|$ along the outer iteration. }
\label{Experiments-DNLA-ex1-DNLA-DeepRitz-error}
\end{subfigure}
\vspace{-0.5cm}
\caption{Numerical results of \eqref{Experiments-DNLA-ex1} using DNLA (deep Ritz) on the test dataset.}
\label{Experiments-DNLA-ex1-DNLA-DeepRitz}
\vspace{-0.5cm}
\end{figure}

\begin{figure}[t!]
    \centering
    \begin{subfigure}[t]{0.48\textwidth}
        \centering
        \includegraphics[width=0.3\textwidth]{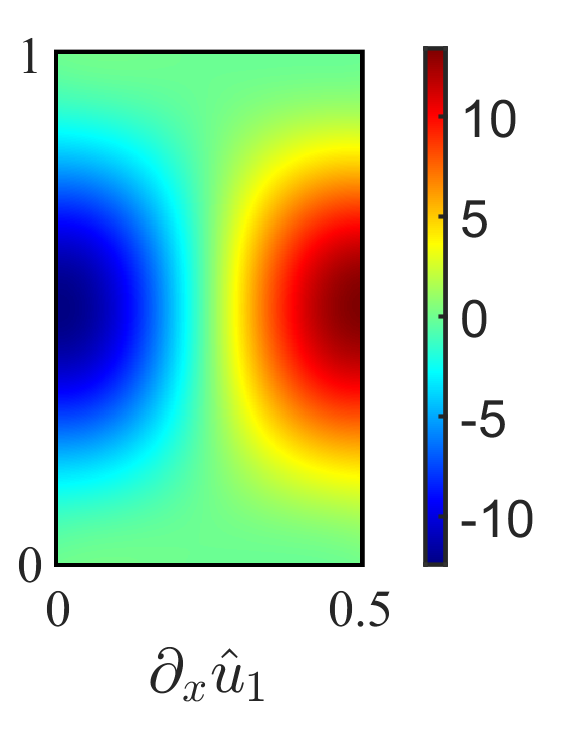}
	   \includegraphics[width=0.3\textwidth]{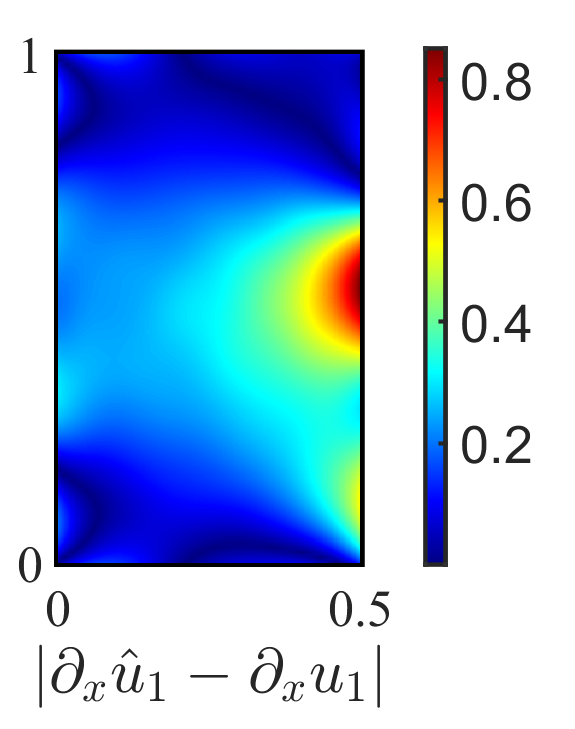}
	   \vspace{-0.1cm}
        \caption{DNLA (PINN)}
    \end{subfigure}%
    ~ 
    \begin{subfigure}[t]{0.48\textwidth}
        \centering
        \includegraphics[width=0.3\textwidth]{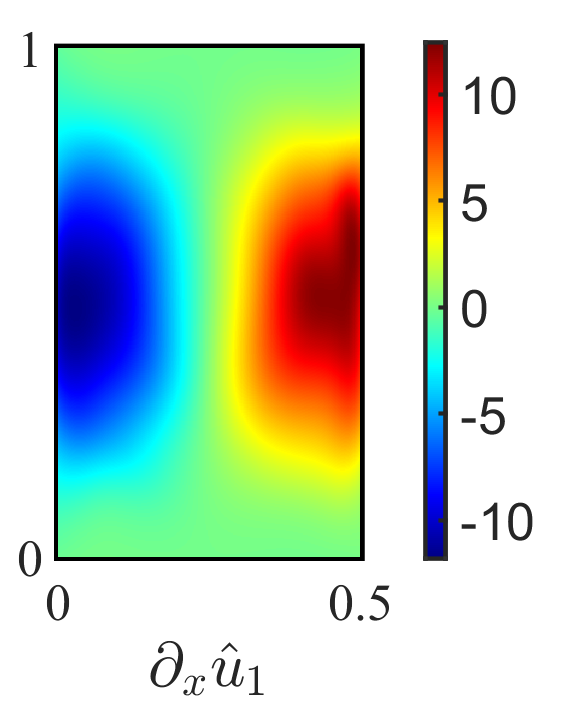}
	   \includegraphics[width=0.3\textwidth]{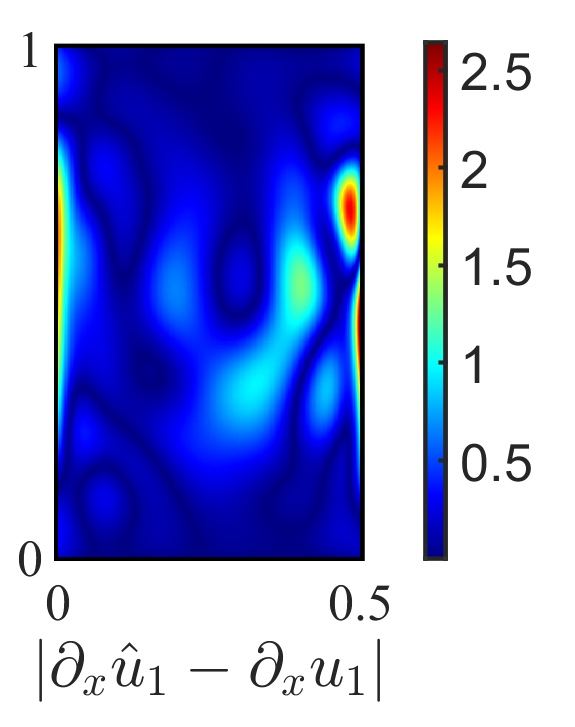}
	   \vspace{-0.1cm}
        \caption{DNLA (deep Ritz)}
    \end{subfigure}
\vspace{-0.5cm}
\caption{Erroneous Dirichlet-to-Neumann maps for \eqref{Experiments-DNLA-ex1}: $\partial_x\hat{u}_1^{[9]}$ and $|\partial_x (\hat{u}_1^{[9]} - u_1)|$.}
\label{Experiments-DNLA-ex1-Overfit-Dirichlet-Subproblem}
\vspace{-0.8cm}
\end{figure}

In contrast, although the Dirichlet-to-Neumann map through the trained solution of Dirichlet subproblem is usually of unacceptable low accuracy, our proposed method doesn't need to explicitly enforce the flux continuity along subdomain interfaces, thereby enabling the effectiveness of convergence in the presence of erroneous interface conditions (see \autoref{Experiments-DNLA-ex1-Overfit-Dirichlet-Subproblem}). To validate our statements, we show in \autoref{Experiments-DNLA-ex1-DNLA-PINN} and \autoref{Experiments-DNLA-ex1-DNLA-DeepRitz} the numerical results using Algorithm \ref{Algorithm-DN-Learning-2Subdomains}, where PINNs \cite{raissi2019physics} and deep Ritz methods \cite{yu2018deep} are employed for solving the Dirichlet subproblem, respectively. 

As can be observed from \autoref{Experiments-DNLA-ex1-DNLA-PINN} and \autoref{Experiments-DNLA-ex1-DNLA-DeepRitz}, the predicted solution using our proposed learning algorithms is in agreement with the true solution, while the Neumann traces shown in \autoref{Experiments-DNLA-ex1-Overfit-Dirichlet-Subproblem} indicate that the network solution of Dirichlet subproblem learns to fit the given Dirichlet boundary condition with erroneous Neumann traces. More quantitatively, we run the simulations for 5 times to calculate the relative-$L_2$ errors, and the results (mean value $\pm$ standard deviation) are reported in \autoref{Experiments-DNLA-ex1-Err-Table}. By employing our proposed compensated deep Ritz method for solving the Neumann subproblem, it can be observed that our learning algorithms work reasonably well, while the DN-PINNs is typically divergent due to the lack of accurate flux transmission across the interface. Moreover, as the solution of \eqref{Experiments-DNLA-ex1} is rather smooth on each subdomain, it can be found in \autoref{Experiments-DNLA-ex1-Err-Table} that DNLA (PINNs) performs better than DNLA (deep Ritz). This is because that second-order derivatives are explicitly involved during the training process, leading to better estimates of the solution's gradient inside the subdomain (see \autoref{Experiments-DNLA-ex1-Overfit-Dirichlet-Subproblem}).

Moreover, by employing the DNLA (PINNs) for solving \eqref{Experiments-DNLA-ex1}, we report in \autoref{Experiments-DNLA-ex1-NetArch-Table} the relative-$L_2$ error and the corresponding number of outer iterations for different architectures, which indicates that the number of outer iterations required to achieve a comparable accuracy remain approximately constant as the width and depth of the network vary across a certain range of values. When the depth goes further deeper, such an observation may no longer be valid due to the vanishing gradient problem.

\begin{table}[t!]
\small
\caption{Relative-$L_2$ errors of DNLA (PINNs) with the number of outer iterations in bracket for different network architectures for example \eqref{Experiments-DNLA-ex1}.}
\vspace{-0.1cm}
\centering
\begin{tabular}{ | c || c | c | c | c | c |    }
\hline
\multicolumn{2}{|c|}{ \diagbox[width=8em]{Depth}{Width} } & 30  & 40 & 50 & 60 \\
\hline
\multicolumn{2}{|c|}{6} & 0.0644 (9) & 0.0316 (9) & 0.0580 (9) & 0.0592 (8) \\
\hline
\multicolumn{2}{|c|}{8} & 0.0482 (9) & 0.0860 (9) & 0.0449 (9) & 0.0490 (9) \\
\hline 
\end{tabular}
\label{Experiments-DNLA-ex1-NetArch-Table}
\vspace{-0.3cm}
\end{table}

\subsubsection{Poisson's Equation with Zigzag Interface}
To demonstrate the advantage of mesh-free property over traditional mesh-based numerical methods \cite{toselli2004domain}, we consider the previous example but with a more complex interface geometry,
\begin{equation}
\begin{array}{cl}
-\Delta u(x,y)  = 4\pi^2 \sin(2 \pi y)  (2 \cos(2 \pi x) - 1)  \ & \text{in}\ \Omega=(0,1)^2,\\
u(x,y) = 0\ \ & \text{on}\ \partial \Omega,
\end{array}
\label{Experiments-DNLA-ex2}
\end{equation}
where the exact solution $u(x,y) = \sin(2\pi y)(\cos(2\pi x)-1)$ and the interface is a curved zigzag line as depicted in \autoref{fig-domain-decomposition}. More precisely, the zigzag function reads
\begin{equation*}
	x = c ( a(20y-\text{floor}(20y)) + b ) + 0.5
\end{equation*}
where coefficients $ a = 0.05 (-1 + 2\times \text{mod}(\text{floor}(20y), 2))$, $b=-0.05\times \text{mod}(\text{floor}(20x),2)$ and $c=-2\times\text{mod}(\text{floor}(10x),2)+1$, therefore enabling the sample generation inside each subdomain and at its boundary. Our proposed learning algorithm can easily handle such irregular boundary shapes, while the finite difference or finite element methods \cite{brenner2008mathematical} requires careful treatment of edges and corners.

\begin{figure}[b!]
\centering
\vspace{-0.2cm}
\begin{subfigure}[htp]{\textwidth}
\centering
\includegraphics[width=0.192\textwidth]{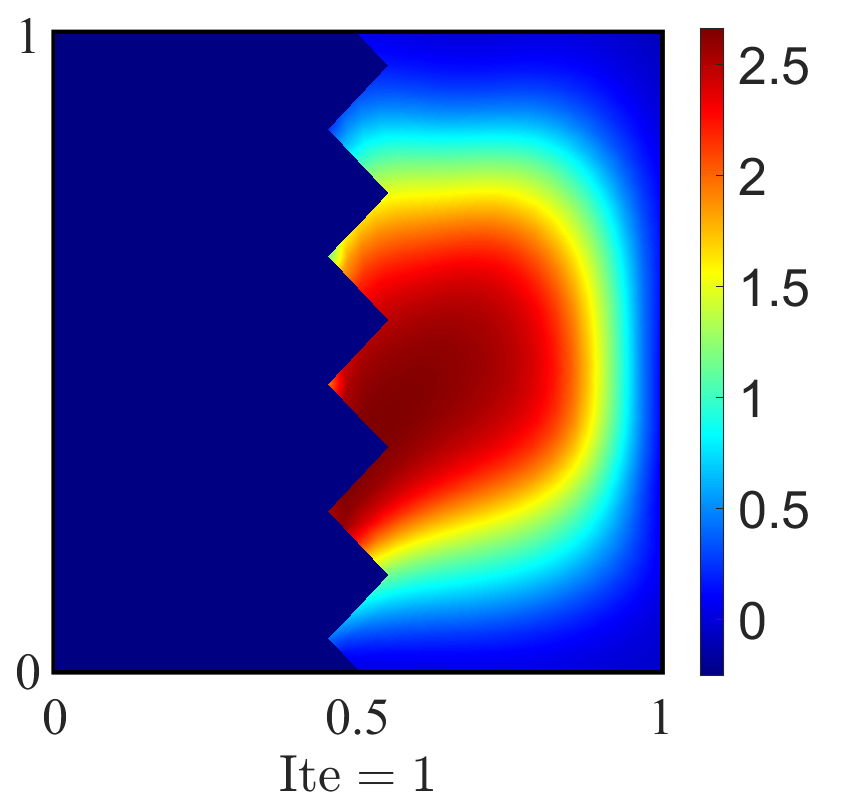}
\includegraphics[width=0.192\textwidth]{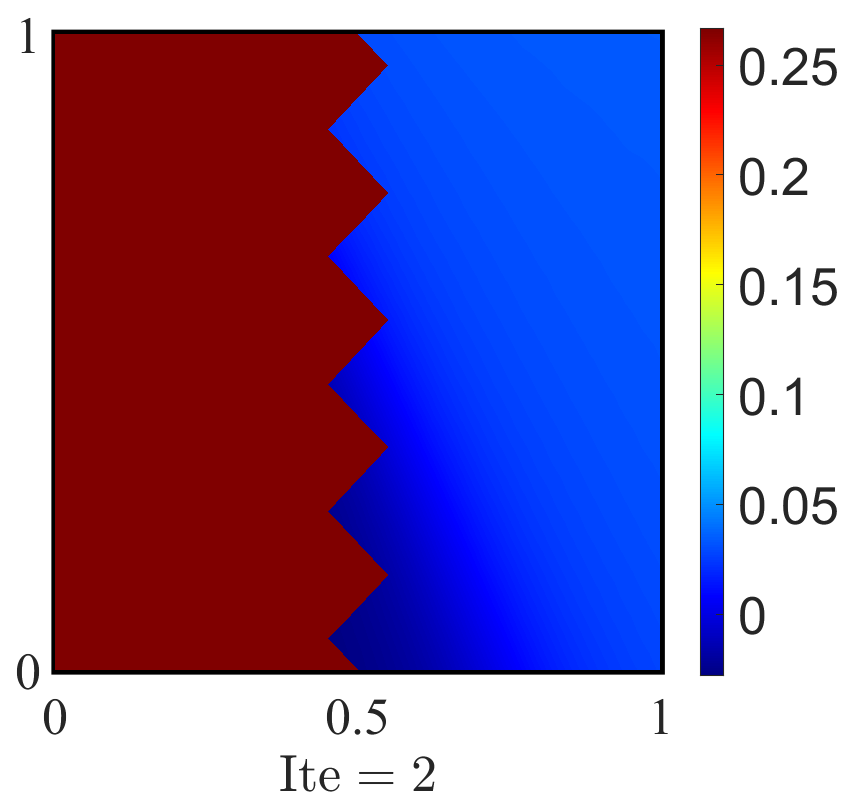}
\includegraphics[width=0.192\textwidth]{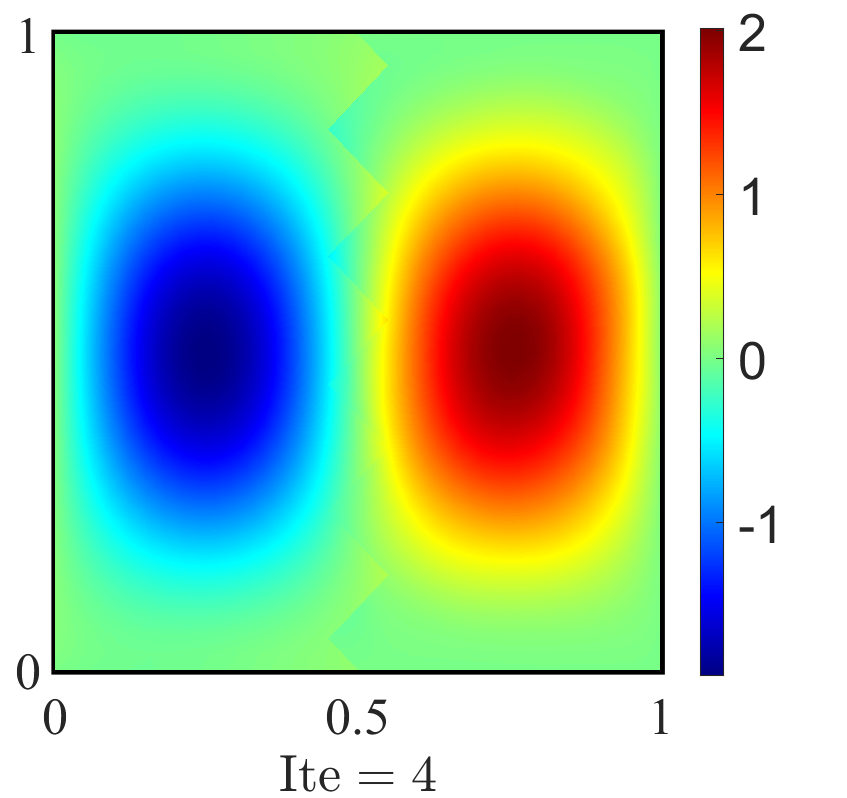}
\includegraphics[width=0.192\textwidth]{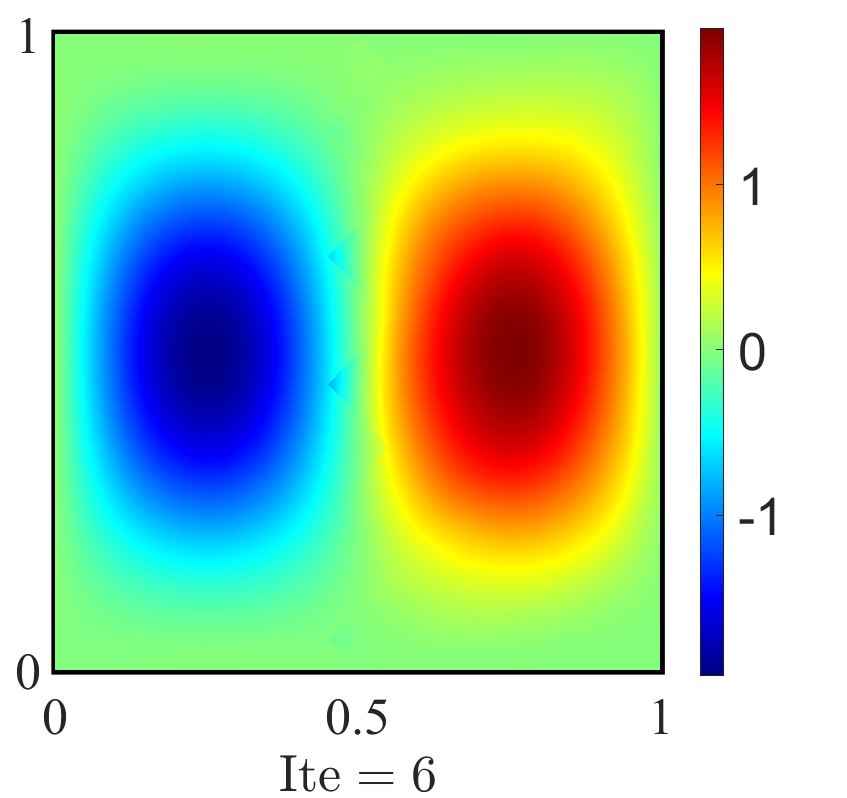}
\vspace{-0.1cm}
\caption{Iterative solutions $\hat{u}^{[k]}(x,y)$ using DNLA (PINNs) along the outer iteration.}
\label{Experiments-DNLA-ex2-DNLA-PINN-solution}
\vspace{-0.2cm}
\end{subfigure}
\begin{subfigure}[htp]{\textwidth}
\centering
\includegraphics[width=0.192\textwidth]{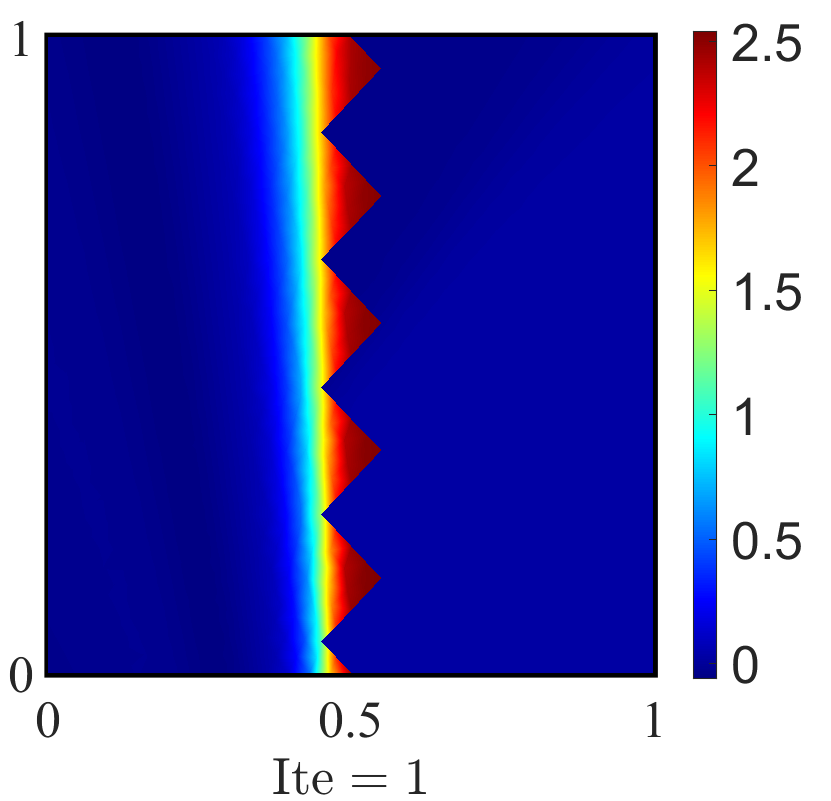}
\includegraphics[width=0.192\textwidth]{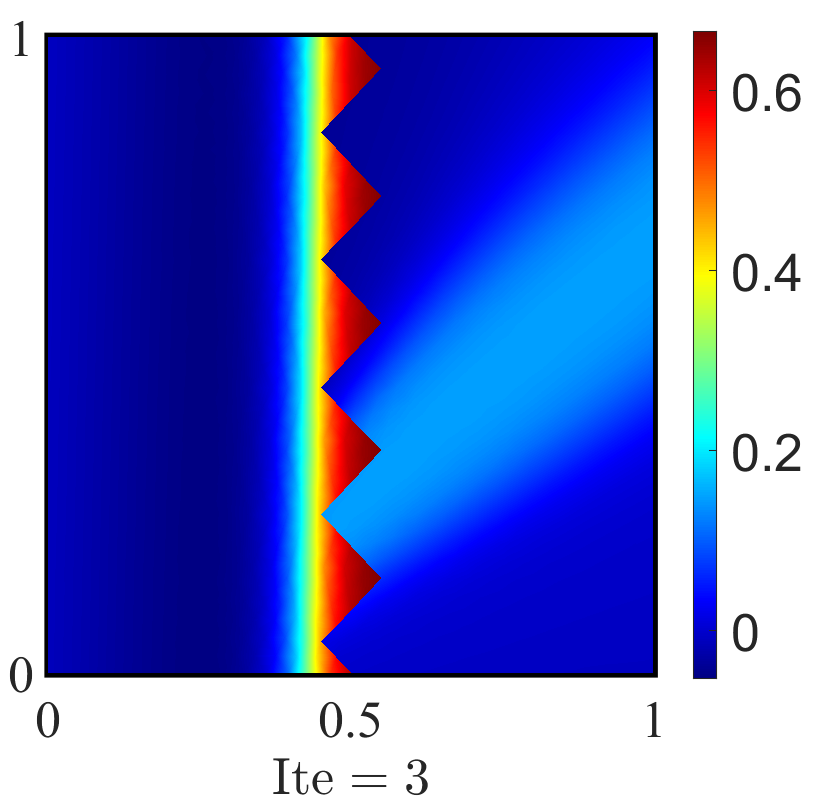}
\includegraphics[width=0.192\textwidth]{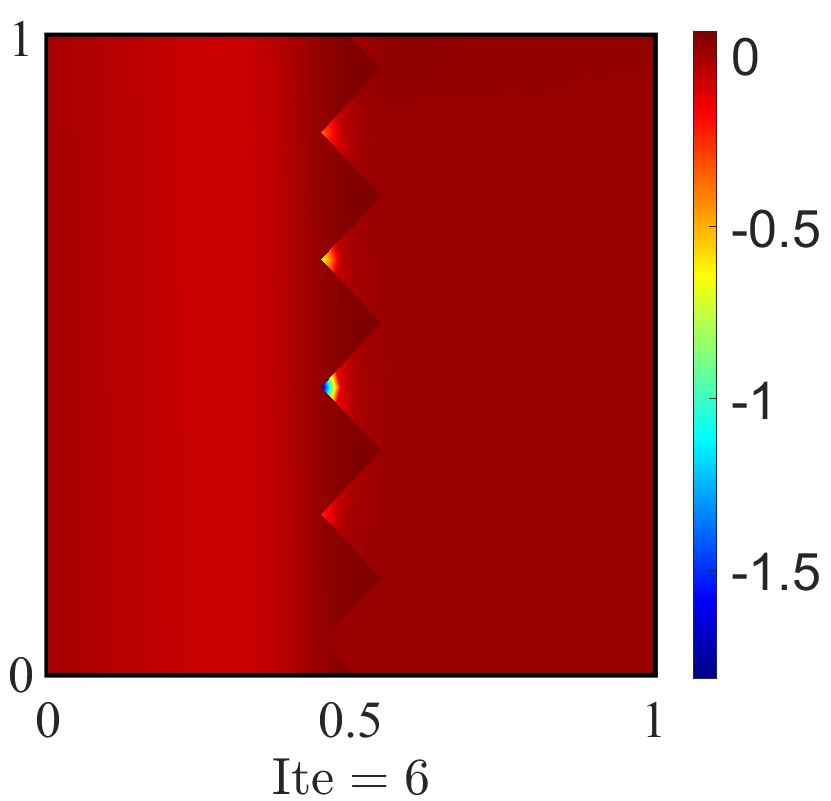}
\includegraphics[width=0.192\textwidth]{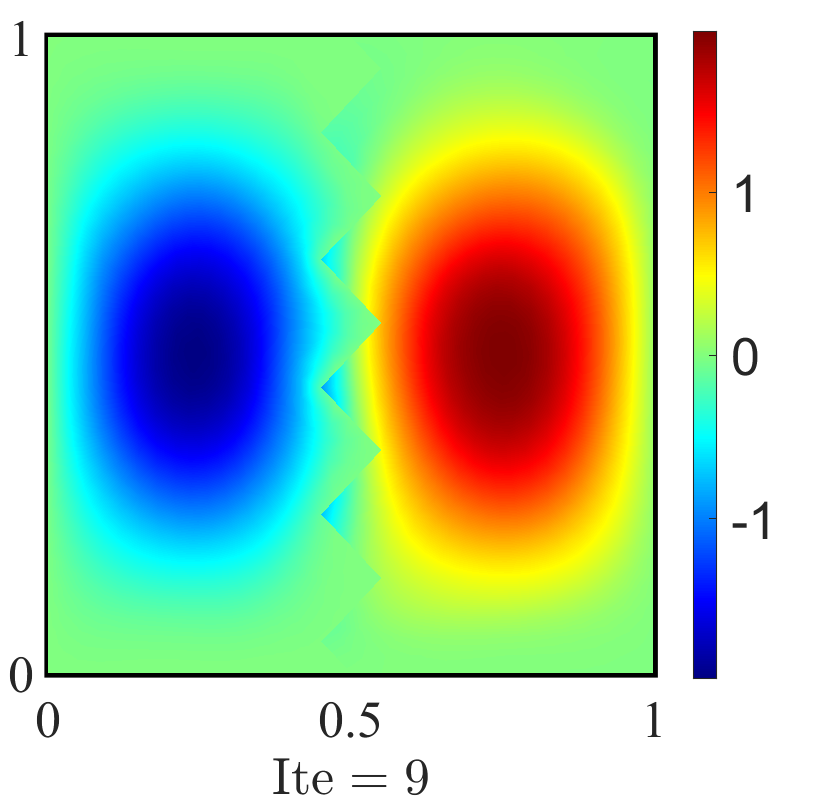}
\vspace{-0.1cm}
\caption{Iterative solutions $\hat{u}^{[k]}(x,y)$ using DNLA (deep Ritz) along the outer iteration.}
\label{Experiments-DNLA-ex2-DNLA-PINN-error}
\end{subfigure}
\vspace{-0.5cm}
\caption{Numerical results of \eqref{Experiments-DNLA-ex2} using DNLA on the test dataset.}
\vspace{-0.55cm}
\label{Experiments-DNLA-ex2-DNLA-PINN}
\end{figure}

\begin{figure}[t!]
\centering
\begin{subfigure}[htp]{\textwidth}
\centering
\includegraphics[width=0.14\textwidth]{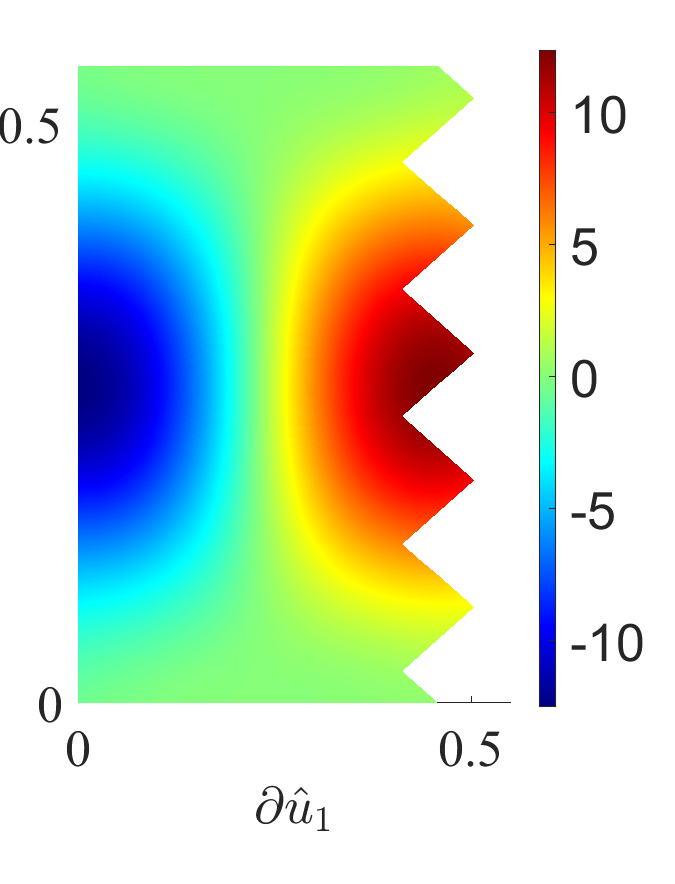}
\includegraphics[width=0.14\textwidth]{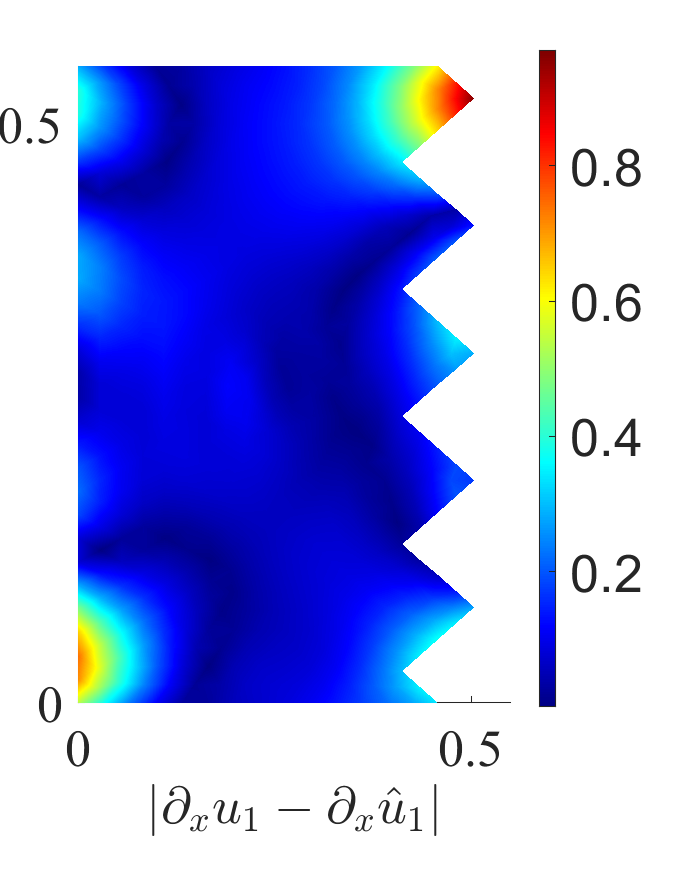}
\includegraphics[width=0.14\textwidth]{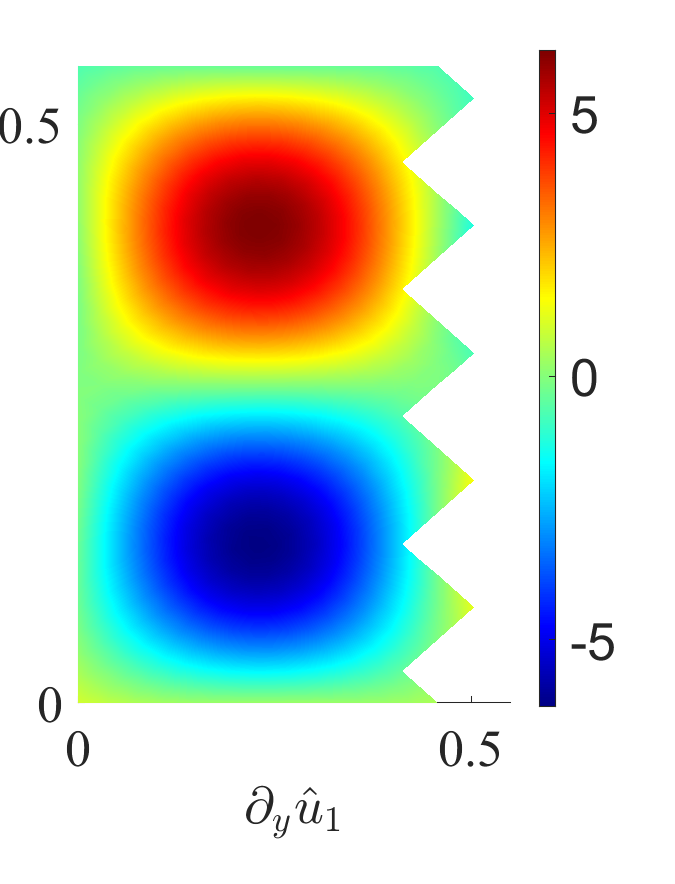}
\includegraphics[width=0.14\textwidth]{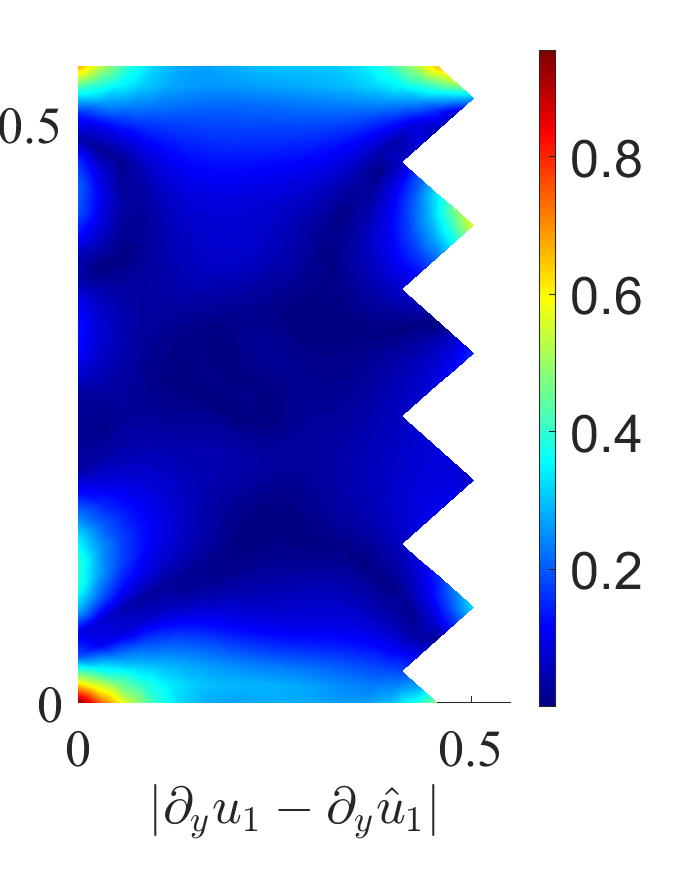}
\vspace{-0.15cm}
\caption{$\partial_x\hat{u}_1^{[4]}$, $\partial_y\hat{u}_1^{[4]}$ and error profiles $|\partial_x (\hat{u}_1^{[4]} - u_1)|$, $|\partial_y (\hat{u}_1^{[4]} - u_1)|$ using DNLA (PINNs). }
\vspace{-0.15cm}
\end{subfigure}
\begin{subfigure}[htp]{\textwidth}
\centering
\includegraphics[width=0.14\textwidth]{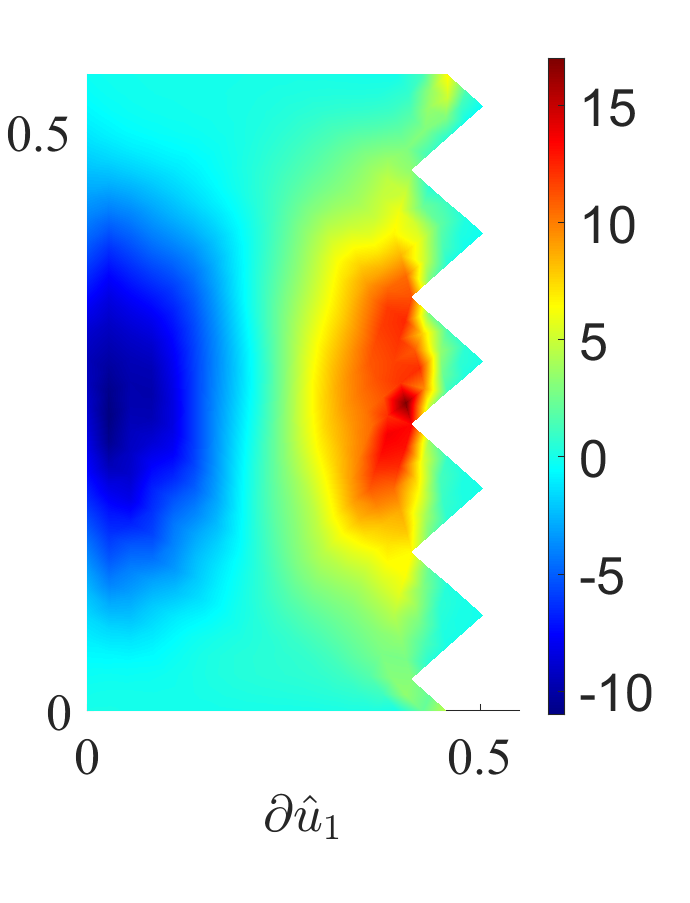}
\includegraphics[width=0.14\textwidth]{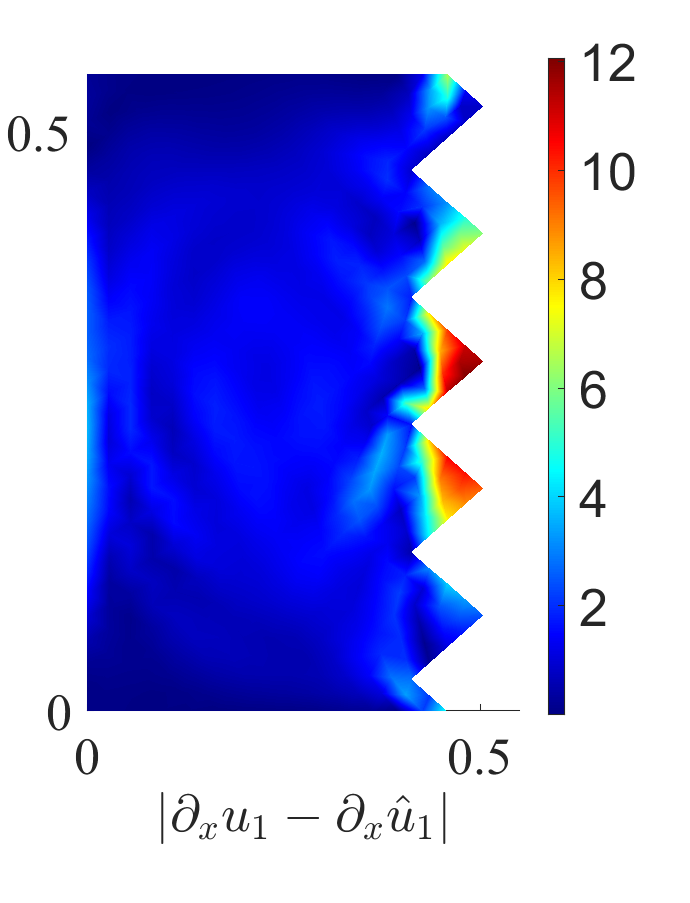}
\includegraphics[width=0.14\textwidth]{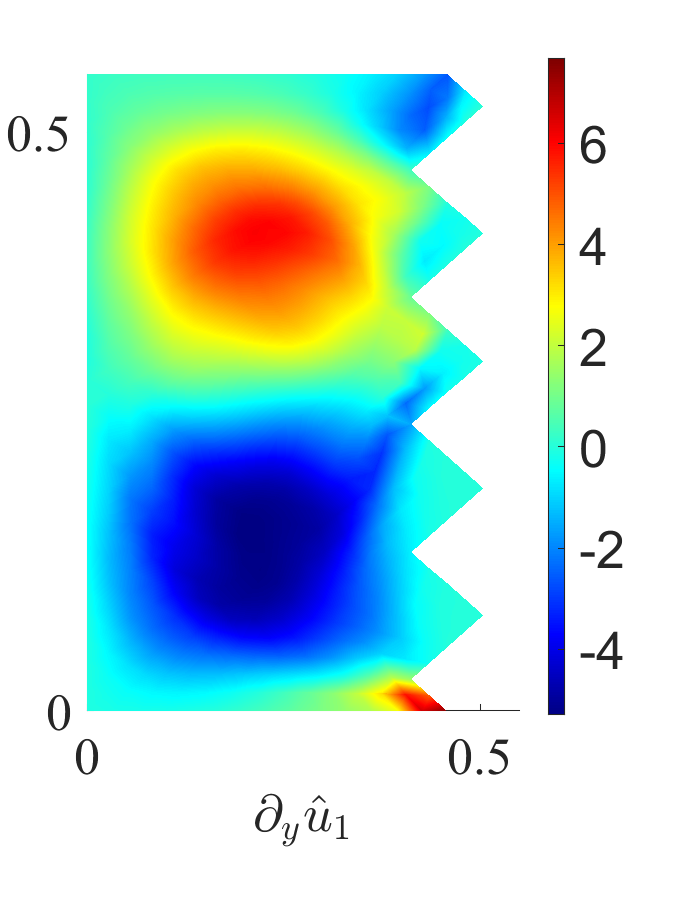}
\includegraphics[width=0.14\textwidth]{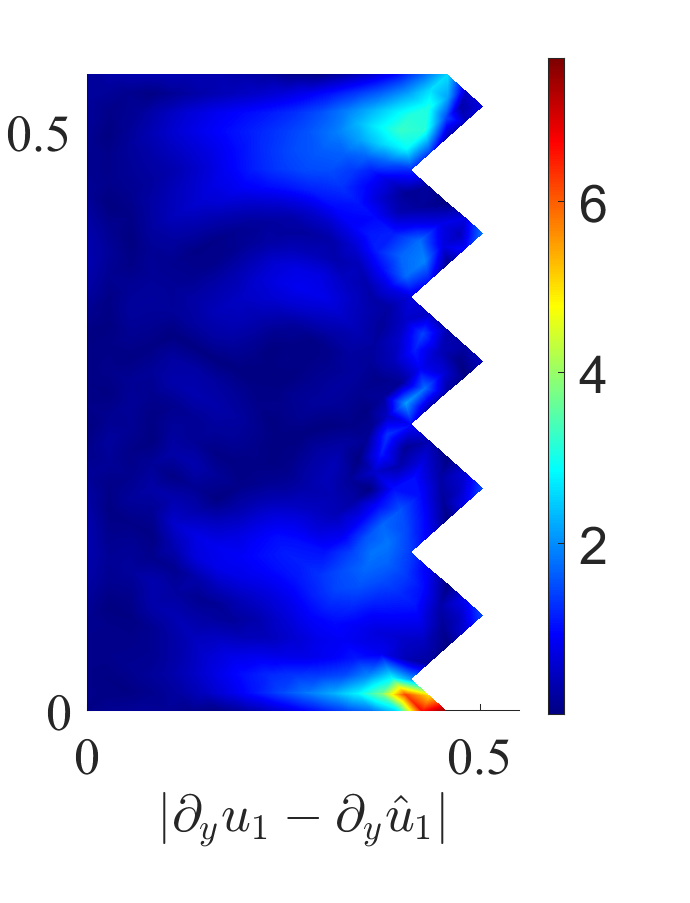}
\vspace{-0.15cm}
\caption{$\partial_x\hat{u}_1^{[9]}$, $\partial_y\hat{u}_1^{[9]}$ and error profiles $|\partial_x (\hat{u}_1^{[9]} - u_1)|$, $|\partial_y (\hat{u}_1^{[9]} - u_1)|$ using DNLA (deep Ritz). }
\end{subfigure}
\vspace{-0.5cm}
\caption{Erroneous Dirichlet-to-Neumann maps for \eqref{Experiments-DNLA-ex2}.}
\label{Experiments-DNLA-ex2-Overfit-Dirichlet-Subproblem}
\vspace{-0.75cm}
\end{figure}

\begin{table}[t!]
\vspace{-0.05cm}
\small
\caption{Relative-$L_2$ errors of the network solution along the outer iteration for \eqref{Experiments-DNLA-ex2}, with mean value ($\pm$ standard deviation) being reported over 5 independent runs.}
\vspace{-0.1cm}
\centering
\renewcommand{\arraystretch}{1.1}
\begin{tabular}{ | c || c | c | c | c | c | c |  }
\hline
\multicolumn{2}{|c|}{ \diagbox[width=16em]{Relative Errors}{Outer Iterations} } & 1  & 2 & 4 & 6 & 11  \\
\hline	
\hline
\multirow{5}{*}{$ \displaystyle \!\! \frac{ \lVert \hat{u}^{[k]} - u \rVert_{L_2} } { \lVert u \rVert_{L_2} }\!\!$} & DN-PINNs & \makecell{2.39 \\ \!($\pm$\! 1.46)\!} & \makecell{1.55 \\ \!($\pm$\! 0.93)\!} &  \makecell{1.47 \\ \!($\pm$\! 0.53)\!} &  \makecell{1.34 \\ \!($\pm$\! 0.60)\!} &  \makecell{1.34 \\ \!($\pm$\! 0.53)\!}  \\ 
\cline{2-7}
& DNLA (PINNs) &  \makecell{1.62 \\ \!($\pm$\! 0.48)\!} &  \makecell{0.98 \\ \!($\pm$\! 0.38)\!} &  \makecell{0.24 \\ \!($\pm$\! 0.13)\!} & \makecell{0.07 \\ \!($\pm$\! 0.02)\!} & - \\ 
\cline{2-7}
& \!\!\! DNLA (deep Ritz)\! &  \makecell{1.66 \\ \!($\pm$\! 0.26)\!} &  \makecell{1.52 \\ \!($\pm$\! 0.19)\!} &  \makecell{1.74 \\ \!($\pm$\! 0.38)\!} &  \makecell{1.37 \\ \!($\pm$\! 0.02)\!} &  \makecell{0.20 \\ \!($\pm$\! 0.05)\!} \\ 
\hline			                                                     
\end{tabular}
\label{Experiments-DNLA-ex2-Err-Table}
\vspace{-0.4cm}
\end{table}

Here, the initial guess $h^{[0]}(x,y) = \sin(2\pi x)(\cos(2\pi y)-1) - 1000\sin(2\pi y)^2\sin(2\pi x)$ is adopted for all methods tested below. As before, the DN-PINNs approach fails to converge to the true solution (see supplementary materials) as the network solution of Dirichlet subproblem is prone to return erroneous Neumann traces at interface. In contrast, by solving the Neumann subproblem through our compensated deep Ritz method, the numerical results in \autoref{Experiments-DNLA-ex2-DNLA-PINN} demonstrate that our DNLA (PINNs) can obtain a satisfactory approximation to the exact solution of \eqref{Experiments-DNLA-ex2}, which also avoids the meshing procedure that is often challenging for problems with complex interfaces. Importantly, DNLA (PINNs) remains effective in the presence of inaccurate flux estimations (see \autoref{Experiments-DNLA-ex2-Overfit-Dirichlet-Subproblem}), making it highly desirable in practice since the erroneous Dirichlet-to-Neumann map always occurs to some extent.

However, when the deep Ritz method \cite{yu2018deep} is used to solve the Dirichlet subproblem, the accuracy of approximate gradients within the subdomain $\Omega_1$ is no longer comparable to that of PINNs \cite{chen2020comparison}. This situation can become even worse for irregular domains (see \autoref{Experiments-DNLA-ex2-Overfit-Dirichlet-Subproblem} and \autoref{Experiments-DNLA-ex2-DNLA-PINN}). To further validate our claims, we present in \autoref{Experiments-DNLA-ex2-Err-Table} the quantitative results from 5 runs, which reveals that DNLA (PINNs) outperforms DN-PINNs and DNLA (deep Ritz) in terms of accuracy.

\subsubsection{Poisson's Equation with Four Subdomains}
Next, we consider the Poisson problem that is divided into four subproblems in two-dimension, i.e.,
\begin{equation}
\begin{array}{cl}
-\Delta u(x,y)  = f(x,y)\ & \text{in}\ \Omega=(0,1)^2, \\
u(x,y) = 0\ \ & \text{on}\ \partial \Omega,
\end{array}
\label{Experiments-DNLA-ex3}
\end{equation}
where $u(x,y) = \sin(2\pi x)\sin(8\pi y) $ and $f(x,y)=$ $65 \pi^2 \sin(2\pi x)\sin(8\pi y)$. Here, the domain is decomposed using the red-black partition \cite{toselli2004domain}, while the multidomains are categorized into two sets \cite{toselli2004domain} as depicted in \autoref{fig-domain-decomposition}. Then, the deep learning-based algorithms are deployed, with the initial guess at interface chosen as $h^{[0]}(x,y)=u(x,y)-100x(x-1)y(y-1)$ in what follows. Due to the high frequency of exact solution, the number of epochs here is $5k$ and the initial learning rate is $0.001$, which will decay at the $2k$-th and $4k$-th epoch. 

\begin{figure}[t!]
\centering
\begin{subfigure}[htp]{\textwidth}
\centering
\includegraphics[width=0.192\textwidth]{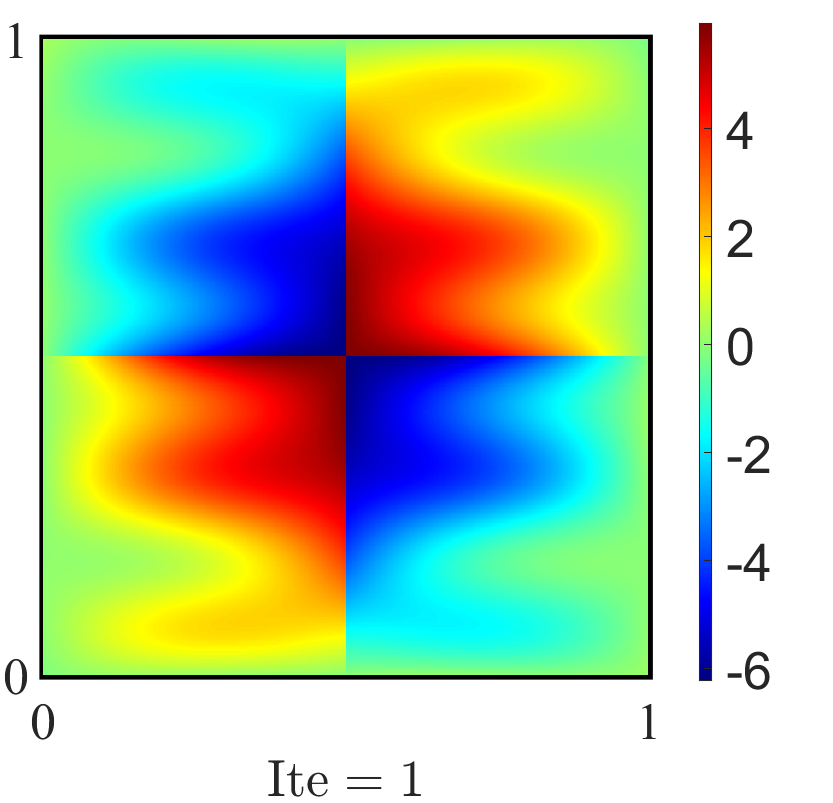}
\includegraphics[width=0.192\textwidth]{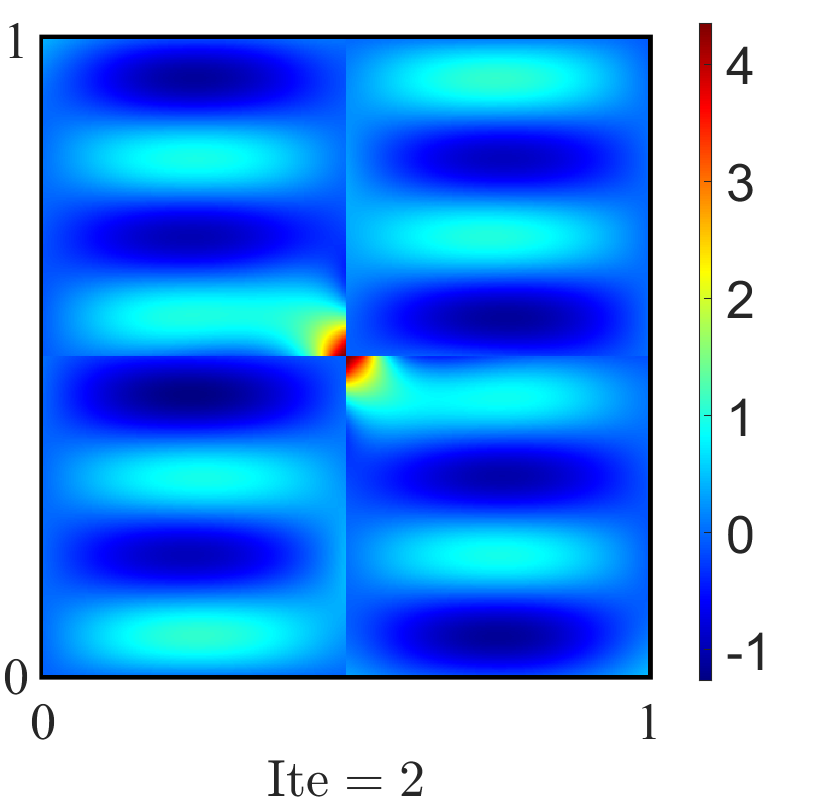}
\includegraphics[width=0.192\textwidth]{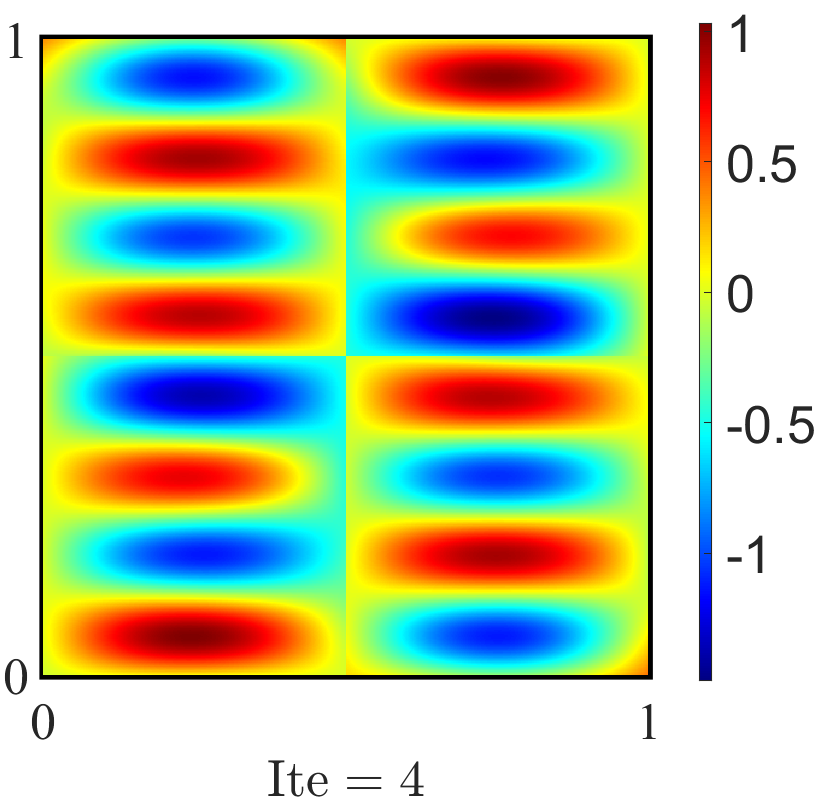}
\includegraphics[width=0.192\textwidth]{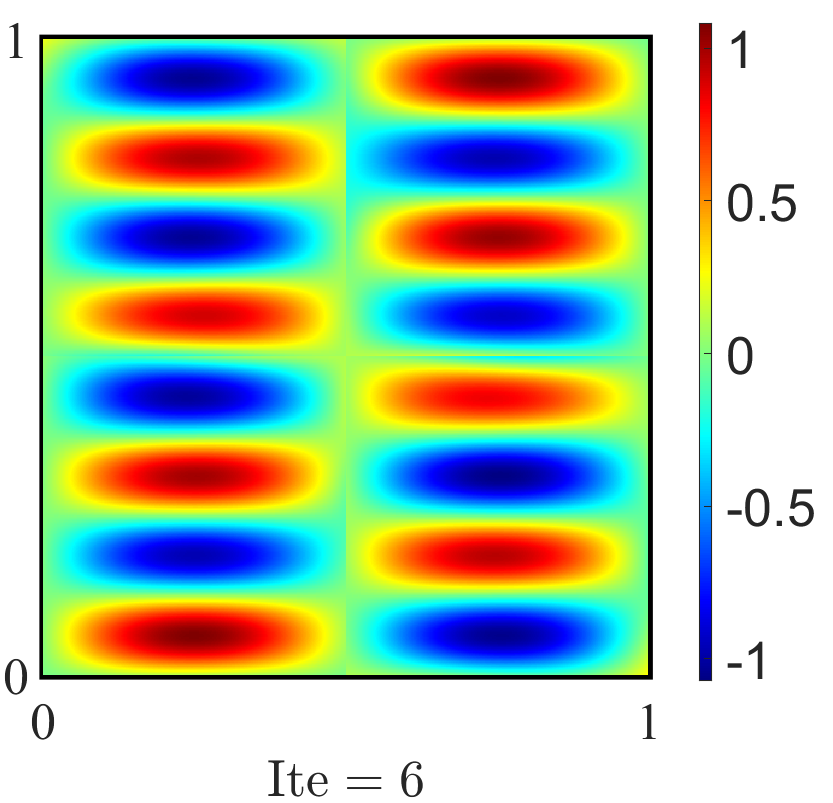}
\vspace{-0.1cm}
\caption{Iterative solutions $\hat{u}^{[k]}(x,y)$ using DNLA (PINNs) along the outer iteration. }
\label{Experiments-DNLA-ex3-DNLA-PINN-solution}
\vspace{-0.15cm}
\end{subfigure}
\begin{subfigure}[htp]{\textwidth}
\centering
\includegraphics[width=0.192\textwidth]{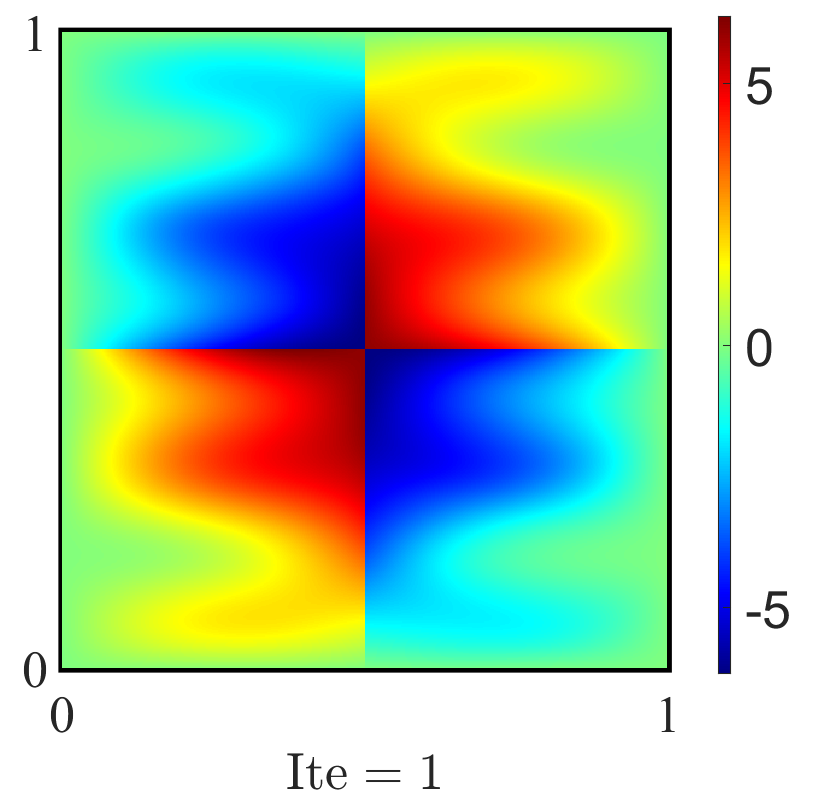}
\includegraphics[width=0.192\textwidth]{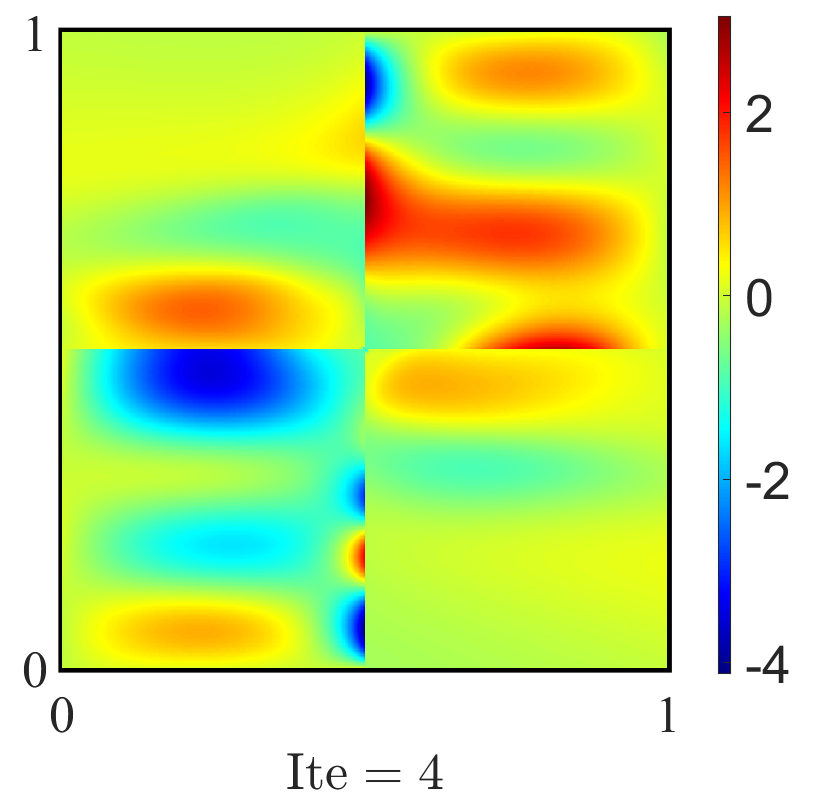}
\includegraphics[width=0.192\textwidth]{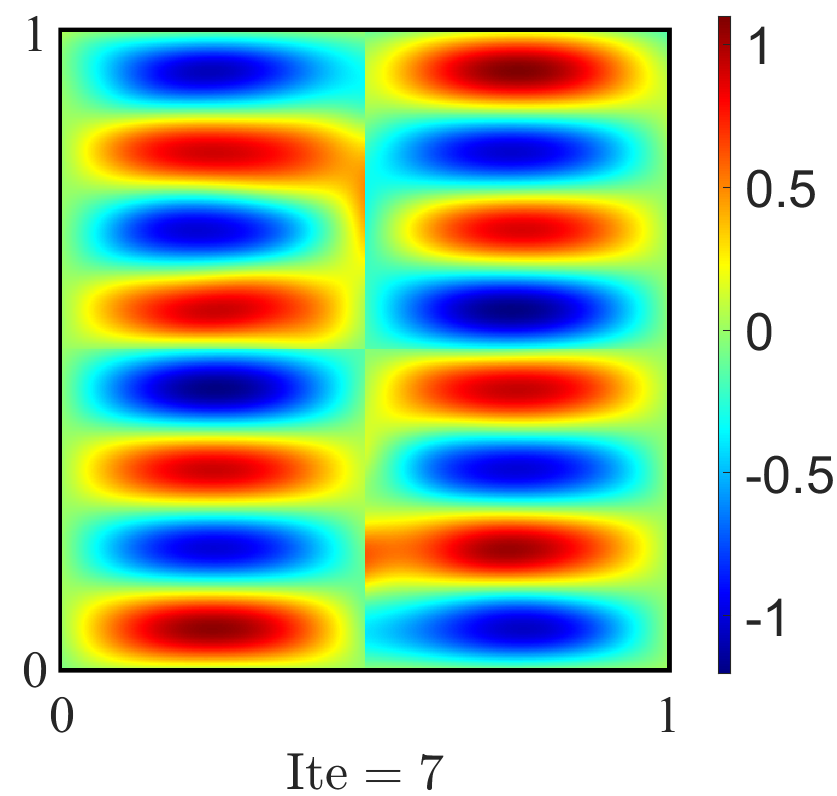}
\includegraphics[width=0.192\textwidth]{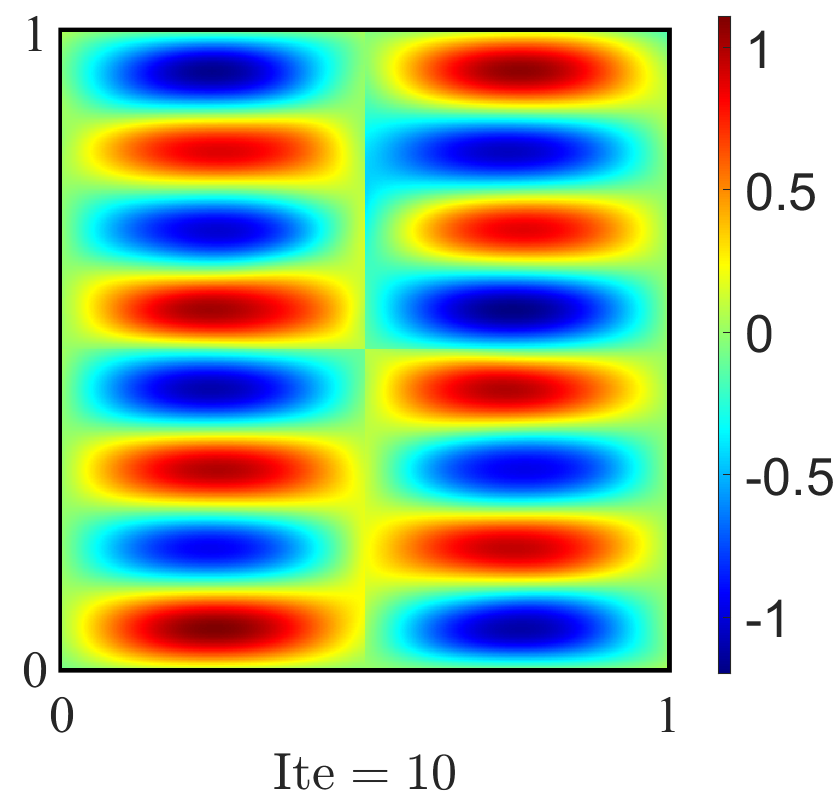}
\caption{Iterative solutions $\hat{u}^{[k]}(x,y)$ using DNLA (deep Ritz) along the outer iteration. }
\vspace{-0.1cm}
\label{Experiments-DNLA-ex3-DNLA-PINN-error}
\end{subfigure}
\vspace{-0.35cm}
\caption{Numerical results of \eqref{Experiments-DNLA-ex3} using DNLA on the test dataset.}
\label{Experiments-DNLA-ex3-DNLA-PINN}
\vspace{-0.6cm}
\end{figure}

\begin{table}[t!]
\vspace{-0.1cm}
\small
\caption{Relative-$L_2$ errors of the network solution along the outer iteration for \eqref{Experiments-DNLA-ex3}, with mean value ($\pm$ standard deviation) being reported over 5 independent runs.}
\vspace{-0.1cm}
\centering
\renewcommand{\arraystretch}{1.1}
\begin{tabular}{ | c || c | c | c | c | c | c |  }
\hline
\multicolumn{2}{|c|}{ \diagbox[width=16em]{Relative Errors}{Outer Iterations} } & 1  & 2 & 4 & 6 & 9  \\
\hline	
\hline
\multirow{5}{*}{$ \displaystyle \!\! \frac{ \lVert \hat{u}^{[k]} - u \rVert_{L_2} } { \lVert u \rVert_{L_2} }\!\!$} & DN-PINNs & \makecell{7.99 \\ \!($\pm$\! 0.08)\!} & \makecell{3.56 \\ \!($\pm$\! 1.44)\!} &  \makecell{5.90 \\ \!($\pm$\! 3.37)\!} &  \makecell{5.21 \\ \!($\pm$\! 2.90)\!} &  \makecell{2.76 \\ \!($\pm$\! 1.49)\!}  \\ 
\cline{2-7}
& DNLA (PINNs) &  \makecell{6.16 \\ \!($\pm$\! 0.2.05)\!} &  \makecell{1.35 \\ \!($\pm$\! 0.61)\!} &  \makecell{0.53 \\ \!($\pm$\! 0.24)\!} &  \makecell{0.31 \\ \!($\pm$\! 0.10)\!} &  \makecell{0.35 \\ \!($\pm$\! 0.10)\!} \\ 
\cline{2-7}
& \!\!\! DNLA (deep Ritz)\! &  \makecell{7.98 \\ \!($\pm$\! 0.11)\!} &  \makecell{3.01 \\ \!($\pm$\! 0.85)\!} &  \makecell{1.51 \\ \!($\pm$\! 0.92)\!} &  \makecell{1.04 \\ \!($\pm$\! 0.93)\!} &  \makecell{0.36 \\ \!($\pm$\! 0.10)\!} \\ 
\hline			                                                     
\end{tabular}
\label{Experiments-DNLA-ex3-Err-Table}
\vspace{-0.2cm}
\end{table}

\begin{figure}[t!]
\centering
\begin{subfigure}[htp]{\textwidth}
\centering
\includegraphics[width=0.176\textwidth]{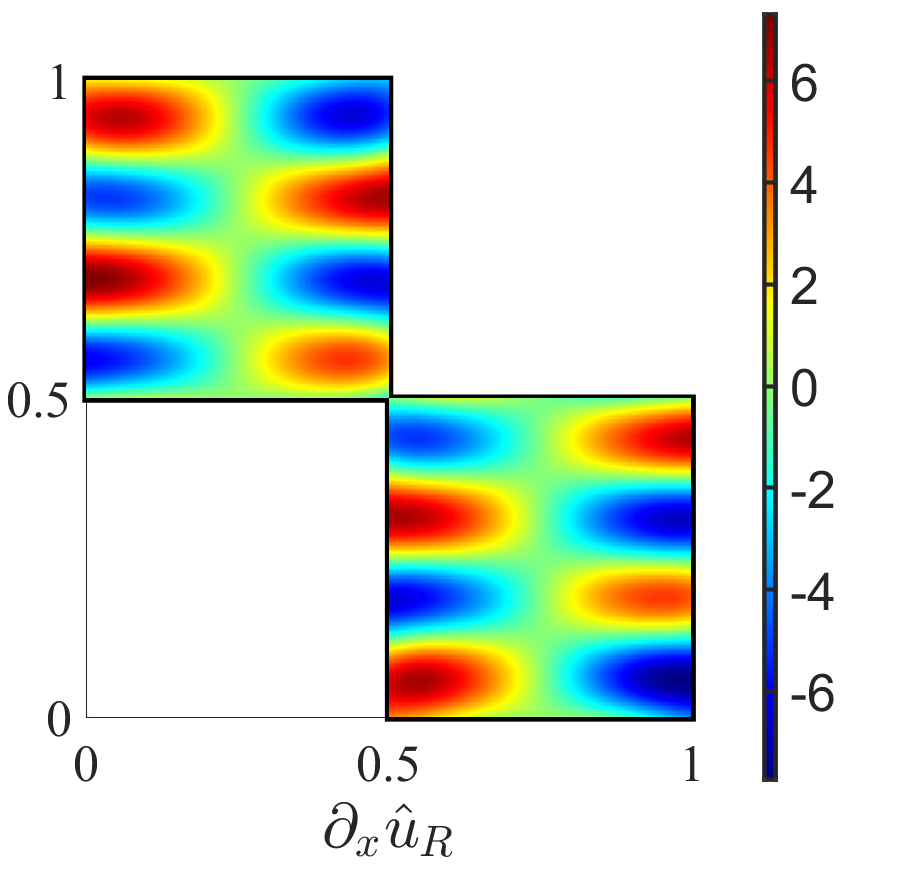}
\includegraphics[width=0.176\textwidth]{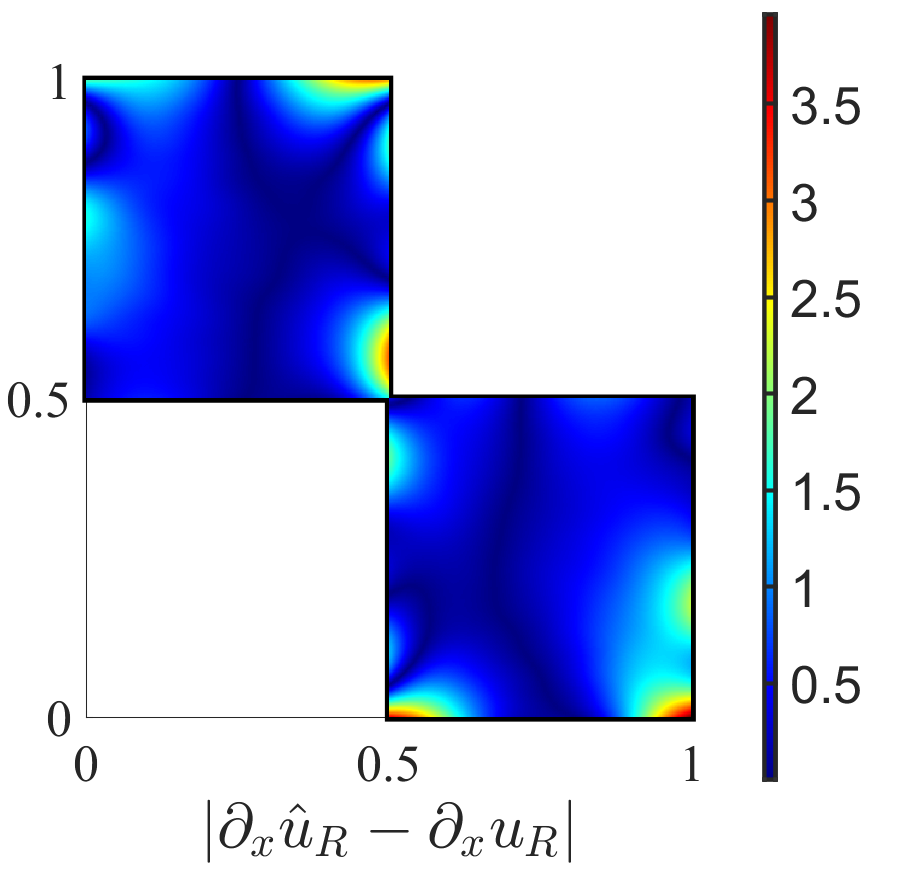}
\includegraphics[width=0.176\textwidth]{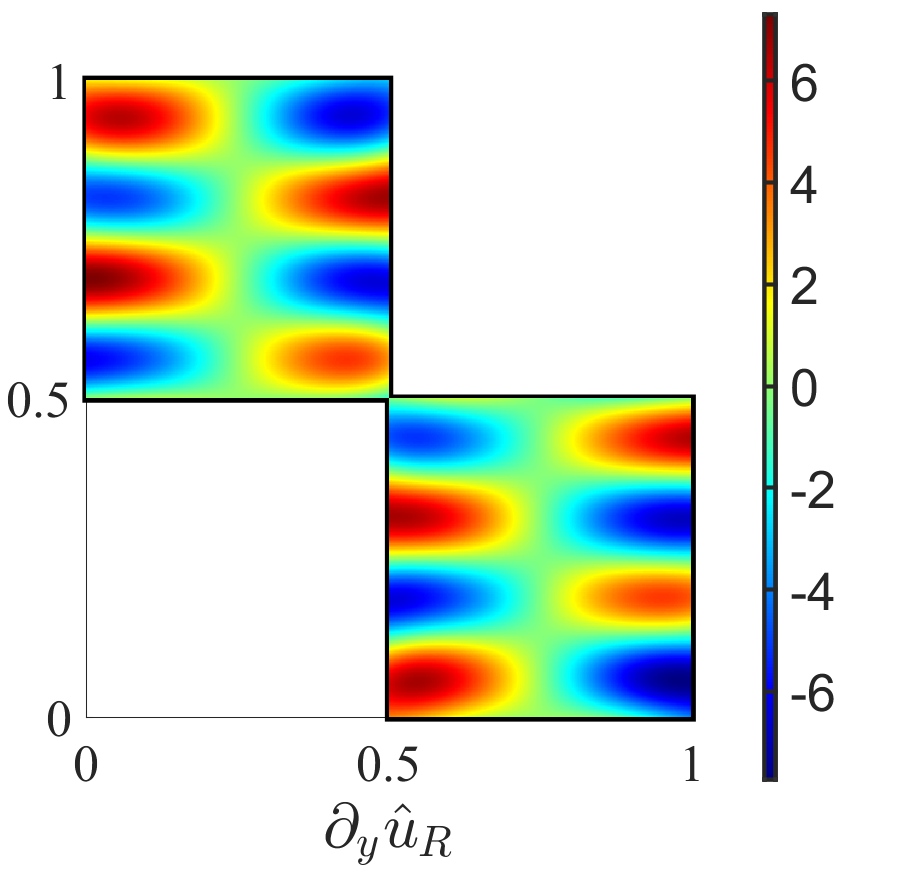}
\includegraphics[width=0.176\textwidth]{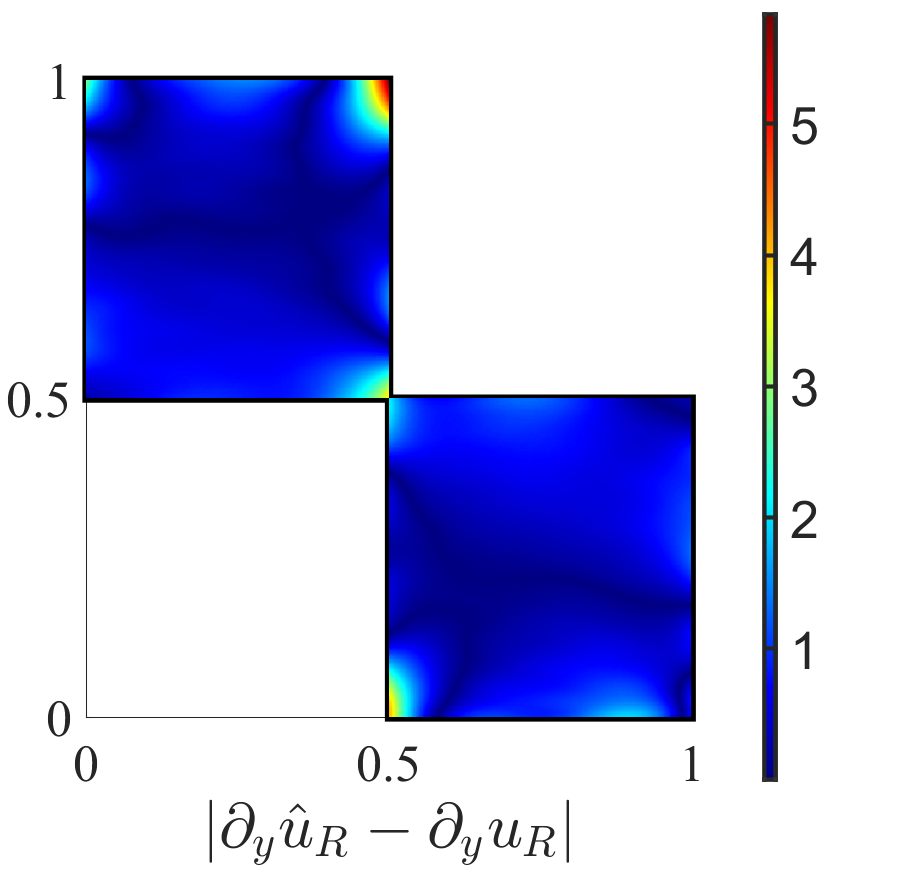}
\vspace{-0.15cm}
\caption{$\partial_x\hat{u}_R^{[3]}$, $\partial_y\hat{u}_R^{[3]}$ and error profiles $|\partial_x (\hat{u}_R^{[3]} - u_R)|$, $|\partial_y (\hat{u}_R^{[3]} - u_R)|$ using DNLA (PINNs). }
\end{subfigure}
\vspace{-0.2cm}
\begin{subfigure}[htp]{\textwidth}
\centering
\includegraphics[width=0.176\textwidth]{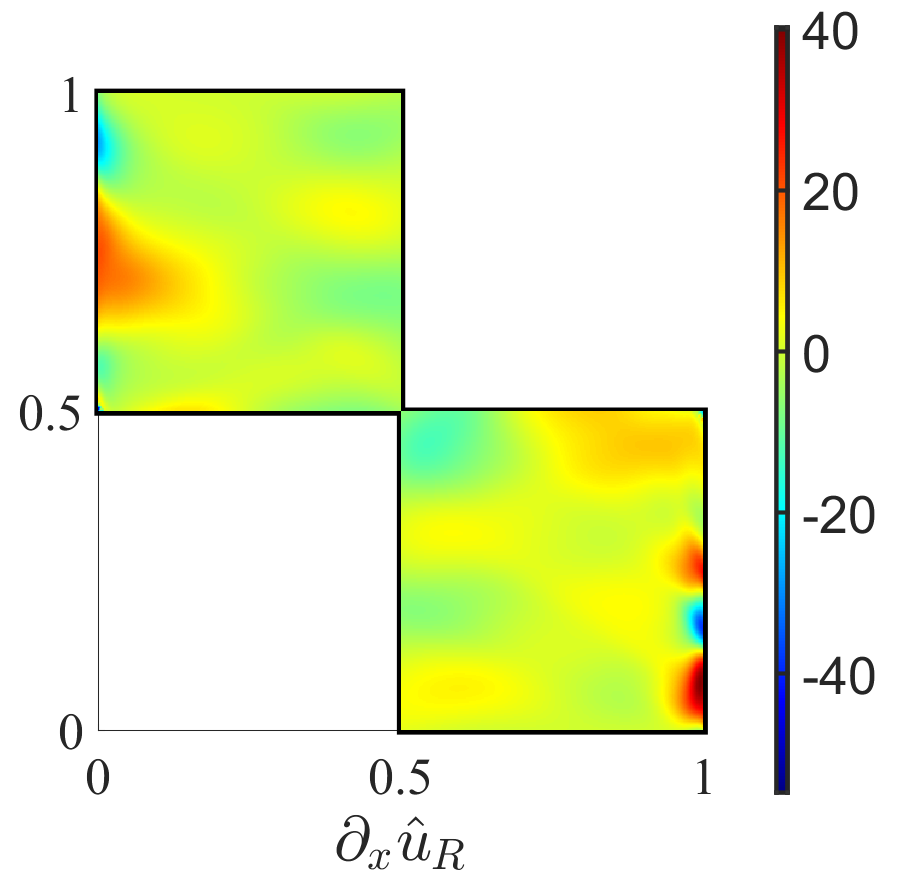}
\includegraphics[width=0.176\textwidth]{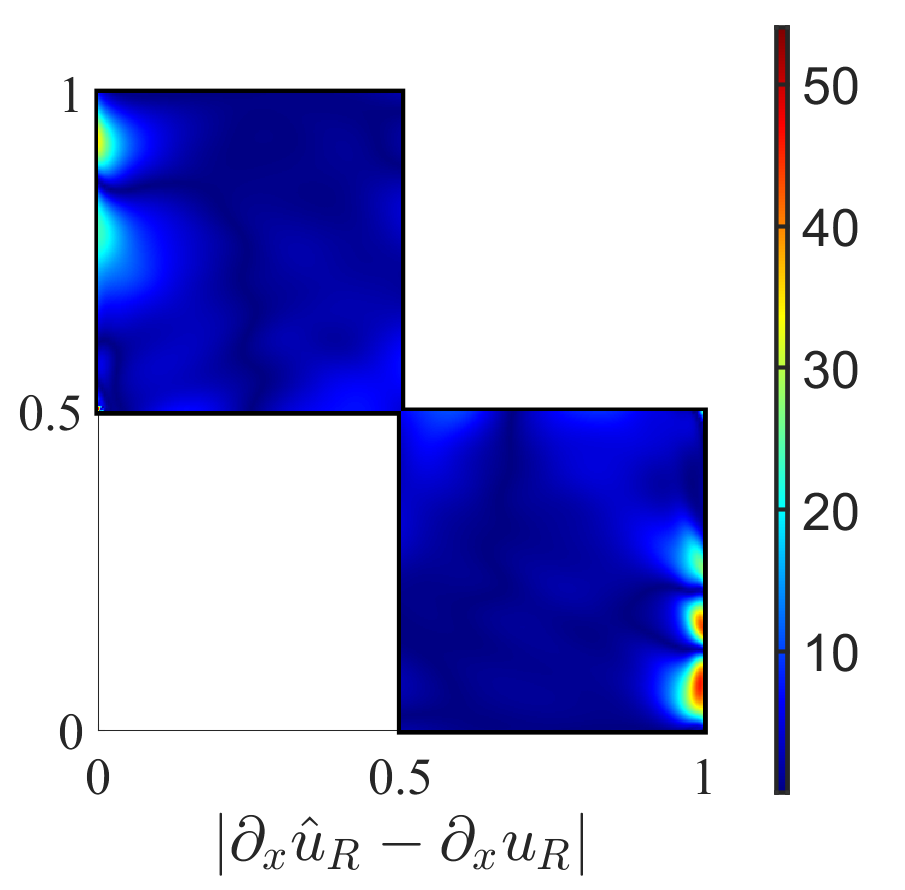}
\includegraphics[width=0.176\textwidth]{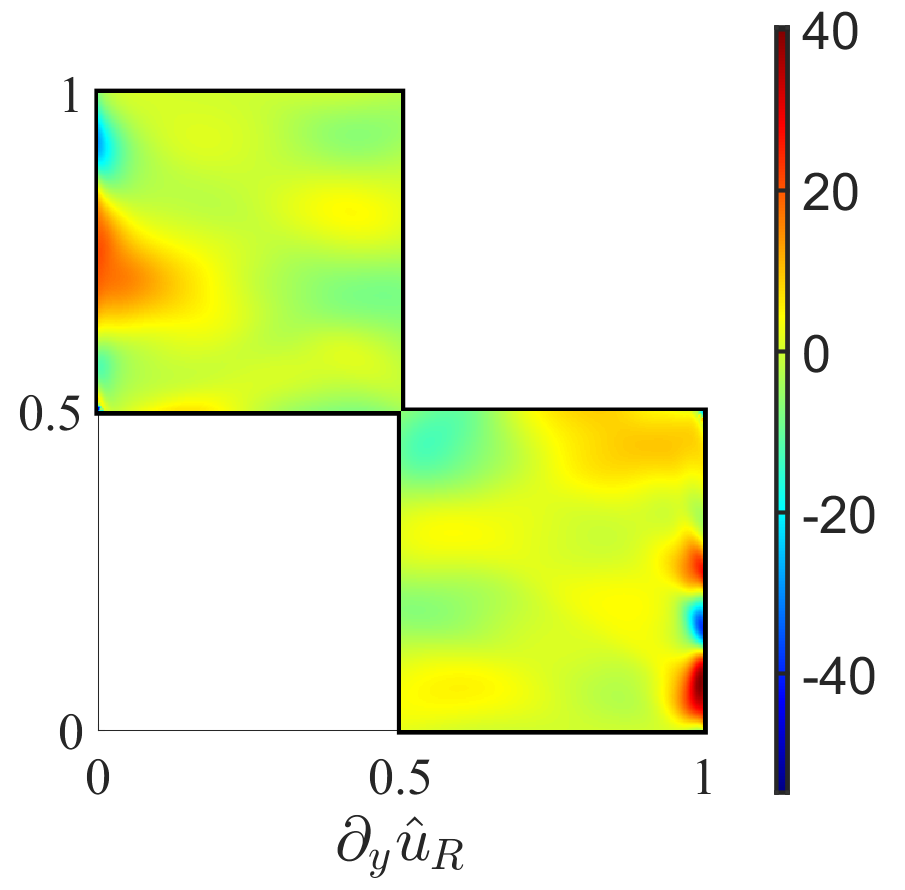}
\includegraphics[width=0.176\textwidth]{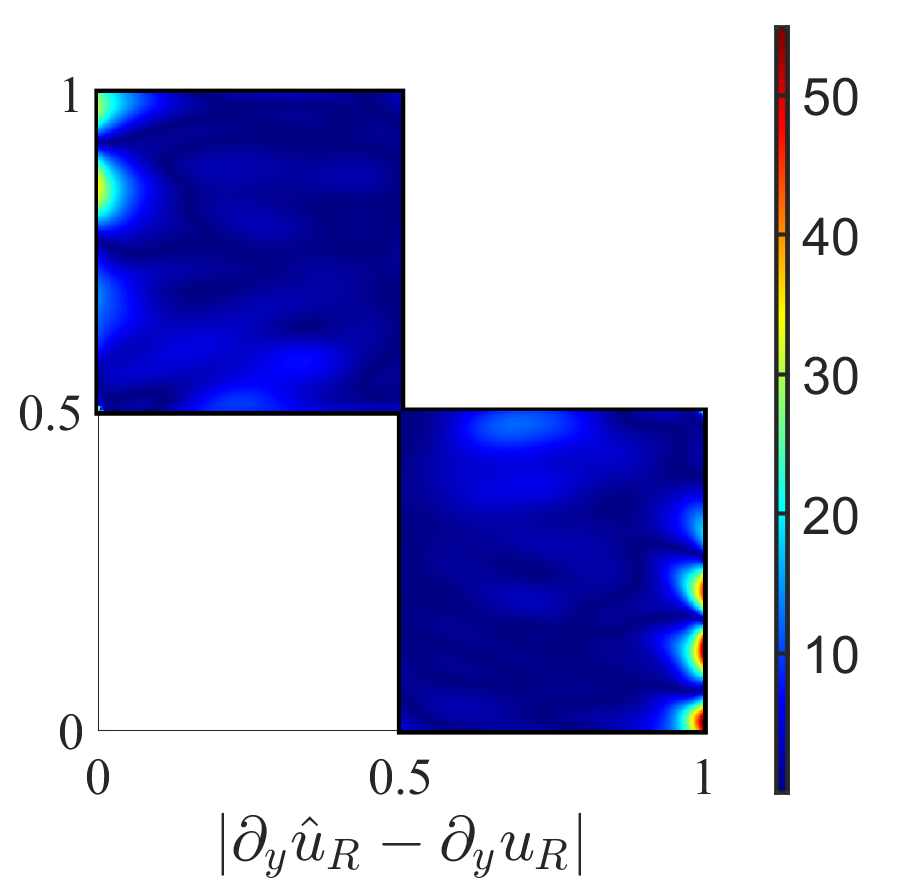}
\vspace{-0.1cm}
\caption{$\partial_x\hat{u}_R^{[5]}$, $\partial_y\hat{u}_R^{[5]}$ and error profiles $|\partial_x (\hat{u}_R^{[5]} - u_R)|$, $|\partial_y (\hat{u}_R^{[5]} - u_R)|$ using DNLA (deep Ritz). }
\end{subfigure}
\vspace{-0.25cm}
\caption{Erroneous Dirichlet-to-Neumann maps for \eqref{Experiments-DNLA-ex3}.}
\label{Experiments-DNLA-ex3-Overfit-Dirichlet-Subproblem}
\vspace{-0.7cm}
\end{figure}

For problem \eqref{Experiments-DNLA-ex3} with non-trivial flux functions along the interface, it is not guaranteed that iterative solutions using DN-PINNs will converge to the true solution due to issue of erroneous Dirichlet-to-Neumann map (see supplementary materials). However, even though the inaccurate flux predicition on subdomain interfaces remains unresolved when using our methods (see \autoref{Experiments-DNLA-ex3-Overfit-Dirichlet-Subproblem}), the compensated deep Ritz method has enabled the Neumann subproblem to be solved with acceptable accuracy. Moreover, we execute the simulation for 5 runs and report the statistical results in \autoref{Experiments-DNLA-ex3-Err-Table} to further demonstrate that DNLA (PINNs) can outperform other methods in terms of accuracy. Notably, as neural networks often fit functions from low to high frequency during the training process \cite{xu2022overview}, the relative-$L_2$ errors for problem \eqref{Experiments-DNLA-ex3} are larger than previous examples and can be further reduced using more sophisticated network architectures \cite{xu2022overview}.

\subsubsection{Poisson's Equation in High Dimension}
As is well known, another key and desirable advantage of using deep learning solvers is that they can tackle difficulties induced by the curse of dimensionality. To this end, we consider a Poisson problem in five dimension, i.e.,
\begin{equation}
\begin{array}{cl}
\displaystyle -\Delta u(x_1,\cdots,x_5)  =  4\pi^2\sum\limits_{i=1}^5 \sin (x_i)\ & \text{in}\ \Omega = (0,1)^5, \\
u(x_1,\cdots,x_5) = 0\ \ & \text{on}\ \partial \Omega,
\end{array}
\label{Experiments-DNLA-ex4}
\end{equation}
where the exact solution $u(x_1,\cdots,x_5) = \sum\limits_{i=1}^5 \sin (x_i)$, and the domain is decomposed into two subdomains $\Omega_1= \big\{(x_1,\cdots,x_5)\in\Omega \,\big|\, x_1<0.5 \big\}$ and $\Omega_2= \big\{(x_1,\cdots,x_5)\in\Omega \,\big|\, x_1>0.5 \big\}$. Here, the initial guess of the Dirichlet data at interface is chosen as $\displaystyle h^{[0]}(\mathbf{x})=u(\mathbf{x})-5000\left(x_1\prod\limits_{i=2}^5 x_i(x_i-1)\right)$, and the fully-connected neural network employed here has 8 hidden layers of 50 neurons each. The computational results using DN-PINNs, DNLA (PINNs), and DNLA (deep Ritz) approaches are shown in \autoref{Experiments-DNLA-ex4-Err-Table}, which implies that our proposed learning algorithms can achieve better performance to the existing DN-PINNs approach.

\begin{table}[t!]
\small
\caption{Relative-$L_2$ errors of the network solution along the outer iteration for \eqref{Experiments-DNLA-ex4}, with mean value ($\pm$ standard deviation) being reported over 5 runs.}
\vspace{-0.1cm}
\centering
\renewcommand{\arraystretch}{1.1}
\begin{tabular}{ | c || c | c | c | c | c |}
\hline
\multicolumn{2}{|c|}{ \diagbox[width=17em]{Relative Errors}{Outer Iterations} } & 1  & 2 & 4 & 6 \\
\hline	
\hline
\multirow{5}{*}{$ \displaystyle \!\! \frac{ \lVert \hat{u}^{[k]} - u \rVert_{L_2} } { \lVert u \rVert_{L_2} }\!\!$} & DN-PINNs & \makecell{0.980 \\ \!($\pm$\! 0.075)\!} & \makecell{0.610 \\ \!($\pm$\! 0.121)\!} &  \makecell{0.454 \\ \!($\pm$\! 0.203)\!} &  \makecell{0.377 \\ \!($\pm$\! 0.281)\!} \\ 
\cline{2-6}
& DNLA (PINNs) &  \makecell{0.310 \\ \!($\pm$\! 0.016)\!} &  \makecell{0.120 \\ \!($\pm$\! 0.009)\!} &  \makecell{0.058 \\ \!($\pm$\! 0.008)\!} &  \makecell{0.041 \\ \!($\pm$\! 0.001)\!} \\ 
\cline{2-6}
& \!\!\! DNLA (deep Ritz)\! &  \makecell{0.287 \\ \!($\pm$\! 0.006)\!} &  \makecell{0.144 \\ \!($\pm$\! 0.008)\!} &  \makecell{0.072 \\ \!($\pm$\! 0.003)\!} &  \makecell{0.064\\ \!($\pm$\! 0.004)\!} \\
\hline	                                                     
\end{tabular}
\label{Experiments-DNLA-ex4-Err-Table}
\vspace{-0.4cm}
\end{table}

\subsubsection{High-Contrast Elliptic Equation}
\begin{figure}[htp]
\centering
\begin{subfigure}[htp]{\textwidth}
\centering
\includegraphics[width=0.192\textwidth]{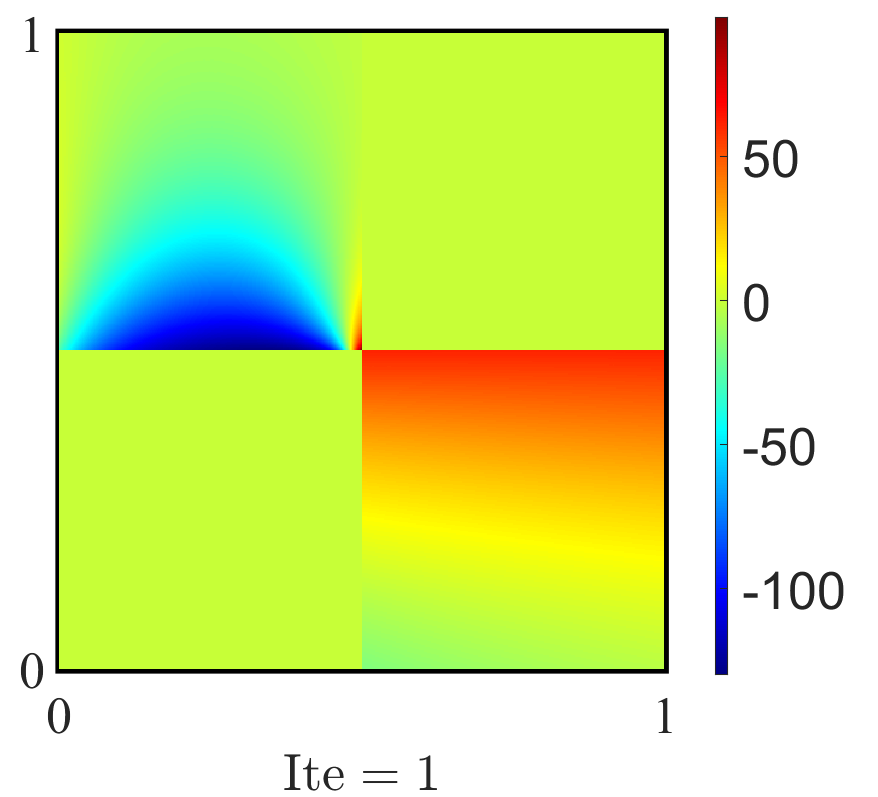}
\includegraphics[width=0.192\textwidth]{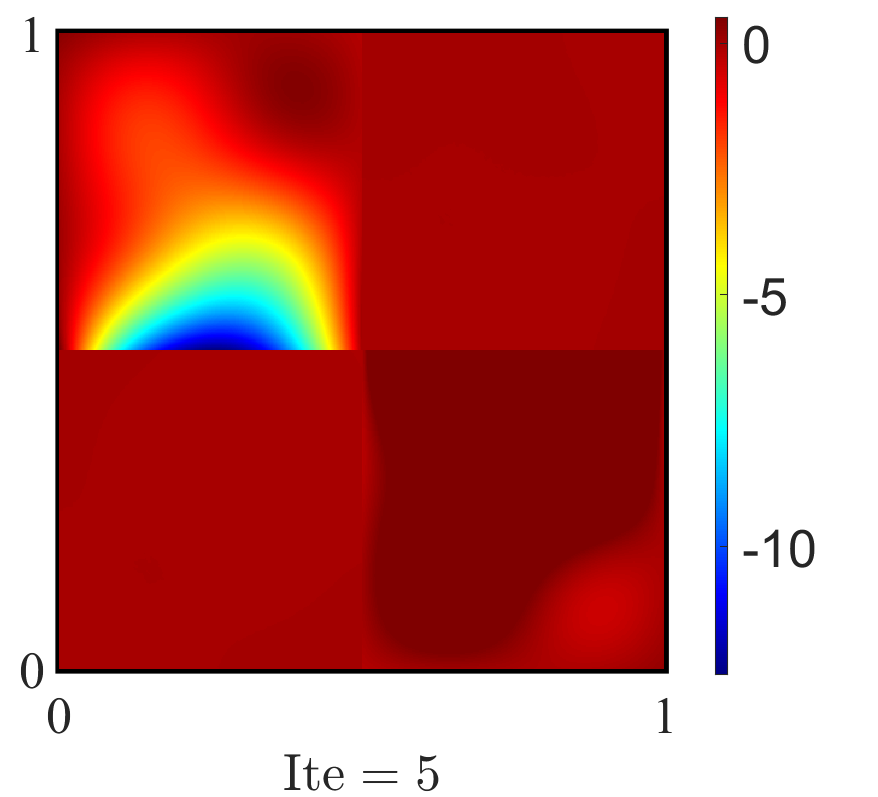}
\includegraphics[width=0.192\textwidth]{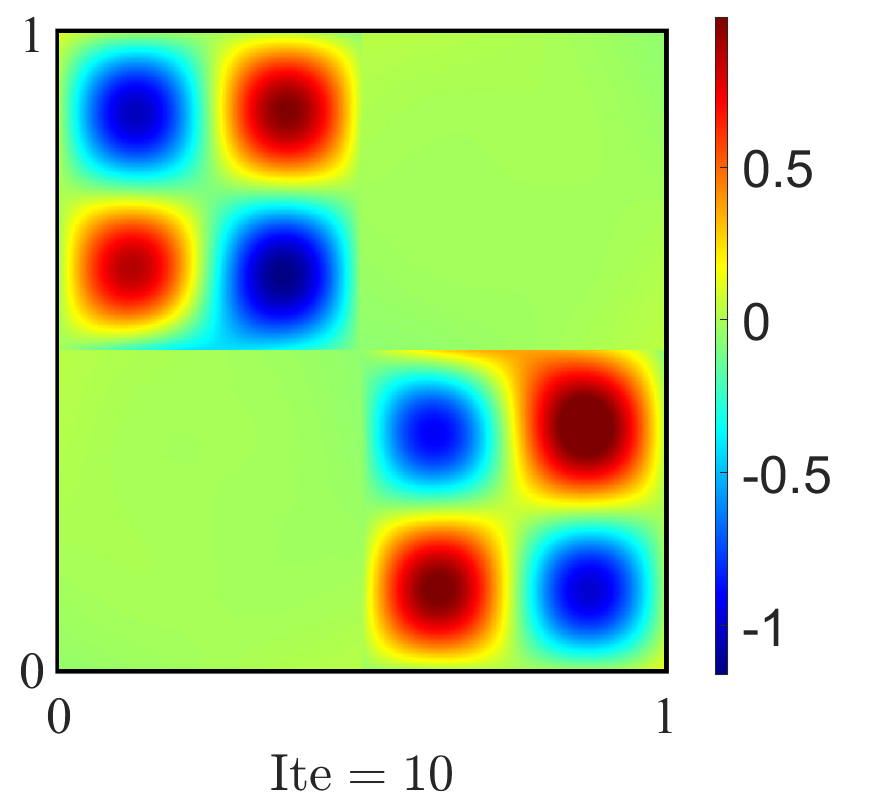}
\includegraphics[width=0.192\textwidth]{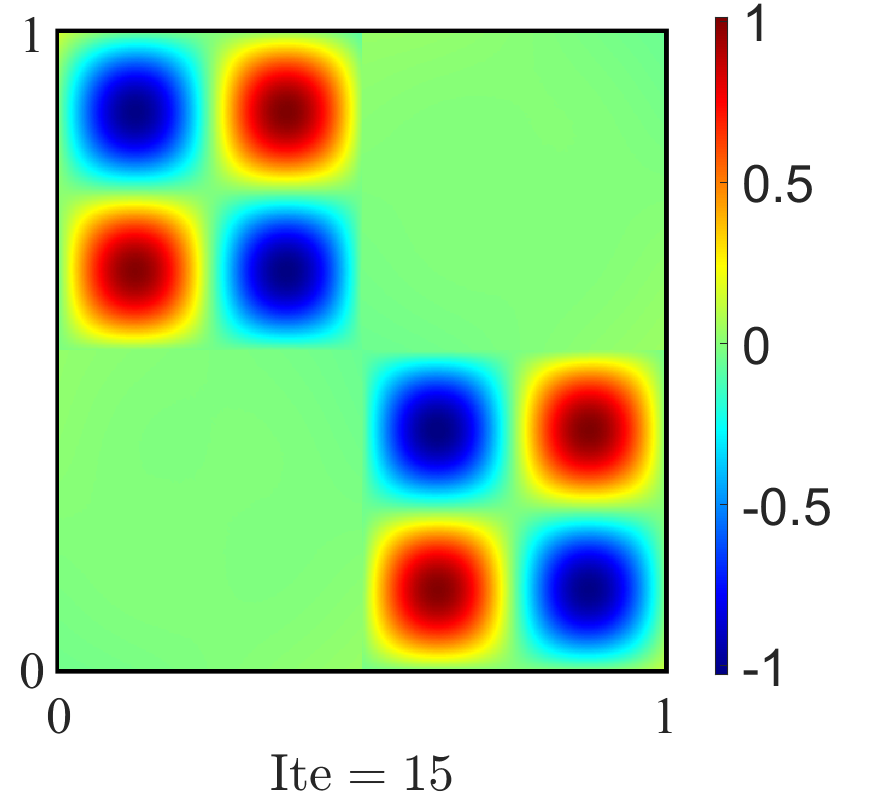}
\vspace{-0.1cm}
\caption{Iterative solutions $\hat{u}^{[k]}(x,y)$ using DNLA (PINNs) along the outer iteration. }
\label{Experiments-DNLA-ex5-DNLA-PINN-solution}
\vspace{-0.2cm}
\end{subfigure}
\begin{subfigure}[htp]{\textwidth}
\centering
\includegraphics[width=0.192\textwidth]{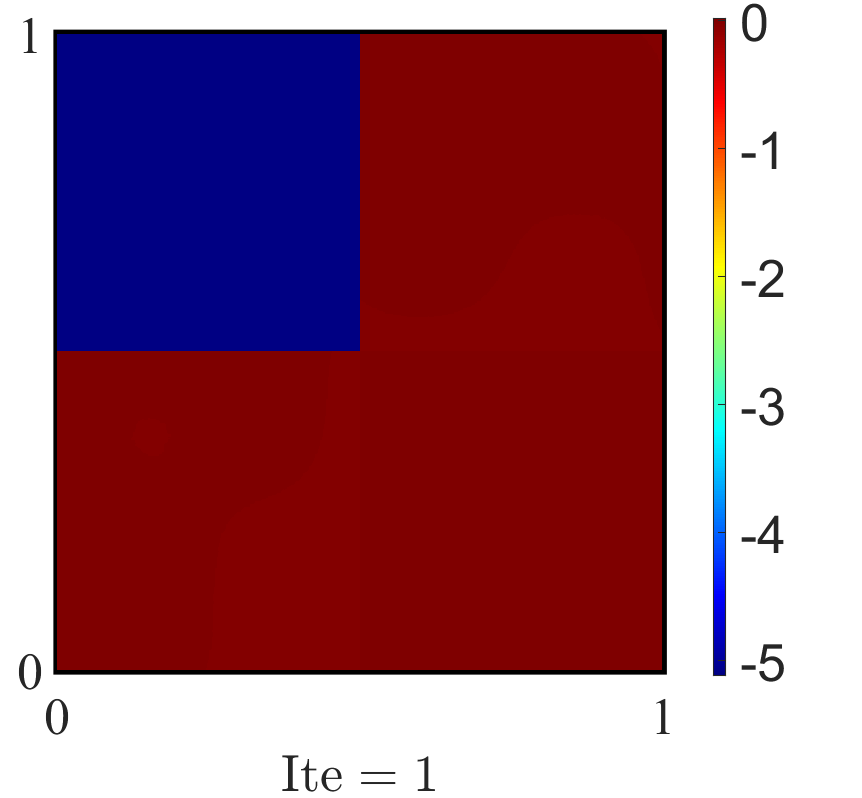}
\includegraphics[width=0.192\textwidth]{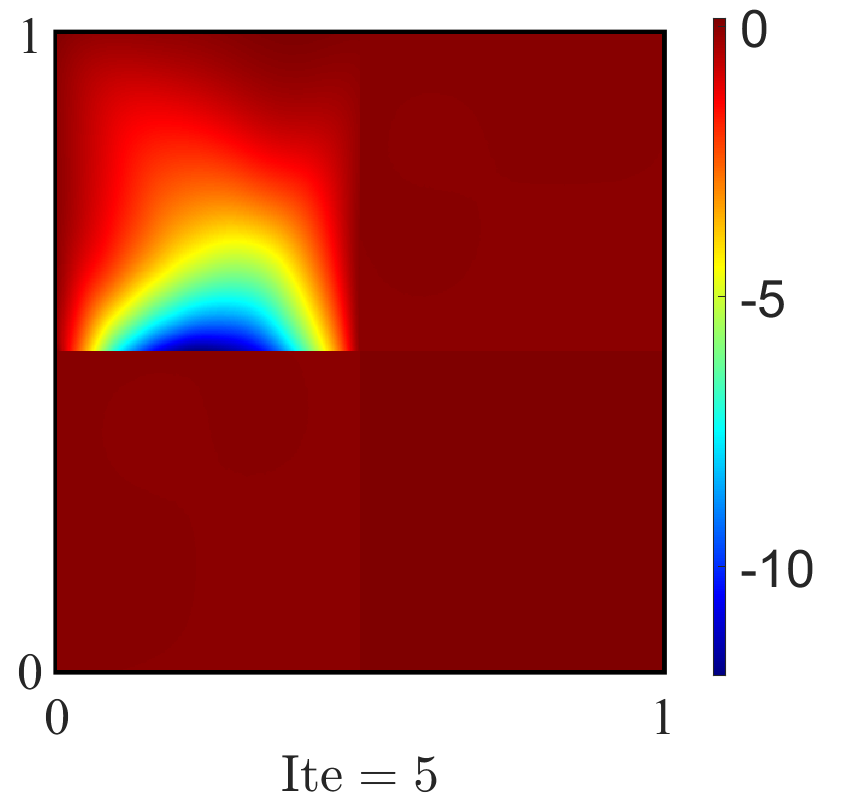}
\includegraphics[width=0.192\textwidth]{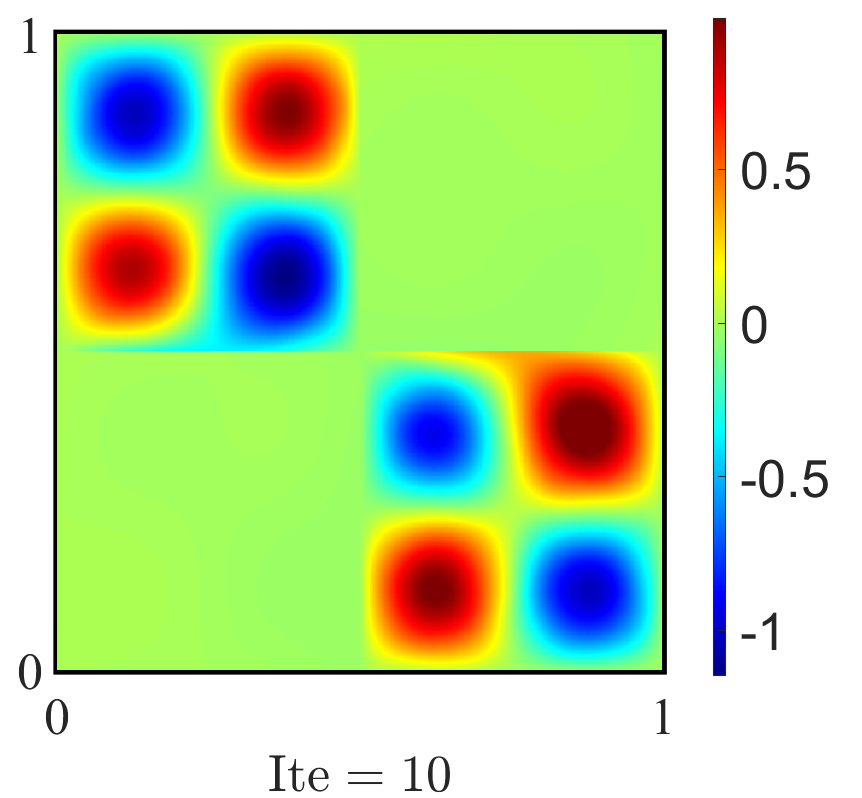}
\includegraphics[width=0.192\textwidth]{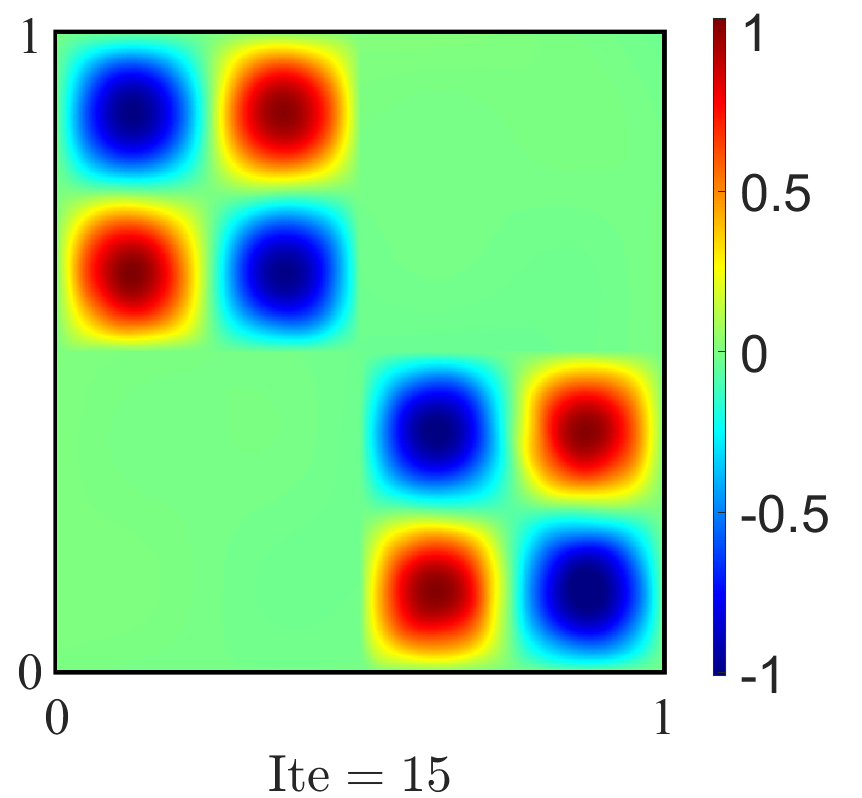}
\vspace{-0.1cm}
\caption{Iterative solutions $\hat{u}^{[k]}(x,y)$ using DNLA (deep Ritz) along the outer iteration. }
\label{Experiments-DNLA-ex5-DNLA-PINN-error}
\end{subfigure}
\vspace{-0.45cm}
\caption{Numerical results of example \eqref{Experiments-DNLA-ex5} using DNLA on the test dataset.}
\label{Experiments-DNLA-ex5-DNLA-PINN}
\vspace{-0.6cm}
\end{figure}

Note that as mentioned in Remark \ref{Remark-High-Contrast}, our proposed Dirichlet-Neumann learning algorithm can also be used to solve the more challenging interface problem with high-contrast coefficients. As such, we consider an elliptic interface problem in two dimension, 
\begin{equation}
\begin{array}{cl}
-\nabla \cdot \left( c(x,y) \nabla u(x,y)  \right) = 32 \pi^2 \sin(4\pi x)\cos(4\pi y)\ \ & \text{in}\ \Omega=(0,1)^2,\\
u(x,y) = 0\ \ & \text{on}\ \partial \Omega, 
\end{array}
\label{Experiments-DNLA-ex5}
\end{equation}
where the computational domain is decomposed into four isolated subdomains as shown in \autoref{fig-extended-4-area}, the exact solution is given by $u(x,y) = \sin(4\pi x) \sin(4\pi y) / c(x,y)$, and the coefficient $c(x,y)$ is piecewise constant with respect to the partition of domain
\begin{equation*}
c(x,y) = \left\{
\begin{array}{cl}
1 \ & \text{in}\ \Omega_1\cup\Omega_3,\\
100\ \ & \text{in}\ \Omega_2\cup\Omega_4.
\end{array}\right.
\end{equation*}
Here, we choose $h^{[0]}=100\cos(100\pi x)\cos(100\pi y)+100xy$ as the initial guess, and the numerical results using DNLA are depicted in \autoref{Experiments-DNLA-ex5-DNLA-PINN}.  Clearly, our method can facilitate the convergence of outer iterations in the presence of erroneous flux estimations (see supplementary materials for more details).

\subsection{Robin-Robin Learning Algorithm}
To demonstrate the effectiveness of our compensated deep Ritz method for realizing the non-overlapping Robin-Robin algorithm, we consider the following Poisson equation in two-dimension
\begin{equation}
\begin{array}{cl}
-\Delta u(x,y)  = 4 \pi^2 \sin(2 \pi x)  (2 \cos(2 \pi y) - 1)  \ & \text{in}\ \Omega=(0,1)^2,\\
u(x,y) = 0\ \ & \text{on}\ \partial \Omega,
\end{array}
\label{Experiments-RRLM-ex1}
\end{equation}
where the exact solution $u(x,y) = \sin(2\pi x)(\cos(2\pi y)-1)$, and the interface $\Gamma=\partial\Omega_1\cap\partial\Omega_2$ is a straight line segment from $(0.5,0)$ to $(0.5,1)$ as depicted in \autoref{Experiments-DNLA-ex1-exact-solution}. 

\begin{figure}[t!]
\centering
\includegraphics[width=0.192\textwidth]{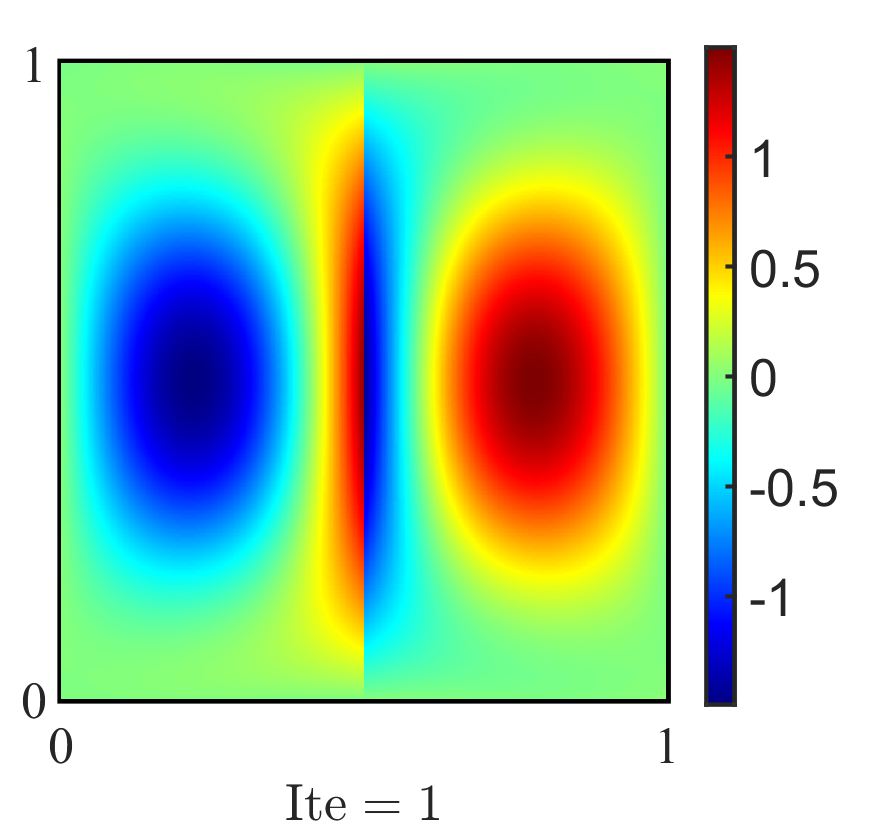}
\includegraphics[width=0.192\textwidth]{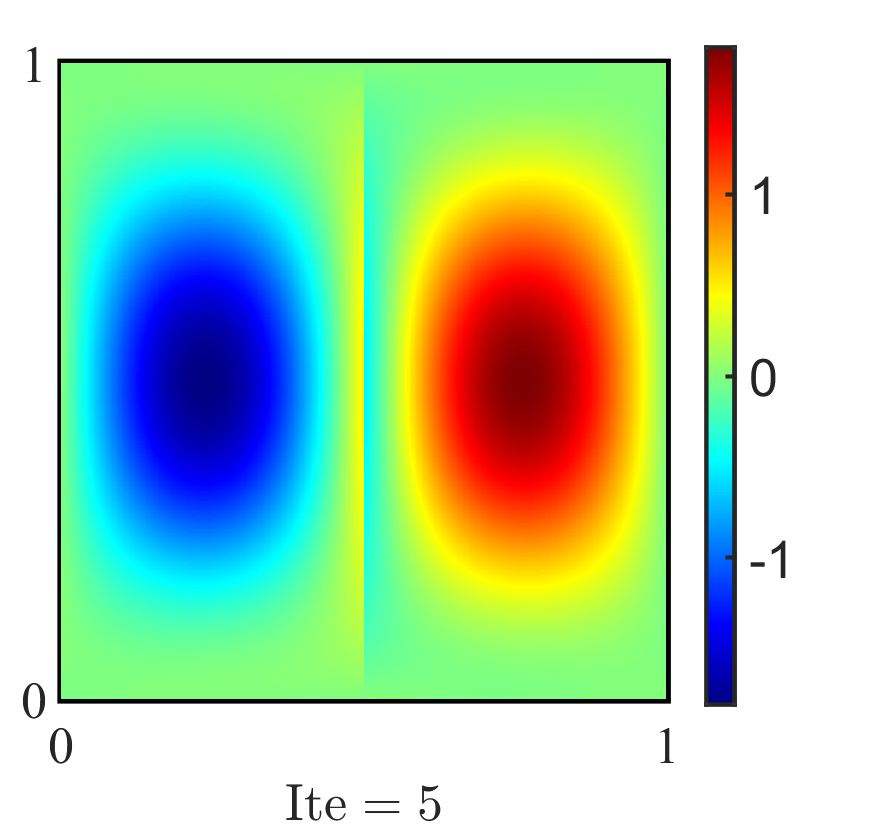}
\includegraphics[width=0.176\textwidth]{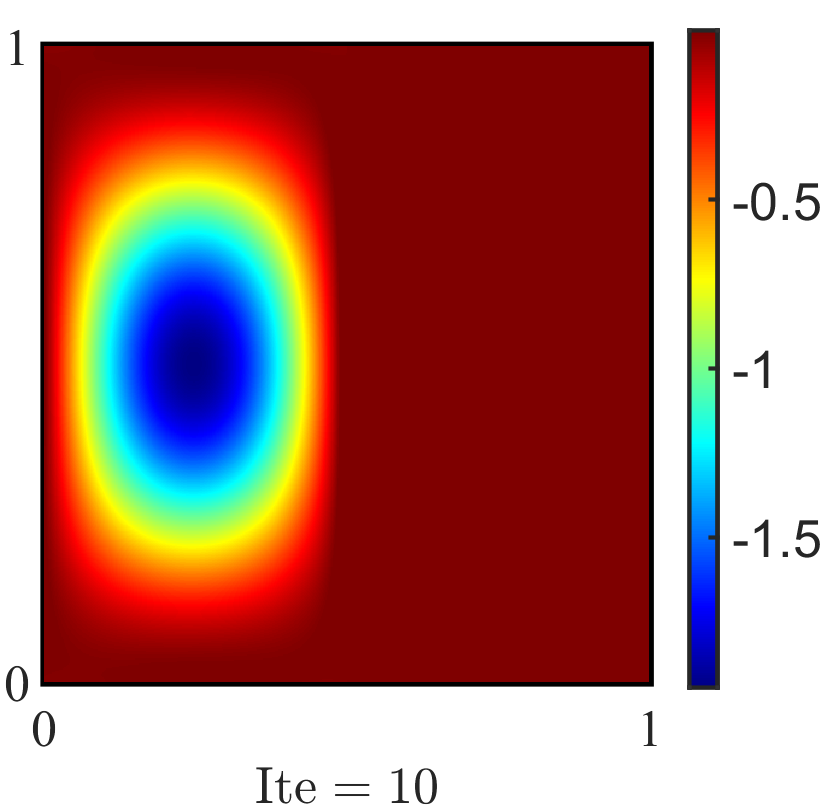}
\includegraphics[width=0.192\textwidth]{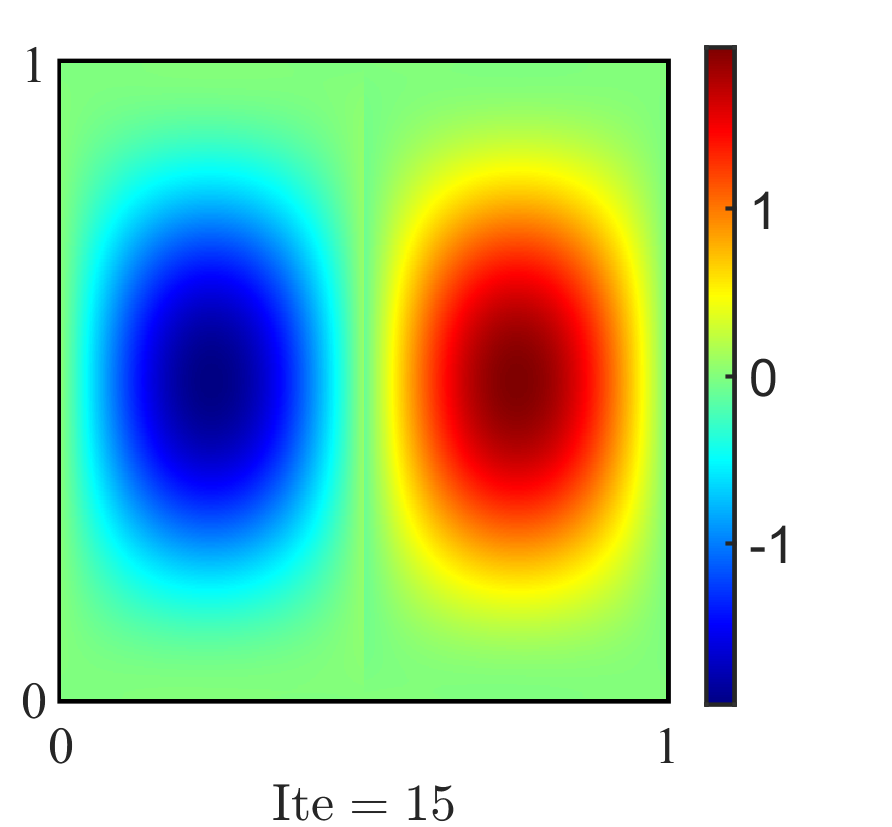}
\vspace{-0.2cm}
\caption{Iterative solutions $\hat{u}^{[k]}$ of \eqref{Experiments-RRLM-ex1} using RR-PINNs for $(\kappa_1,\kappa_2)=(1,1)$.}
\label{Experiments-RRLM-ex1-RR-PINNs}
\vspace{-0.5cm}
\end{figure}

By choosing $(\kappa_1,\kappa_2)=(1,0.01)$, the computational results using RR-PINNs, i.e., Algorithm \ref{Algorithm-RR-Learning-2Subdomains}, in a typical simulation is depicted \autoref{Experiments-RRLM-ex1-RR-PINNs}, which can converge to the true solution but requires extra outer iterations when compared to the DNLA (PINNs) or DNLA (deep Ritz) approach (see \autoref{Experiments-DNLA-ex1-DNLA-PINN} or \autoref{Experiments-DNLA-ex1-DNLA-DeepRitz}). 

\begin{figure}[t!]
\centering
\includegraphics[width=0.192\textwidth]{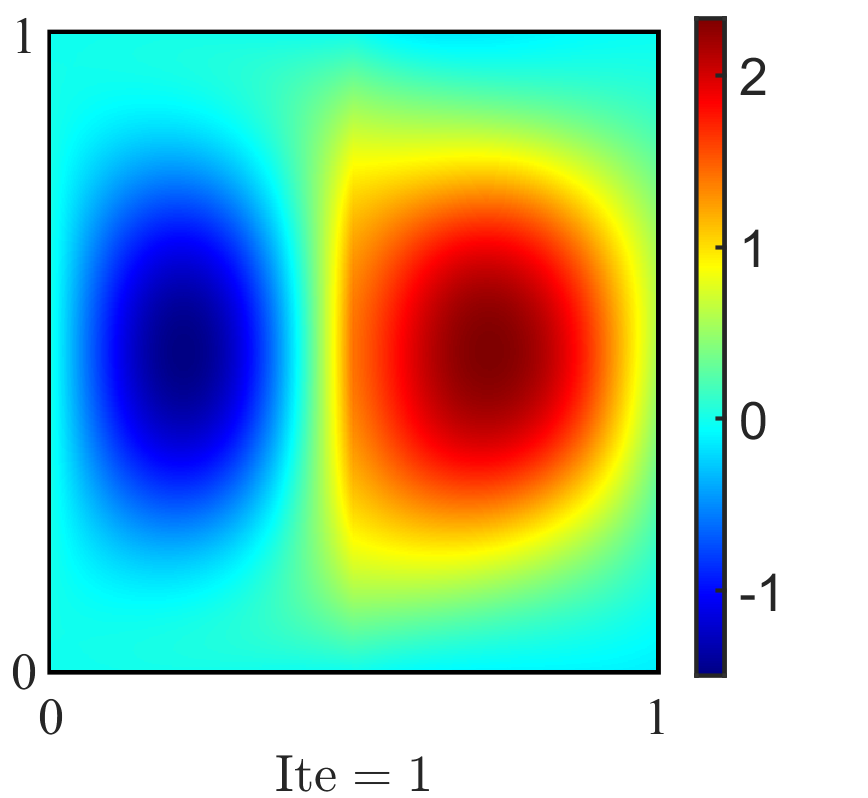}
\includegraphics[width=0.192\textwidth]{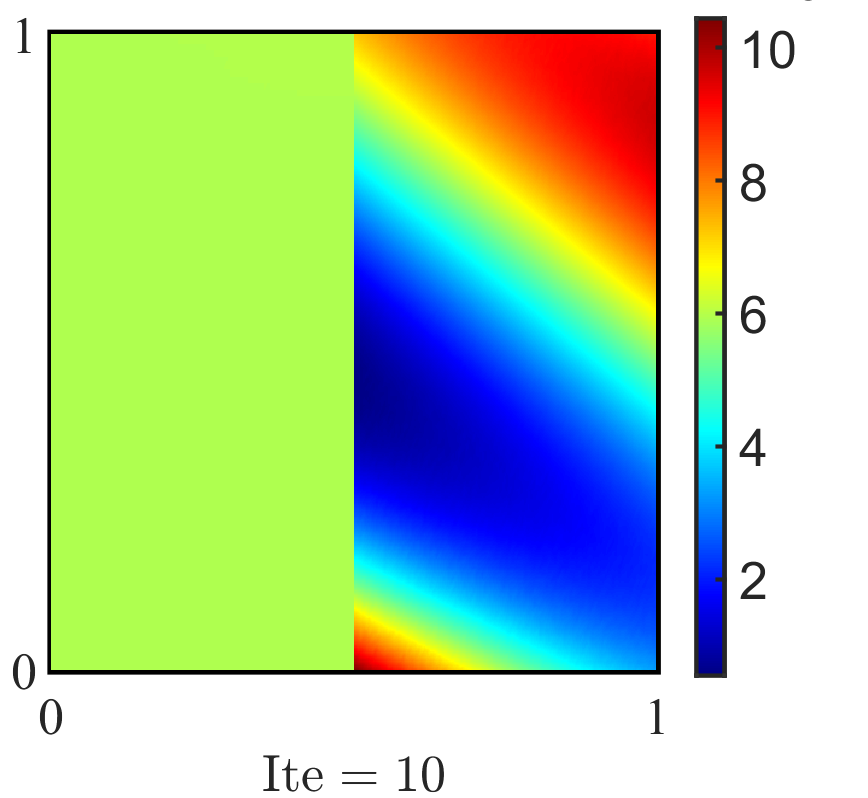}
\includegraphics[width=0.192\textwidth]{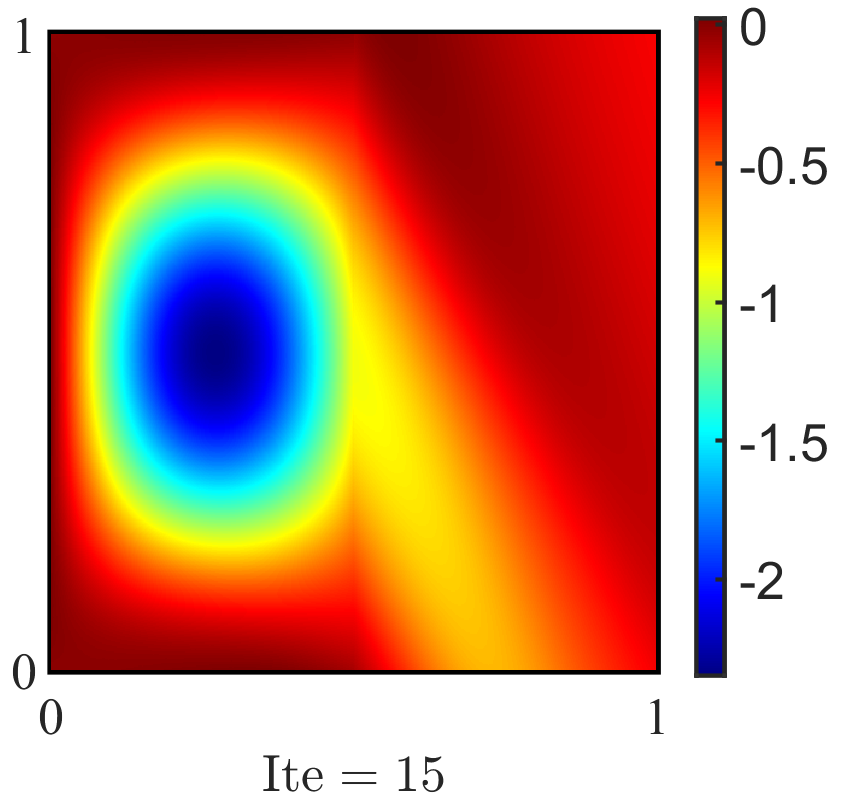}
\includegraphics[width=0.192\textwidth]{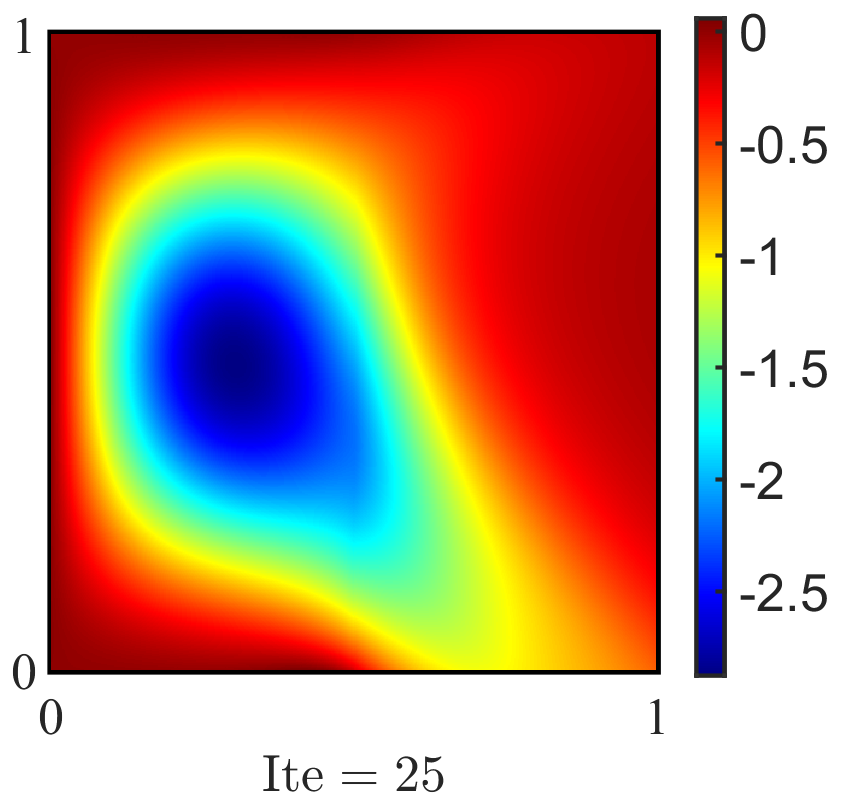}
\vspace{-0.2cm}
\caption{Iterative solutions $\hat{u}^{[k]}$ of \eqref{Experiments-RRLM-ex1} using RR-PINNs for $(\kappa_1,\kappa_2)=(1,1000)$.}
\label{Experiments-RRLM-ex1-RR-PINNs-1-1000}
\vspace{-0.5cm}
\end{figure}

\begin{figure}[t!]
\begin{subfigure}[htp]{\textwidth}
\centering
\includegraphics[width=0.192\textwidth]{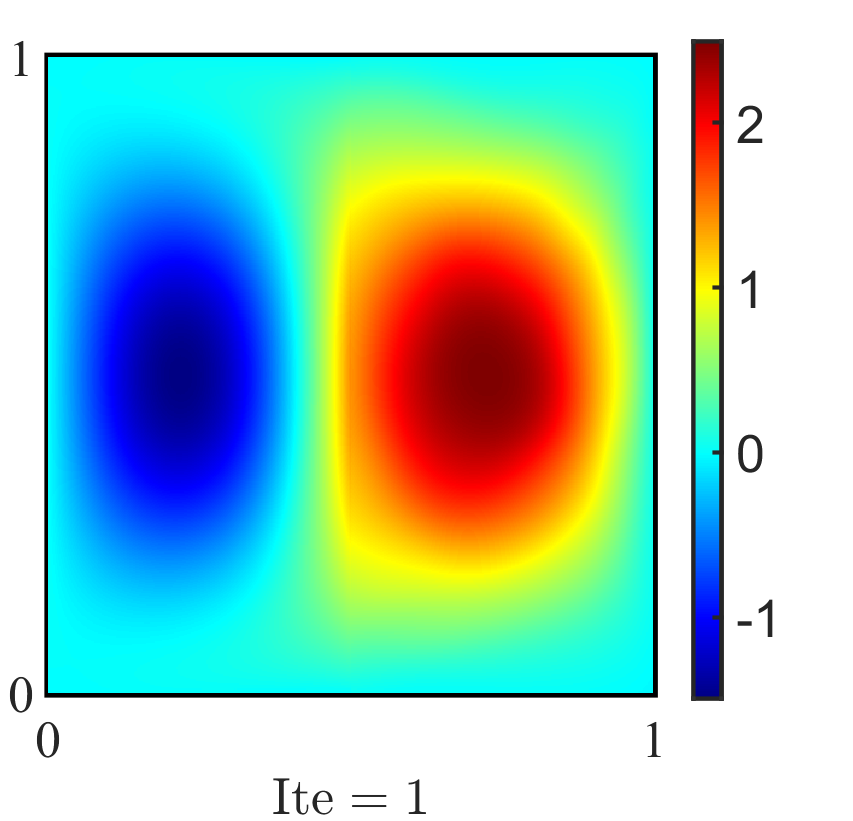}
\includegraphics[width=0.192\textwidth]{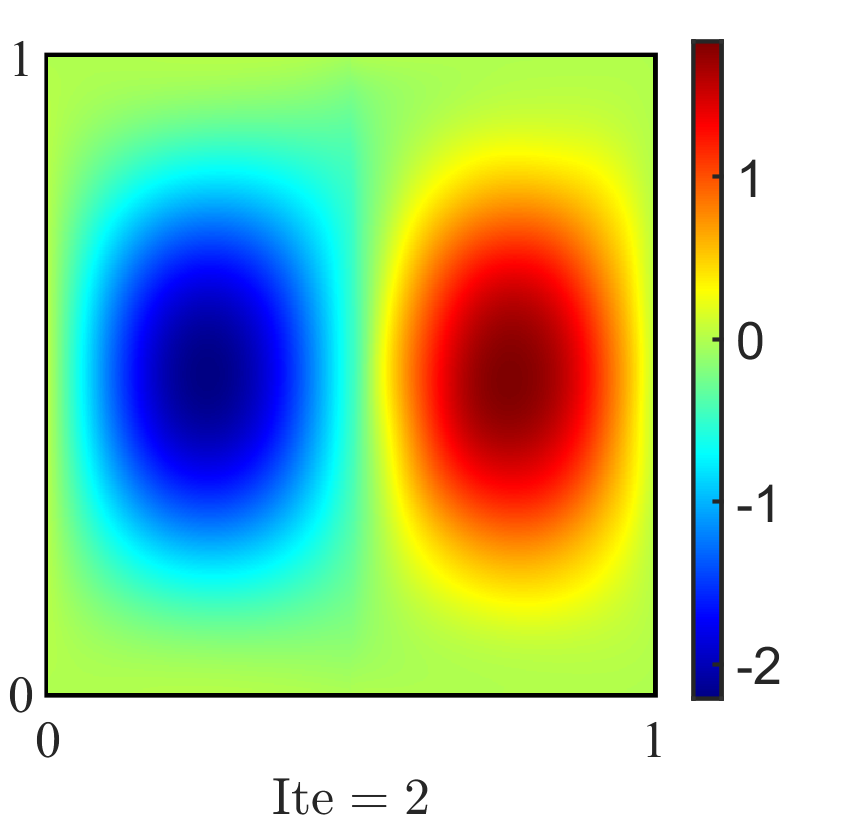}
\includegraphics[width=0.192\textwidth]{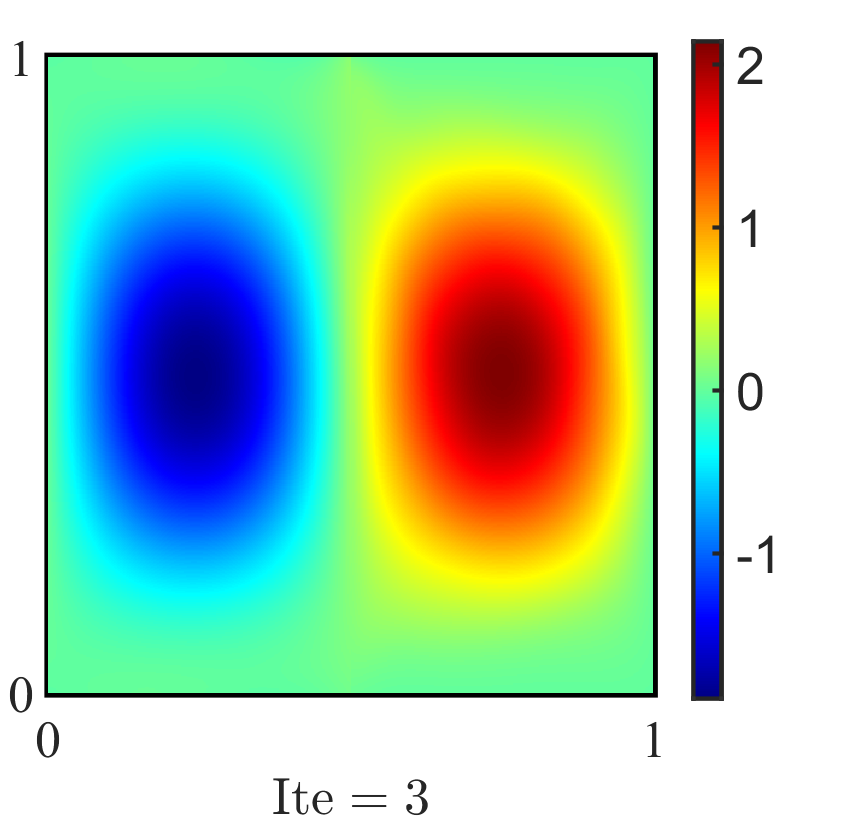}
\includegraphics[width=0.192\textwidth]{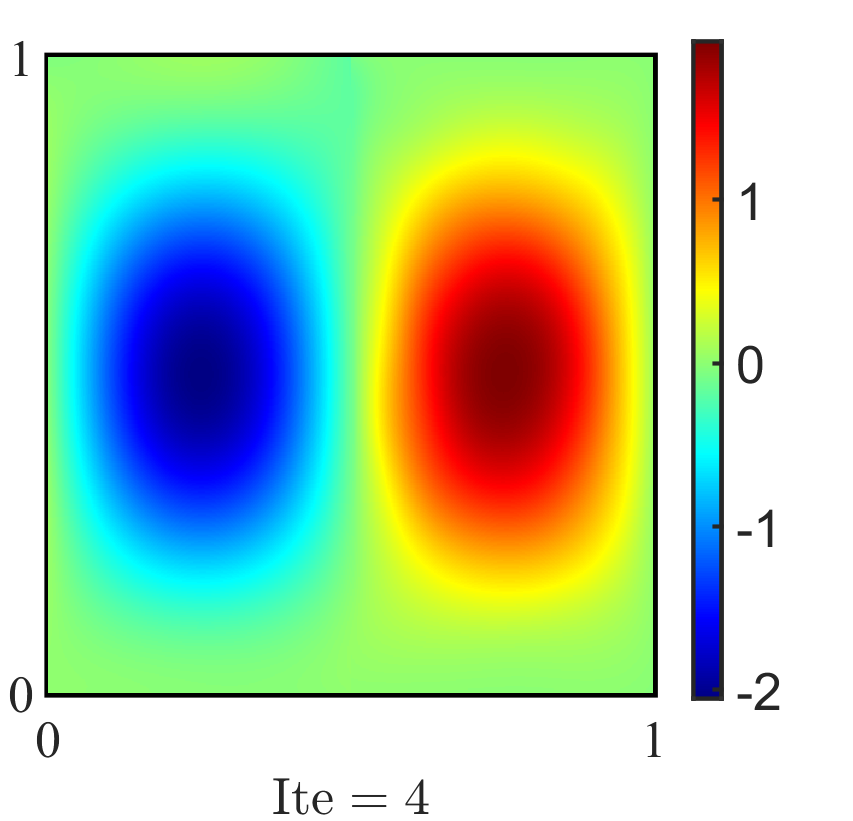}
\caption{Iterative solutions $\hat{u}^{[k]}(x,y)$ using RRLM (PINNs) along the outer iteration. }
\label{Experiments-RRLM-ex2-RRLM-PINN-solution}
\end{subfigure}
\begin{subfigure}[htp]{\textwidth}
\centering
\includegraphics[width=0.192\textwidth]{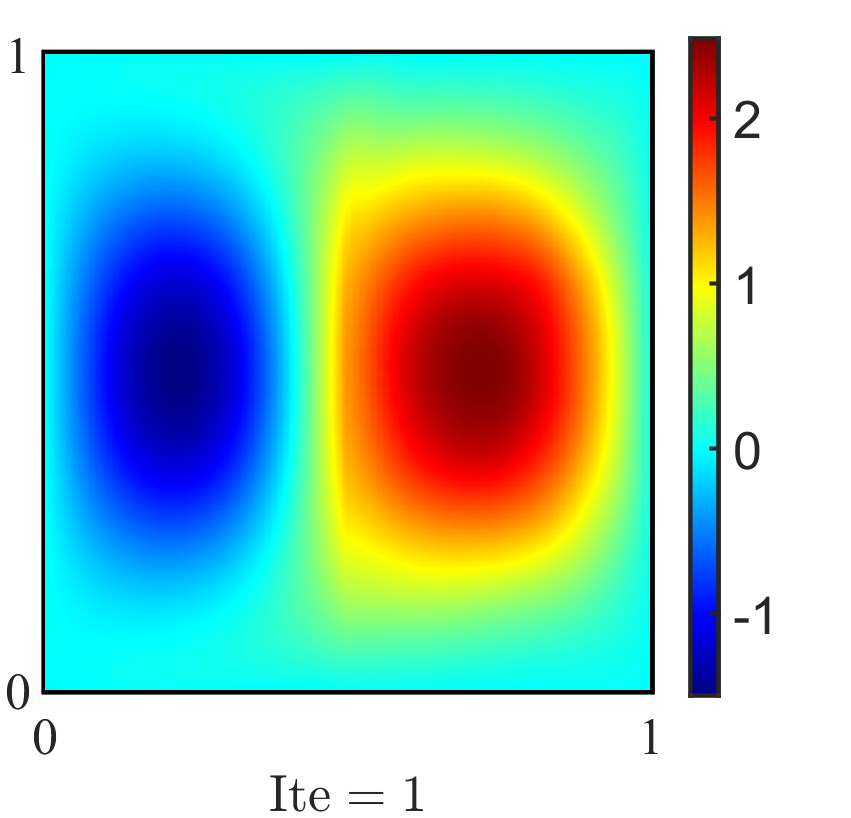}
\includegraphics[width=0.192\textwidth]{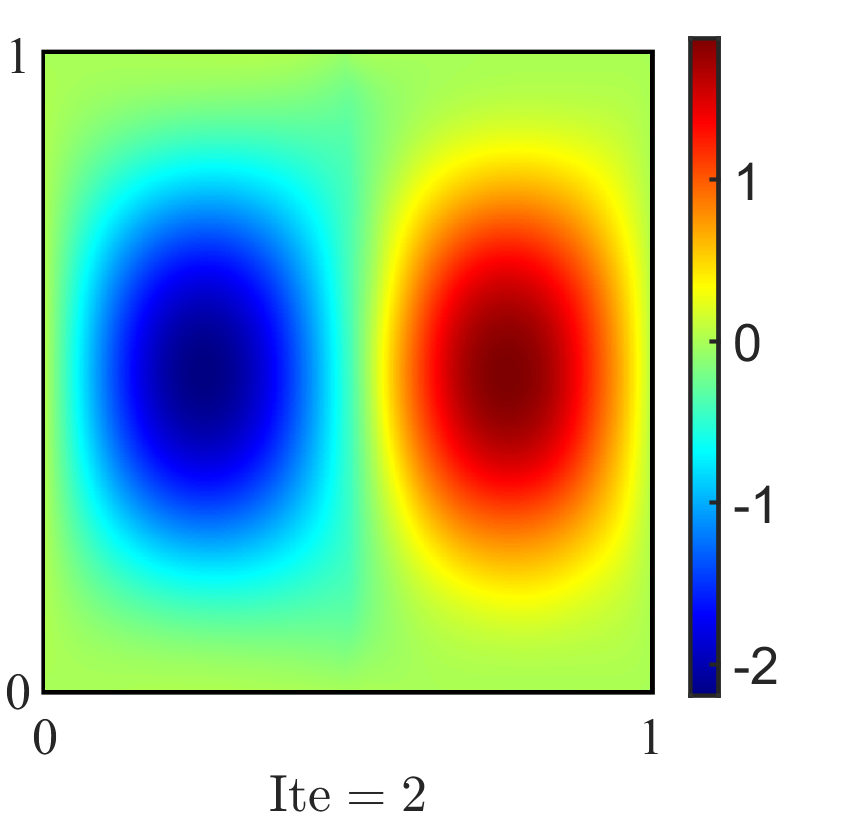}
\includegraphics[width=0.192\textwidth]{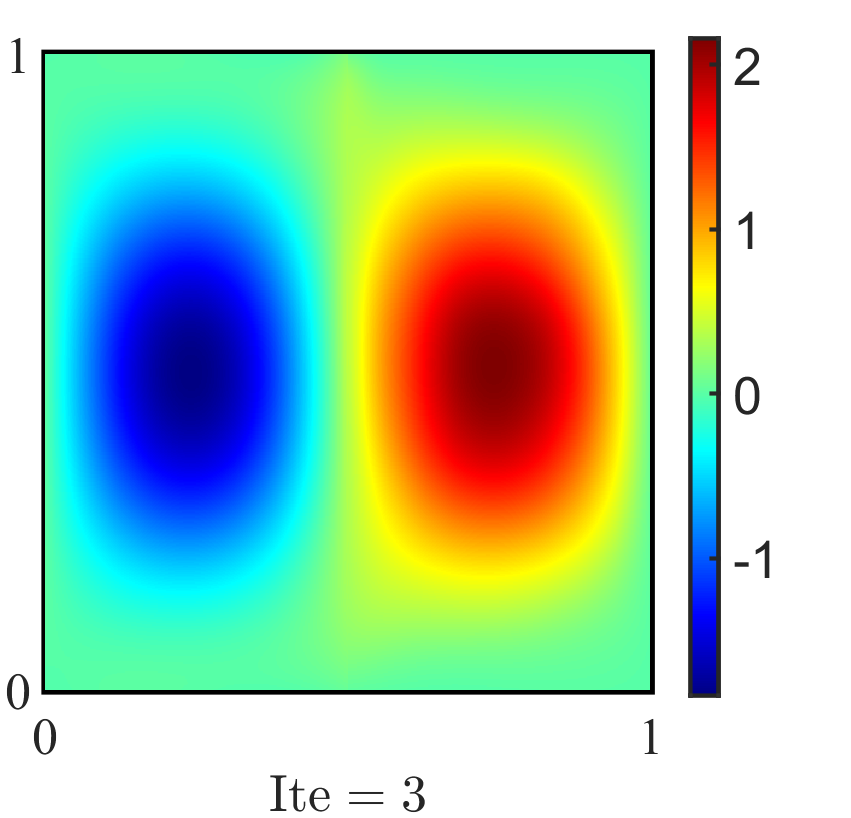}
\includegraphics[width=0.192\textwidth]{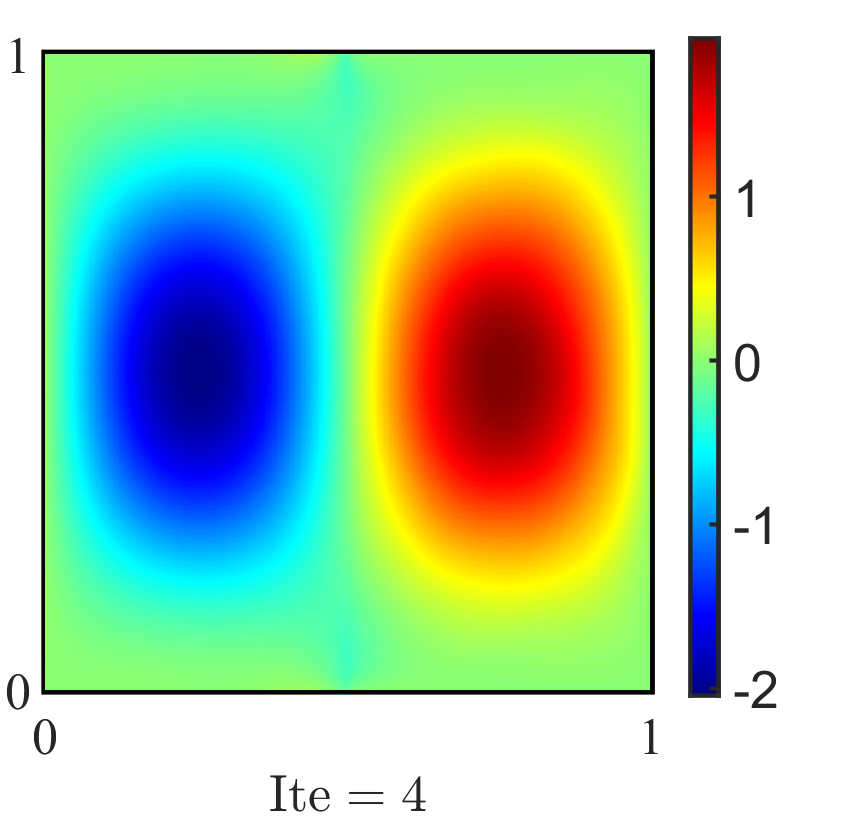}
\caption{Iterative solutions $\hat{u}^{[k]}(x,y)$ using RRLM (deep Ritz) along the outer iteration. }
\label{Experiments-RRLM-ex2-RRLM-PINN-error}
\end{subfigure}
\vspace{-0.45cm}
\caption{Numerical results of example \eqref{Experiments-RRLM-ex1} using RRLM for $(\kappa_1,\kappa_2)=(1,1000)$.}
\label{Experiments-RRLM-ex2-RRLM-PINN}
\vspace{-0.7cm}
\end{figure}

To accelerate the convergence of outer iterations, we set $(\kappa_1,\kappa_2)=(1,1000)$ in the following experiments. Unfortunately, the RR-PINNs approach suffers from the issue of weight imbalance and therefore fails to work (see \autoref{Experiments-RRLM-ex1-RR-PINNs-1-1000}). On the contrary, our compensated deep Ritz method (see \autoref{Experiments-RRLM-ex2-RRLM-PINN}) can handle the issue of weight imbalance and converge effectively, which only requires the replacement of the second subproblem solver with our proposed learning approach.

\section{Conclusion}
In this paper, a systematic study is presented for realizing classical non-overlapping DDMs through the use of artificial neural networks, which is based on the information exchange between neighbouring subproblems rather than domain partition strategies. For methods that rely on a direct flux exchange across subdomain interfaces, a key difficulty of deploying deep learning approaches as decomposed subproblem solvers is the issue of erroneous Dirichlet-to-Neumann map that always occurs to a greater or lesser extent in practice. To deal with the inaccurate flux estimation at interface, we develop a novel learning approach, i.e., the compensated deep Ritz method using neural network extension operators, to enable reliable flux transmission in the presence of erroneous interface conditions. As an immediate result, it allows us to construct effective learning approaches for realizing classical Dirichlet-Neumann, Neumann-Neumann, and Dirichlet-Dirichlet algorithms. On the other hand, the Robin-Robin algorithm, which only requires the exchange of Dirichlet traces but may suffer from the issue of weight imbalance, can also benefit from our compensated deep Ritz method. Finally, we conduct numerical experiments on a series of elliptic boundary value problems to demonstrate the effectiveness of our proposed learning algorithms. Possible future explorations would involve the coarse space acceleration \cite{mercier2021coarse}, adaptive sampling techniques \cite{he2022mesh},  efficient parallel iteration, and improvements of network architecture that could potentially further accelerate the convergence at a reduced cost.
\section{Acknowledgement}

We are grateful to the anonymous reviewers for their valuable feedback, which helped us improve the manuscript. This research was  conducted using computational resources and
services at the HPC center, School of Mathematical Sciences, Tongji University.




\bibliographystyle{siamplain}
\bibliography{references}
\end{document}


\maketitle

\section{Supplement to Section 5.1.2}

\begin{figure}[H]
\centering
\begin{subfigure}[htp]{\textwidth}
\centering
\includegraphics[width=0.192\textwidth]{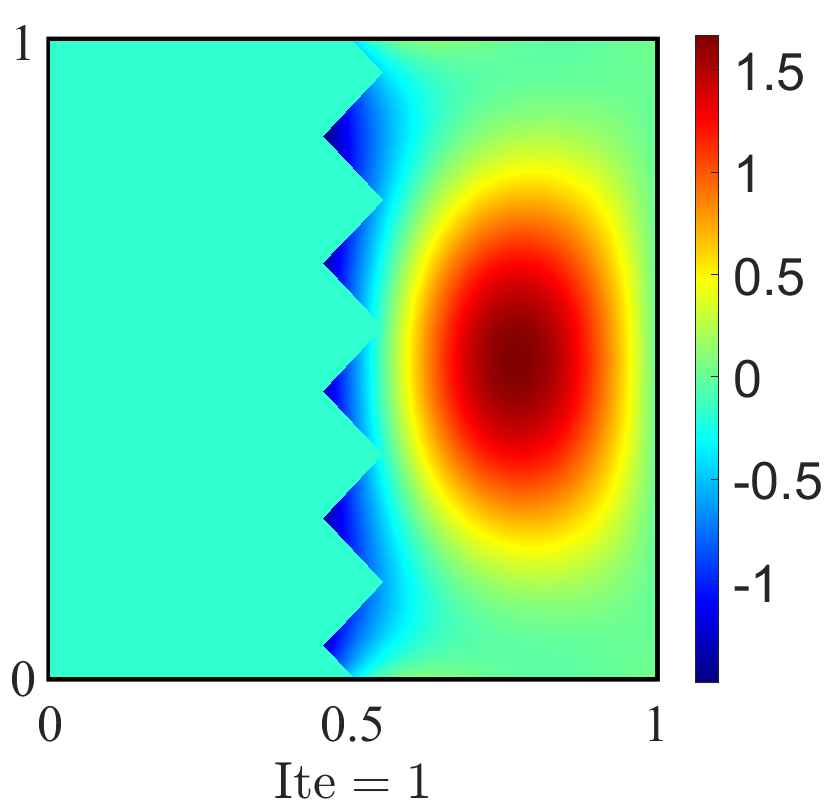}
\includegraphics[width=0.192\textwidth]{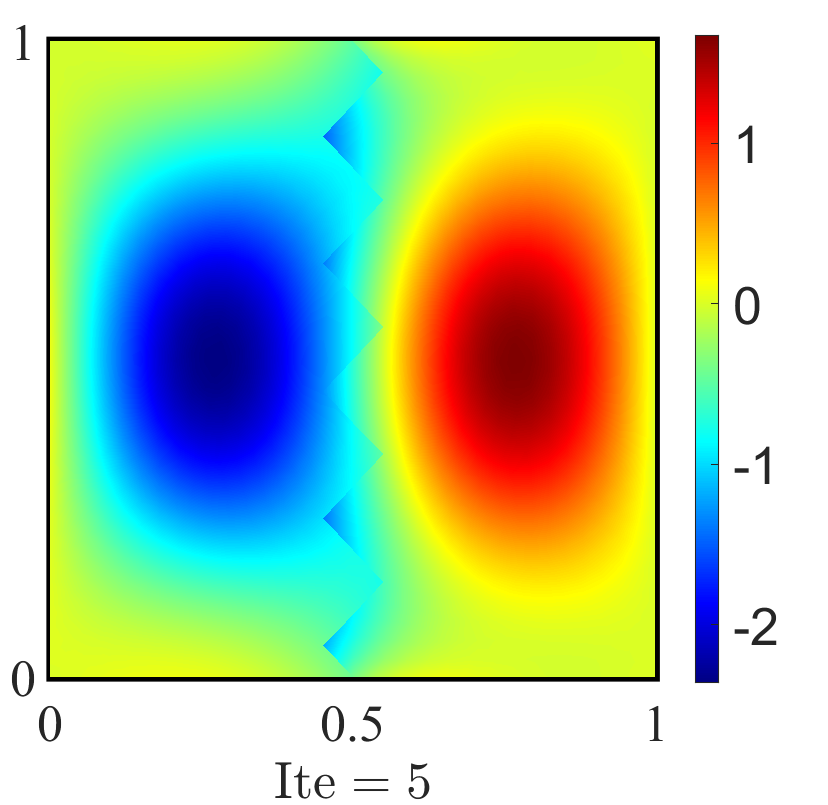}
\includegraphics[width=0.192\textwidth]{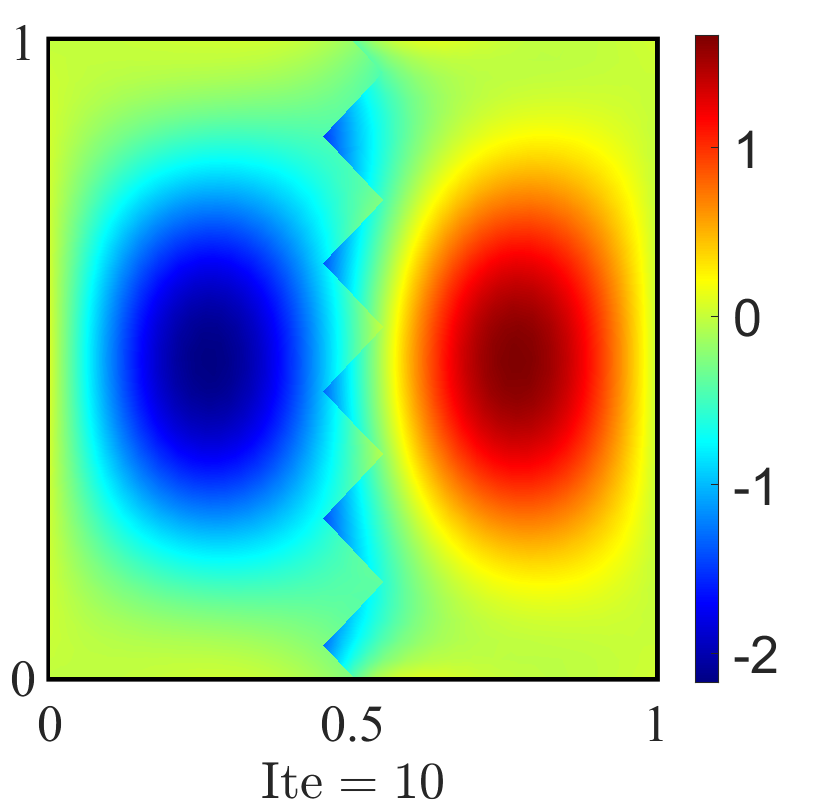}
\includegraphics[width=0.192\textwidth]{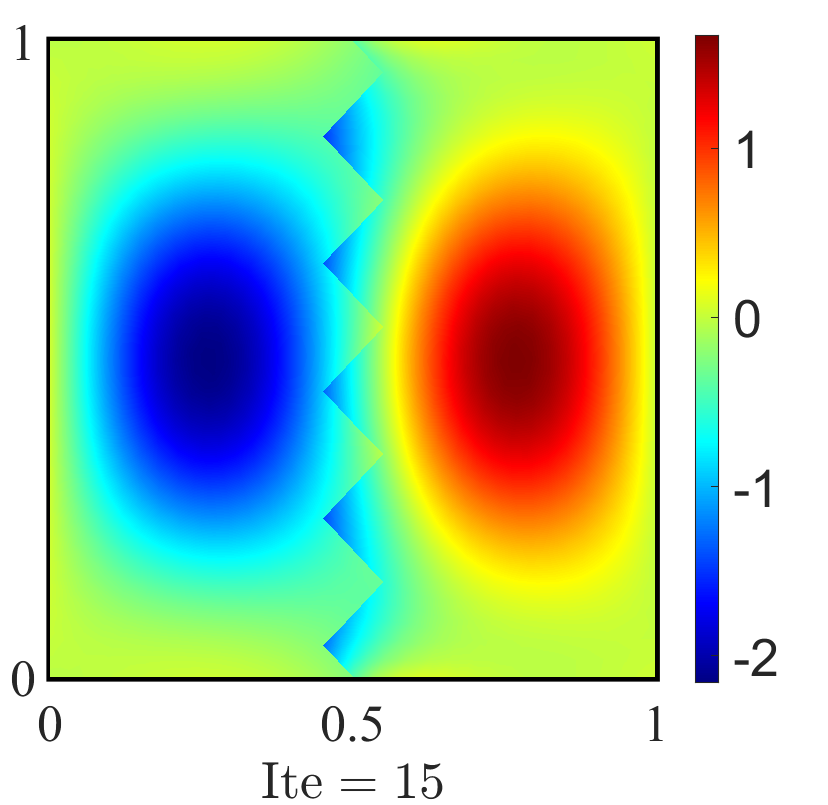}
\vspace{-0.1cm}
\caption{The numerical solutions $\hat{u}^{[k]}(x,y)$ along the outer iterations. }
\vspace{-0.2cm}
\end{subfigure}
\begin{subfigure}[htp]{\textwidth}
\centering
\includegraphics[width=0.192\textwidth]{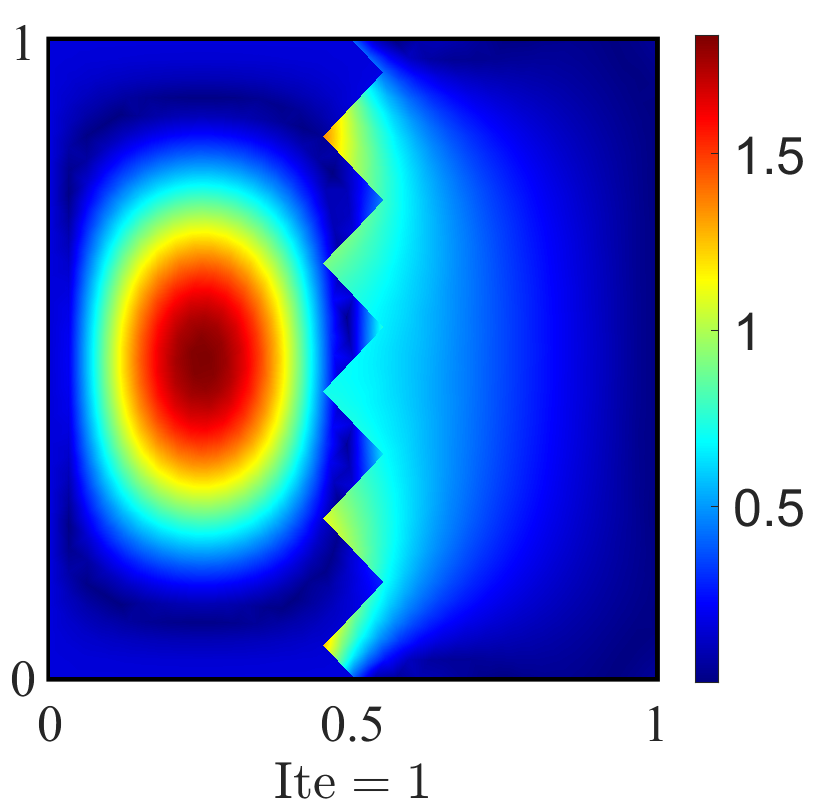}
\includegraphics[width=0.192\textwidth]{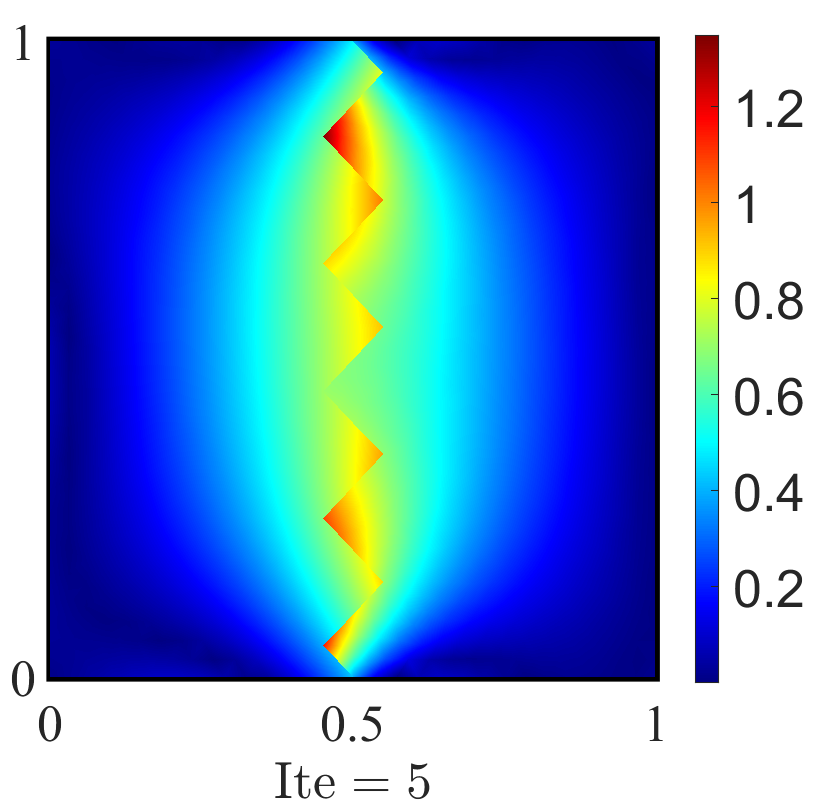}
\includegraphics[width=0.192\textwidth]{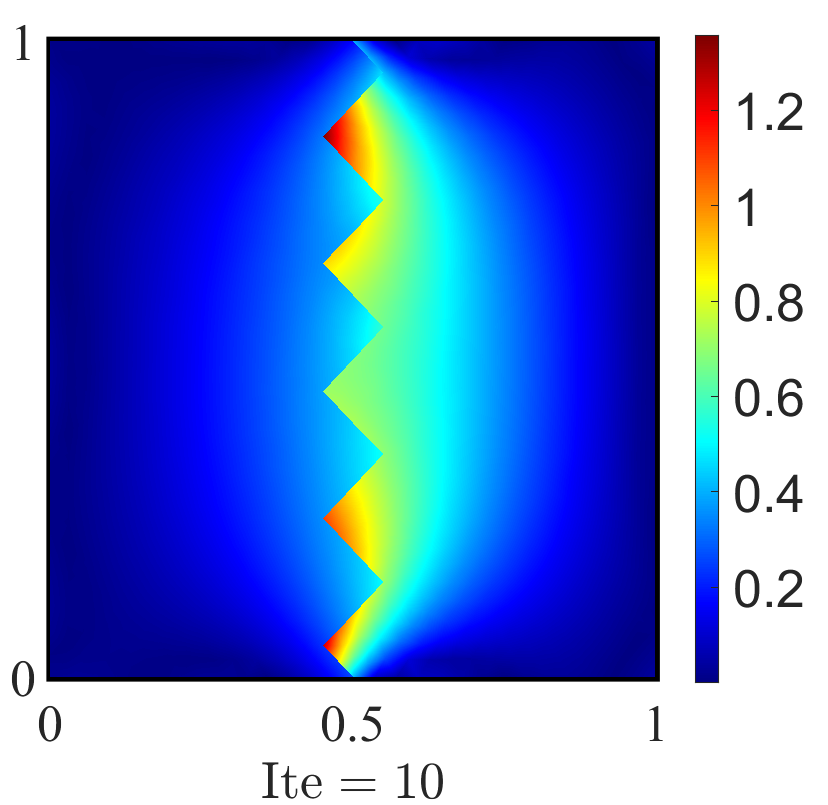}
\includegraphics[width=0.192\textwidth]{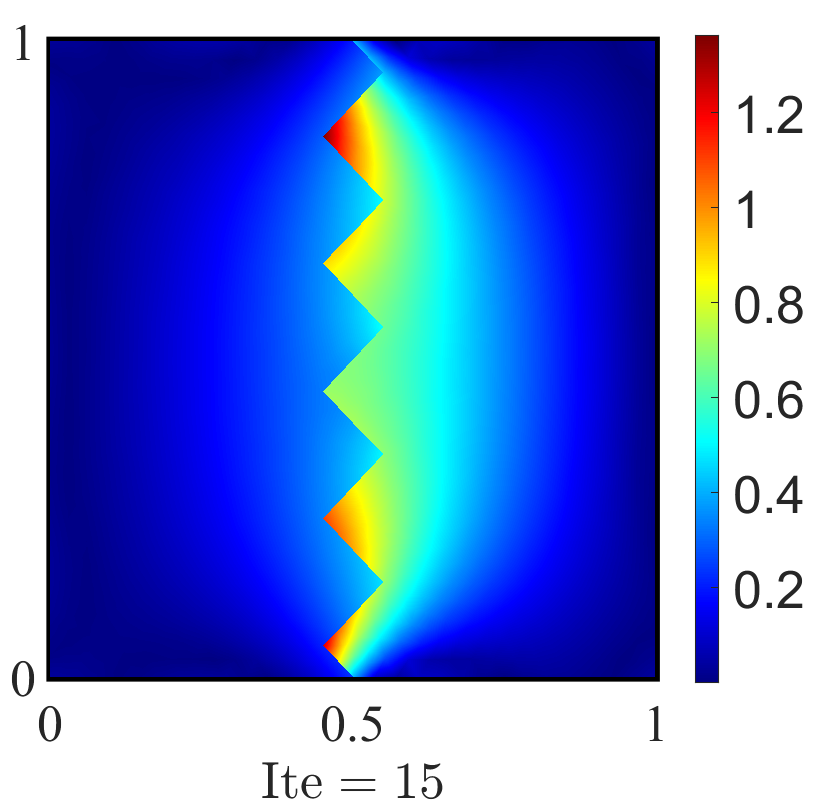}
\vspace{-0.1cm}
\caption{The pointwise absolute errors $|\hat{u}^{[k]}(x,y) - u(x,y)|$ along the outer iterations. }
\end{subfigure}
\vspace{-0.45cm}
\caption{Numerical results of example (5.2) using the DN-PINNs on testdata.}
\label{Experiments-DNLM-ex2-DN-PINNs}
\vspace{-0.6cm}
\end{figure}

\begin{figure}[H]
\centering
\begin{subfigure}[htp]{\textwidth}
\centering
\includegraphics[width=0.192\textwidth]{figure-DNLM//fig-DN-ex2-DNLM-PINN-u-NN-ite-1.png}
\includegraphics[width=0.192\textwidth]{figure-DNLM//fig-DN-ex2-DNLM-PINN-u-NN-ite-2.png}
\includegraphics[width=0.192\textwidth]{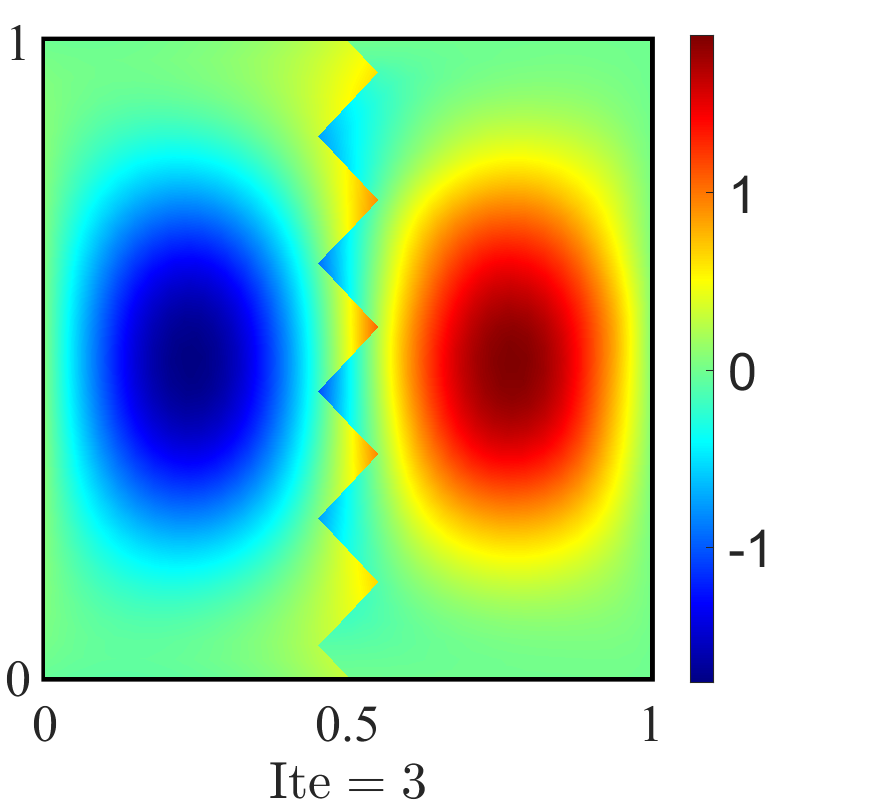}
\includegraphics[width=0.192\textwidth]{figure-DNLM//fig-DN-ex2-DNLM-PINN-u-NN-ite-4.png}
\vspace{-0.1cm}
\caption{The numerical solutions $\hat{u}^{[k]}(x,y)$ along the outer iterations.}
\vspace{-0.2cm}
\end{subfigure}
\begin{subfigure}[htp]{\textwidth}
\centering
\includegraphics[width=0.192\textwidth]{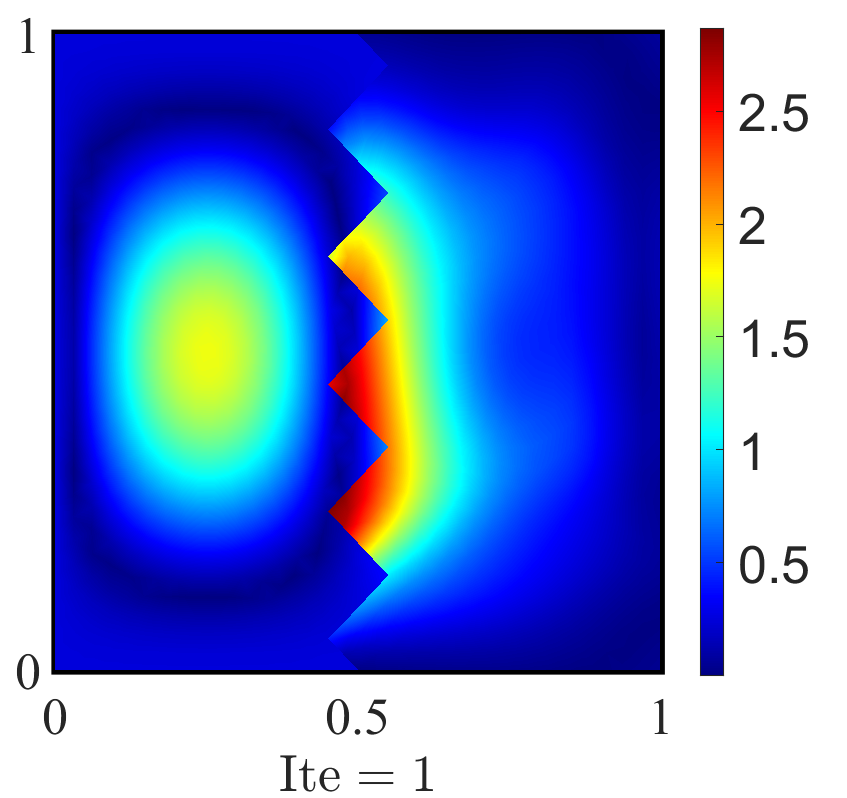}
\includegraphics[width=0.192\textwidth]{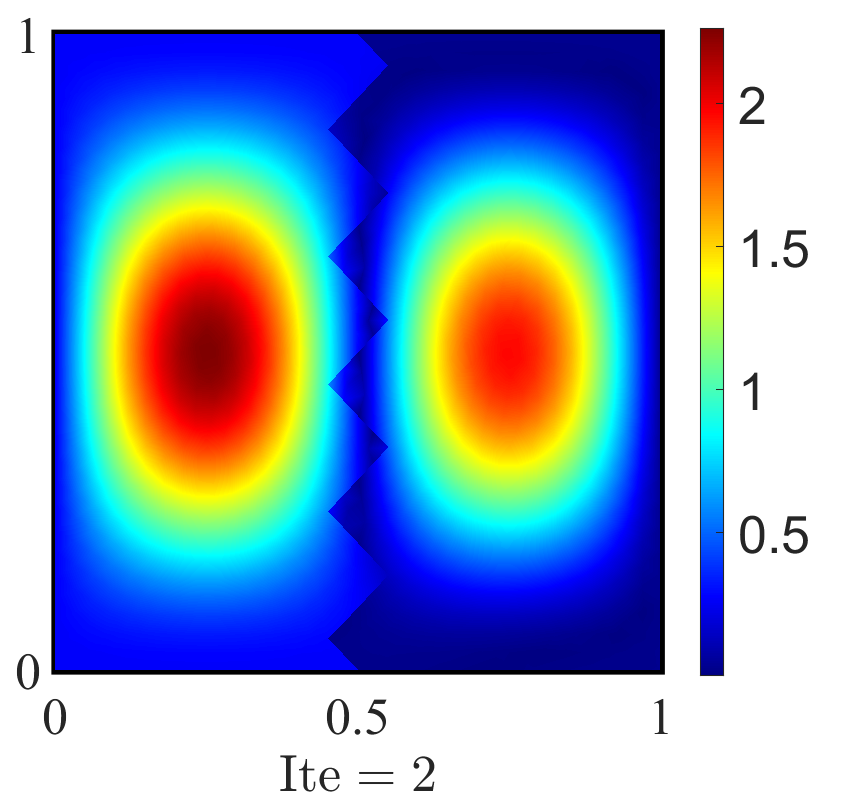}
\includegraphics[width=0.192\textwidth]{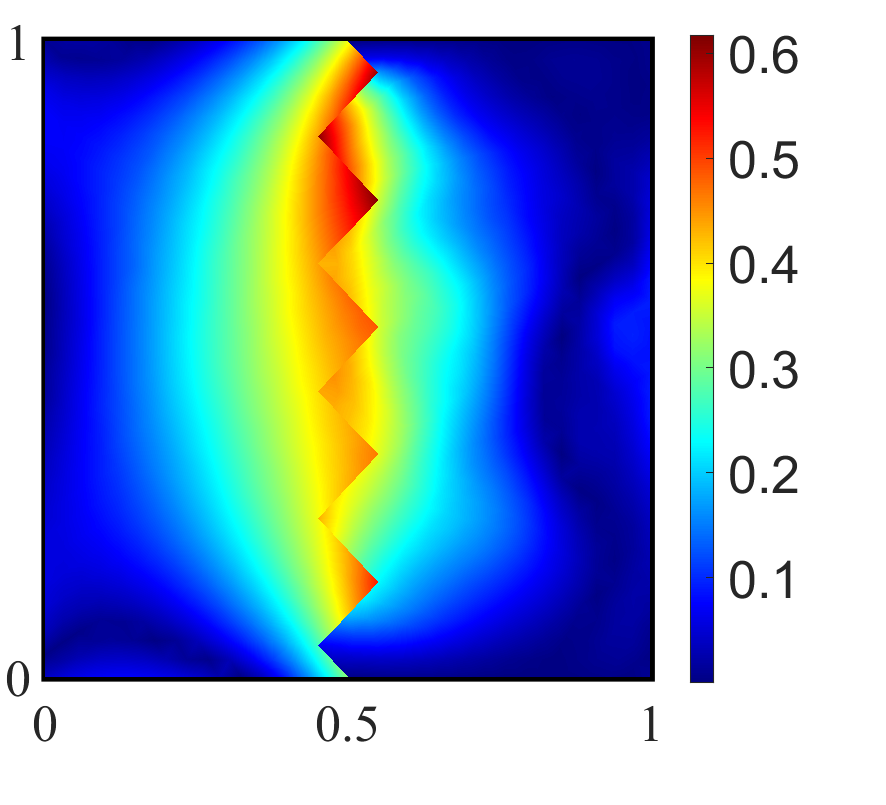}
\includegraphics[width=0.192\textwidth]{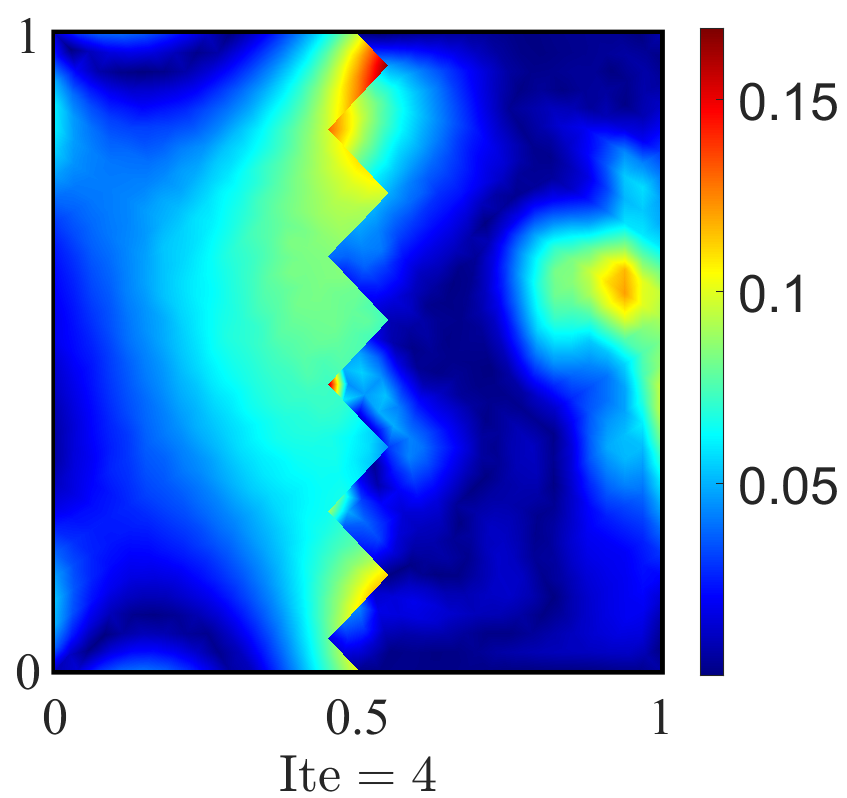}
\vspace{-0.1cm}
\caption{The pointwise absolute errors $|\hat{u}^{[k]}(x,y) - u(x,y)|$ along the outer iterations. }
\end{subfigure}
\vspace{-0.45cm}
\caption{Numerical results of example (5.2) using our DNLM (PINN) on testdata.}
\label{Experiments-DNLM-ex2-DNLM-PINN}
\end{figure}

The numerical results using the DN-PINNs scheme in a typical simulation are depicted in \autoref{Experiments-DNLM-ex2-DN-PINNs}, which fails to converge to the exact solution. Our proposed DNLM (PINN) can facilitate the convergence of outer iteration in the presence of interface overfitting (see \autoref{Experiments-DNLM-ex2-DNLM-PINN}).

\section{Supplement to Section 5.1.3}

\begin{figure}[htp]
\centering
\includegraphics[width=0.178\textwidth]{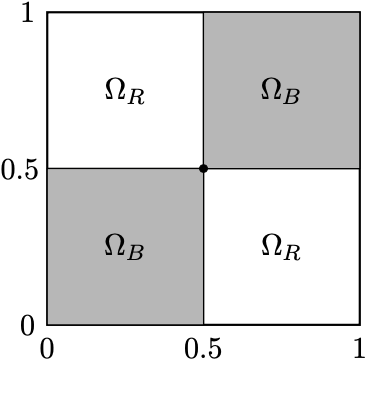}
\hspace{0.15cm}
\includegraphics[width=0.2\textwidth]{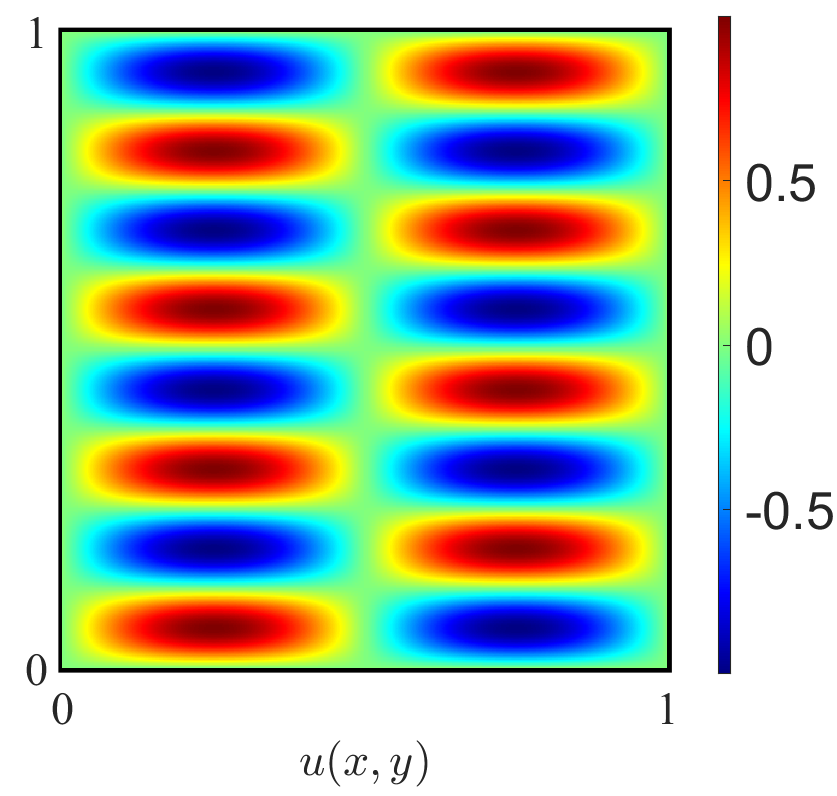}
\includegraphics[width=0.2\textwidth]{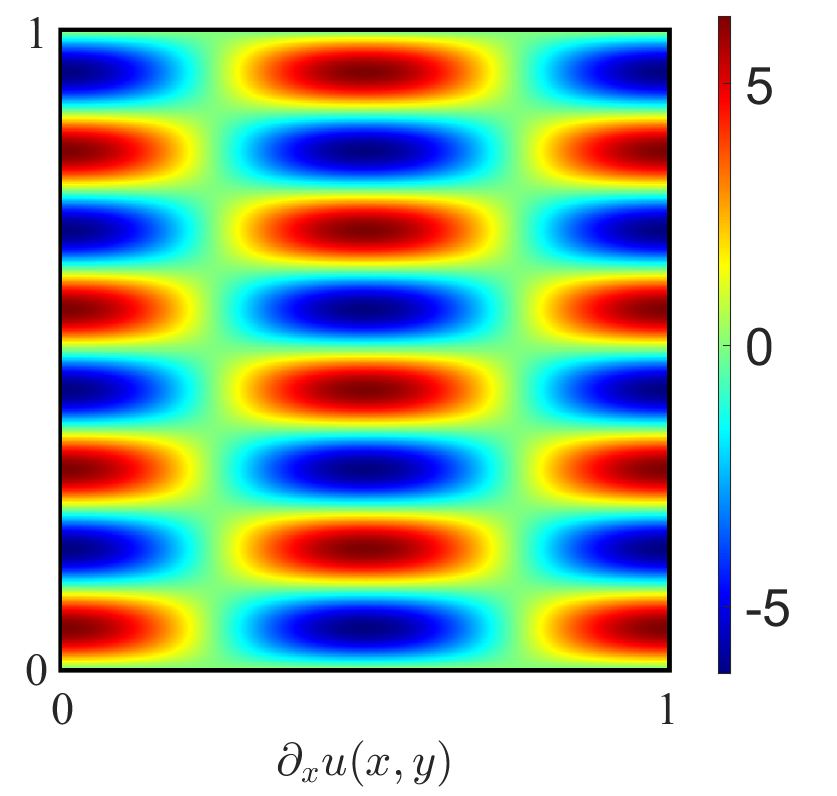}
\includegraphics[width=0.2\textwidth]{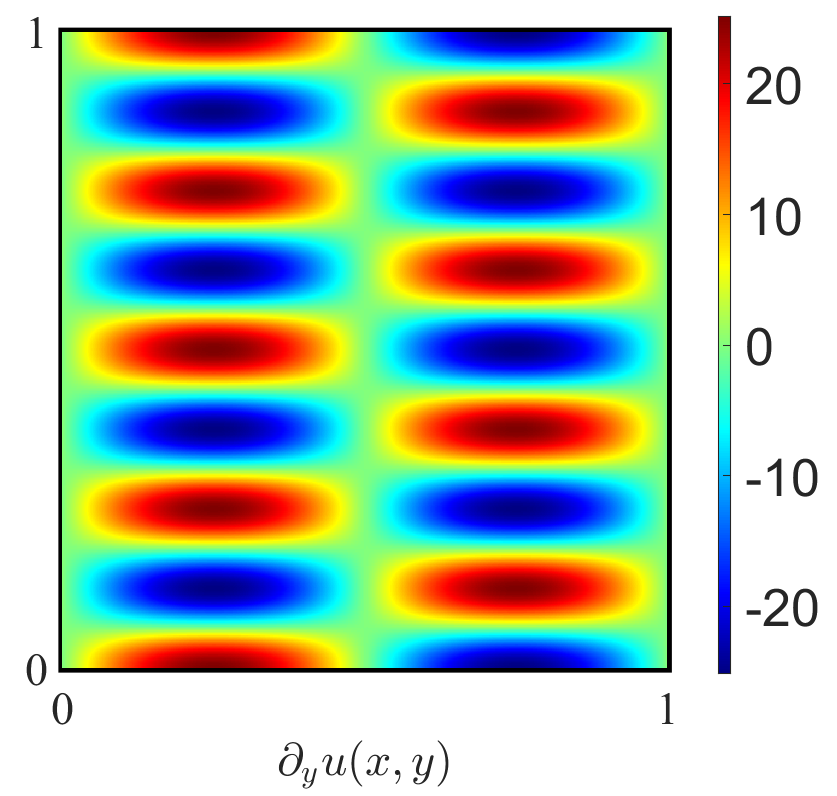}
\caption{From left to right: decomposition of domain into two subregions, exact solution $u(x,y)$ and its partial derivatives $\partial_x u(x,y)$, $\partial_y u(x,y)$ for example (5.3).}
\label{Experiments-DNLM-ex3-exact-solution}
\vspace{-0.5cm}
\end{figure}

The decomposition of domain and the exact solution for example (5.3) are shown in \autoref{Experiments-DNLM-ex3-exact-solution}. 

\begin{figure}[H]
\centering
\begin{subfigure}[htp]{\textwidth}
\centering
\includegraphics[width=0.192\textwidth]{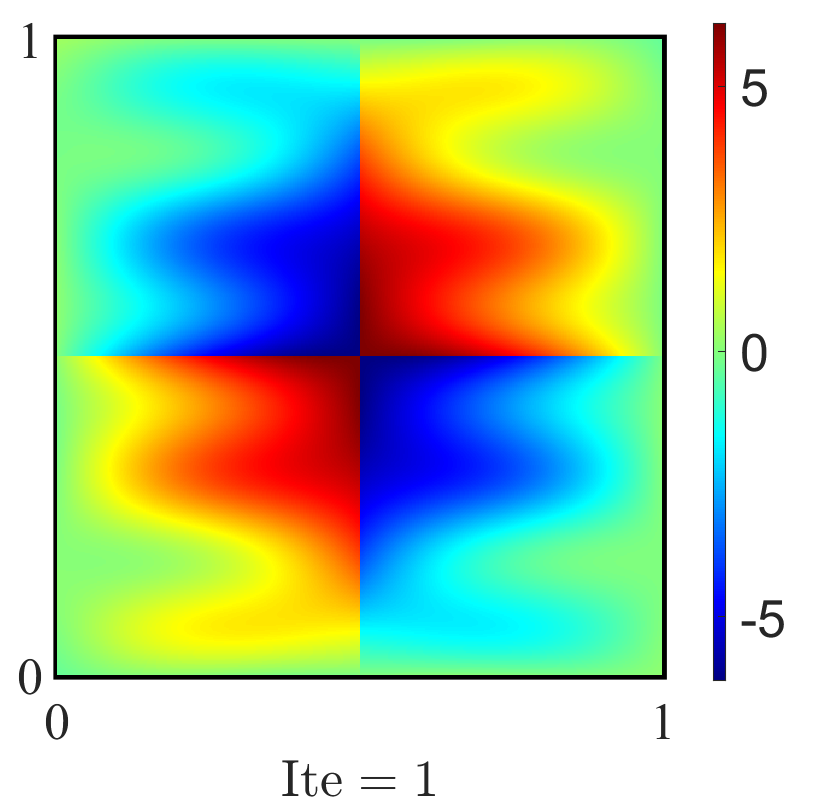}
\includegraphics[width=0.192\textwidth]{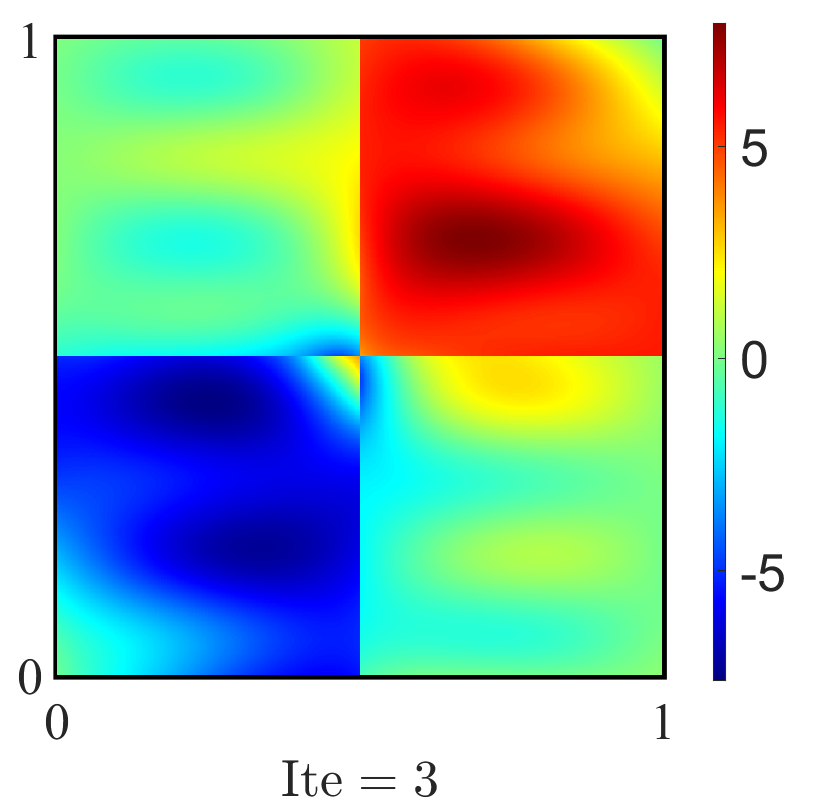}
\includegraphics[width=0.192\textwidth]{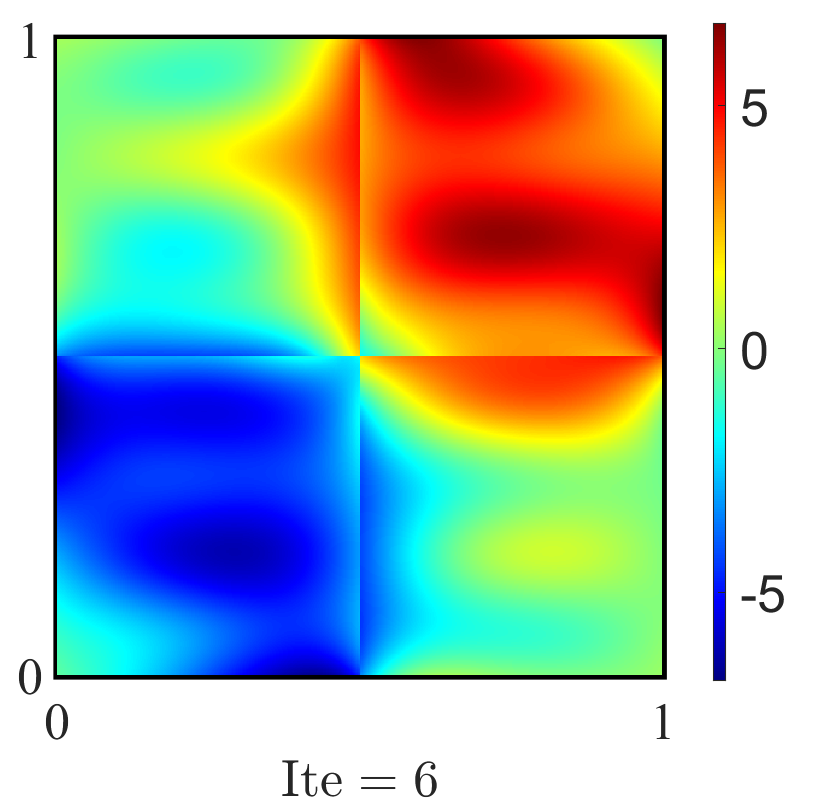}
\includegraphics[width=0.192\textwidth]{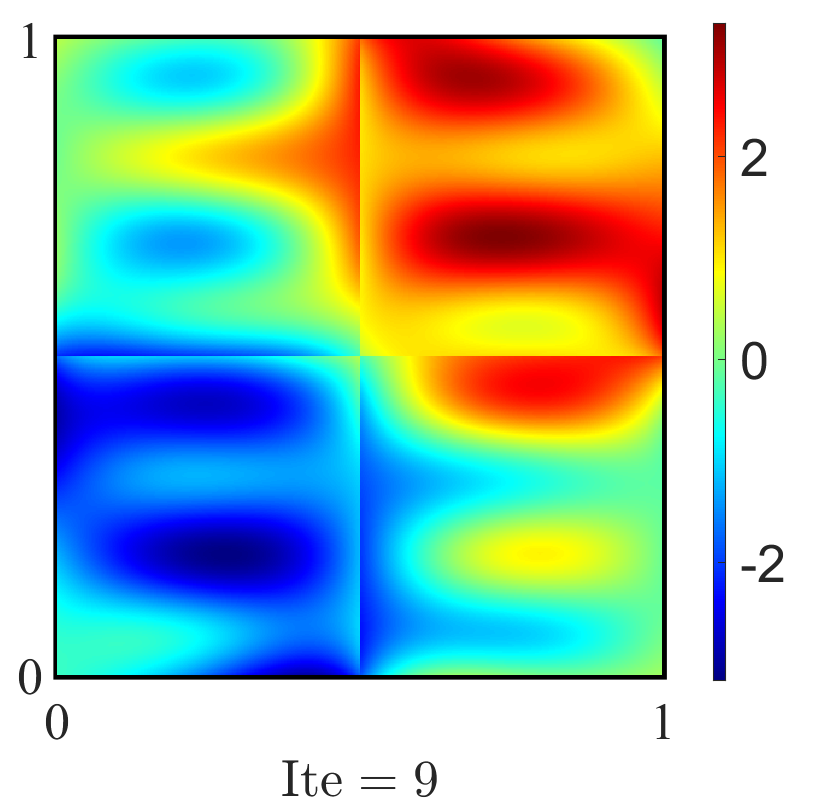}
\vspace{-0.1cm}
\caption{The numerical solutions $\hat{u}^{[k]}(x,y)$ along the outer iterations. }
\vspace{-0.2cm}
\end{subfigure}
\begin{subfigure}[htp]{\textwidth}
\centering
\includegraphics[width=0.192\textwidth]{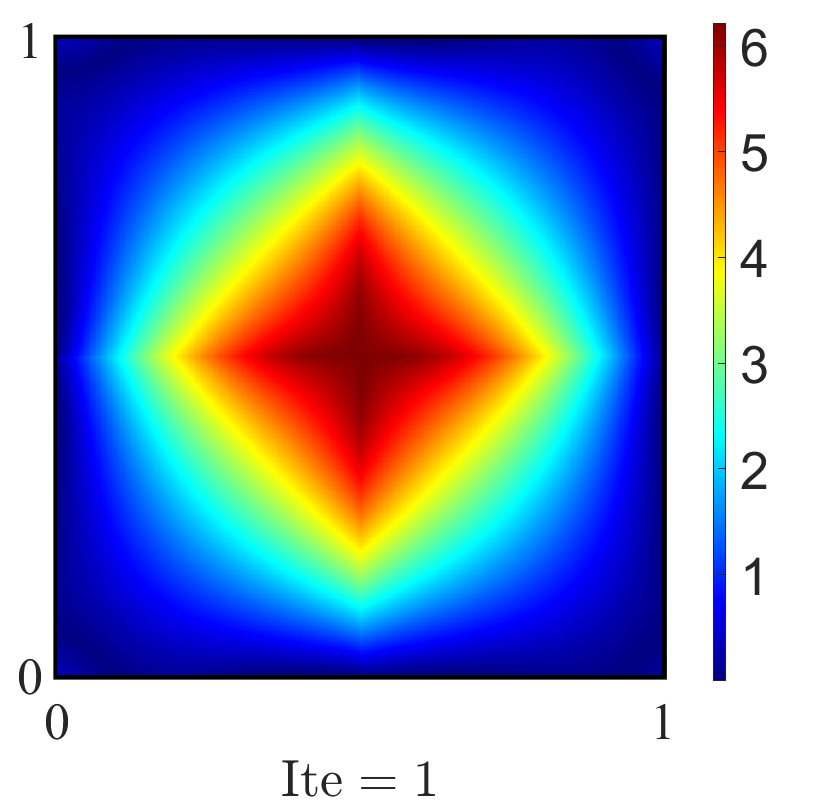}
\includegraphics[width=0.192\textwidth]{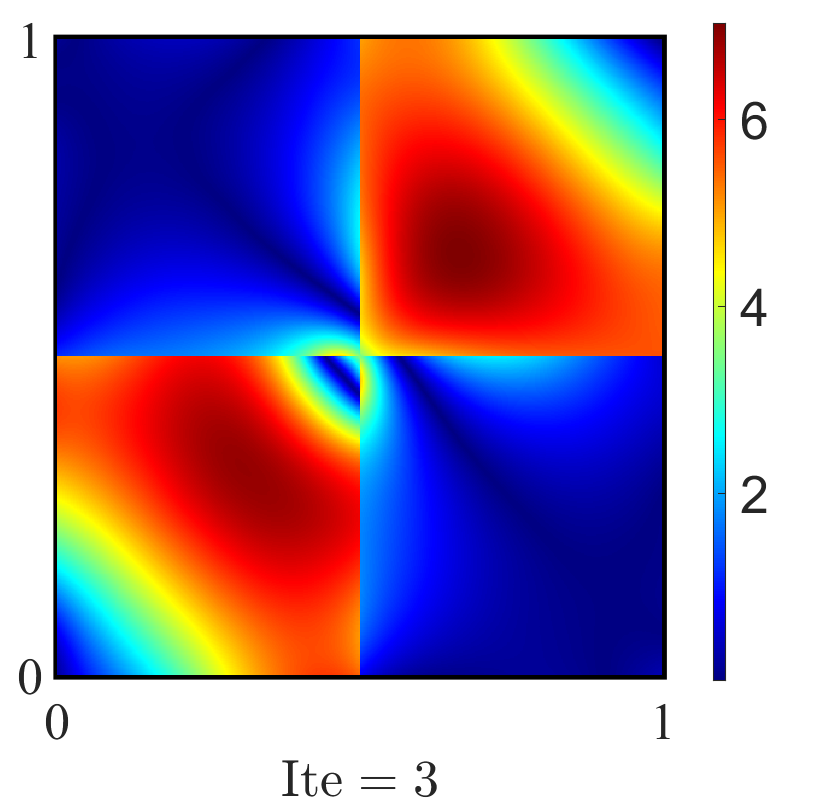}
\includegraphics[width=0.192\textwidth]{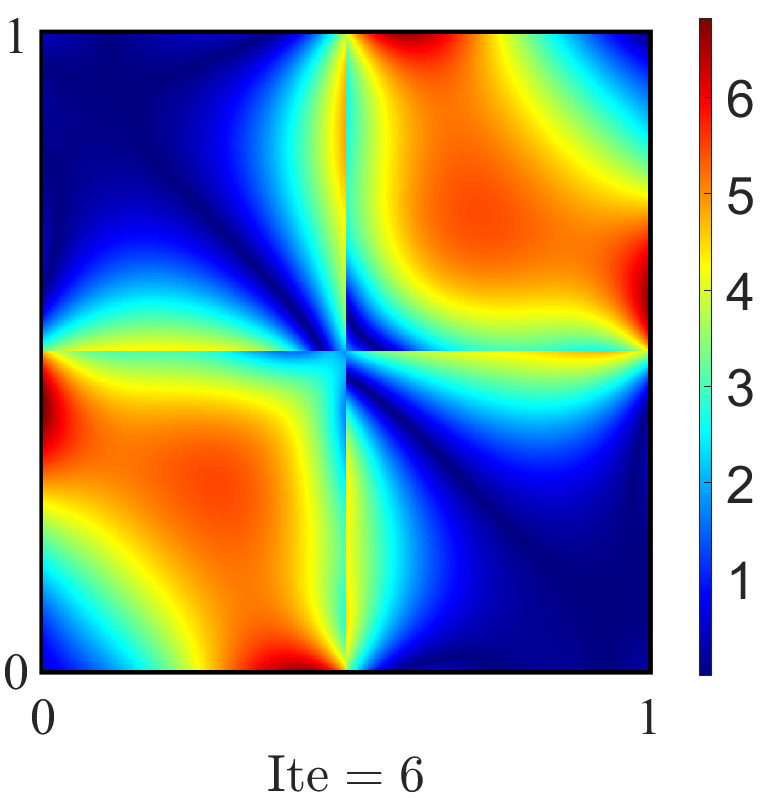}
\includegraphics[width=0.192\textwidth]{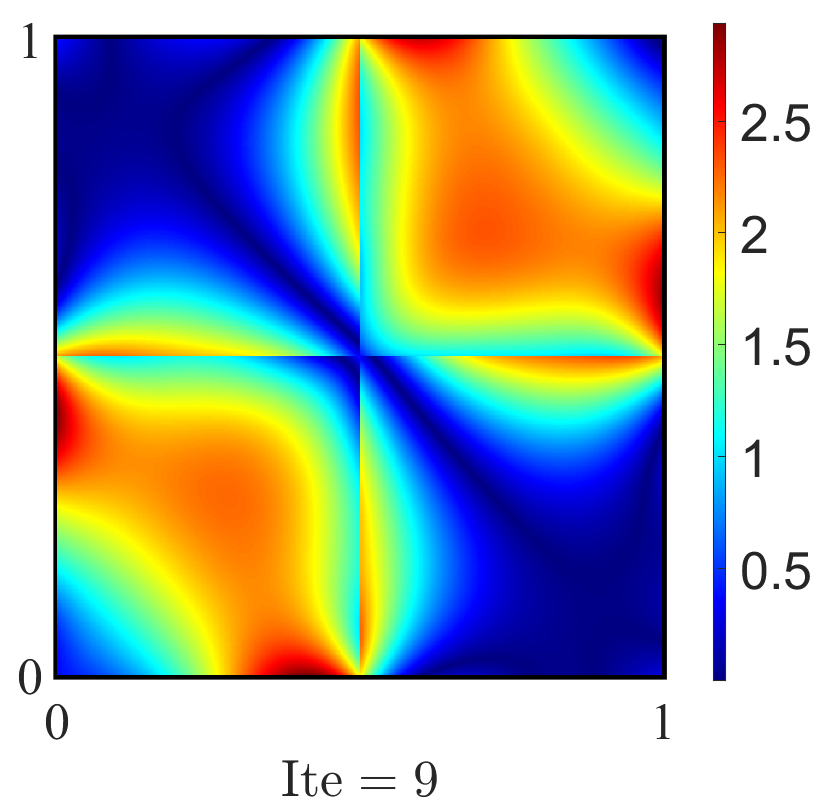}
\vspace{-0.1cm}
\caption{The pointwise absolute errors $|\hat{u}^{[k]}(x,y) - u(x,y)|$ along the outer iterations. }
\end{subfigure}
\vspace{-0.45cm}
\caption{Numerical results of example (5.3) using the DN-PINNs on testdata.}
\label{Experiments-DNLM-ex3-DN-PINNs}
\vspace{-0.7cm}
\end{figure}

\begin{figure}[H]
\centering
\begin{subfigure}[htp]{\textwidth}
\centering
\includegraphics[width=0.192\textwidth]{figure-DNLM//fig-DN-ex3-DNLM-PINN-u-NN-ite-1.png}
\includegraphics[width=0.192\textwidth]{figure-DNLM//fig-DN-ex3-DNLM-PINN-u-NN-ite-2.png}
\includegraphics[width=0.192\textwidth]{figure-DNLM//fig-DN-ex3-DNLM-PINN-u-NN-ite-4.png}
\includegraphics[width=0.192\textwidth]{figure-DNLM//fig-DN-ex3-DNLM-PINN-u-NN-ite-6.png}
\vspace{-0.1cm}
\caption{The numerical solutions $\hat{u}^{[k]}(x,y)$ along the outer iterations. }
\vspace{-0.2cm}
\end{subfigure}
\begin{subfigure}[htp]{\textwidth}
\centering
\includegraphics[width=0.192\textwidth]{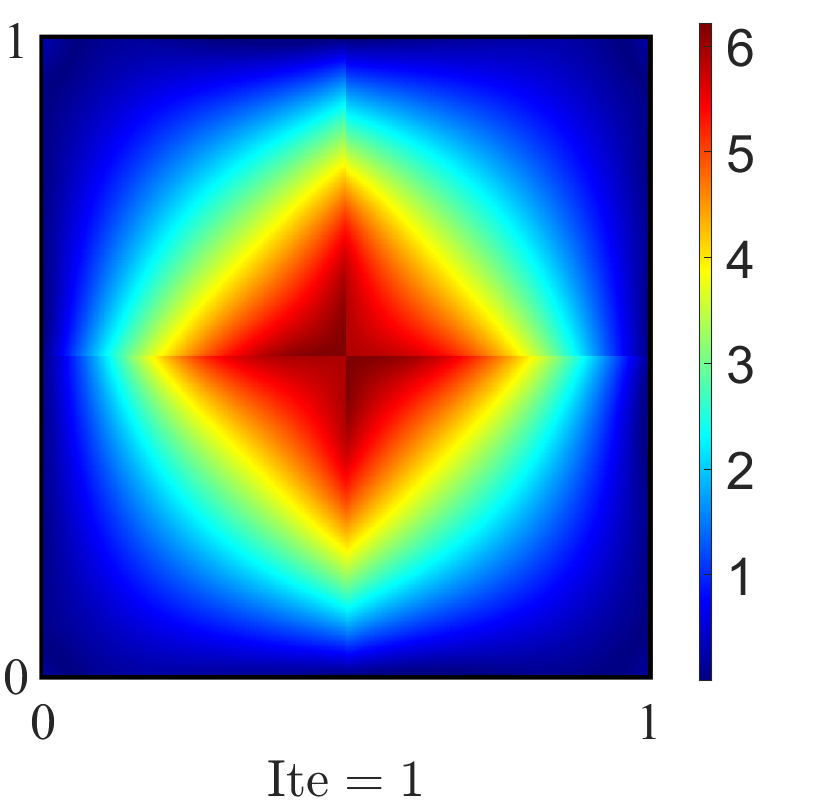}
\includegraphics[width=0.192\textwidth]{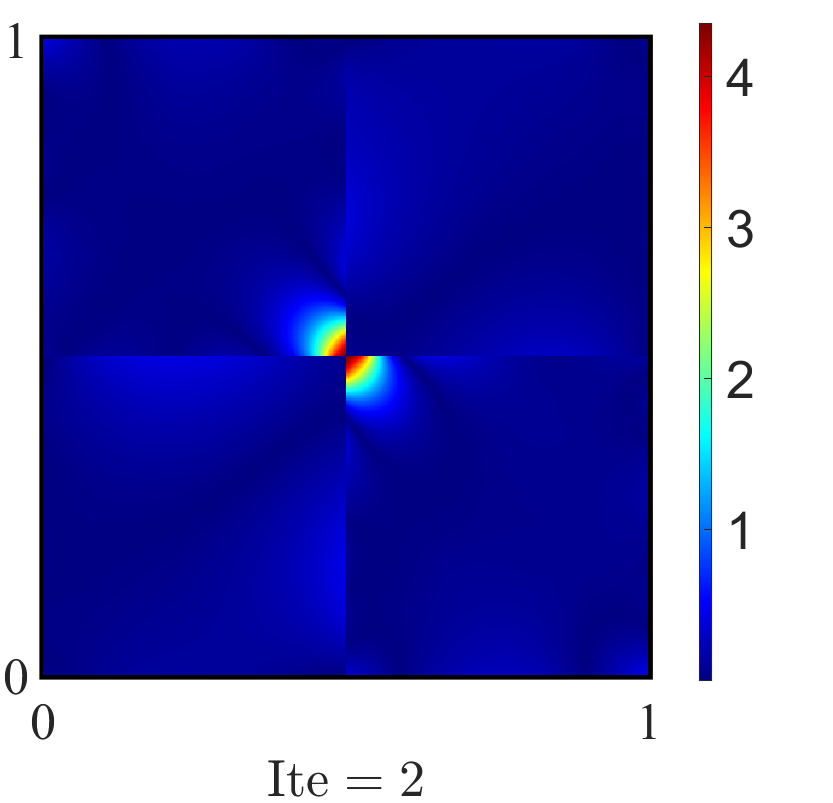}
\includegraphics[width=0.192\textwidth]{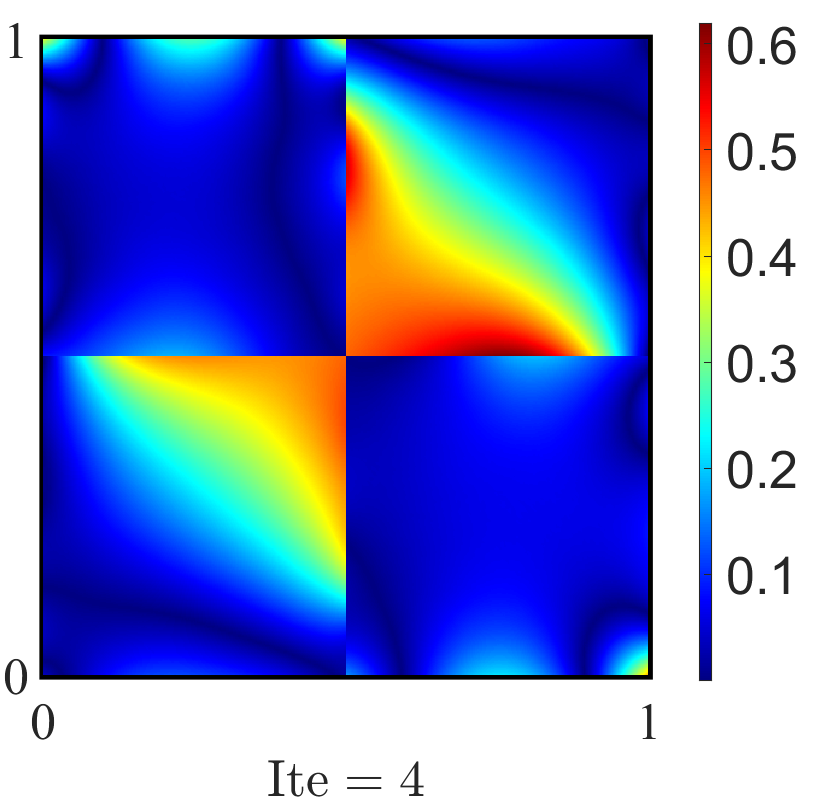}
\includegraphics[width=0.192\textwidth]{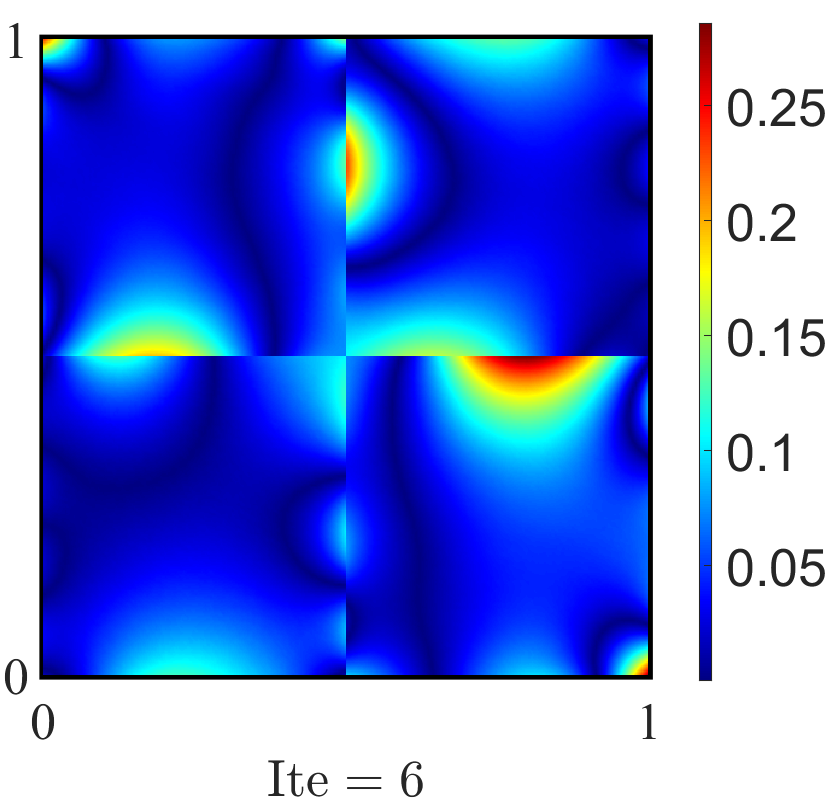}
\vspace{-0.1cm}
\caption{The pointwise absolute errors $|\hat{u}^{[k]}(x,y) - u(x,y)|$ along the outer iterations. }
\end{subfigure}
\vspace{-0.45cm}
\caption{Numerical results of example (5.3) using our DNLM (PINN) on testdata.}
\label{Experiments-DNLM-ex3-DNLM-PINN}
\vspace{-0.5cm}
\end{figure}

The numerical results using DN-PINNs and DNLM (PINN) are depicted in \autoref{Experiments-DNLM-ex3-DN-PINNs} and \autoref{Experiments-DNLM-ex3-DNLM-PINN}, which implies that our proposed method can converge to the exact solution while the DN-PINNs scheme fails.

\section{Supplement to Section 5.1.5}

\begin{figure}[htp]
\centering
\includegraphics[width=0.170\textwidth]{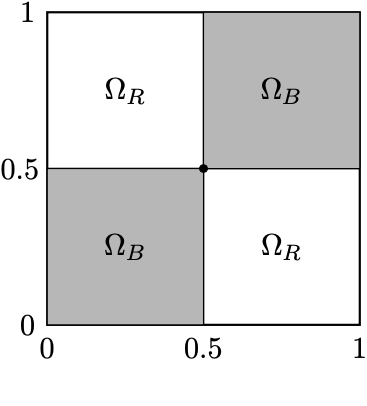}
\hspace{0.15cm}
\includegraphics[width=0.196\textwidth]{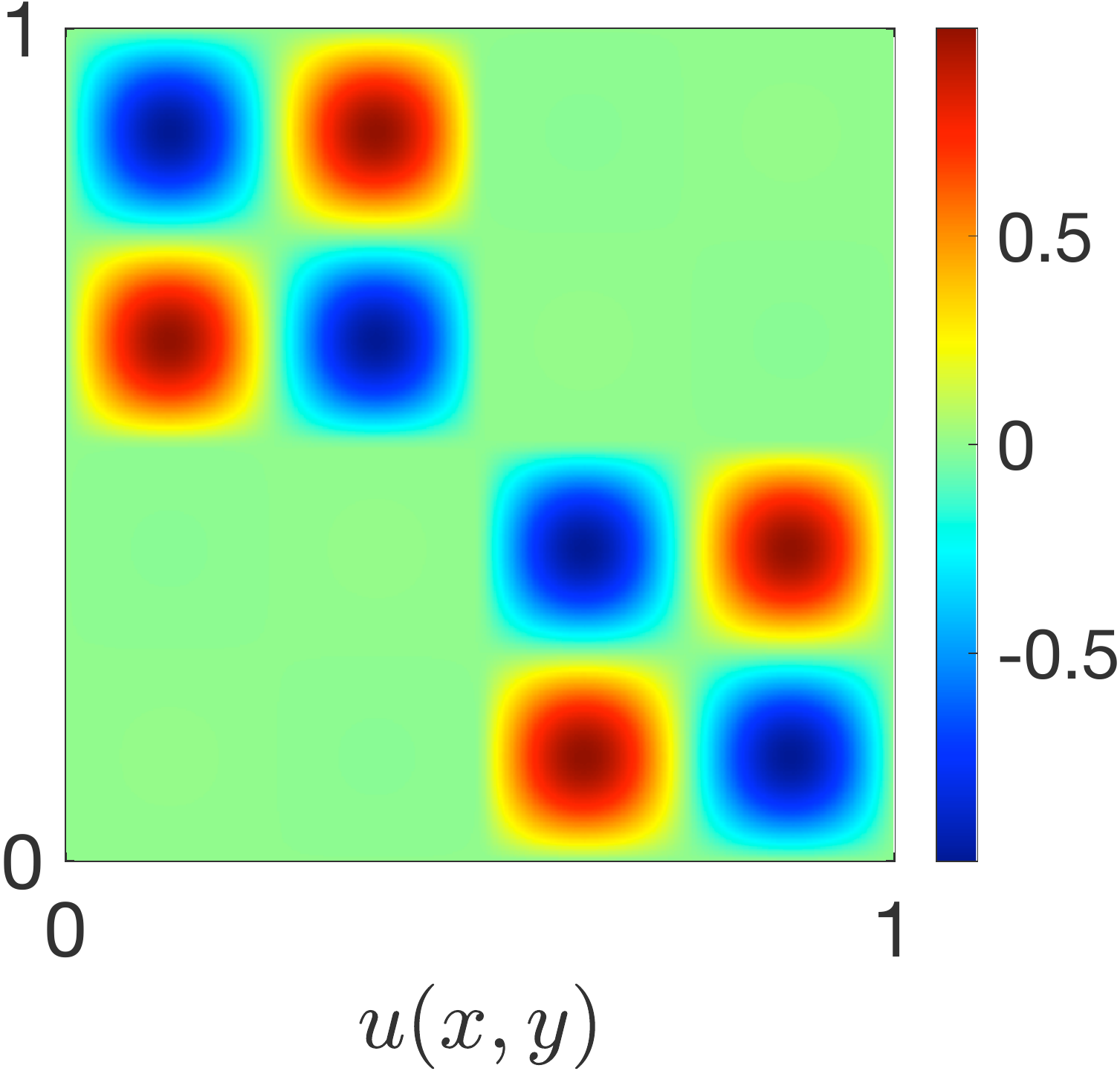}
\includegraphics[width=0.192\textwidth]{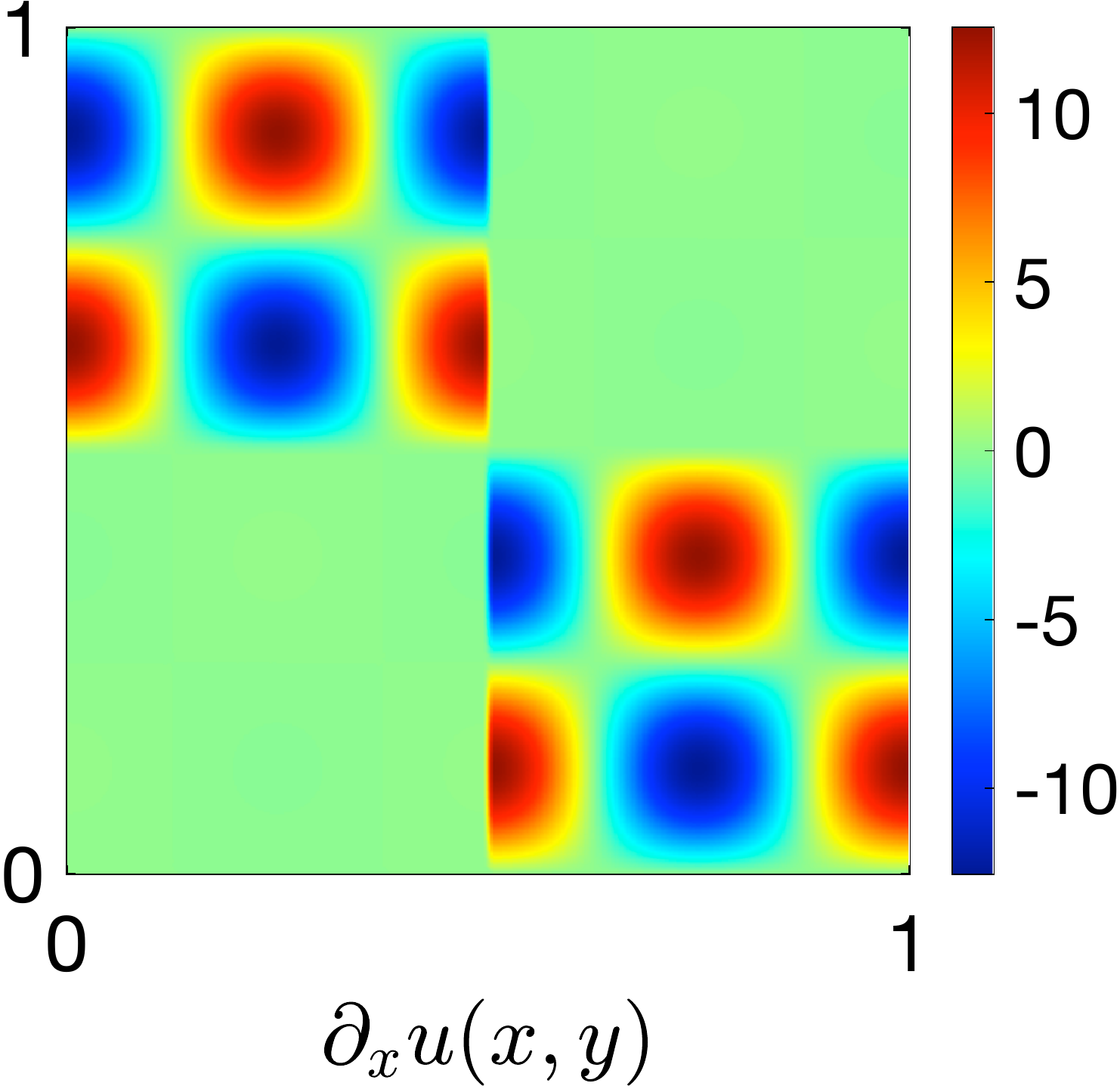}
\includegraphics[width=0.192\textwidth]{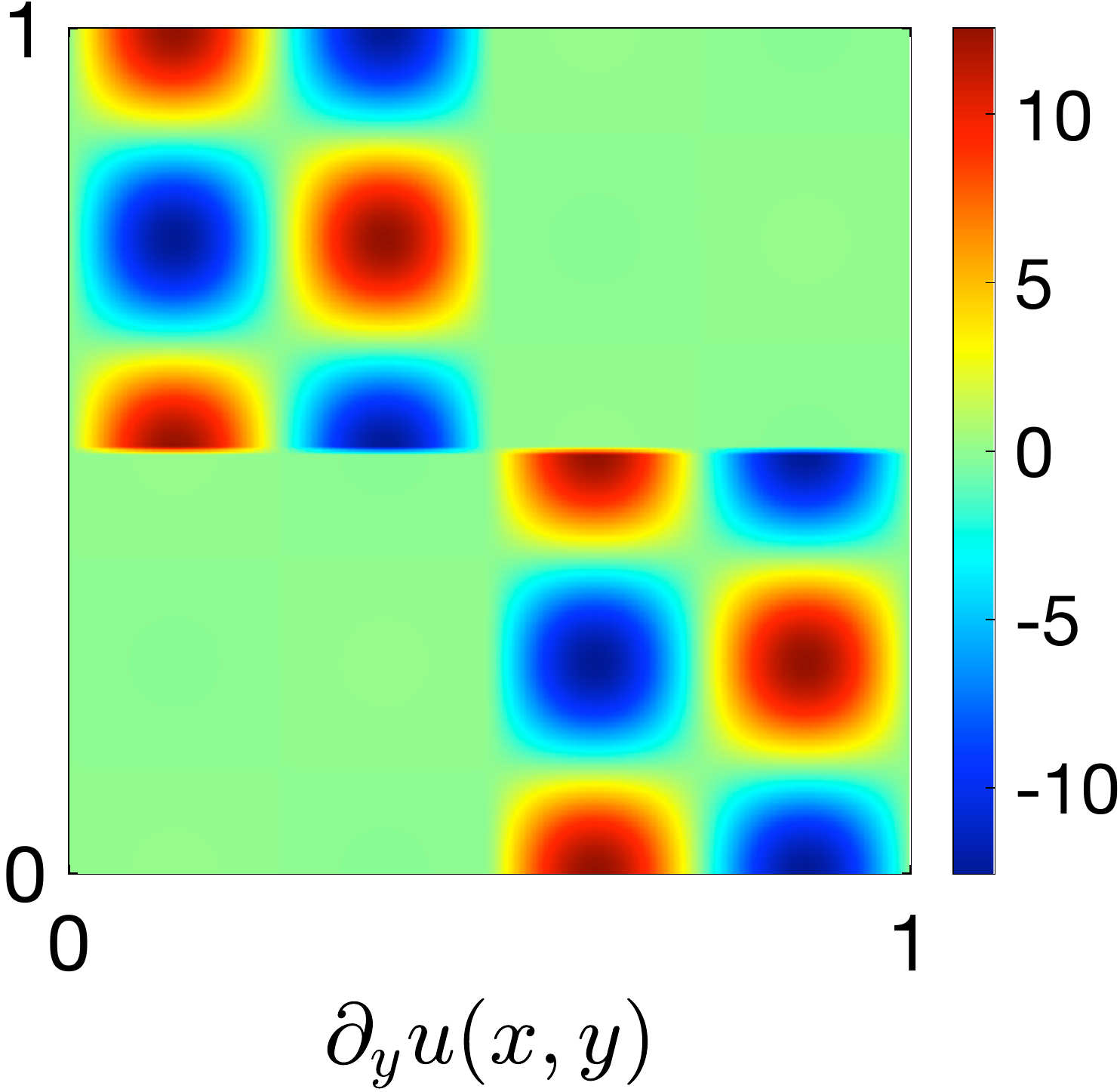}
\caption{From left to right: decomposition of domain into two subregions, exact solution $u(x,y)$ and its partial derivatives $\partial_x u(x,y)$, $\partial_y u(x,y)$ for example (5.5).}
\label{Experiments-DNLM-ex5-exact-solution}
\vspace{-0.5cm}
\end{figure}

The decomposition of domain and the exact solution for example (5.5) are shown in \autoref{Experiments-DNLM-ex5-exact-solution}. 

\begin{figure}[H]
\centering
\begin{subfigure}[htp]{\textwidth}
\centering
\includegraphics[width=0.192\textwidth]{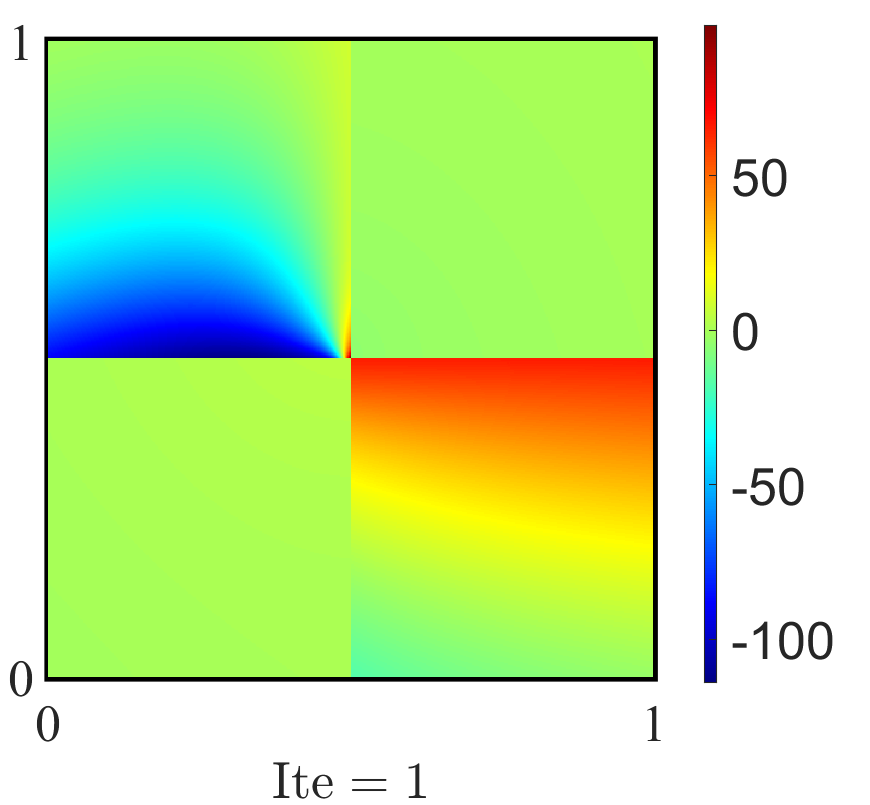}
\includegraphics[width=0.192\textwidth]{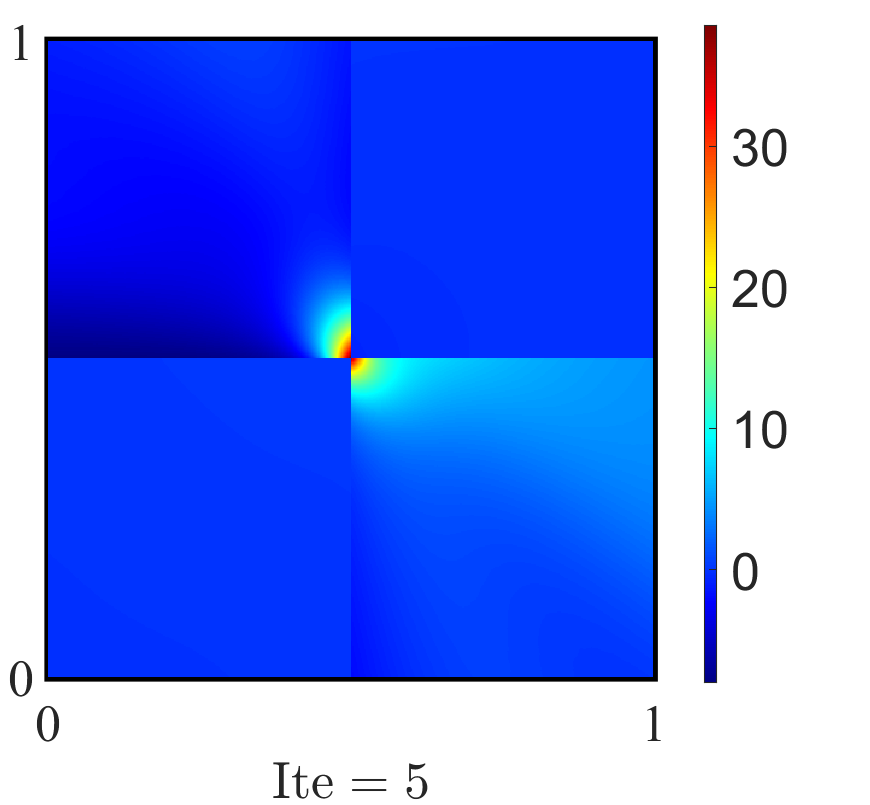}
\includegraphics[width=0.192\textwidth]{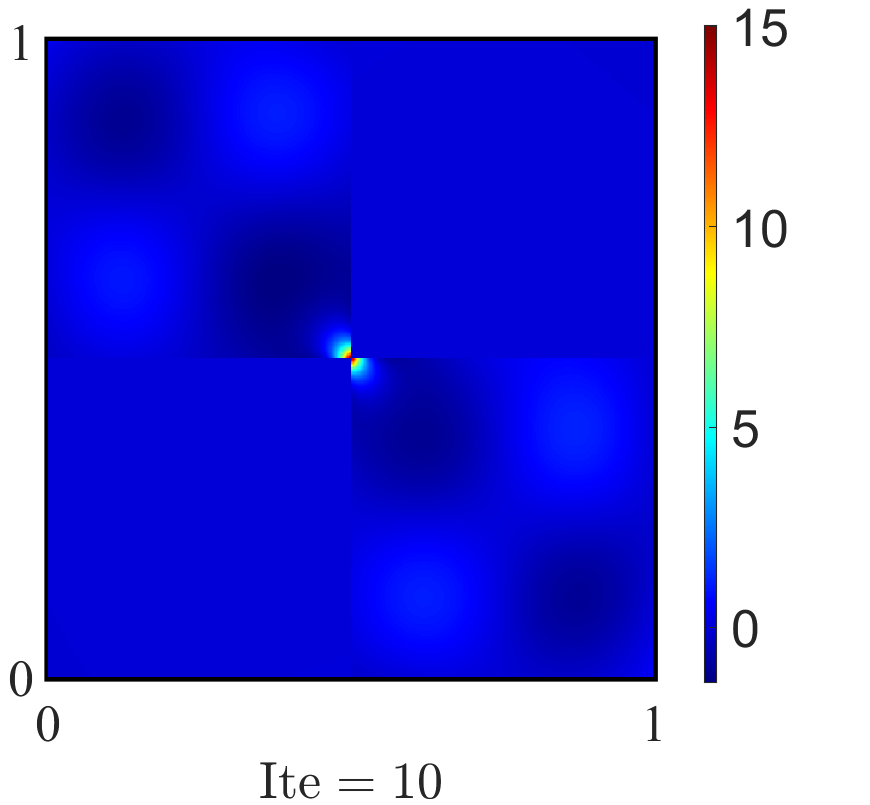}
\includegraphics[width=0.192\textwidth]{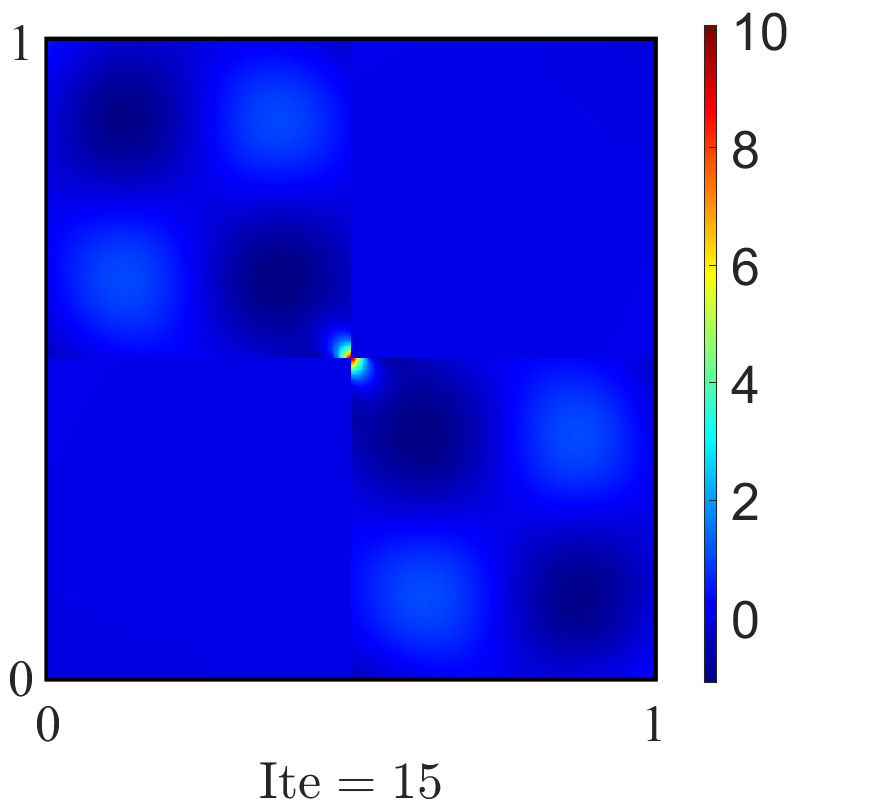}
\vspace{-0.1cm}
\caption{The numerical solutions $\hat{u}^{[k]}(x,y)$ along the outer iterations. }
\vspace{-0.2cm}
\end{subfigure}
\begin{subfigure}[htp]{\textwidth}
\centering
\includegraphics[width=0.192\textwidth]{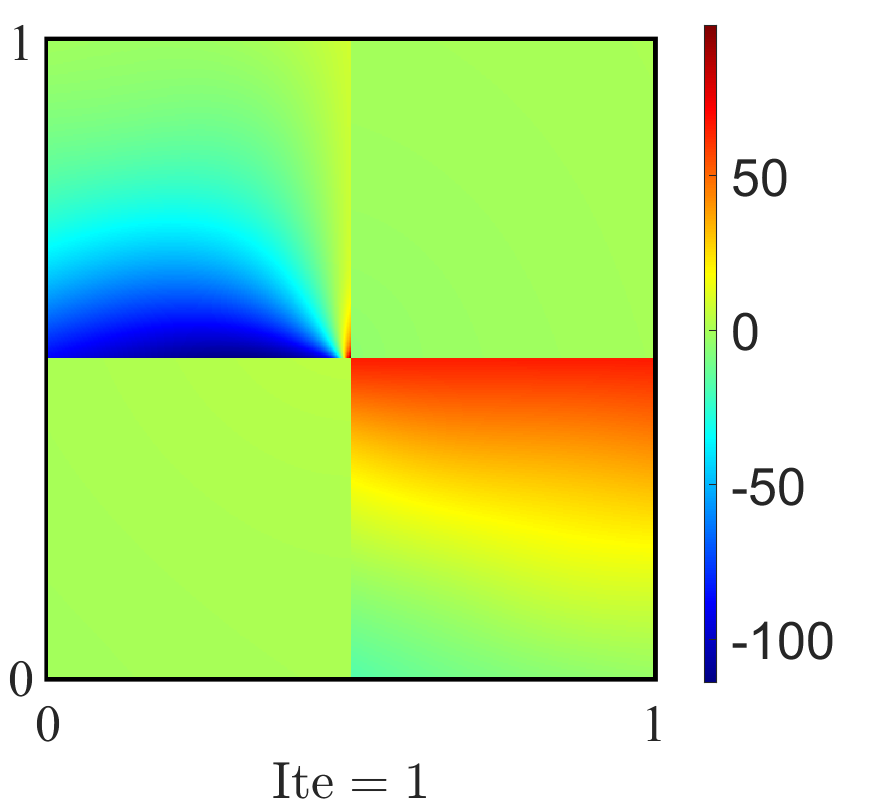}
\includegraphics[width=0.192\textwidth]{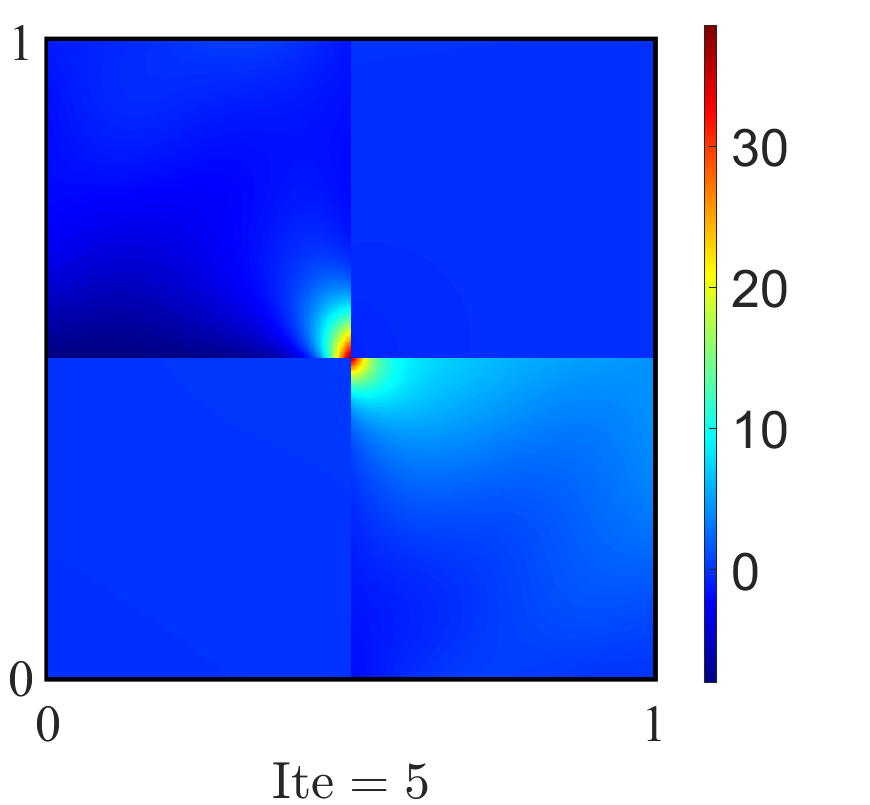}
\includegraphics[width=0.192\textwidth]{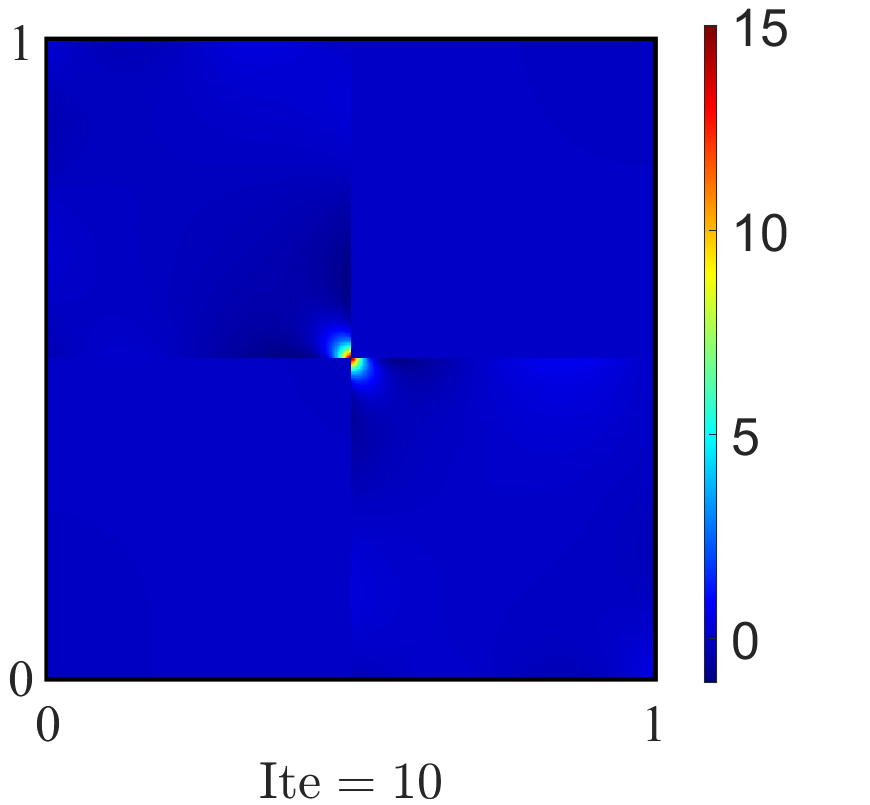}
\includegraphics[width=0.192\textwidth]{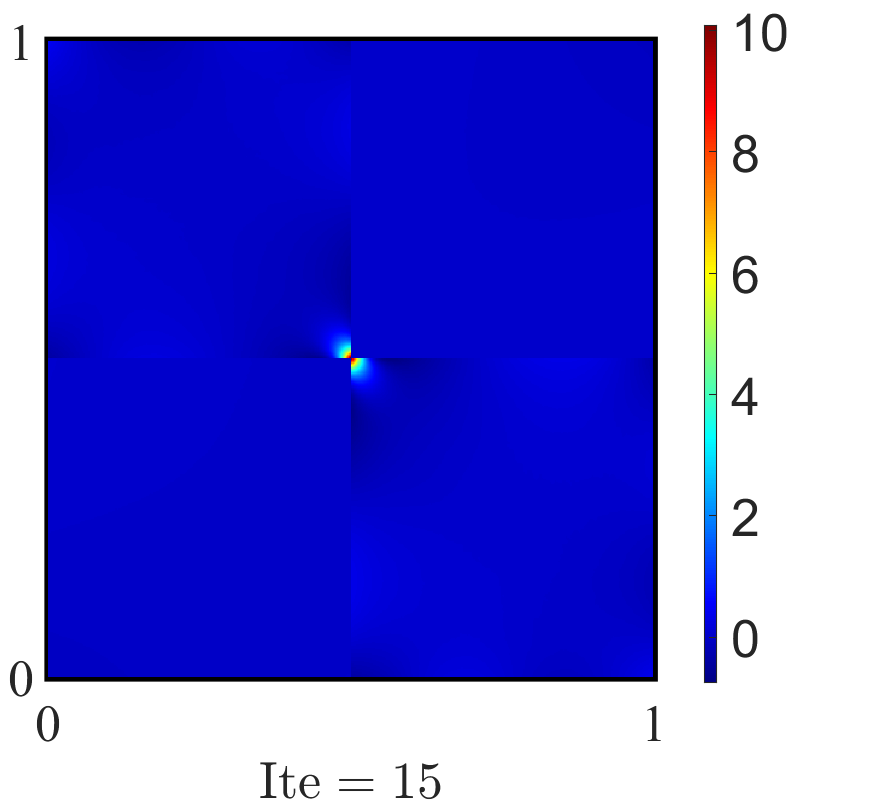}
\vspace{-0.1cm}
\caption{The pointwise absolute errors $|\hat{u}^{[k]}(x,y) - u(x,y)|$ along the outer iterations. }
\end{subfigure}
\vspace{-0.45cm}
\caption{Numerical results of example (5.5) using the DN-PINNs on testdata.}
\label{Experiments-DNLM-ex5-DN-PINNs}
\vspace{-0.7cm}
\end{figure}

\begin{figure}[H]
\centering
\begin{subfigure}[htp]{\textwidth}
\centering
\includegraphics[width=0.192\textwidth]{figure-DNLM//fig-DN-ex5-DNLM-PINN-u-NN-ite-1.png}
\includegraphics[width=0.192\textwidth]{figure-DNLM//fig-DN-ex5-DNLM-PINN-u-NN-ite-5.png}
\includegraphics[width=0.192\textwidth]{figure-DNLM//fig-DN-ex5-DNLM-PINN-u-NN-ite-10.png}
\includegraphics[width=0.192\textwidth]{figure-DNLM//fig-DN-ex5-DNLM-PINN-u-NN-ite-15.png}
\vspace{-0.1cm}
\caption{The numerical solutions $\hat{u}^{[k]}(x,y)$ along the outer iterations. }
\vspace{-0.2cm}
\end{subfigure}
\begin{subfigure}[htp]{\textwidth}
\centering
\includegraphics[width=0.192\textwidth]{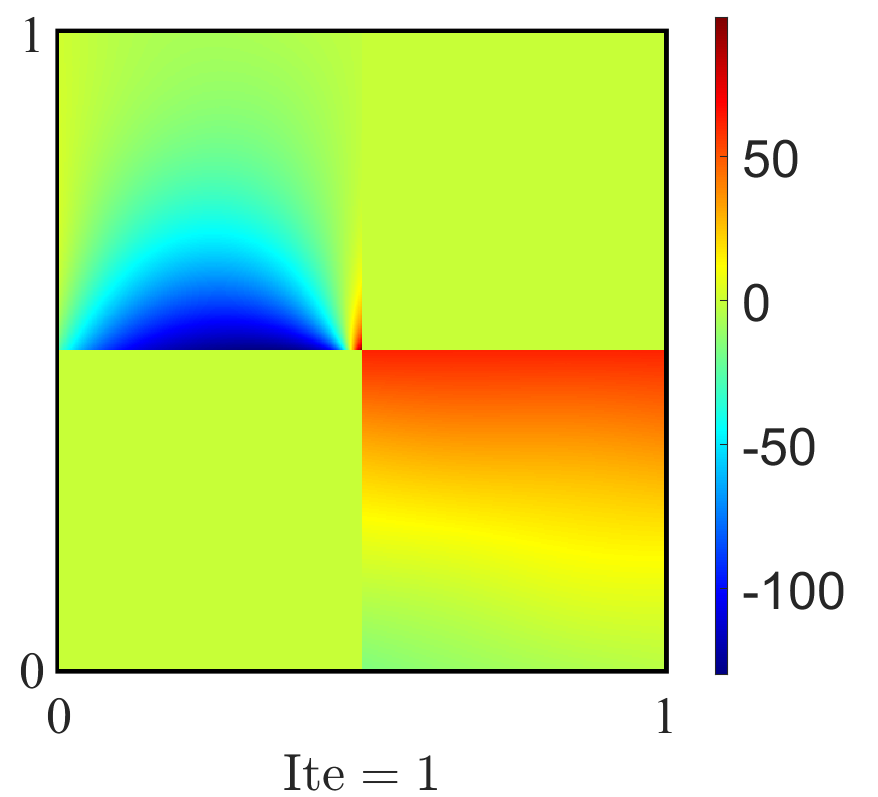}
\includegraphics[width=0.192\textwidth]{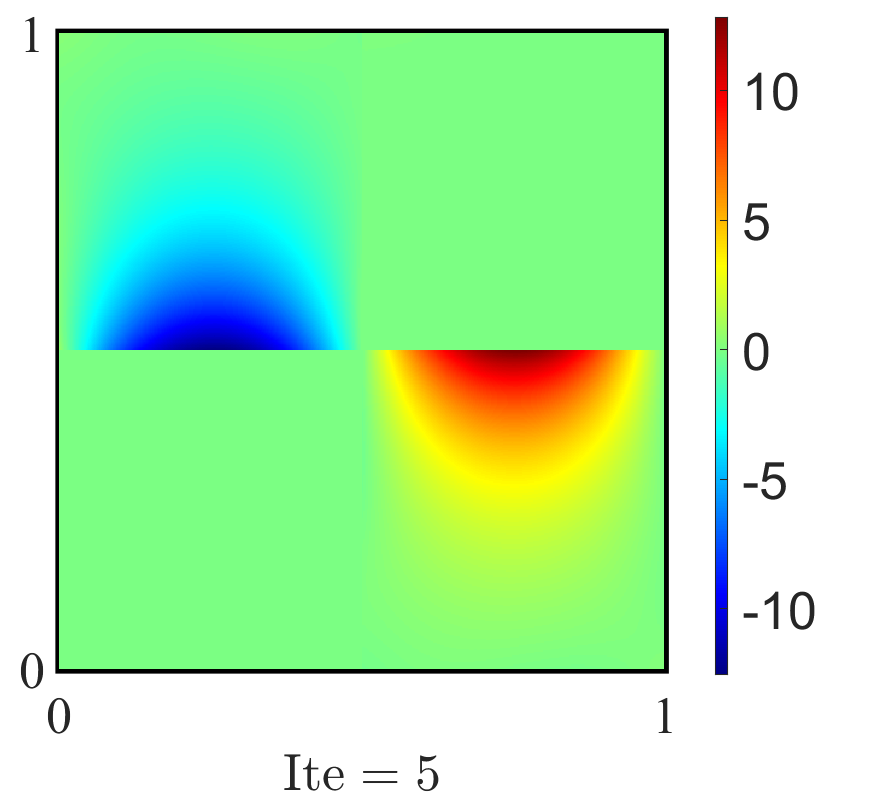}
\includegraphics[width=0.192\textwidth]{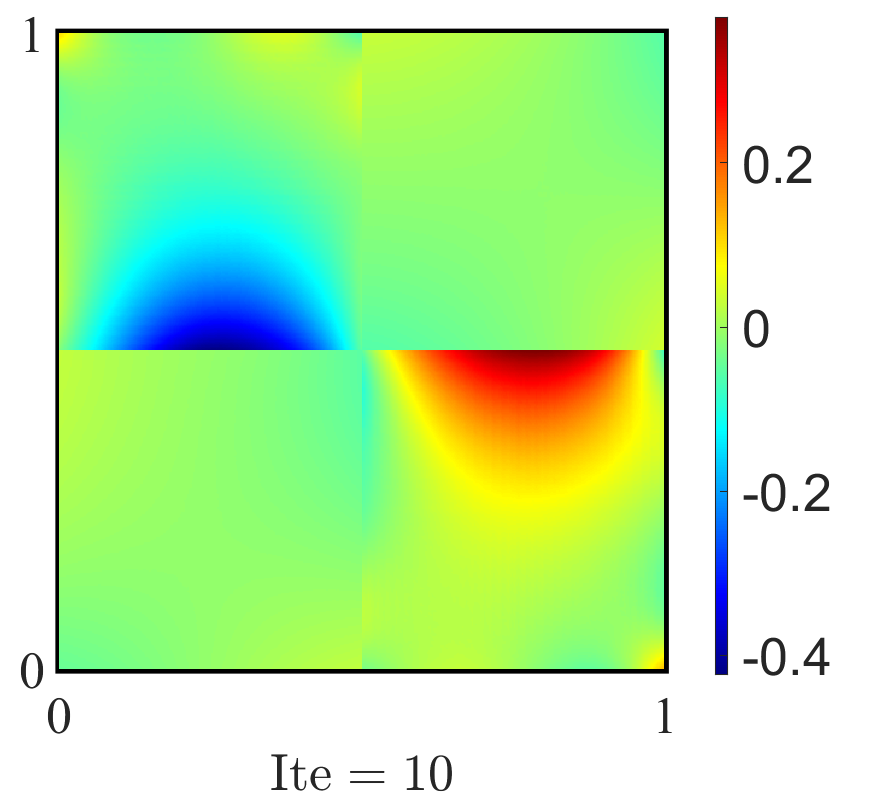}
\includegraphics[width=0.192\textwidth]{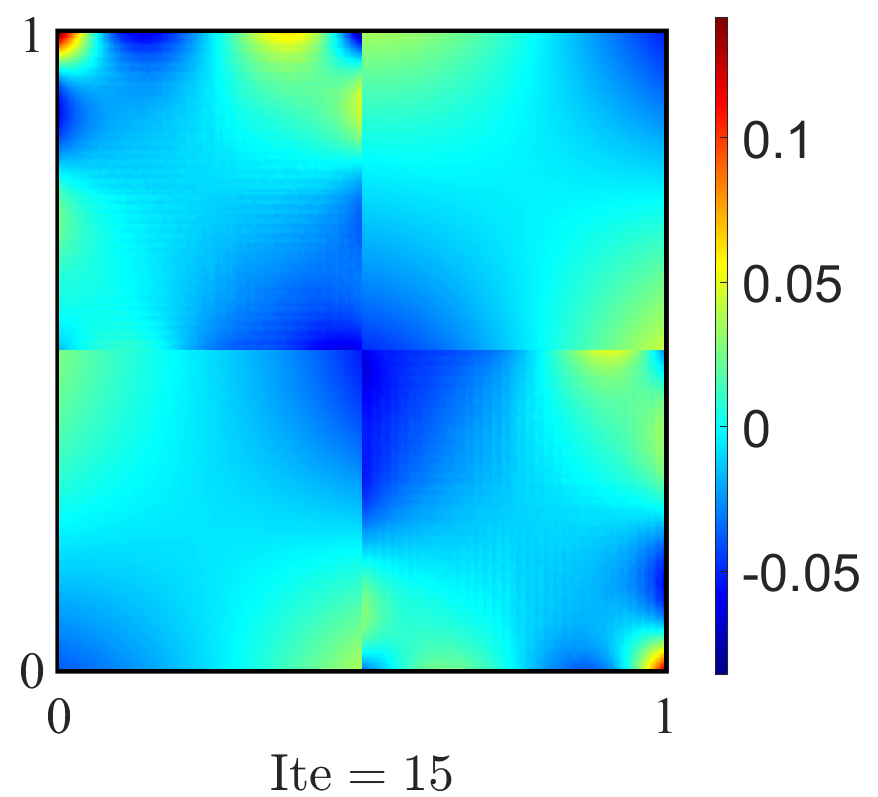}
\vspace{-0.1cm}
\caption{The pointwise absolute errors $|\hat{u}^{[k]}(x,y) - u(x,y)|$ along the outer iterations. }
\end{subfigure}
\vspace{-0.45cm}
\caption{Numerical results of example (5.5) using our DNLM (PINN) on testdata.}
\label{Experiments-DNLM-ex5-DNLM-PINN}
\vspace{-0.5cm}
\end{figure}

The numerical results using DN-PINNs and DNLM (PINN) are depicted in \autoref{Experiments-DNLM-ex5-DN-PINNs} and \autoref{Experiments-DNLM-ex5-DNLM-PINN}, which implies that our proposed method can converge to the exact solution while the DN-PINNs scheme fails.

\begin{figure}[H]
\centering
\begin{subfigure}[htp]{\textwidth}
\centering
\includegraphics[width=0.192\textwidth]{figure-DNLM//fig-DN-ex5-DNLM-DeepRitz-u-NN-ite-1.png}
\includegraphics[width=0.192\textwidth]{figure-DNLM//fig-DN-ex5-DNLM-DeepRitz-u-NN-ite-5.png}
\includegraphics[width=0.192\textwidth]{figure-DNLM//fig-DN-ex5-DNLM-DeepRitz-u-NN-ite-10.png}
\includegraphics[width=0.192\textwidth]{figure-DNLM//fig-DN-ex5-DNLM-DeepRitz-u-NN-ite-15.png}
\vspace{-0.1cm}
\caption{The numerical solutions $\hat{u}^{[k]}(x,y)$ along the outer iterations. }
\vspace{-0.2cm}
\end{subfigure}
\begin{subfigure}[htp]{\textwidth}
\centering
\includegraphics[width=0.192\textwidth]{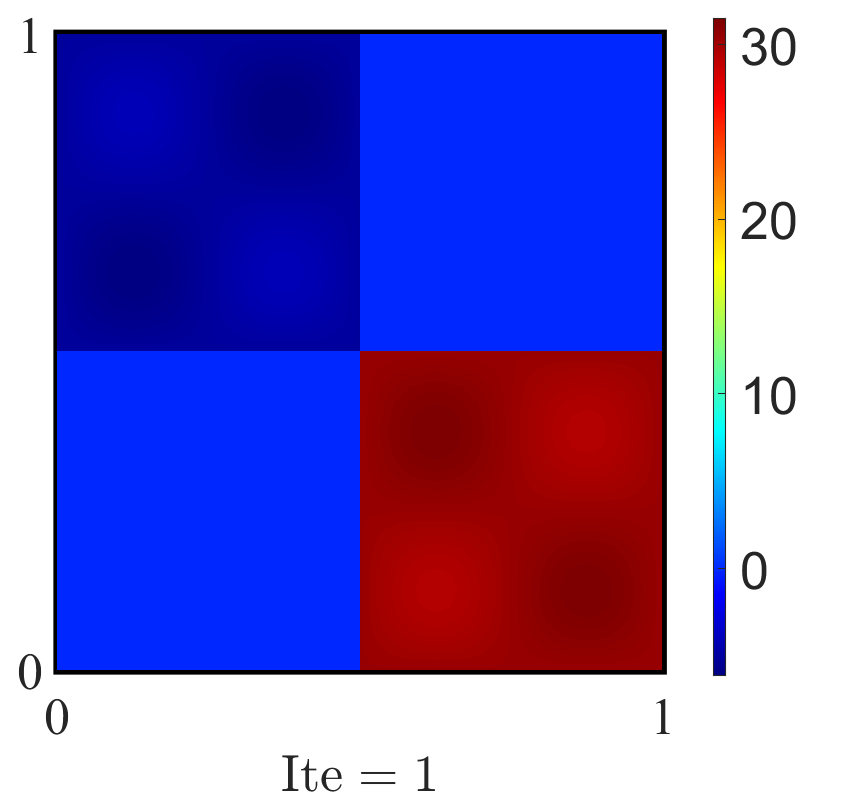}
\includegraphics[width=0.192\textwidth]{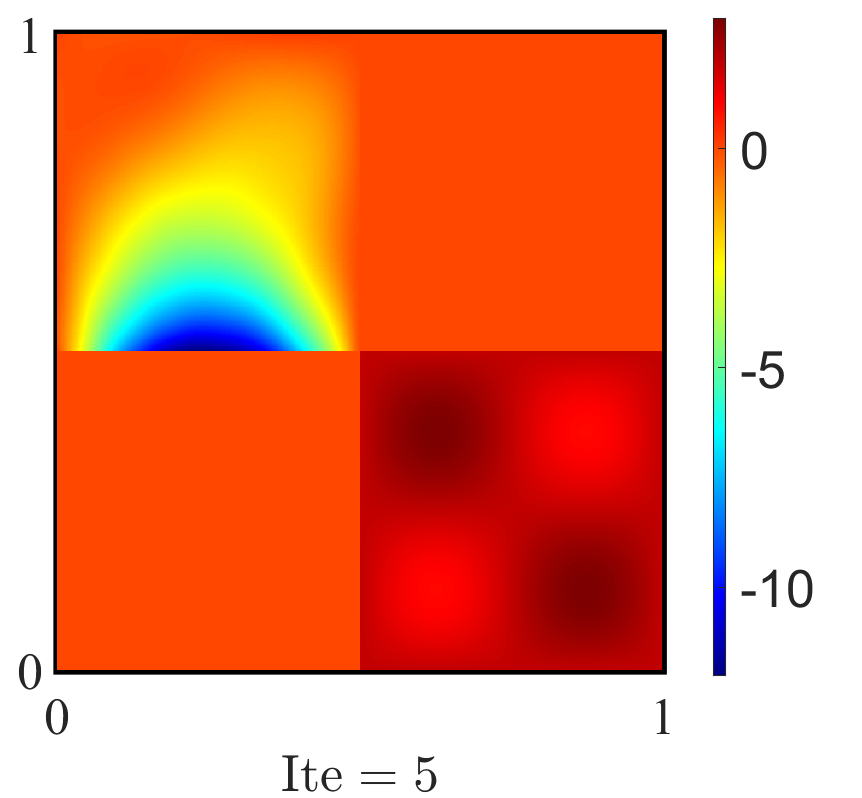}
\includegraphics[width=0.192\textwidth]{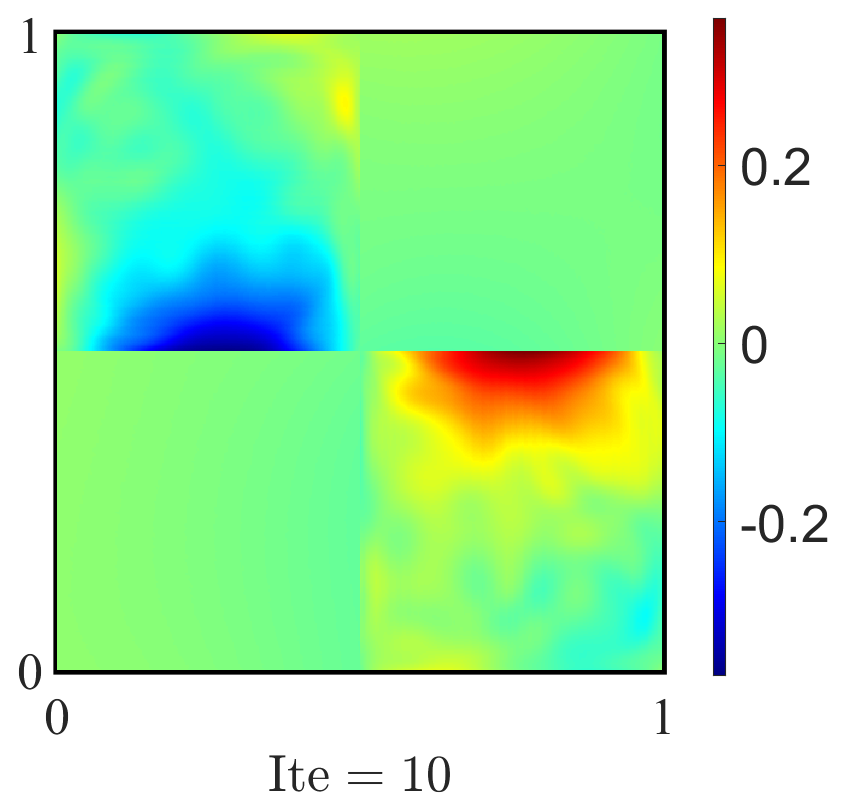}
\includegraphics[width=0.192\textwidth]{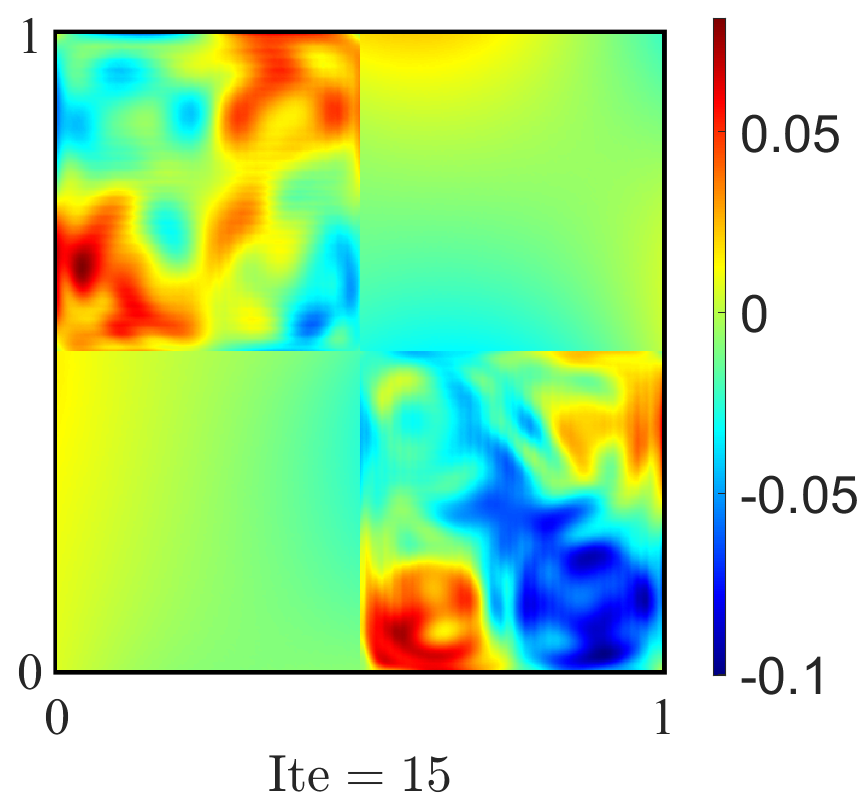}
\vspace{-0.1cm}
\caption{The pointwise absolute errors $|\hat{u}^{[k]}(x,y) - u(x,y)|$ along the outer iterations. }
\end{subfigure}
\vspace{-0.45cm}
\caption{Numerical results of (5.5) using our DNLM (deep Ritz) on testdata.}
\label{Experiments-DNLM-ex5-DNLM-DeepRitz}
\vspace{-0.5cm}
\end{figure}

\begin{figure}[H]
\centering
\begin{subfigure}[htp]{\textwidth}
\centering
\includegraphics[width=0.192\textwidth]{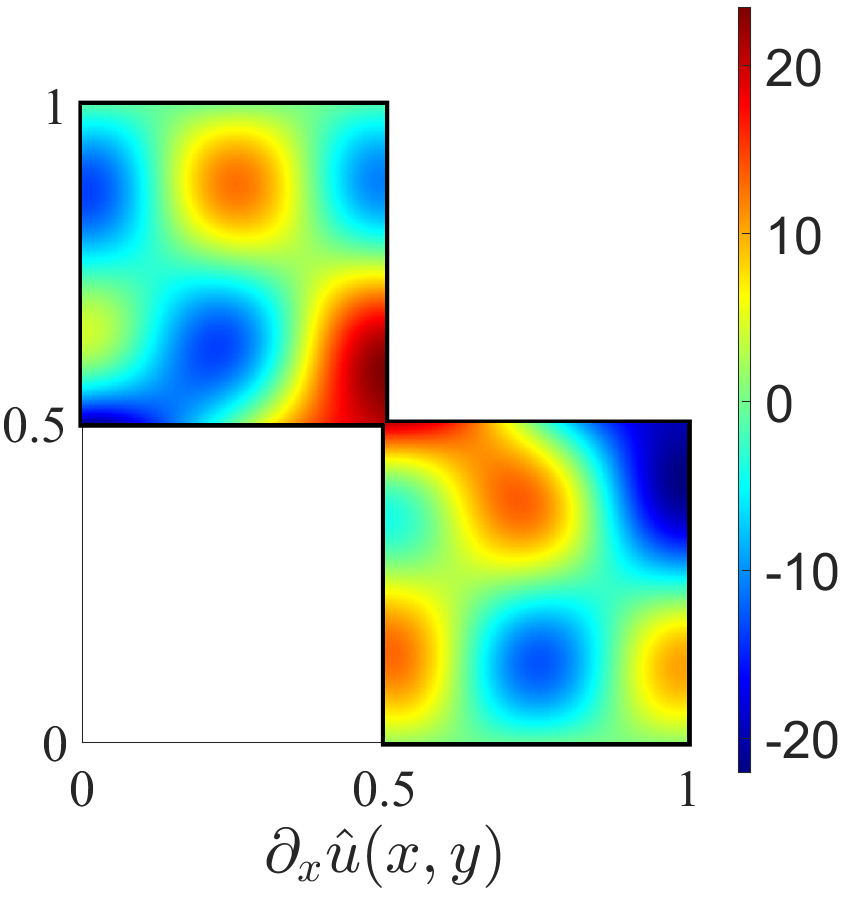}
\includegraphics[width=0.192\textwidth]{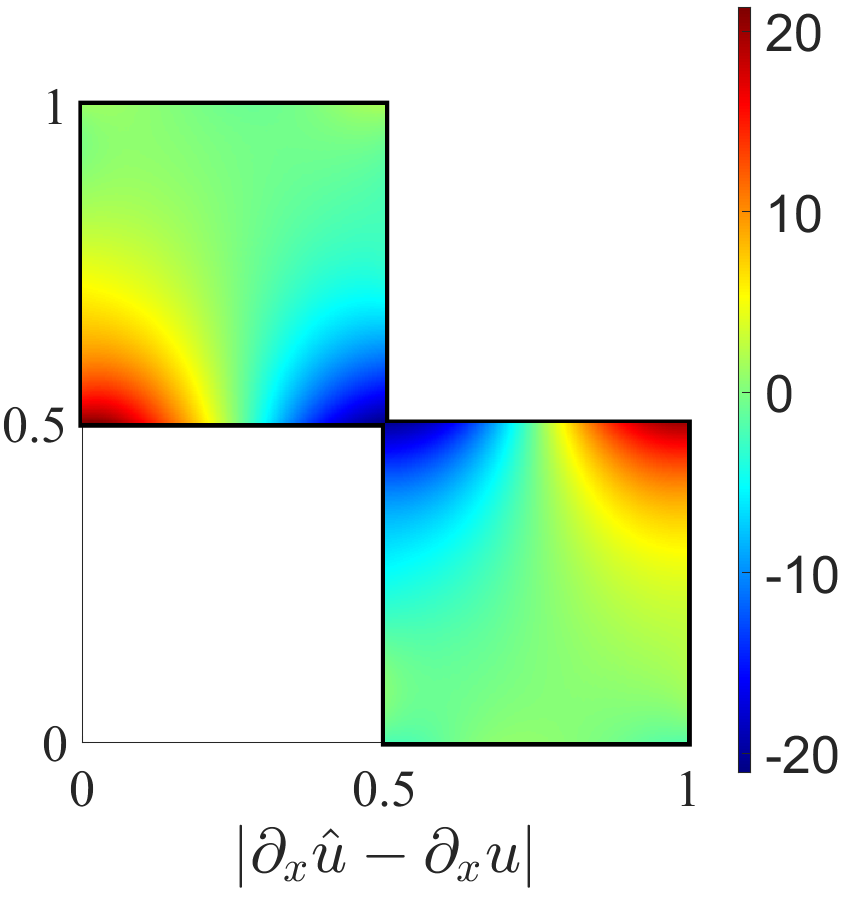}
\includegraphics[width=0.192\textwidth]{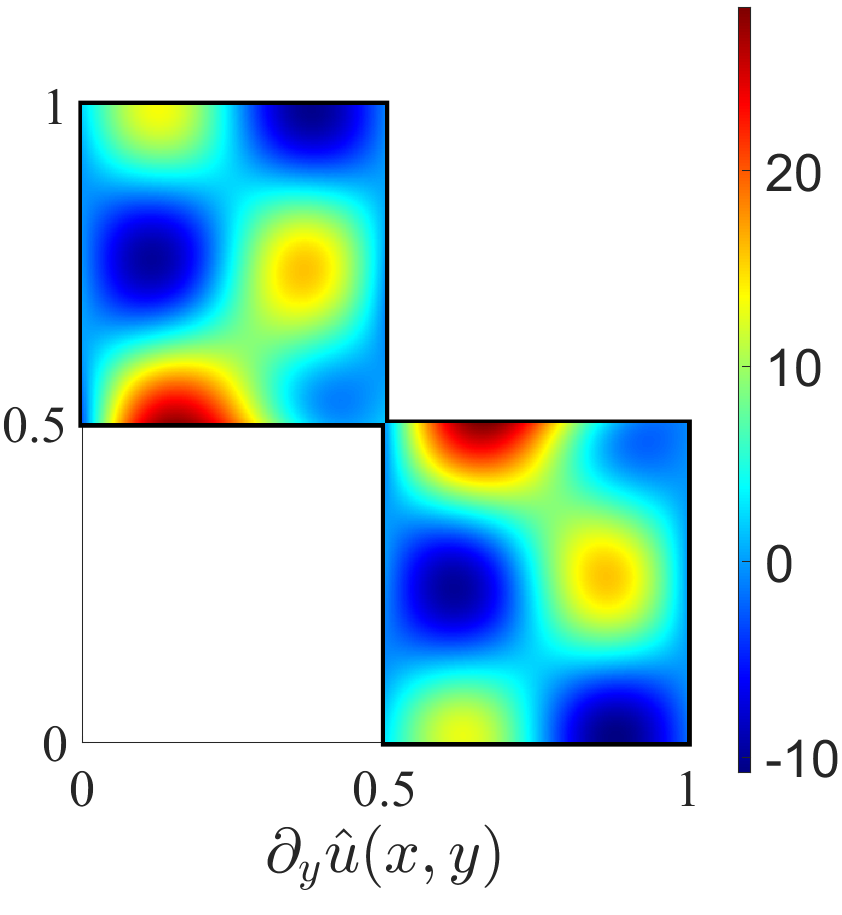}
\includegraphics[width=0.192\textwidth]{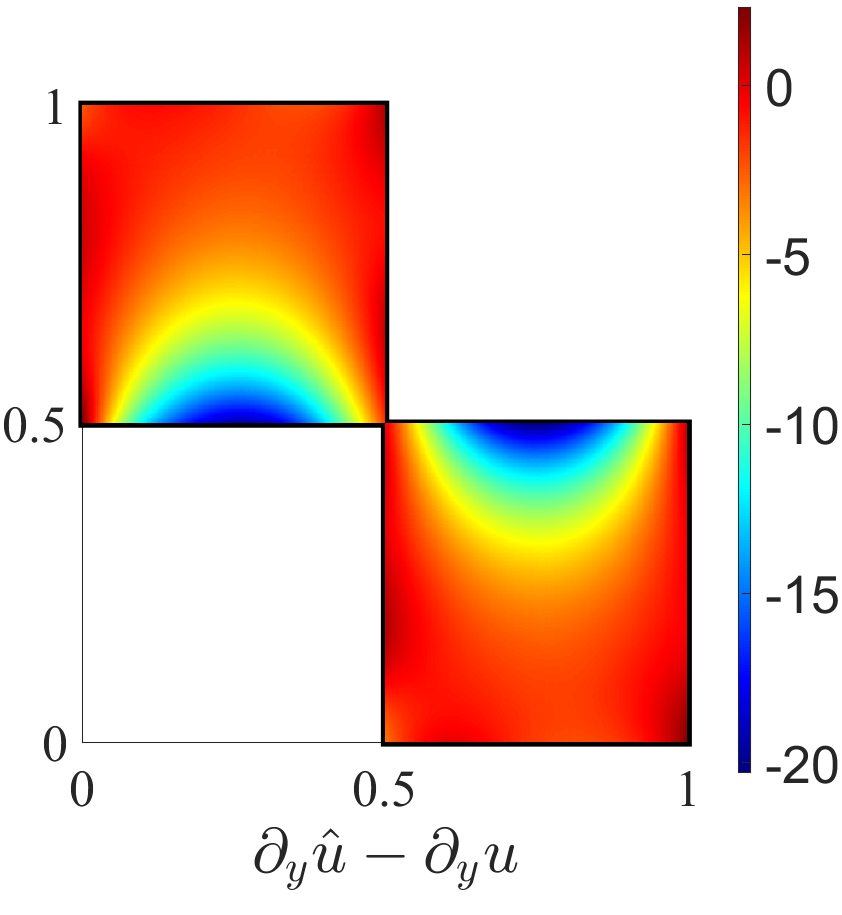}
\vspace{-0.1cm}
\caption{Derivatives $\partial_x\hat{u}^{[7]}_R$, $\partial_y\hat{u}^{[7]}_R$ and errors $|\partial_x \hat{u}^{[7]}_R - \partial_x u_R|$, $|\partial_y \hat{u}^{[7]}_R - \partial_y u_R|$ using DNLM (PINN). }
\vspace{-0.2cm}
\end{subfigure}
\begin{subfigure}[htp]{\textwidth}
\centering
\includegraphics[width=0.192\textwidth]{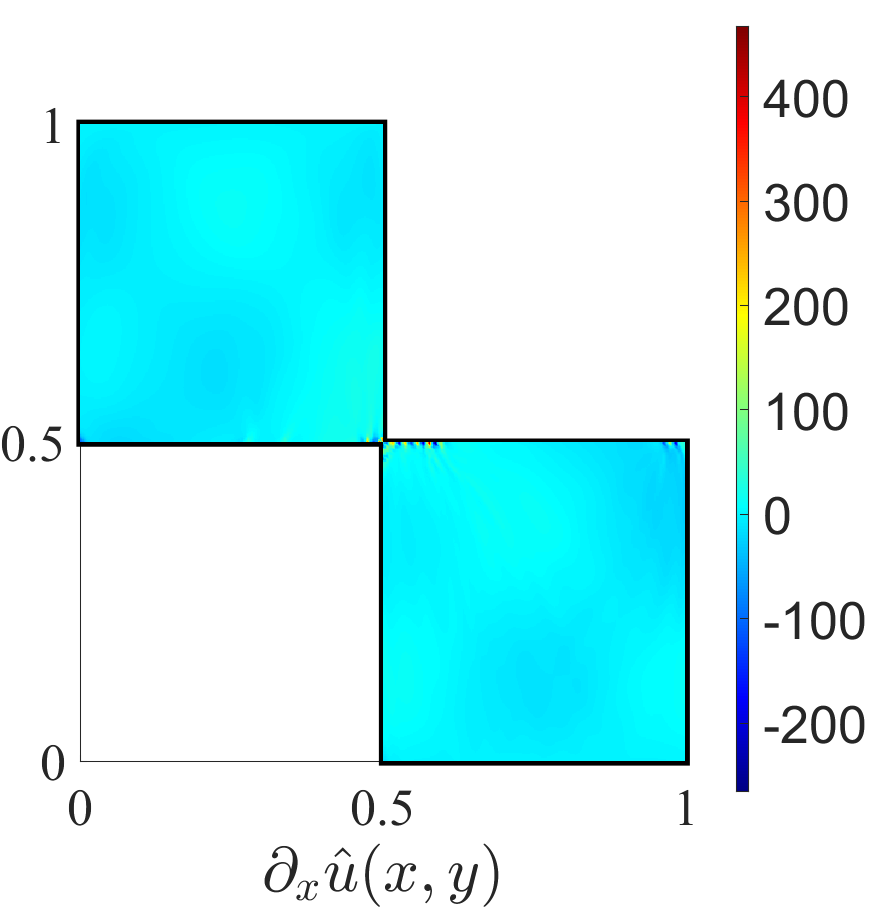}
\includegraphics[width=0.192\textwidth]{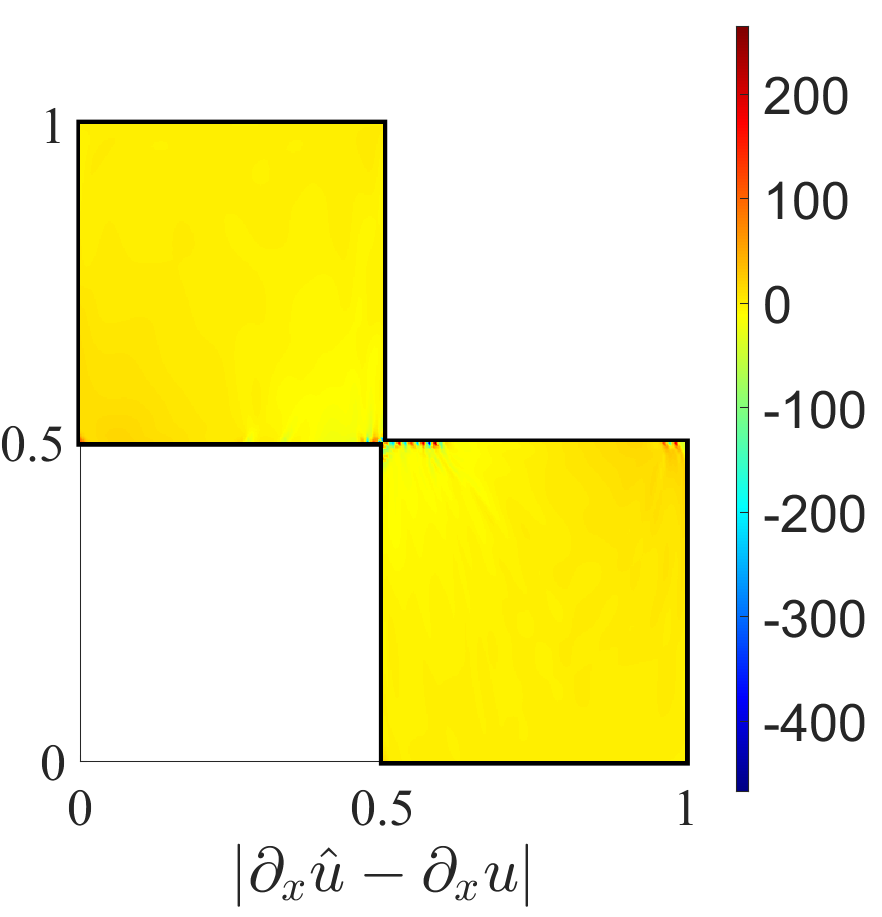}
\includegraphics[width=0.192\textwidth]{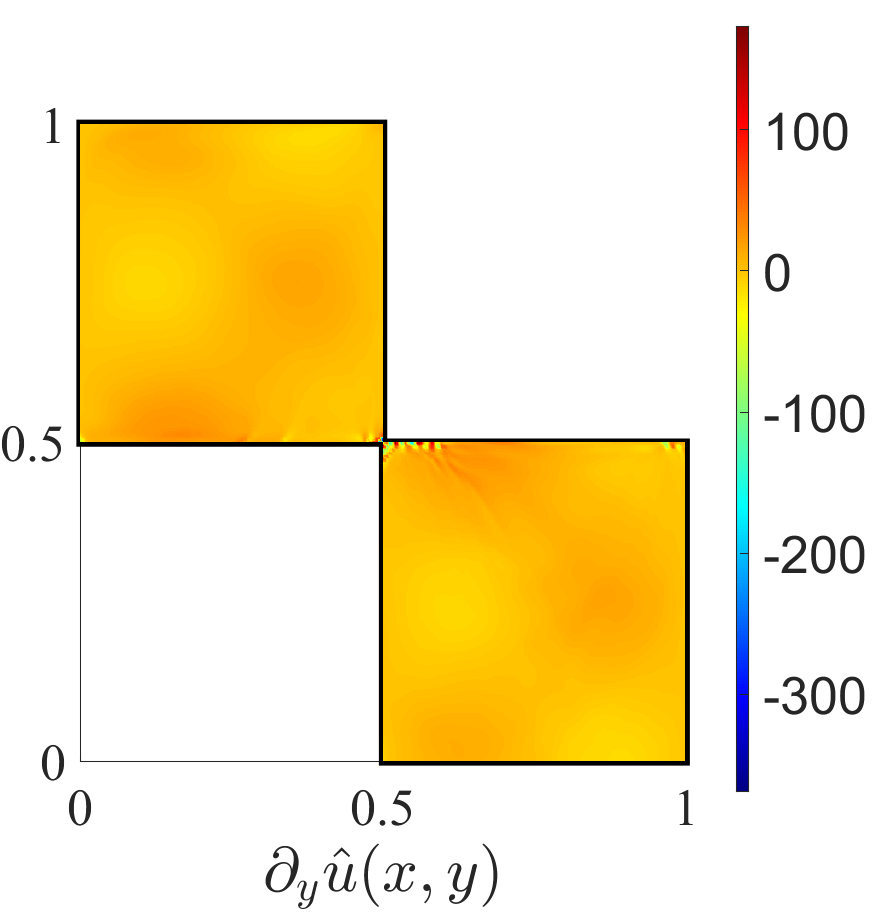}
\includegraphics[width=0.192\textwidth]{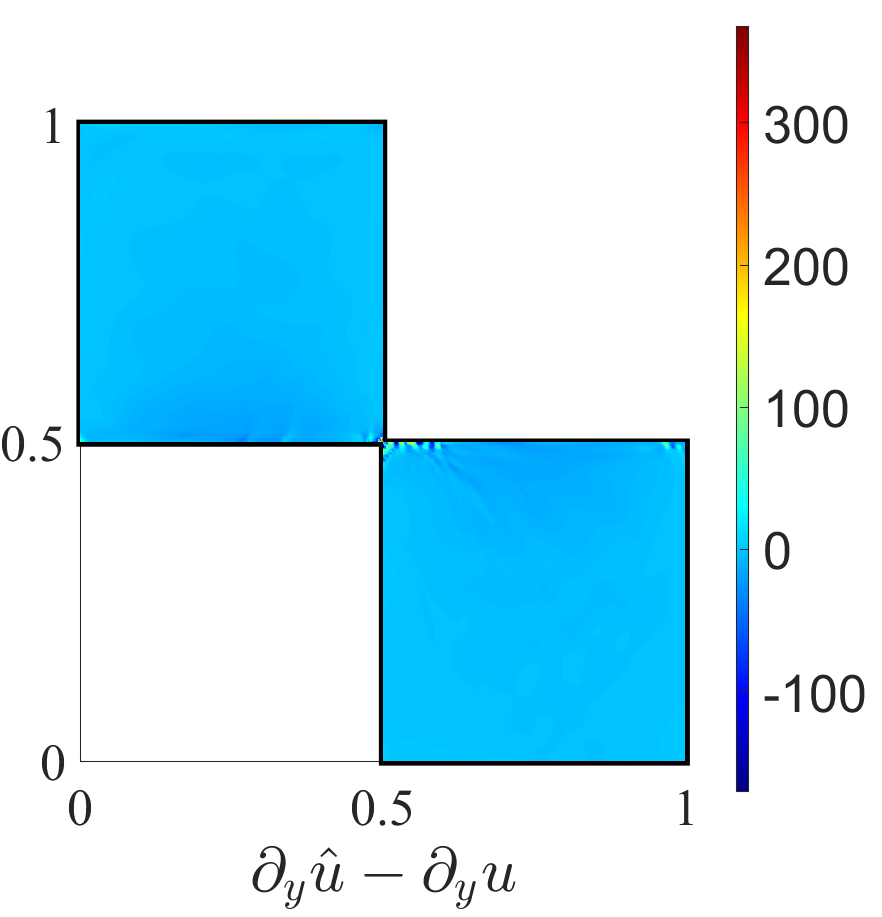}
\vspace{-0.1cm}
\caption{Derivatives $\partial_x\hat{u}^{[7]}_R$, $\partial_y\hat{u}^{[7]}_R$ and $|\partial_x \hat{u}^{[7]}_R - \partial_x u_R|$, $|\partial_y \hat{u}^{[7]}_R - \partial_y u_R|$ using DNLM (deep Ritz). }
\end{subfigure}
\vspace{-0.45cm}
\caption{Overfitting phenomenon in solving the Dirichlet subproblem of (5.5).}
\label{Experiments-DNLM-ex5-Overfit-Dirichlet-Subproblem}
\vspace{-0.7cm}
\end{figure}

Moreover, we show in \autoref{Experiments-DNLM-ex5-DNLM-DeepRitz} the results using DNLM (Deep Ritz), in \autoref{Experiments-DNLM-ex5-Overfit-Dirichlet-Subproblem} the issue of interface overfitting, and in \autoref{Experiments-DNLM-ex5-Err-Table} the statistical results.

\vspace{-0.2cm}
\begin{table}[htp]
\small
\caption{ Relative $L_2$ errors of the predicted solution along the outer iteration $k$ for example (5.5), with mean value ($\pm$ standard deviation) being reported over 5 runs.}
\centering
\renewcommand{\arraystretch}{1.1}
\begin{tabular}{ | c || c | c | c | c | c | c |  }
\hline
\multicolumn{2}{|c|}{ \diagbox[width=16em]{Relative Errors}{Outer Iterations} } & 1  & 4 & 7 & 10 & 14  \\
\hline	
\hline
\multirow{5}{*}{$ \displaystyle \!\! \frac{ \lVert \hat{u}^{[k]} - u \rVert_{L_2} } { \lVert u \rVert_{L_2} }\!\!$} & DN-PINNs & \makecell{112.5 \\ \!($\pm$\! 9.453)\!} & \makecell{20.55 \\ \!($\pm$\! 0.53)\!} &  \makecell{4.55 \\ \!($\pm$\! 0.72)\!} &  \makecell{1.88 \\ \!($\pm$\! 0.77)\!} &  \makecell{1.28 \\ \!($\pm$\! 0.45)\!}  \\ 
\cline{2-7}
& DNLM (PINN) &  \makecell{106.8 \\ \!\!($\pm$\! 2.65)\!\!} &  \makecell{20.22 \\ \!($\pm$\! 0.19)\!} &  \makecell{2.57 \\ \!($\pm$\! 0.02)\!} &  \makecell{0.37 \\ \!($\pm$\! 0.10)\!} &  \makecell{0.12 \\ \!($\pm$\! 0.01)\!} \\ 
\cline{2-7}
& \!\!\! DNLM (Deep Ritz)\! &  \makecell{57.14 \\ \!\!($\pm$\! 24.30)\!\!} &  \makecell{13.38 \\ \!($\pm$\! 2.97)\!} &  \makecell{2.51 \\ \!($\pm$\! 0.03)\!} &  \makecell{0.33 \\ \!($\pm$\! 0.01)\!} &  \makecell{0.12 \\ \!($\pm$\! 0.01)\!} \\ 
\hline
\end{tabular}
\label{Experiments-DNLM-ex5-Err-Table}
\end{table}


\section{Supplement to Section 5.2}

\begin{figure}[H]
\centering
\begin{subfigure}[htp]{\textwidth}
\centering
\includegraphics[width=0.192\textwidth]{figure-RRLM//fig-RR-ex2-RR-PINNs-u-NN-ite-1.png}
\includegraphics[width=0.192\textwidth]{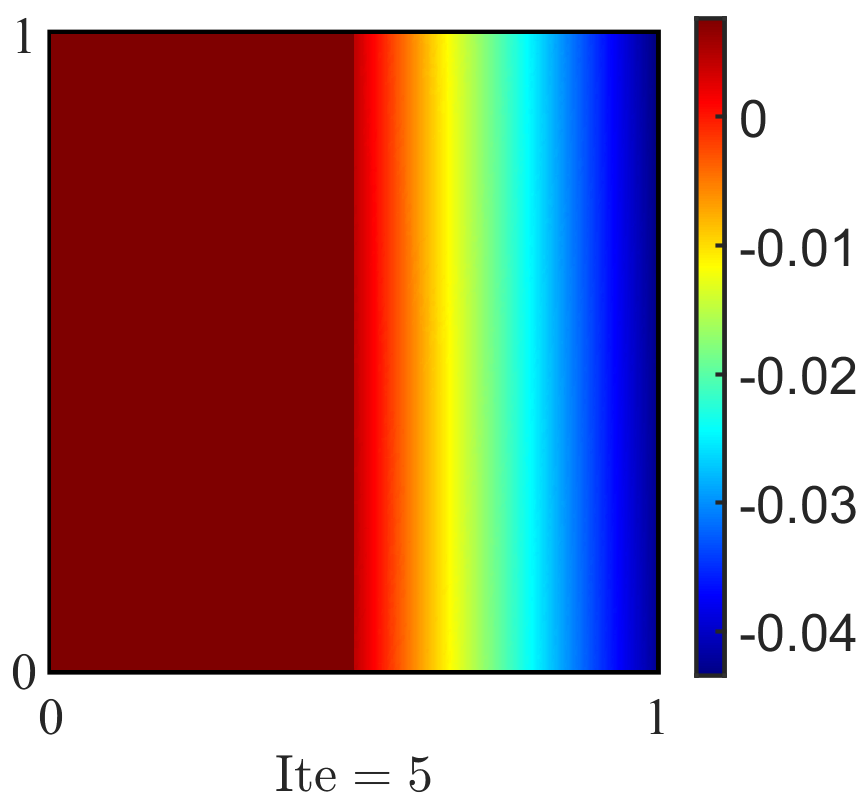}
\includegraphics[width=0.192\textwidth]{figure-RRLM//fig-RR-ex2-RR-PINNs-u-NN-ite-15.png}
\includegraphics[width=0.192\textwidth]{figure-RRLM//fig-RR-ex2-RR-PINNs-u-NN-ite-25.png}
\vspace{-0.1cm}
\caption{The numerical solutions $\hat{u}^{[k]}(x,y)$ along the outer iterations. }
\vspace{-0.2cm}
\end{subfigure}
\begin{subfigure}[htp]{\textwidth}
\centering
\includegraphics[width=0.192\textwidth]{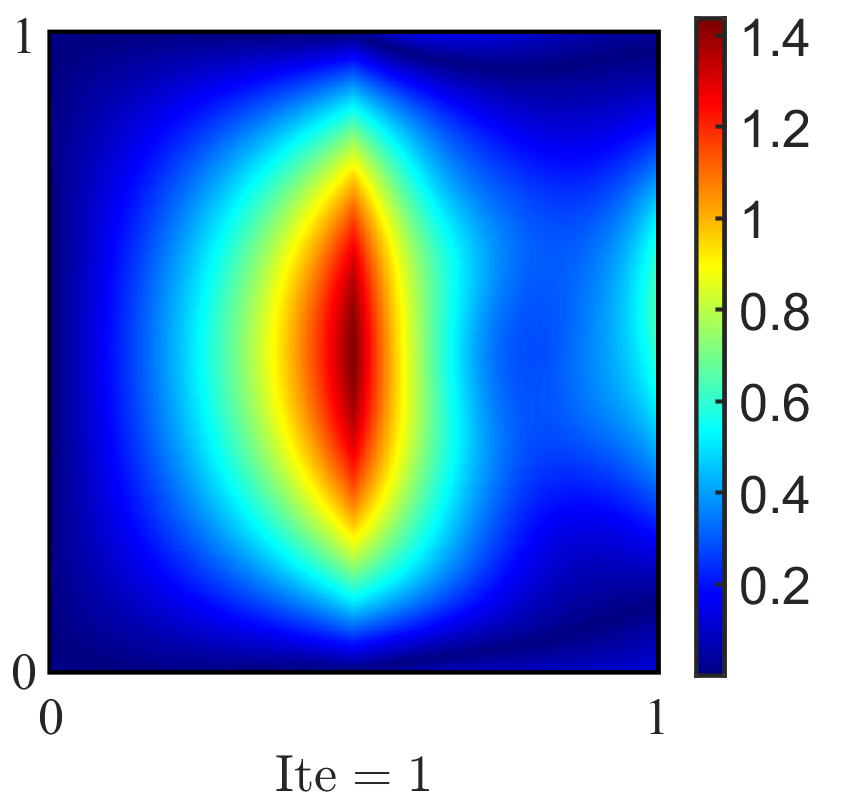}
\includegraphics[width=0.192\textwidth]{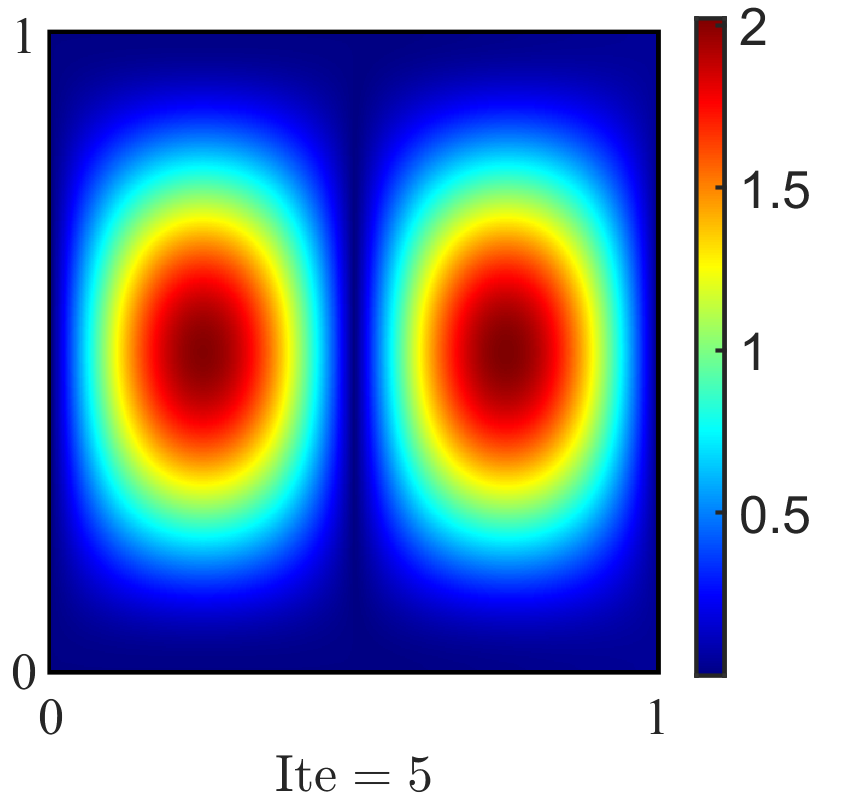}
\includegraphics[width=0.192\textwidth]{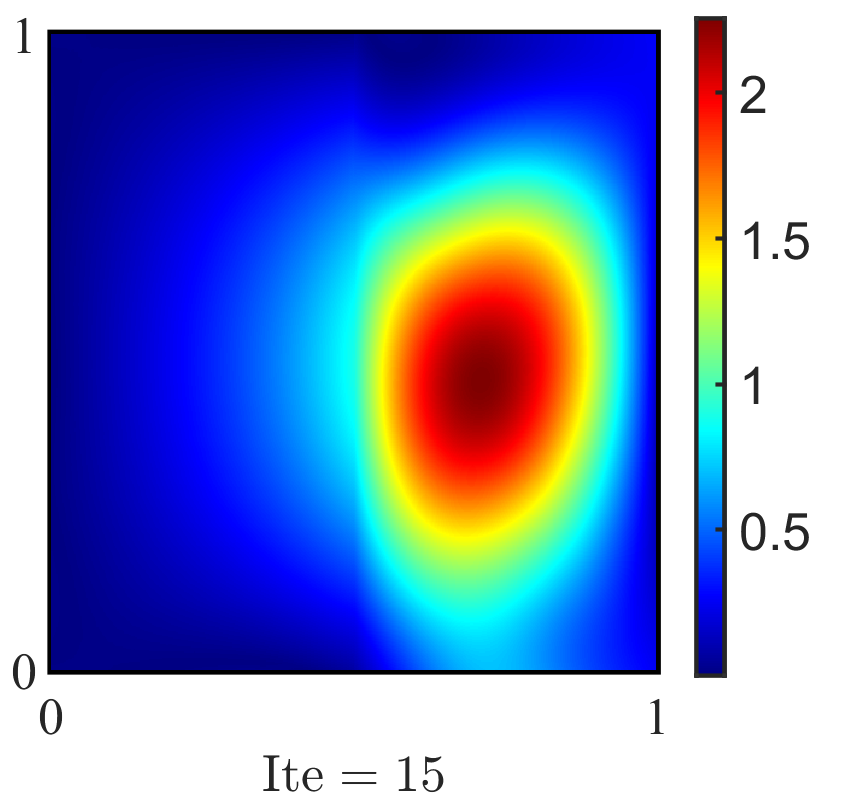}
\includegraphics[width=0.192\textwidth]{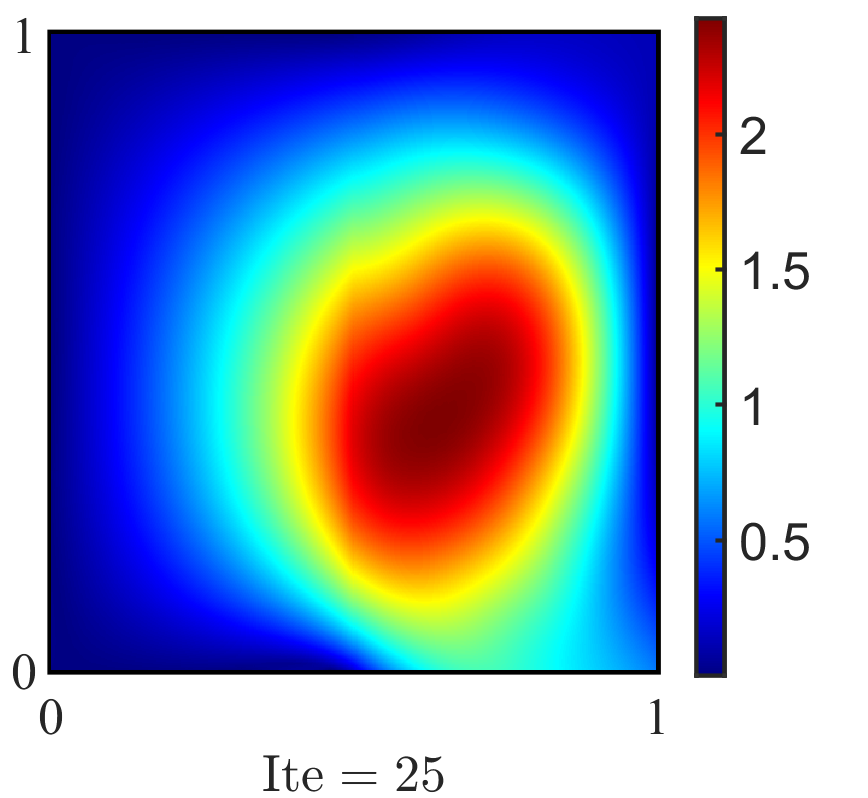}
\vspace{-0.1cm}
\caption{The pointwise absolute errors $|\hat{u}^{[k]}(x,y) - u(x,y)|$ along the outer iterations. }
\end{subfigure}
\vspace{-0.45cm}
\caption{Numerical results of (5.6) using the RR-PINNs on testdata.}
\label{Experiments-RRLM-ex1-RR-PINNs}
\vspace{-0.6cm}
\end{figure}

\begin{figure}[H]
\centering
\begin{subfigure}[htp]{\textwidth}
\centering
\includegraphics[width=0.192\textwidth]{figure-RRLM//fig-RR-ex2-RRLM-PINN-u-NN-ite-1.png}
\includegraphics[width=0.192\textwidth]{figure-RRLM//fig-RR-ex2-RRLM-PINN-u-NN-ite-2.png}
\includegraphics[width=0.192\textwidth]{figure-RRLM//fig-RR-ex2-RRLM-PINN-u-NN-ite-3.png}
\includegraphics[width=0.192\textwidth]{figure-RRLM//fig-RR-ex2-RRLM-PINN-u-NN-ite-4.png}
\vspace{-0.1cm}
\caption{The numerical solutions $\hat{u}^{[k]}(x,y)$ along the outer iterations. }
\vspace{-0.2cm}
\end{subfigure}
\begin{subfigure}[htp]{\textwidth}
\centering
\includegraphics[width=0.192\textwidth]{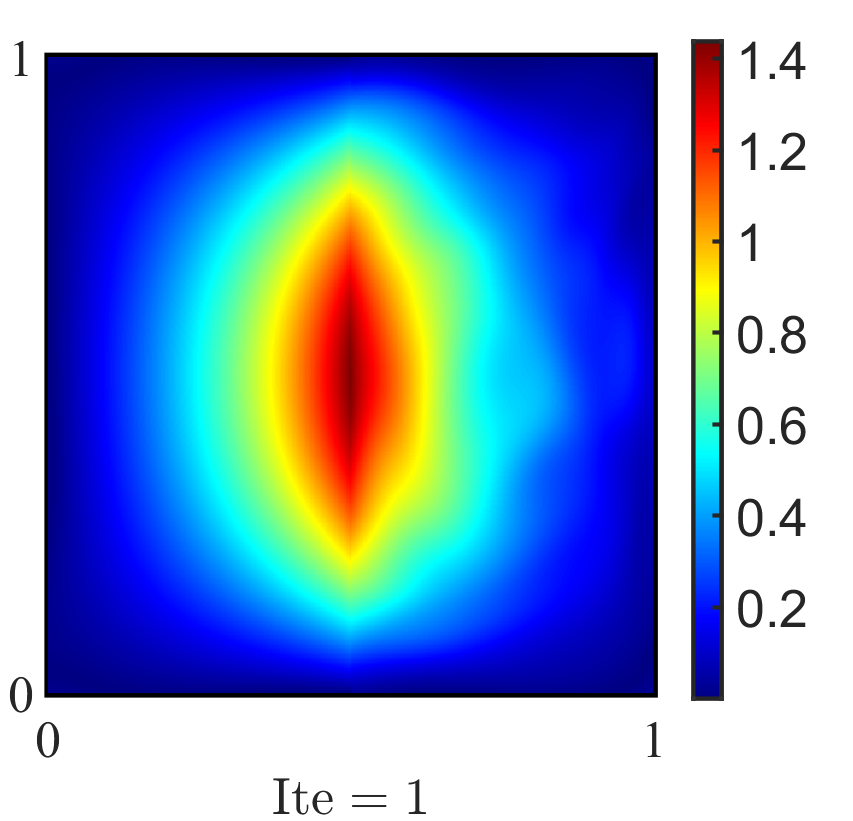}
\includegraphics[width=0.192\textwidth]{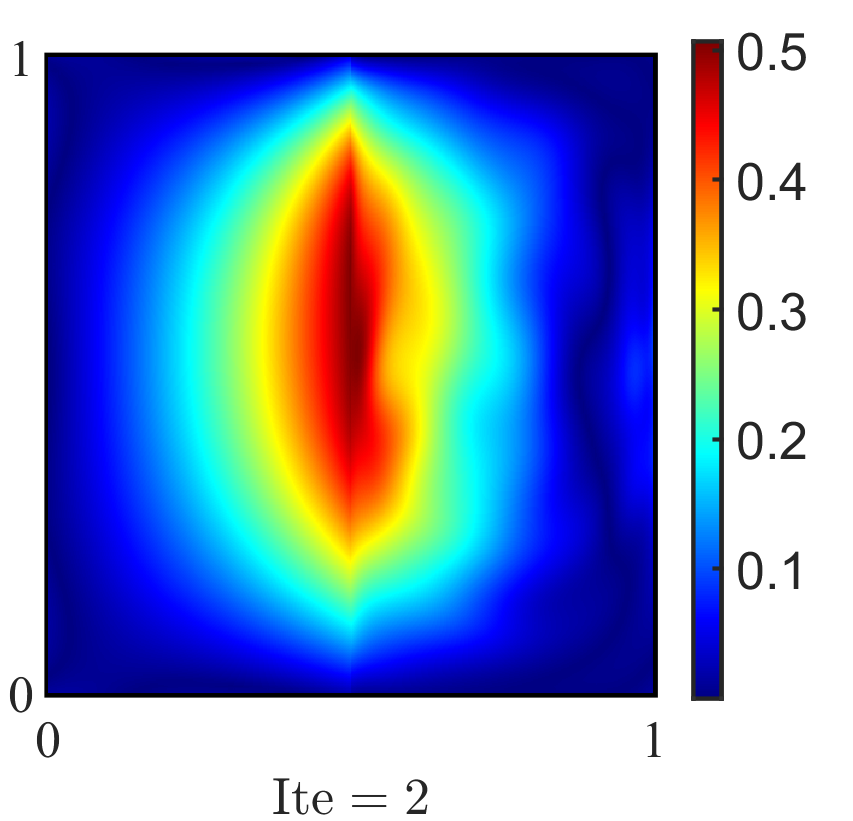}
\includegraphics[width=0.192\textwidth]{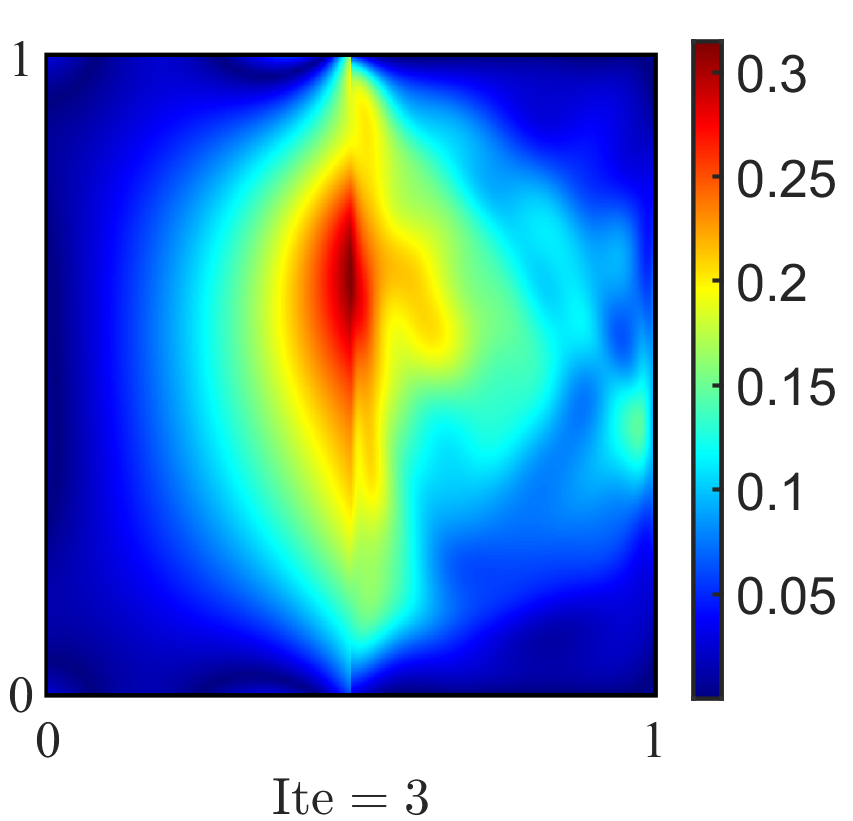}
\includegraphics[width=0.192\textwidth]{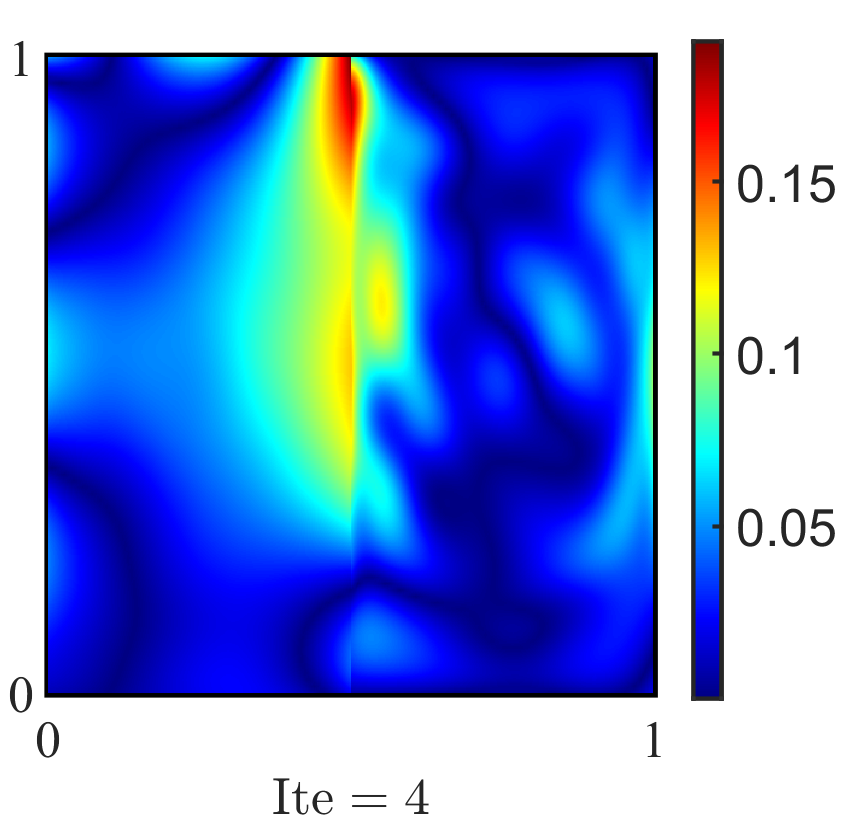}
\vspace{-0.1cm}
\caption{The pointwise absolute errors $|\hat{u}^{[k]}(x,y) - u(x,y)|$ along the outer iterations. }
\end{subfigure}
\vspace{-0.45cm}
\caption{Numerical results of (5.6) using our RRLM (PINN) on testdata.}
\label{Experiments-RRLM-ex2-RRLM-PINN}
\vspace{-0.6cm}
\end{figure}

\begin{figure}[H]
\centering
\begin{subfigure}[htp]{\textwidth}
\centering
\includegraphics[width=0.192\textwidth]{figure-RRLM//fig-RR-ex2-RRLM-DeepRitz-u-NN-ite-1.png}
\includegraphics[width=0.192\textwidth]{figure-RRLM//fig-RR-ex2-RRLM-DeepRitz-u-NN-ite-2.png}
\includegraphics[width=0.192\textwidth]{figure-RRLM//fig-RR-ex2-RRLM-DeepRitz-u-NN-ite-3.png}
\includegraphics[width=0.192\textwidth]{figure-RRLM//fig-RR-ex2-RRLM-DeepRitz-u-NN-ite-4.png}
\vspace{-0.1cm}
\caption{The numerical solutions $\hat{u}^{[k]}(x,y)$ along the outer iterations. }
\vspace{-0.2cm}
\end{subfigure}
\begin{subfigure}[htp]{\textwidth}
\centering
\includegraphics[width=0.192\textwidth]{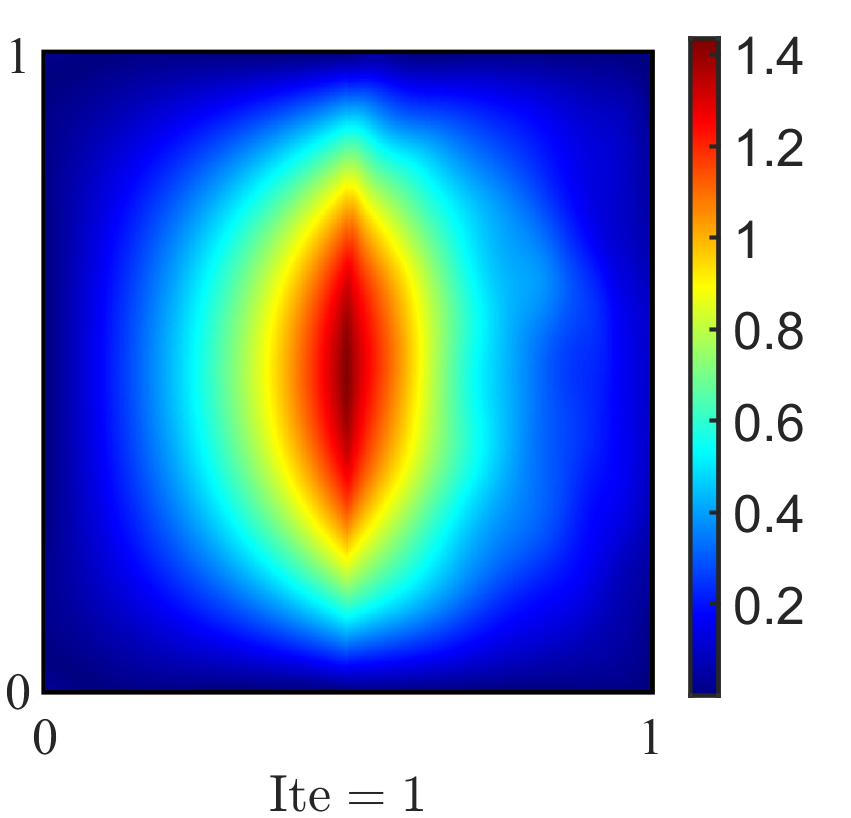}
\includegraphics[width=0.192\textwidth]{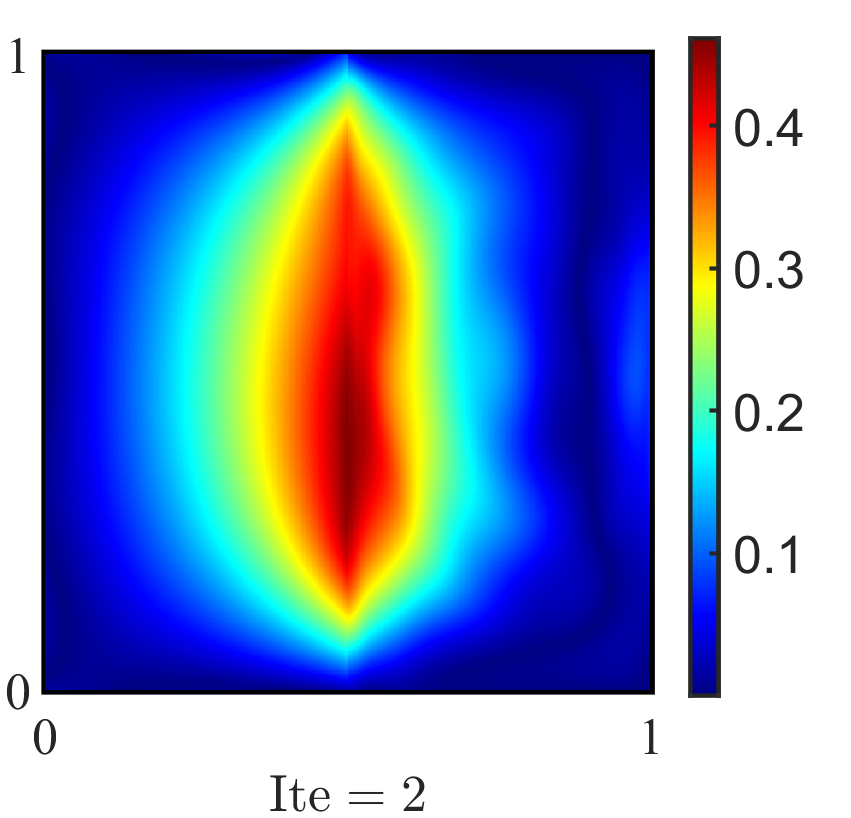}
\includegraphics[width=0.192\textwidth]{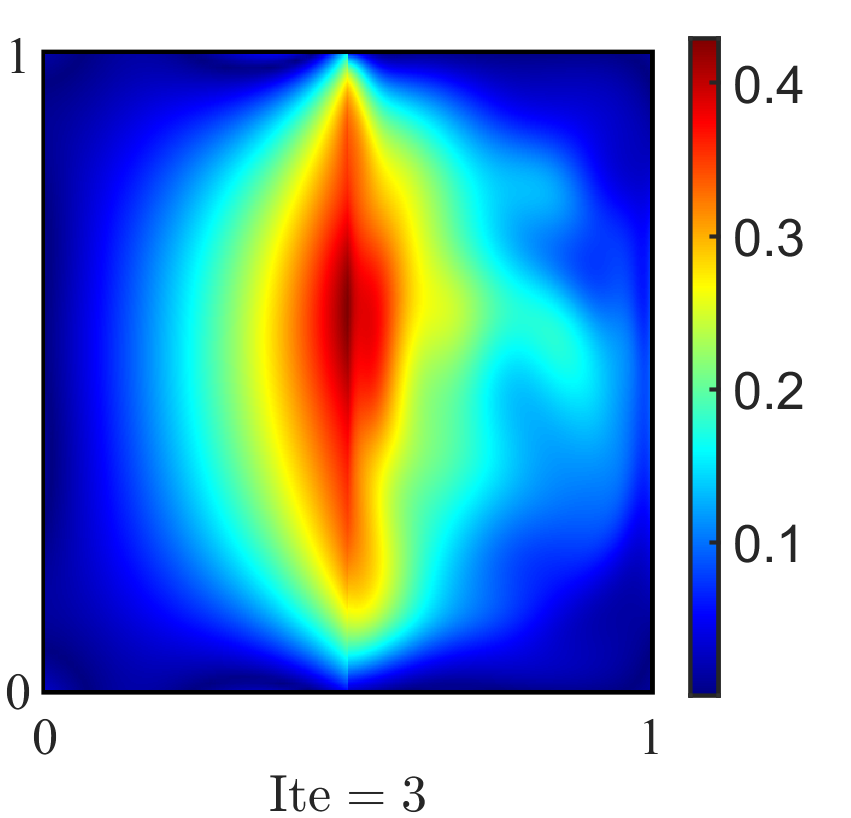}
\includegraphics[width=0.192\textwidth]{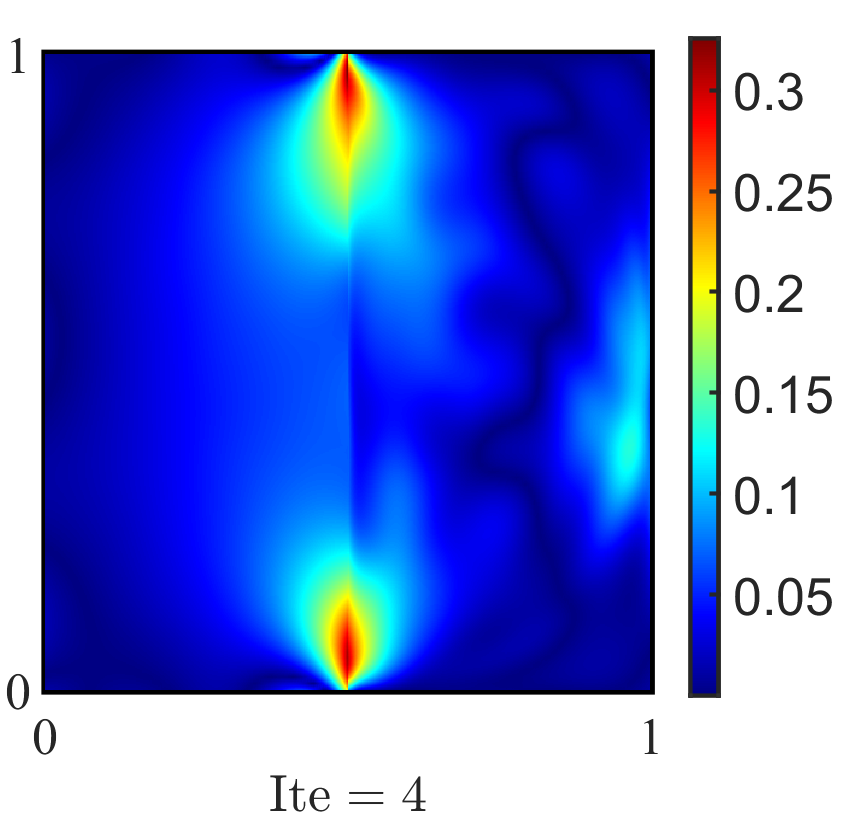}
\vspace{-0.1cm}
\caption{The pointwise absolute errors $|\hat{u}^{[k]}(x,y) - u(x,y)|$ along the outer iterations. }
\end{subfigure}
\vspace{-0.45cm}
\caption{Numerical results of (5.6) using our RRLM (deep Ritz) on testdata.}
\label{Experiments-RRLM-ex2-RRLM-DeepRitz}
\vspace{-0.5cm}
\end{figure}

When solving (5.6) using RR-PINNs for $\kappa_2=1000$, it fails to converge as depicted in \autoref{Experiments-RRLM-ex1-RR-PINNs}. On the other hand, by using our compensated deep Ritz method, the numerical results shown in \autoref{Experiments-RRLM-ex2-RRLM-PINN} and \autoref{Experiments-RRLM-ex2-RRLM-DeepRitz} imply that our methods can converge to the true solution.

